\documentclass[1pt]{amsart}

\usepackage{amsmath,amsthm, amscd, amssymb, amsfonts}

\usepackage{hhline}

\newcommand{\Z}{{\mathbb Z}}

\theoremstyle{plain}

\title{The character tables of  centralizers in  Sporadic Simple Groups of ${\rm HS}$ and ${\rm Co_3}$}
\author{ \small Shouchuan Zhang,    \ \ Jing Cheng,\ \  Jieqiong He
}
\address{ Mathematics Department of Hunan University,\newline \indent Changsha China,
410082, E-mail: z9491@yahoo.com.cn }

\date{}

\begin{document}
\newtheorem{Proposition}{\quad Proposition}[section]
\newtheorem{Theorem}{\quad Theorem}
\newtheorem{Definition}[Proposition]{\quad Definition}
\newtheorem{Corollary}[Proposition]{\quad Corollary}
\newtheorem{Lemma}[Proposition]{\quad Lemma}
\newtheorem{Example}[Proposition]{\quad Example}

\maketitle \addtocounter{section}{-1}

\numberwithin{equation}{section}

\date{}

\begin {abstract}  To classify the finite dimensional pointed Hopf
algebras with   $G= {\rm HS}$ or ${\rm Co3}$ we obtain the
representatives of conjugacy classes of $G$ and all character tables
of centralizers of these representatives by means of software {\rm
GAP}.


keywords: {\rm GAP}, Hopf algebra, sporadic simple group, character.
\end {abstract}

\section{Introduction}\label {s0}

This article is to contribute to the classification of
finite-dimensional complex pointed Hopf algebras  with sporadic
simple group $G= {\rm HS}$ or ${\rm Co3}$

 Many papers are about the classification of finite dimensional
pointed Hopf algebras, for example,  \cite{ AS02, AS00, AS05, He06,
AHS08, AG03, AFZ08, AZ07,Gr00,  Fa07, AF06, AF07, ZZC04, ZCZ08,
ZZWCY08}. In these research  ones need  the centralizers and
character tables of groups. In this paper we obtain   the
representatives of conjugacy classes of   sporadic simple group $
{\rm HS}$ and  ${\rm Co3}$, as well as  all character tables of
centralizers of these representatives by means of software {\rm
GAP}.

\section {${\rm HS}$ }

In this section $G$ denotes the sporadic simple group  ${\rm HS}$.
It is clear that $G$  is a  sub-group of  symmetric group $\mathbb
S_{100}$ (see \cite {GAP}).

\subsection {Program}

 gap$>$ x:=Indeterminate(GF(2));

gap$>$ SetName(x,"x");

gap$>$ f:=x$^{23}$-1;

gap$>$ Factors(f);

gap$>$ f:=First(Factors(f),i-$>$Degree(i)$>$1);

gap$>$ LoadPackage("guava");

gap$>$ cod:=GeneratorPolCode(f,23,GF(2));

gap$>$ IsPerfectCode(cod);

gap$>$ ext:=ExtendedCode(cod);

gap$>$ WeightDistribution(ext);

gap$>$ autext:=AutomorphismGroup(ext);

gap$>$ gens:=SmallGeneratingSet(autext);;

gap$>$ m24:=Group(gens);

gap$>$ SetName(m24,"m24");

gap$>$ Size(m24);

gap$>$ Transitivity(m24,[1..24]);

gap$>$ st:=Stabilizer(m24,[23,24],OnSets);

gap$>$ gens:=SmallGeneratingSet(st);;

gap$>$ st2:=Group(gens);

gap$>$ m22a:=Action(st2,[1..22]);

gap$>$ SetName(m22a,"m22a");

gap$>$ Size(m22a);

gap$>$ s:=SylowSubgroup(m22a,2);;

gap$>$ pciso:=IsomorphismPcGroup(s);

gap$>$ a:=Image(pciso);

gap$>$ n:=Filtered(NormalSubgroups(a),i-$>$Size(i)=16

$>$ and IsElementaryAbelian(i));

gap$>$ n:=List(n,i-$>$PreImage(pciso,i));;

gap$>$ e:=Filtered(n,i-$>$IsRegular(i,MovedPoints(i)));;Length(e);

gap$>$ e:=e[1];;

gap$>$ h:=Normalizer(m22a,e);;

gap$>$ ophom:=ActionHomomorphism(m22a,RightCosets(m22a,h),OnRight);;

gap$>$ mop:=Image(ophom);

gap$>$ DegreeAction(mop);

gap$>$ dp:=DirectProduct(m22a,mop);;

gap$>$ emb1:=Embedding(dp,1);;

gap$>$ emb2:=Embedding(dp,2);;

gap$>$ diag:=List(GeneratorsOfGroup(m22a),

$>$ i-$>$Image(emb1,i)*Image(emb2,Image(ophom,i)));;

gap$>$ diag:=Group(diag,());;

gap$>$ SetName(diag,"M22.2-99");

gap$>$ LoadPackage("grape");

gap$>$ gamma:=NullGraph(diag,100);

gap$>$ AddEdgeOrbit(gamma,[1,100]); AddEdgeOrbit(gamma,[100,1]);

gap$>$ hexad:=First(Orbits(h,[1..22]),i-$>$Length(i)=6);

gap$>$ for i in hexad do

$>$      AddEdgeOrbit(gamma,[i,23]);

$>$      AddEdgeOrbit(gamma,[23,i]);

$>$    od;.

gap$>$ Adjacency(gamma,23);

gap$>$ stab:=Stabilizer(diag,23);;

gap$>$ orbs:=Orbits(stab,[24..99]);;

gap$>$ orbreps:=List(orbs,i-$>$i[1]);

gap$>$ rep1:=orbreps[1];

gap$>$ Adjacency(gamma,rep1);

gap$>$ Intersection(hexad,Adjacency(gamma,rep1));

gap$>$ rep2:=orbreps[2];

gap$>$ Adjacency(gamma,rep2);

gap$>$ Intersection(hexad,Adjacency(gamma,rep2));

gap$>$ AddEdgeOrbit(gamma,[23,rep2]); AddEdgeOrbit(gamma,[rep2,23]);

gap$>$ IsSimpleGraph(gamma);

gap$>$ Adjacency(gamma,23);

gap$>$ IsDistanceRegular(gamma);

gap$>$ aug:=AutGroupGraph(gamma);

gap$>$ Size(aug);

gap$>$ small := SmallGeneratingSet(aug);;

gap$>$ Length(small);

gap$>$ aug2 := Group(small);

gap$>$ DisplayCompositionSeries(aug2);

gap$>$ hs:=DerivedSubgroup(aug2);Display(hs);Display(Order(hs));

ccl:=ConjugacyClasses(hs);Display(ccl);

q:=NrConjugacyClasses(hs);; Display (q);

gap $>$  ccl:=ConjugacyClasses(hs);;

gap $>$ q:=NrConjugacyClasses(hs);; Display (q);

gap $>$  for i in [1..q] do s:=Representative(ccl[i]);;Display(s);
od;

gap $>$  ccl:=ConjugacyClasses(hs);;

gap $>$ q:=NrConjugacyClasses(hs);; Display (q);

gap $>$  for i in [1..q] do
r:=Order(Representative(ccl[i]));Display(r); od;

gap $>$  for i in [1..q] do s:=Representative(ccl[i]);;
cen:=Centralizer(hs,s);; cl := ConjugacyClasses(cen);; t :=
NrConjugacyClasses(cen);;

gap $>$
s1:=Representative(ccl[1]);;cen1:=Centralizer(hs,s1);;Display
(cen1);

gap $>$  cl1:=ConjugacyClasses(cen1);char:=CharacterTable (cen1);;

gap $>$  Display (char);

gap $>$
s1:=Representative(ccl[2]);;cen1:=Centralizer(hs,s1);;Display
(cen1);

gap $>$  cl1:=ConjugacyClasses(cen1);char:=CharacterTable (cen1);;

gap $>$  Display (char);

gap $>$
s1:=Representative(ccl[3]);;cen1:=Centralizer(hs,s1);;Display
(cen1);

gap $>$  cl1:=ConjugacyClasses(cen1);char:=CharacterTable (cen1);;

gap $>$  Display (char);

gap $>$
s1:=Representative(ccl[4]);;cen1:=Centralizer(hs,s1);;Display
(cen1);

gap $>$  cl1:=ConjugacyClasses(cen1);char:=CharacterTable (cen1);;

gap $>$  Display (char);

gap $>$
s1:=Representative(ccl[5]);;cen1:=Centralizer(hs,s1);;Display
(cen1);

gap $>$  cl1:=ConjugacyClasses(cen1);char:=CharacterTable (cen1);;

gap $>$  Display (char);

gap $>$
s1:=Representative(ccl[6]);;cen1:=Centralizer(hs,s1);;Display
(cen1);

gap $>$  cl1:=ConjugacyClasses(cen1);char:=CharacterTable (cen1);;

gap $>$  Display (char);

gap $>$
s1:=Representative(ccl[7]);;cen1:=Centralizer(hs,s1);;Display
(cen1);

gap $>$  cl1:=ConjugacyClasses(cen1);char:=CharacterTable (cen1);;

gap $>$  Display (char);

gap $>$
s1:=Representative(ccl[8]);;cen1:=Centralizer(hs,s1);;Display
(cen1);

gap $>$  cl1:=ConjugacyClasses(cen1);char:=CharacterTable (cen1);;

gap $>$  Display (char);

gap $>$
s1:=Representative(ccl[9]);;cen1:=Centralizer(hs,s1);;Display
(cen1);

gap $>$  cl1:=ConjugacyClasses(cen1);char:=CharacterTable (cen1);;

gap $>$  Display (char);

gap $>$
s1:=Representative(ccl[10]);;cen1:=Centralizer(hs,s1);;Display
(cen1);

gap $>$  cl1:=ConjugacyClasses(cen1);char:=CharacterTable (cen1);;

gap $>$  Display (char);

gap $>$
s1:=Representative(ccl[11]);;cen1:=Centralizer(hs,s1);;Display
(cen1);

gap $>$  cl1:=ConjugacyClasses(cen1);char:=CharacterTable (cen1);;

gap $>$  Display (char);

gap $>$
s1:=Representative(ccl[12]);;cen1:=Centralizer(hs,s1);;Display
(cen1);

gap $>$  cl1:=ConjugacyClasses(cen1);char:=CharacterTable (cen1);;

gap $>$  Display (char);

gap $>$
s1:=Representative(ccl[13]);;cen1:=Centralizer(hs,s1);;Display
(cen1);

gap $>$  cl1:=ConjugacyClasses(cen1);char:=CharacterTable (cen1);;

gap $>$  Display (char);

gap $>$
s1:=Representative(ccl[14]);;cen1:=Centralizer(hs,s1);;Display
(cen1);

gap $>$  cl1:=ConjugacyClasses(cen1);char:=CharacterTable (cen1);;

gap $>$  Display (char);

gap $>$
s1:=Representative(ccl[15]);;cen1:=Centralizer(hs,s1);;Display
(cen1);

 gap $>$  cl1:=ConjugacyClasses(cen1);char:=CharacterTable
(cen1);;

gap $>$  Display (char);

gap $>$
s1:=Representative(ccl[16]);;cen1:=Centralizer(hs,s1);;Display
(cen1);

gap $>$  cl1:=ConjugacyClasses(cen1);char:=CharacterTable (cen1);;

gap $>$  Display (char);

gap $>$
s1:=Representative(ccl[17]);;cen1:=Centralizer(hs,s1);;Display
(cen1);

 gap $>$  cl1:=ConjugacyClasses(cen1);char:=CharacterTable
(cen1);;

gap $>$  Display (char);

gap $>$
s1:=Representative(ccl[18]);;cen1:=Centralizer(hs,s1);;Display
(cen1);

gap $>$  cl1:=ConjugacyClasses(cen1);char:=CharacterTable (cen1);;

gap $>$  Display (char);

gap $>$
s1:=Representative(ccl[19]);;cen1:=Centralizer(hs,s1);;Display
(cen1);

gap $>$  cl1:=ConjugacyClasses(cen1);char:=CharacterTable (cen1);;

gap $>$  Display (char);

gap $>$
s1:=Representative(ccl[20]);;cen1:=Centralizer(hs,s1);;Display
(cen1);

 gap $>$  cl1:=ConjugacyClasses(cen1);char:=CharacterTable
(cen1);;

gap $>$  Display (char);

gap $>$
s1:=Representative(ccl[21]);;cen1:=Centralizer(hs,s1);;Display
(cen1);

 gap $>$  cl1:=ConjugacyClasses(cen1);char:=CharacterTable
(cen1);;

 gap $>$  Display (char);

gap $>$
s1:=Representative(ccl[22]);;cen1:=Centralizer(hs,s1);;Display
(cen1);

gap $>$  cl1:=ConjugacyClasses(cen1);char:=CharacterTable (cen1);;

gap $>$  Display (char);

gap $>$
s1:=Representative(ccl[23]);;cen1:=Centralizer(hs,s1);;Display
(cen1);

gap $>$  cl1:=ConjugacyClasses(cen1);char:=CharacterTable (cen1);;

gap $>$  Display (char);

gap $>$
s1:=Representative(ccl[24]);;cen1:=Centralizer(hs,s1);;Display
(cen1);

gap $>$  cl1:=ConjugacyClasses(cen1);char:=CharacterTable (cen1);;

gap $>$  Display (char);
\subsection {The character tables }

The order of $G$ is 44352000.

The generators of $G$ are:\\
$(1,85,8,96)(2,98,39,68,20,93,56,30)(3,7,47,29,4,16,38,62)(5,45,46,89,12,53,52,80)(6,28)(9,24,76,10,\\
13,60,25,18)(11,86,37,73,22,92,32,63)(14,35,43,97,19,54,59,81)(17,66)(23,41,55,95,77,31,40,84)(26,64)\\
(27,42,58,88,69,57,50,82)(33,75,78,99,44,34,65,87)(36,71,61,90,51,74,70,91)(48,94,67,83)(49,100,72,79)$,

$(1,34)(2,97,39,84,66,83,8,96)(3,15,19,62,73,67,24,69)(4,56)(5,88,76,89,21,92,77,82)(6,44,71,46,23,60,\\
16,12)(7,87,18,94,40,100,70,98)(9,42,72,29,45,61,25,14)(10,63,43,36,74,78,27,17)(11,33,75,58,30,59,28,53)\\
(13,54,65,32)(20,37,64,68)(22,38,41,55,35,49,26,57)(31,50)(47,79,52,99,51,85,48,95)(81,86,93,90)$.

The representatives of conjugacy classes of   $G$ are:\\
$s_1$=(1) , $s_2$=(1, 33, 64, 44, 42, 63, 9, 78) (2, 15, 7, 14, 47,
57, 18, 87) (3, 8, 69, 39, 96, 91, 68, 52) (4, 37,
    74, 29, 26, 72, 62, 88) (5, 25, 51, 82, 76, 12, 80, 60) (6, 10, 50, 75) (11, 59, 13, 99, 41, 94, 85,
    22) (16, 48, 53, 86) (17, 81) (19, 92, 45, 35, 31, 55, 58, 65) (20, 28, 49, 73, 38, 79, 43, 54) (21, 98,
    70, 32, 90, 89, 93, 83) (23, 36, 100, 71, 97, 61, 67, 24) (27, 77, 56, 66) (30, 95, 84, 46) (34, 40) ,
$s_3$=(1, 64, 42, 9) (2, 7, 47, 18) (3, 69, 96, 68) (4, 74, 26, 62)
(5, 51, 76, 80) (6, 50) (8, 39, 91, 52) (10,
    75) (11, 13, 41, 85) (12, 60, 25, 82) (14, 57, 87, 15) (16, 53) (19, 45, 31, 58) (20, 49, 38, 43) (21,
    70, 90, 93) (22, 59, 99, 94) (23, 100, 97, 67) (24, 36, 71, 61) (27, 56) (28, 73, 79, 54) (29, 72, 88,
    37) (30, 84) (32, 89, 83, 98) (33, 44, 63, 78) (35, 55, 65, 92) (46, 95) (48, 86) (66, 77) ,
$s_4$=(1, 42) (2, 47) (3, 96) (4, 26) (5, 76) (7, 18) (8, 91) (9,
64) (11, 41) (12, 25) (13, 85) (14, 87) (15,
    57) (19, 31) (20, 38) (21, 90) (22, 99) (23, 97) (24, 71) (28, 79) (29, 88) (32, 83) (33, 63) (35,
    65) (36, 61) (37, 72) (39, 52) (43, 49) (44, 78) (45, 58) (51, 80) (54, 73) (55, 92) (59, 94) (60,
    82) (62, 74) (67, 100) (68, 69) (70, 93) (89, 98) ,  $s_5$=(1, 70, 64, 90, 42, 93, 9, 21) (2, 36, 7, 71, 47,
    61, 18, 24) (3, 73, 69, 79, 96, 54, 68, 28) (4, 13, 74, 41, 26, 85, 62, 11) (5, 45, 51, 31, 76, 58, 80,
    19) (6, 46, 50, 95) (8, 20, 39, 49, 91, 38, 52, 43) (10, 30, 75, 84) (12, 35, 60, 55, 25, 65, 82,
    92) (14, 23, 57, 100, 87, 97, 15, 67) (16, 66, 53, 77) (17, 40) (22, 29, 59, 72, 99, 88, 94, 37) (27,
    86, 56, 48) (32, 78, 89, 33, 83, 44, 98, 63) (34, 81) ,
$s_6$=(1, 57, 42, 15) (2, 78, 47, 44) (3, 31, 96, 19) (4, 59, 26,
94) (5, 73, 76, 54) (7, 33, 18, 63) (8, 55, 91, 92) (9, 14, 64, 87)
(11, 29, 41, 88) (12,
    20, 25, 38) (13, 72, 85, 37) (16, 56) (17, 34) (21, 23, 90, 97) (22, 74, 99, 62) (24, 32, 71, 83) (27,
    53) (28, 51, 79, 80) (35, 39, 65, 52) (36, 89, 61, 98) (40, 81) (43, 60, 49, 82) (45, 69, 58, 68) (48,
    66) (67, 70, 100, 93) (77, 86) ,  $s_7$=(1, 14) (2, 63) (3, 45) (4, 22) (5, 28) (6, 50) (7, 78) (8, 35) (9,
    15) (10, 75) (11, 37) (12, 43) (13, 29) (16, 27) (17, 34) (18, 44) (19, 68) (20, 60) (21, 67) (23,
    70) (24, 98) (25, 49) (26, 99) (30, 84) (31, 69) (32, 36) (33, 47) (38, 82) (39, 55) (40, 81) (41,
    72) (42, 87) (46, 95) (48, 77) (51, 73) (52, 92) (53, 56) (54, 80) (57, 64) (58, 96) (59, 74) (61,
    83) (62, 94) (65, 91) (66, 86) (71, 89) (76, 79) (85, 88) (90, 100) (93, 97) ,
$s_8$=(1, 16, 61, 41, 26, 85, 62, 47, 66, 21) (2, 36, 7, 71, 11, 64,
93, 74, 18, 24) (3, 43, 69, 91, 30, 6, 58, 92,
    51, 35) (4, 17, 42, 90, 40, 13, 9, 53, 77, 70) (5, 20, 19, 49, 79, 95, 75, 82, 80, 39) (8, 45, 12, 31, 65,
    84, 50, 96, 52, 73) (10, 38, 54, 55, 28, 60, 68, 25, 76, 46) (14, 27, 83, 72, 99, 88, 94, 33, 86,
    67) (15, 56, 48, 23, 22, 34, 87, 100, 81, 29) (32, 78, 89, 37, 57, 97, 59, 44, 98, 63) ,
$s_9$=(1, 61, 26, 62, 66) (2, 7, 11, 93, 18) (3, 69, 30, 58, 51) (4,
42, 40, 9, 77) (5, 19, 79, 75, 80) (6, 92, 35,
    43, 91) (8, 12, 65, 50, 52) (10, 54, 28, 68, 76) (13, 53, 70, 17, 90) (14, 83, 99, 94, 86) (15, 48, 22,
    87, 81) (16, 41, 85, 47, 21) (20, 49, 95, 82, 39) (23, 34, 100, 29, 56) (24, 36, 71, 64, 74) (25, 46,
    38, 55, 60) (27, 72, 88, 33, 67) (31, 84, 96, 73, 45) (32, 89, 57, 59, 98) (37, 97, 44, 63, 78) ,
$s_{10}$=(1, 43, 80, 66, 35, 75, 62, 92, 79, 26, 6, 19, 61, 91, 5)
(2, 36, 78, 18, 24, 63, 93, 74, 44, 11, 64, 97, 7,
    71, 37) (3, 38, 83, 51, 46, 14, 58, 25, 86, 30, 60, 94, 69, 55, 99) (4, 9, 42, 77, 40) (8, 20, 27, 52, 39,
    67, 50, 82, 33, 65, 95, 88, 12, 49, 72) (10, 16, 31, 76, 21, 45, 68, 47, 73, 28, 85, 96, 54, 41, 84) (13,
    34, 87, 90, 23, 22, 17, 56, 48, 70, 29, 15, 53, 100, 81) (32, 59, 89, 98, 57) ,
$s_{11}$=(1, 75, 6) (2, 63, 64) (3, 14, 60) (5, 35, 26) (7, 78, 74)
(8, 67, 95) (10, 45, 85) (11, 37, 24) (12, 27,
    82) (13, 22, 29) (15, 34, 17) (16, 68, 96) (18, 44, 71) (19, 43, 62) (20, 50, 88) (21, 28, 84) (23, 70,
    81) (25, 69, 83) (30, 99, 46) (31, 47, 54) (33, 49, 52) (36, 93, 97) (38, 58, 94) (39, 65, 72) (41, 76,
    73) (48, 100, 90) (51, 86, 55) (53, 87, 56) (61, 80, 92) (66, 79, 91) ,
$s_{12}$=(1, 19, 23, 84, 81, 82, 52, 53, 38, 45, 10, 4, 72, 77, 32,
86, 98, 85, 50, 42) (2, 12, 71, 51, 96, 29, 80,
    25, 43, 48, 16, 22, 79, 26, 75, 33, 24, 59, 92, 46) (3, 95, 100, 47, 15, 20, 61, 41, 44, 5) (6, 67, 90,
    60, 94, 66, 17, 63, 36, 76, 14, 97, 70, 87, 35, 58, 7, 34, 40, 69) (8, 83, 55, 93, 57, 9, 30, 64, 37,
    73) (11, 65, 49, 54, 13, 88, 74, 62, 39, 56, 21, 89, 91, 68, 18, 27, 31, 28, 99, 78) ,
$s_{13}$=(1, 4, 23, 77, 81, 86, 52, 85, 38, 42, 10, 19, 72, 84, 32,
82, 98, 53, 50, 45) (2, 22, 71, 26, 96, 33, 80,
    59, 43, 46, 16, 12, 79, 51, 75, 29, 24, 25, 92, 48) (3, 95, 100, 47, 15, 20, 61, 41, 44, 5) (6, 97, 90,
    87, 94, 58, 17, 34, 36, 69, 14, 67, 70, 60, 35, 66, 7, 63, 40, 76) (8, 83, 55, 93, 57, 9, 30, 64, 37,
    73) (11, 89, 49, 68, 13, 27, 74, 28, 39, 78, 21, 65, 91, 54, 18, 88, 31, 62, 99, 56) ,
$s_{14}$=(1, 23, 81, 52, 38, 10, 72, 32, 98, 50) (2, 71, 96, 80, 43,
16, 79, 75, 24, 92) (3, 100, 15, 61, 44) (4, 77,
    86, 85, 42, 19, 84, 82, 53, 45) (5, 95, 47, 20, 41) (6, 90, 94, 17, 36, 14, 70, 35, 7, 40) (8, 55, 57, 30,
    37) (9, 64, 73, 83, 93) (11, 49, 13, 74, 39, 21, 91, 18, 31, 99) (12, 51, 29, 25, 48, 22, 26, 33, 59,
    46) (27, 28, 78, 65, 54, 88, 62, 56, 89, 68) (34, 69, 67, 60, 66, 63, 76, 97, 87, 58) ,
$s_{15}$=(1, 81, 38, 72, 98) (2, 96, 43, 79, 24) (3, 15, 44, 100,
61) (4, 86, 42, 84, 53) (5, 47, 41, 95, 20) (6, 94,
    36, 70, 7) (8, 57, 37, 55, 30) (9, 73, 93, 64, 83) (10, 32, 50, 23, 52) (11, 13, 39, 91, 31) (12, 29, 48,
    26, 59) (14, 35, 40, 90, 17) (16, 75, 92, 71, 80) (18, 99, 49, 74, 21) (19, 82, 45, 77, 85) (22, 33, 46,
    51, 25) (27, 78, 54, 62, 89) (28, 65, 88, 56, 68) (34, 67, 66, 76, 87) (58, 69, 60, 63, 97) ,
$s_{16}$=(1, 82, 10, 86) (2, 29, 16, 33) (3, 20) (4, 98, 19, 52) (5,
15) (6, 66, 14, 58) (7, 67, 17, 97) (8, 9) (11,
    88, 21, 27) (12, 80, 22, 24) (13, 56, 18, 78) (23, 53, 72, 85) (25, 79, 59, 71) (26, 92, 51, 43) (28,
    49, 62, 91) (30, 83) (31, 65, 74, 89) (32, 42, 81, 45) (34, 90, 63, 70) (35, 69, 94, 76) (36, 87, 40,
    60) (37, 93) (38, 77, 50, 84) (39, 68, 99, 54) (41, 100) (44, 47) (46, 96, 48, 75) (55, 64) (57,
    73) (61, 95) ,
$s_{17}$=(1, 37, 46, 32, 43, 24, 76) (2, 33, 91, 13, 87, 20, 96) (3,
63, 74, 81, 55, 95, 85) (4,
    54, 72, 35, 84, 65, 31) (5, 44, 49, 99, 16, 19, 42) (7, 10, 12, 94, 82, 86, 64) (8, 38, 17, 15, 56, 39,
    47) (9, 88, 57, 14, 92, 90, 18) (11, 97, 71, 75, 78, 23, 41) (21, 62, 52, 29, 60, 59, 25) (22, 34, 83,
    73, 77, 40, 26) (27, 45, 69, 48, 36, 93, 58) (28, 61, 30, 68, 67, 100, 70) (50, 98, 89, 53, 80, 66,
    51) ,
$s_{18}$=(1, 87, 46, 96, 20, 51, 78, 32, 42, 89, 74, 86) (2, 60, 6,
38, 10, 98, 5, 27, 72, 100, 73,
    53) (3, 52, 76, 25, 90, 66) (4, 95, 40, 28, 11, 93) (7, 8, 26, 97, 30, 80, 70, 77, 65, 62, 85, 21) (9, 22,
    37, 99, 36, 79, 67, 75, 43, 41, 92, 13) (12, 88, 29, 31, 48, 69, 82, 71, 63, 59, 47, 16) (14, 49, 44, 34,
    94, 45) (15, 17, 91, 23, 33, 58, 68, 83, 39, 64, 55, 54) (18, 61, 81, 24) (19, 50, 84, 35) (56, 57) ,
$s_{19}$=(1, 46, 20, 78, 42, 74) (2, 6, 10, 5, 72, 73) (3, 76, 90)
(4, 40, 11) (7, 26, 30, 70, 65, 85) (8, 97, 80, 77,
    62, 21) (9, 37, 36, 67, 43, 92) (12, 29, 48, 82, 63, 47) (13, 22, 99, 79, 75, 41) (14, 44, 94) (15, 91,
    33, 68, 39, 55) (16, 88, 31, 69, 71, 59) (17, 23, 58, 83, 64, 54) (18, 81) (19, 84) (24, 61) (25, 66,
    52) (27, 100, 53, 60, 38, 98) (28, 93, 95) (32, 89, 86, 87, 96, 51) (34, 45, 49) (35, 50) ,
$s_{20}$=(1, 53, 42, 27, 20, 38) (2, 96, 72, 86, 10, 32) (3, 94, 76,
14, 90, 44) (4, 95, 40, 28, 11, 93) (5, 89, 6,
    51, 73, 87) (7, 36, 65, 9, 30, 43) (8, 13, 62, 75, 80, 99) (12, 71, 63, 31, 48, 16) (15, 17, 39, 64, 33,
    58) (18, 19) (21, 41, 77, 79, 97, 22) (23, 55, 54, 68, 83, 91) (24, 50) (25, 49, 66, 34, 52, 45) (26,
    67, 85, 37, 70, 92) (29, 59, 47, 69, 82, 88) (35, 61) (46, 60, 74, 100, 78, 98) (56, 57) (81, 84) ,
$s_{21}$=(1, 5, 22, 79, 61, 90, 59, 33) (2, 88, 65, 69, 8, 13, 67,
12) (3, 49) (4, 93, 89, 92, 62, 39, 15, 44) (6, 75,
    27, 16, 84, 95, 53, 56) (7, 58, 23, 11, 51, 20, 54, 72) (9, 71, 83, 87, 64, 24, 32, 14) (10, 40, 30,
    48) (17, 86, 34, 77) (18, 45, 97, 41, 80, 38, 73, 37) (21, 94, 63, 42, 76, 99, 28, 36) (25, 47, 29, 35,
    68, 91, 85, 100) (26, 70, 98, 55, 74, 52, 57, 78) (43, 96) (46, 81, 50, 66) (60, 82) ,
$s_{22}$=(1, 42, 17, 71, 24) (2, 96, 70, 54, 25) (3, 66, 47, 41, 83)
(4, 5, 26, 46, 6) (7, 68, 90, 60, 73) (8, 44, 99,
    22, 97) (9, 36, 49, 34, 20) (10, 91, 84, 23, 78) (11, 12, 82, 13, 35) (14, 79, 98, 55, 92) (15, 18, 93,
    32, 65) (19, 27, 58, 38, 43) (21, 77, 69, 57, 85) (28, 48, 87, 89, 86) (29, 40, 37, 72, 88) (30, 52, 62,
    39, 75) (31, 61, 56, 64, 45) (33, 50, 94, 95, 67) (51, 100, 59, 63, 80) ,
$s_{23}$=(1, 44, 84, 95, 20, 57, 68, 28, 92, 90, 52) (2, 26, 43, 23,
94, 3, 50, 16, 100, 21, 33) (4, 99, 83, 30, 15,
    69, 93, 24, 39, 5, 75) (6, 72, 22, 34, 89, 98, 53, 66, 47, 77, 40) (7, 42, 76, 37, 19, 88, 17, 27, 45, 96,
    80) (8, 49, 71, 29, 14, 79, 85, 25, 48, 54, 62) (9, 64, 91, 31, 63, 10, 87, 61, 38, 56, 51) (11, 65, 46,
    32, 67, 13, 82, 70, 35, 73, 81) (12, 58, 60, 36, 86, 74, 18, 78, 41, 97, 59) ,
$s_{24}$=(1, 84, 20, 68, 92, 52, 44, 95, 57, 28, 90) (2, 43, 94, 50,
100, 33, 26, 23, 3, 16, 21) (4, 83, 15, 93, 39,
    75, 99, 30, 69, 24, 5) (6, 22, 89, 53, 47, 40, 72, 34, 98, 66, 77) (7, 76, 19, 17, 45, 80, 42, 37, 88, 27,
    96) (8, 71, 14, 85, 48, 62, 49, 29, 79, 25, 54) (9, 91, 63, 87, 38, 51, 64, 31, 10, 61, 56) (11, 46, 67,
    82, 35, 81, 65, 32, 13, 70, 73) (12, 60, 86, 18, 41, 59, 58, 36, 74, 78,
    97).

 Obviously,  $s_1$ is the unity element
and $G= G^{s_1}$.

The character table of $G^{s_1} =G$:\\
\begin{tabular}{c|cccccccccccccccccccccc}
  & & & & & & & & & &10 & & & & & & & & & &20 & &\\\hline
$\chi_1^{(1)}$&1&1&1&1&1&1&1&1&1&1&1&1&1&1&1&1&1&1&1&1&1&1
\\$\chi_1^{(2)}$&22&.&2&6&.&2&-2&-2&2&-1&4&-1&-1&1&-3&-6&1&.&.&-2&.&2
\\$\chi_1^{(3)}$&77&-1&1&13&-1&5&1&1&-3&.&5&.&.&-2&2&5&.&-1&1&1&1&2
\\$\chi_1^{(4)}$&154&.&-2&10&.&6&10&.&4&1&1&-2&-2&.&4&-2&.&1&1&1&.&-1
\\$\chi_1^{(5)}$&154&-2&2&10&2&-2&-10&.&4&1&1&.&.&.&4&-10&.&-1&1&-1&.&-1
\\$\chi_1^{(6)}$&154&2&2&10&-2&-2&-10&.&4&1&1&.&.&.&4&-10&.&-1&1&-1&.&-1
\\$\chi_1^{(7)}$&175&1&3&15&1&-1&11&1&5&-1&4&.&.&.&.&15&.&.&.&2&-1&.
\\$\chi_1^{(8)}$&231&-1&-1&7&-1&-1&-9&1&1&1&6&.&.&2&6&15&.&.&-2&.&-1&1
\\$\chi_1^{(9)}$&693&-1&1&21&-1&5&9&-1&3&.&.&1&1&1&-7&21&.&.&.&.&1&-2
\\$\chi_1^{(10)}$&770&.&-2&34&.&2&-10&.&.&.&5&1&1&-1&-5&-14&.&1&1&-1&-2&.
\\$\chi_1^{(11)}$&770&.&-2&-14&.&-2&10&.&.&.&5&A&-A&1&-5&-10&.&-1&1&1&.&.
\\$\chi_1^{(12)}$&770&.&-2&-14&.&-2&10&.&.&.&5&-A&A&1&-5&-10&.&-1&1&1&.&.
\\$\chi_1^{(13)}$&825&1&1&25&1&1&9&-1&-5&1&6&.&.&.&.&-15&-1&.&-2&.&1&.
\\$\chi_1^{(14)}$&896&.&.&.&.&.&16&1&1&1&-4&.&.&.&-4&.&.&.&.&-2&.&1
\\$\chi_1^{(15)}$&896&.&.&.&.&.&16&1&1&1&-4&.&.&.&-4&.&.&.&.&-2&.&1
\\$\chi_1^{(16)}$&1056&.&.&32&.&.&.&.&-4&-1&-6&.&.&2&6&.&-1&.&2&.&.&1
\\$\chi_1^{(17)}$&1386&.&-2&-6&.&-2&18&-2&6&.&.&1&1&-1&11&6&.&.&.&.&.&1
\\$\chi_1^{(18)}$&1408&.&.&.&.&.&16&1&-7&-1&4&.&.&.&8&.&1&.&.&-2&.&-2
\\$\chi_1^{(19)}$&1750&.&2&-10&.&6&10&.&.&.&-5&.&.&.&.&-10&.&-1&-1&1&-2&.
\\$\chi_1^{(20)}$&1925&-1&1&5&-1&-3&1&1&5&-1&-1&.&.&.&.&-35&.&1&-1&1&1&.
\\$\chi_1^{(21)}$&1925&1&-3&5&1&5&-19&1&5&-1&-1&.&.&.&.&5&.&-1&-1&-1&1&.
\\$\chi_1^{(22)}$&2520&.&.&24&.&-8&.&.&.&.&.&-1&-1&-1&-5&24&.&.&.&.&.&.
\\$\chi_1^{(23)}$&2750&.&2&-50&.&2&-10&.&.&.&5&.&.&.&.&10&-1&1&1&-1&.&.
\\$\chi_1^{(24)}$&3200&.&.&.&.&.&-16&-1&-5&1&-4&.&.&.&.&.&1&.&.&2&.&.
\end{tabular}

\begin{tabular}{c|cc}
  & &\\\hline
$\chi_1^{(1)}$&1&1
\\$\chi_1^{(2)}$&.&.
\\$\chi_1^{(3)}$&.&.
\\$\chi_1^{(4)}$&.&.
\\$\chi_1^{(5)}$&.&.
\\$\chi_1^{(6)}$&.&.
\\$\chi_1^{(7)}$&-1&-1
\\$\chi_1^{(8)}$&.&.
\\$\chi_1^{(9)}$&.&.
\\$\chi_1^{(10)}$&.&.
\\$\chi_1^{(11)}$&.&.
\\$\chi_1^{(12)}$&.&.
\\$\chi_1^{(13)}$&.&.
\\$\chi_1^{(14)}$&B&/B
\\$\chi_1^{(15)}$&/B&B
\\$\chi_1^{(16)}$&.&.
\\$\chi_1^{(17)}$&.&.
\\$\chi_1^{(18)}$&.&.
\\$\chi_1^{(19)}$&1&1
\\$\chi_1^{(20)}$&.&.
\\$\chi_1^{(21)}$&.&.
\\$\chi_1^{(22)}$&1&1
\\$\chi_1^{(23)}$&.&.
\\$\chi_1^{(24)}$&-1&-1
\end{tabular}

\noindent \noindent where   A =
-E(20)-E(20)$^9$+E(20)$^{13}$+E(20)$^{17}$
  = -ER(-5) = -i5;
B = E(11)$^2$+E(11)$^6$+E(11)$^7$+E(11)$^8$+E(11)$^{10}$
  = (-1-ER(-11))/2 = -1-b11.

The generators of $G^{s_2}$ are:\\
(  1, 33, 64, 44, 42, 63,  9, 78)(  2, 15,  7, 14, 47, 57, 18, 87)
    (  3,  8, 69, 39, 96, 91, 68, 52)(  4, 37, 74, 29, 26, 72, 62, 88)
    (  5, 25, 51, 82, 76, 12, 80, 60)(  6, 10, 50, 75)( 11, 59, 13, 99, 41, 94, 85, 22)( 16, 48, 53, 86)( 17, 81)( 19, 92, 45, 35, 31, 55, 58, 65)
    ( 20, 28, 49, 73, 38, 79, 43, 54)( 21, 98, 70, 32, 90, 89, 93, 83)
    ( 23, 36,100, 71, 97, 61, 67, 24)( 27, 77, 56, 66)( 30, 95, 84, 46)( 34, 40),

  (  1,  2, 64,  7, 42, 47,  9, 18)(  3, 65, 69, 92, 96, 35, 68, 55)
    (  4, 41, 74, 85, 26, 11, 62, 13)(  5, 43, 51, 20, 76, 49, 80, 38)
    (  6, 75, 50, 10)(  8, 19, 39, 45, 91, 31, 52, 58)( 12, 73, 60, 79, 25, 54, 82,  28)( 14, 63, 57, 78, 87, 33, 15, 44)( 16, 77, 53, 66)( 17, 40)
    ( 21, 71, 70, 61, 90, 24, 93, 36)( 22, 72, 59, 88, 99, 37, 94, 29)
    ( 23, 83,100, 98, 97, 32, 67, 89)( 27, 48, 56, 86)( 30, 46, 84, 95)( 34,
    81).

The representatives of conjugacy classes of   $G^{s_2}$ are:\\
 (1),
 (1,2,64,7,42,47,9,18)(3,65,69,92,96,35,68,55)(4,41,74,85,26,11,62,13)(5,43,
    51,20,76,49,80,38)(6,75,50,10)(8,\\19,39,45,91,31,52,58)(12,73,60,79,25,54,82,
    28)(14,63,57,78,87,33,15,44)(16,77,53,66)(17,40)(21,71,70,61,90,24,\\93,36)(22,72,
    59,88,99,37,94,29)(23,83,100,98,97,32,67,89)(27,48,56,86)(30,46,84,95)(34,81),

  (1,7,9,2,42,18,64,47)(3,92,68,65,96,55,69,35)(4,85,62,41,26,13,74,11)(5,20,80,43,
    76,38,51,49)(6,10,50,75)(8,\\45,52,19,91,58,39,31)(12,79,82,73,25,28,60,54)(14,78,
    15,63,87,44,57,33)(16,66,53,77)(17,40)(21,61,93,71,90,36,\\70,24)(22,88,94,72,99,
    29,59,37)(23,98,67,83,97,89,100,32)(27,86,56,48)(30,95,84,46)(34,81),

  (1,9,42,64)(2,18,47,7)(3,68,96,69)(4,62,26,74)(5,80,76,51)(6,50)(8,52,91,39)(10,
    75)(11,85,41,13)(12,82,25,\\60)(14,15,87,57)(16,53)(19,58,31,45)(20,43,38,49)(21,
    93,90,70)(22,94,99,59)(23,67,97,100)(24,61,71,36)(27,\\56)(28,54,79,73)(29,37,88,
    72)(30,84)(32,98,83,89)(33,78,63,44)(35,92,65,55)(46,95)(48,86)(66,77),

  (1,14)(2,63)(3,45)(4,22)(5,28)(6,50)(7,78)(8,35)(9,15)(10,75)(11,37)(12,43)(13,
    29)(16,27)(17,34)(18,44)\\(19,68)(20,60)(21,67)(23,70)(24,98)(25,49)(26,99)(30,
    84)(31,69)(32,36)(33,47)(38,82)(39,55)(40,81)(41,72)\\(42,87)(46,95)(48,77)(51,
    73)(52,92)(53,56)(54,80)(57,64)(58,96)(59,74)(61,83)(62,94)(65,91)(66,86)(71,
    89)\\(76,79)(85,88)(90,100)(93,97),

     (1,15,42,57)(2,44,47,78)(3,19,96,31)(4,94,
    26,59)(5,54,76,73)(7,63,18,33)(8,92,91,55)(9,87,64,14)(11,88,41,\\29)(12,38,25,
    20)(13,37,85,72)(16,56)(17,34)(21,97,90,23)(22,62,99,74)(24,83,71,32)(27,53)(28,
    80,79,51)(35,\\52,65,39)(36,98,61,89)(40,81)(43,82,49,60)(45,68,58,69)(48,66)(67,
    93,100,70)(77,86),

    (1,18,9,47,42,7,64,2)(3,55,68,35,96,92,69,65)(4,13,62,11,
    26,85,74,41)(5,38,80,49,76,20,51,43)(6,10,50,75)\\(8,58,52,31,91,45,39,19)(12,28,
    82,54,25,79,60,73)(14,44,15,33,87,78,57,63)(16,66,53,77)(17,40)(21,36,93,24,90,
    61,70,71)(22,29,94,37,99,88,59,72)(23,89,67,32,97,98,100,83)(27,86,56,48)(30,95,
    84,46)(34,81),

     (1,33,64,44,42,63,9,78)(2,15,7,14,47,57,18,87)(3,8,69,39,96,91,
    68,52)(4,37,74,29,26,72,62,88)(5,25,51,82,76,\\12,80,60)(6,10,50,75)(11,59,13,99,
    41,94,85,22)(16,48,53,86)(17,81)(19,92,45,35,31,55,58,65)(20,28,49,73,38,79,\\43,
    54)(21,98,70,32,90,89,93,83)(23,36,100,71,97,61,67,24)(27,77,56,66)(30,95,84,
    46)(34,40),

     (1,42)(2,47)(3,96)(4,26)(5,76)(7,18)(8,91)(9,64)(11,41)(12,25)(13,
    85)(14,87)(15,57)(19,31)(20,38)(21,90)\\(22,99)(23,97)(24,71)(28,79)(29,88)(32,
    83)(33,63)(35,65)(36,61)(37,72)(39,52)(43,49)(44,78)(45,58)(51,80)(54,\\73)(55,
    92)(59,94)(60,82)(62,74)(67,100)(68,69)(70,93)(89,98),

  (1,44,9,33,42,78,64,63)(2,14,18,15,47,87,7,57)(3,39,68,8,96,52,69,91)(4,29,62,37,
    26,88,74,72)(5,82,80,25,76,\\60,51,12)(6,75,50,10)(11,99,85,59,41,22,13,94)(16,86,
    53,48)(17,81)(19,35,58,92,31,65,45,55)(20,73,43,28,38,54,\\49,79)(21,32,93,98,90,
    83,70,89)(23,71,67,36,97,24,100,61)(27,66,56,77)(30,46,84,95)(34,40),

  (1,47,64,18,42,2,9,7)(3,35,69,55,96,65,68,92)(4,11,74,13,26,41,62,85)(5,49,51,38,
    76,43,80,20)(6,75,50,10)\\(8,31,39,58,91,19,52,45)(12,54,60,28,25,73,82,79)(14,33,
    57,44,87,63,15,78)(16,77,53,66)(17,40)(21,24,70,36,90,\\71,93,61)(22,37,59,29,99,
    72,94,88)(23,32,100,89,97,83,67,98)(27,48,56,86)(30,46,84,95)(34,81),

  (1,57,42,15)(2,78,47,44)(3,31,96,19)(4,59,26,94)(5,73,76,54)(7,33,18,63)(8,55,91,
    92)(9,14,64,87)(11,29,41,\\88)(12,20,25,38)(13,72,85,37)(16,56)(17,34)(21,23,90,
    97)(22,74,99,62)(24,32,71,83)(27,53)(28,51,79,80)(35,39,\\65,52)(36,89,61,98)(40,
    81)(43,60,49,82)(45,69,58,68)(48,66)(67,70,100,93)(77,86),

  (1,63,64,78,42,33,9,44)(2,57,7,87,47,15,18,14)(3,91,69,52,96,8,68,39)(4,72,74,88,
    26,37,62,29)(5,12,51,60,76,\\25,80,82)(6,10,50,75)(11,94,13,22,41,59,85,99)(16,48,
    53,86)(17,81)(19,55,45,65,31,92,58,35)(20,79,49,54,38,28,\\43,73)(21,89,70,83,90,
    98,93,32)(23,61,100,24,97,36,67,71)(27,77,56,66)(30,95,84,46)(34,40),

  (1,64,42,9)(2,7,47,18)(3,69,96,68)(4,74,26,62)(5,51,76,80)(6,50)(8,39,91,52)(10,
    75)(11,13,41,85)(12,60,25,\\82)(14,57,87,15)(16,53)(19,45,31,58)(20,49,38,43)(21,
    70,90,93)(22,59,99,94)(23,100,97,67)(24,36,71,61)(27,56)\\(28,73,79,54)(29,72,88,
    37)(30,84)(32,89,83,98)(33,44,63,78)(35,55,65,92)(46,95)(48,86)(66,77),

  (1,78,9,63,42,44,64,33)(2,87,18,57,47,14,7,15)(3,52,68,91,96,39,69,8)(4,88,62,72,
    26,29,74,37)(5,60,80,12,76,\\82,51,25)(6,75,50,10)(11,22,85,94,41,99,13,59)(16,86,
    53,48)(17,81)(19,65,58,55,31,35,45,92)(20,54,43,79,38,73,\\49,28)(21,83,93,89,90,
    32,70,98)(23,24,67,61,97,71,100,36)(27,66,56,77)(30,46,84,95)(34,40),

  (1,87)(2,33)(3,58)(4,99)(5,79)(6,50)(7,44)(8,65)(9,57)(10,75)(11,72)(12,49)(13,
    88)(14,42)(15,64)(16,27)(17,\\34)(18,78)(19,69)(20,82)(21,100)(22,26)(23,93)(24,
    89)(25,43)(28,76)(29,85)(30,84)(31,68)(32,61)(35,91)(36,83)\\(37,41)(38,60)(39,
    92)(40,81)(45,96)(46,95)(47,63)(48,77)(51,54)(52,55)(53,56)(59,62)(66,86)(67,
    90)(70,97)(71,\\98)(73,80)(74,94).

The character table of $G^{s_2}$:\\
\begin{tabular}{c|cccccccccccccccc}
  & & & & & & & & & &10 & & & & & &\\\hline
$\chi_2^{(1)}$&1&1&1&1&1&1&1&1&1&1&1&1&1&1&1&1
\\$\chi_2^{(2)}$&1&-1&-1&1&-1&-1&-1&1&1&1&-1&-1&1&1&1&-1
\\$\chi_2^{(3)}$&1&-1&-1&1&1&1&-1&-1&1&-1&-1&1&-1&1&-1&1
\\$\chi_2^{(4)}$&1&1&1&1&-1&-1&1&-1&1&-1&1&-1&-1&1&-1&-1
\\$\chi_2^{(5)}$&1&A&-A&-1&-1&1&-A&-A&1&A&A&1&-A&-1&A&-1
\\$\chi_2^{(6)}$&1&-A&A&-1&-1&1&A&A&1&-A&-A&1&A&-1&-A&-1
\\$\chi_2^{(7)}$&1&A&-A&-1&1&-1&-A&A&1&-A&A&-1&A&-1&-A&1
\\$\chi_2^{(8)}$&1&-A&A&-1&1&-1&A&-A&1&A&-A&-1&-A&-1&A&1
\\$\chi_2^{(9)}$&1&B&-/B&A&-1&-A&/B&B&-1&-/B&-B&A&-B&-A&/B&1
\\$\chi_2^{(10)}$&1&-/B&B&-A&-1&A&-B&-/B&-1&B&/B&-A&/B&A&-B&1
\\$\chi_2^{(11)}$&1&/B&-B&-A&-1&A&B&/B&-1&-B&-/B&-A&-/B&A&B&1
\\$\chi_2^{(12)}$&1&-B&/B&A&-1&-A&-/B&-B&-1&/B&B&A&B&-A&-/B&1
\\$\chi_2^{(13)}$&1&B&-/B&A&1&A&/B&-B&-1&/B&-B&-A&B&-A&-/B&-1
\\$\chi_2^{(14)}$&1&-/B&B&-A&1&-A&-B&/B&-1&-B&/B&A&-/B&A&B&-1
\\$\chi_2^{(15)}$&1&/B&-B&-A&1&-A&B&-/B&-1&B&-/B&A&/B&A&-B&-1
\\$\chi_2^{(16)}$&1&-B&/B&A&1&A&-/B&B&-1&-/B&B&-A&-B&-A&/B&-1
\end{tabular}

\noindent \noindent where   A = -E(4)
  = -ER(-1) = -i;
B = -E(8).

The generators of $G^{s_3}$ are:\\
 (  1, 64, 42,  9)(  2,  7, 47, 18)(  3, 69, 96, 68)(  4, 74, 26, 62)
    (  5, 51, 76, 80)(  6, 50)(  8, 39, 91, 52)( 10, 75)( 11, 13, 41, 85)
    ( 12, 60, 25, 82)( 14, 57, 87, 15)( 16, 53)( 19, 45, 31, 58)( 20, 49, 38, 43)
    ( 21, 70, 90, 93)( 22, 59, 99, 94)( 23,100, 97, 67)( 24, 36, 71, 61)( 27, 56)
    ( 28, 73, 79, 54)( 29, 72, 88, 37)( 30, 84)( 32, 89, 83, 98)( 33, 44, 63, 78)
    ( 35, 55, 65, 92)( 46, 95)( 48, 86)( 66, 77), (  1,  2, 64,  7, 42, 47,  9, 18)
    (  3, 65, 69, 92, 96, 35, 68, 55)(  4, 41, 74, 85, 26, 11, 62, 13)
    (  5, 43, 51, 20, 76, 49, 80, 38)(  6, 75, 50, 10)(  8, 19, 39, 45, 91, 31, 52,
      58)( 12, 73, 60, 79, 25, 54, 82, 28)( 14, 63, 57, 78, 87, 33, 15, 44)
    ( 16, 77, 53, 66)( 17, 40)( 21, 71, 70, 61, 90, 24, 93, 36)( 22, 72, 59, 88, 99,
     37, 94, 29)( 23, 83,100, 98, 97, 32, 67, 89)( 27, 48, 56, 86)( 30, 46, 84, 95)
    ( 34, 81), (  1,  3, 61, 60)(  2, 54, 21, 55)(  4, 26)(  5, 97, 91, 44)
    (  6, 10, 84, 46)(  7, 28, 70, 65)(  8, 78, 76, 23)(  9, 68, 71, 12)
    ( 11, 85, 41, 13)( 14, 58, 83, 38)( 15, 31, 89, 49)( 17, 34)( 18, 79, 93, 35)
    ( 19, 98, 43, 57)( 20, 87, 45, 32)( 24, 25, 64, 69)( 29, 72, 88, 37)
    ( 30, 95, 50, 75)( 33, 80,100, 39)( 36, 82, 42, 96)( 40, 81)( 47, 73, 90, 92)
    ( 48, 86)( 51, 67, 52, 63)( 62, 74)( 66, 77), (  1,  5, 42, 76)(  2, 58, 47, 45)
    (  3, 33, 96, 63)(  4, 72, 26, 37)(  6, 30)(  7, 19, 18, 31)(  8, 36, 91, 61)
    (  9, 80, 64, 51)( 11, 99, 41, 22)( 12, 97, 25, 23)( 13, 94, 85, 59)
    ( 14, 79, 87, 28)( 15, 73, 57, 54)( 16, 48)( 17, 40)( 20, 90, 38, 21)
    ( 24, 39, 71, 52)( 27, 77)( 29, 74, 88, 62)( 32, 35, 83, 65)( 34, 81)
    ( 43, 70, 49, 93)( 44, 68, 78, 69)( 50, 84)( 53, 86)( 55, 98, 92, 89)( 56, 66)
    ( 60, 67, 82,100).

The representatives of conjugacy classes of   $G^{s_3}$ are:\\
 (1),
  (1,2,64,7,42,47,9,18)(3,65,69,92,96,35,68,55)(4,41,74,85,26,11,62,13)(5,43,
    51,20,76,49,80,38)(6,75,50,10)(8,19,\\39,45,91,31,52,58)(12,73,60,79,25,54,82,
    28)(14,63,57,78,87,33,15,44)(16,77,53,66)(17,40)(21,71,70,61,90,24,93,\\36)(22,72,
    59,88,99,37,94,29)(23,83,100,98,97,32,67,89)(27,48,56,86)(30,46,84,95)(34,81),

  (1,3,61,60)(2,54,21,55)(4,26)(5,97,91,44)(6,10,84,46)(7,28,70,65)(8,78,76,23)(9,
    68,71,12)(11,85,41,13)(14,58,\\83,38)(15,31,89,49)(17,34)(18,79,93,35)(19,98,43,
    57)(20,87,45,32)(24,25,64,69)(29,72,88,37)(30,95,50,75)(33,80,\\100,39)(36,82,42,
    96)(40,81)(47,73,90,92)(48,86)(51,67,52,63)(62,74)(66,77),

  (1,5,42,76)(2,58,47,45)(3,33,96,63)(4,72,26,37)(6,30)(7,19,18,31)(8,36,91,61)(9,
    80,64,51)(11,99,41,22)(12,97,\\25,23)(13,94,85,59)(14,79,87,28)(15,73,57,54)(16,
    48)(17,40)(20,90,38,21)(24,39,71,52)(27,77)(29,74,88,62)(32,35,\\83,65)(34,81)(43,
    70,49,93)(44,68,78,69)(50,84)(53,86)(55,98,92,89)(56,66)(60,67,82,100),

  (1,7,9,2,42,18,64,47)(3,92,68,65,96,55,69,35)(4,85,62,41,26,13,74,11)(5,20,80,43,
    76,38,51,49)(6,10,50,75)(8,45,\\52,19,91,58,39,31)(12,79,82,73,25,28,60,54)(14,78,
    15,63,87,44,57,33)(16,66,53,77)(17,40)(21,61,93,71,90,36,70,24)\\(22,88,94,72,99,
    29,59,37)(23,98,67,83,97,89,100,32)(27,86,56,48)(30,95,84,46)(34,81),

  (1,8)(2,20)(3,67)(4,72)(5,61)(6,50)(7,49)(9,52)(10,46)(11,22)(12,78)(13,59)(14,
    65)(15,55)(16,48)(17,40)(18,43)\\(19,70)(21,58)(23,69)(24,51)(25,44)(26,37)(27,
    77)(28,32)(29,62)(30,84)(31,93)(33,60)(34,81)(35,87)(36,76)(38,47)\\(39,64)(41,
    99)(42,91)(45,90)(53,86)(54,98)(56,66)(57,92)(63,82)(68,97)(71,80)(73,89)(74,
    88)(75,95)(79,83)(85,94)\\(96,100),

     (1,9,42,64)(2,18,47,7)(3,68,96,69)(4,62,26,
    74)(5,80,76,51)(6,50)(8,52,91,39)(10,75)(11,85,41,13)(12,82,25,60)\\(14,15,87,
    57)(16,53)(19,58,31,45)(20,43,38,49)(21,93,90,70)(22,94,99,59)(23,67,97,100)(24,
    61,71,36)(27,56)(28,54,\\79,73)(29,37,88,72)(30,84)(32,98,83,89)(33,78,63,44)(35,
    92,65,55)(46,95)(48,86)(66,77),

    (1,14)(2,63)(3,45)(4,22)(5,28)(6,50)(7,78)(8,
    35)(9,15)(10,75)(11,37)(12,43)(13,29)(16,27)(17,34)(18,44)(19,68)\\(20,60)(21,
    67)(23,70)(24,98)(25,49)(26,99)(30,84)(31,69)(32,36)(33,47)(38,82)(39,55)(40,
    81)(41,72)(42,87)(46,95)\\(48,77)(51,73)(52,92)(53,56)(54,80)(57,64)(58,96)(59,
    74)(61,83)(62,94)(65,91)(66,86)(71,89)(76,79)(85,88)(90,100)\\(93,97),

  (1,15,42,57)(2,44,47,78)(3,19,96,31)(4,94,26,59)(5,54,76,73)(7,63,18,33)(8,92,91,
    55)(9,87,64,14)(11,88,41,29)\\(12,38,25,20)(13,37,85,72)(16,56)(17,34)(21,97,90,
    23)(22,62,99,74)(24,83,71,32)(27,53)(28,80,79,51)(35,52,65,39)\\(36,98,61,89)(40,
    81)(43,82,49,60)(45,68,58,69)(48,66)(67,93,100,70)(77,86),

  (1,19,61,43)(2,76,21,8)(3,89,60,15)(4,59)(5,90,91,47)(6,10,84,46)(7,80,70,39)(9,
    58,71,38)(11,72,41,37)(12,87,68,\\32)(13,88,85,29)(14,69,83,25)(16,56)(18,51,93,
    52)(20,64,45,24)(22,62)(23,92,78,73)(26,94)(27,53)(28,67,65,63)(30,\\95,50,75)(31,
    36,49,42)(33,79,100,35)(44,54,97,55)(48,77)(57,96,98,82)(66,86)(74,99),

  (1,23,64,100,42,97,9,67)(2,32,7,89,47,83,18,98)(3,51,69,76,96,80,68,5)(4,29,74,72,
    26,88,62,37)(6,95,50,46)(8,60,\\39,25,91,82,52,12)(10,84,75,30)(11,22,13,59,41,99,
    85,94)(14,70,57,90,87,93,15,21)(16,86,53,48)(17,81)(19,28,45,73,\\31,79,58,54)(20,
    55,49,65,38,92,43,35)(24,63,36,78,71,33,61,44)(27,66,56,77)(34,40),

  (1,24,42,71)(2,70,47,93)(3,25,96,12)(4,74,26,62)(5,52,76,39)(6,30)(7,90,18,21)(8,
    80,91,51)(9,61,64,36)(10,95)\\(11,85,41,13)(14,98,87,89)(15,83,57,32)(16,53)(19,
    20,31,38)(22,59,99,94)(23,33,97,63)(27,56)(28,92,79,55)(29,37,\\88,72)(35,54,65,
    73)(43,45,49,58)(44,67,78,100)(46,75)(48,86)(50,84)(60,69,82,68)(66,77),

  (1,28)(2,96)(3,47)(4,41)(5,87)(6,84)(7,68)(8,32)(9,54)(10,75)(11,26)(12,93)(13,
    62)(14,76)(15,51)(16,77)(17,81)\\(18,69)(19,44)(20,100)(21,60)(22,37)(23,43)(24,
    55)(25,70)(27,48)(29,59)(30,50)(31,78)(33,58)(34,40)(35,61)(36,65)\\(38,67)(39,
    89)(42,79)(45,63)(46,95)(49,97)(52,98)(53,66)(56,86)(57,80)(64,73)(71,92)(72,
    99)(74,85)(82,90)(83,91)\\(88,94),

    (1,32)(2,100)(3,38)(4,99)(5,35)(6,30)(7,
    97)(8,28)(9,98)(10,95)(11,37)(12,31)(13,29)(14,36)(15,24)(16,27)(17,34)\\(18,
    23)(19,25)(20,96)(21,33)(22,26)(39,73)(40,81)(41,72)(42,83)(43,69)(44,70)(45,
    82)(46,75)(47,67)(48,77)(49,68)\\(50,84)(51,55)(52,54)(53,56)(57,71)(58,60)(59,
    62)(61,87)(63,90)(64,89)(65,76)(66,86)(74,94)(78,93)(79,91)(80,92)\\(85,88),

  (1,35,42,65)(2,60,47,82)(3,21,96,90)(4,41,26,11)(5,83,76,32)(7,25,18,12)(8,14,91,
    87)(9,92,64,55)(10,95)(13,74,85,\\62)(15,39,57,52)(16,77)(17,81)(19,23,31,97)(20,
    63,38,33)(22,37,99,72)(24,73,71,54)(27,48)(28,36,79,61)(29,94,88,59)\\(34,40)(43,
    44,49,78)(45,100,58,67)(46,75)(51,98,80,89)(53,66)(56,86)(68,93,69,70),

  (1,36)(2,90)(3,82)(4,26)(5,8)(6,84)(7,93)(9,24)(10,46)(12,69)(14,32)(15,98)(18,
    70)(19,49)(20,58)(21,47)(22,99)\\(23,44)(25,68)(28,35)(30,50)(31,43)(33,67)(38,
    45)(39,51)(42,61)(52,80)(54,92)(55,73)(57,89)(59,94)(60,96)(62,74)\\(63,100)(64,
    71)(65,79)(75,95)(76,91)(78,97)(83,87),

     (1,42)(2,47)(3,96)(4,26)(5,76)(7,
    18)(8,91)(9,64)(11,41)(12,25)(13,85)(14,87)(15,57)(19,31)(20,38)(21,90)(22,
    99)(23,97)(24,71)(28,79)(29,88)(32,83)(33,63)(35,65)(36,61)(37,72)(39,52)(43,
    49)(44,78)(45,58)(51,80)(54,73)(55,\\92)(59,94)(60,82)(62,74)(67,100)(68,69)(70,
    93)(89,98),

    (1,43,61,19)(2,8,21,76)(3,15,60,89)(4,59)(5,47,91,90)(6,46,84,
    10)(7,39,70,80)(9,38,71,58)(11,37,41,72)(12,32,\\68,87)(13,29,85,88)(14,25,83,
    69)(16,56)(18,52,93,51)(20,24,45,64)(22,62)(23,73,78,92)(26,94)(27,53)(28,63,65,
    67)(30,75,50,95)(31,42,49,36)(33,35,100,79)(44,55,97,54)(48,77)(57,82,98,96)(66,
    86)(74,99),

    (1,44,9,33,42,78,64,63)(2,14,18,15,47,87,7,57)(3,39,68,8,96,52,69,
    91)(4,29,62,37,26,88,74,72)(5,82,80,25,76,60,\\51,12)(6,75,50,10)(11,99,85,59,41,
    22,13,94)(16,86,53,48)(17,81)(19,35,58,92,31,65,45,55)(20,73,43,28,38,54,49,
    79)(21,32,93,98,90,83,70,89)(23,71,67,36,97,24,100,61)(27,66,56,77)(30,46,84,
    95)(34,40),

     (1,60,61,3)(2,55,21,54)(4,26)(5,44,91,97)(6,46,84,10)(7,65,70,
    28)(8,23,76,78)(9,12,71,68)(11,13,41,85)(14,38,83,\\58)(15,49,89,31)(17,34)(18,35,
    93,79)(19,57,43,98)(20,32,45,87)(24,69,64,25)(29,37,88,72)(30,75,50,95)(33,39,
    100,80)(36,96,42,82)(40,81)(47,92,90,73)(48,86)(51,63,52,67)(62,74)(66,77),

  (1,64,42,9)(2,7,47,18)(3,69,96,68)(4,74,26,62)(5,51,76,80)(6,50)(8,39,91,52)(10,
    75)(11,13,41,85)(12,60,25,82)\\(14,57,87,15)(16,53)(19,45,31,58)(20,49,38,43)(21,
    70,90,93)(22,59,99,94)(23,100,97,67)(24,36,71,61)(27,56)(28,\\73,79,54)(29,72,88,
    37)(30,84)(32,89,83,98)(33,44,63,78)(35,55,65,92)(46,95)(48,86)(66,77),

  (1,89,42,98)(2,97,47,23)(3,43,96,49)(4,94,26,59)(5,55,76,92)(6,84)(7,67,18,100)(8,
    73,91,54)(9,32,64,83)(10,46)\\(11,29,41,88)(12,58,25,45)(13,72,85,37)(14,71,87,
    24)(15,36,57,61)(16,56)(17,34)(19,82,31,60)(20,68,38,69)(21,44,\\90,78)(22,62,99,
    74)(27,53)(28,39,79,52)(30,50)(33,70,63,93)(35,51,65,80)(40,81)(48,66)(75,
    95)(77,86).

The character table of $G^{s_3}$:\\
\begin{tabular}{c|cccccccccccccccccccccc}
  & & & & & & & & & & 10& & & & & & & & & & &20 &\\\hline
$\chi_3^{(1)}$&1&1&1&1&1&1&1&1&1&1&1&1&1&1&1&1&1&1&1&1&1&1
\\$\chi_3^{(2)}$&1&-1&-1&-1&-1&-1&1&-1&-1&1&1&1&1&-1&1&1&1&1&1&-1&1&-1
\\$\chi_3^{(3)}$&1&-1&-1&1&-1&1&1&1&1&-1&-1&1&1&1&1&1&1&-1&-1&-1&1&1
\\$\chi_3^{(4)}$&1&-1&1&-1&-1&-1&1&1&1&1&-1&1&-1&1&-1&1&1&1&-1&1&1&1
\\$\chi_3^{(5)}$&1&-1&1&1&-1&1&1&-1&-1&-1&1&1&-1&-1&-1&1&1&-1&1&1&1&-1
\\$\chi_3^{(6)}$&1&1&-1&-1&1&-1&1&1&1&-1&1&1&-1&1&-1&1&1&-1&1&-1&1&1
\\$\chi_3^{(7)}$&1&1&-1&1&1&1&1&-1&-1&1&-1&1&-1&-1&-1&1&1&1&-1&-1&1&-1
\\$\chi_3^{(8)}$&1&1&1&-1&1&-1&1&-1&-1&-1&-1&1&1&-1&1&1&1&-1&-1&1&1&-1
\\$\chi_3^{(9)}$&1&A&A&-1&-A&1&-1&-1&1&A&-A&1&1&1&-1&-1&1&-A&A&-A&-1&-1
\\$\chi_3^{(10)}$&1&-A&-A&-1&A&1&-1&-1&1&-A&A&1&1&1&-1&-1&1&A&-A&A&-1&-1
\\$\chi_3^{(11)}$&1&A&A&1&-A&-1&-1&1&-1&-A&A&1&1&-1&-1&-1&1&A&-A&-A&-1&1
\\$\chi_3^{(12)}$&1&-A&-A&1&A&-1&-1&1&-1&A&-A&1&1&-1&-1&-1&1&-A&A&A&-1&1
\\$\chi_3^{(13)}$&1&A&-A&-1&-A&1&-1&1&-1&A&A&1&-1&-1&1&-1&1&-A&-A&A&-1&1
\\$\chi_3^{(14)}$&1&-A&A&-1&A&1&-1&1&-1&-A&-A&1&-1&-1&1&-1&1&A&A&-A&-1&1
\\$\chi_3^{(15)}$&1&A&-A&1&-A&-1&-1&-1&1&-A&-A&1&-1&1&1&-1&1&A&A&A&-1&-1
\\$\chi_3^{(16)}$&1&-A&A&1&A&-1&-1&-1&1&A&A&1&-1&1&1&-1&1&-A&-A&-A&-1&-1
\\$\chi_3^{(17)}$&2&.&.&.&.&.&-2&-2&2&.&.&-2&.&-2&.&2&2&.&.&.&-2&2
\\$\chi_3^{(18)}$&2&.&.&.&.&.&-2&2&-2&.&.&-2&.&2&.&2&2&.&.&.&-2&-2
\\$\chi_3^{(19)}$&2&.&.&.&.&.&2&-2&-2&.&.&-2&.&2&.&-2&2&.&.&.&2&2
\\$\chi_3^{(20)}$&2&.&.&.&.&.&2&2&2&.&.&-2&.&-2&.&-2&2&.&.&.&2&-2
\\$\chi_3^{(21)}$&4&.&.&.&.&.&B&.&.&.&.&.&.&.&.&.&-4&.&.&.&-B&.
\\$\chi_3^{(22)}$&4&.&.&.&.&.&-B&.&.&.&.&.&.&.&.&.&-4&.&.&.&B&.
\end{tabular}

\noindent \noindent where   A = -E(4)
  = -ER(-1) = -i;
B = -4*E(4)
  = -4*ER(-1) = -4i.

The generators of $G^{s_4}$ are:\\
 (  1, 42)(  2, 47)(  3, 96)(  4, 26)(  5, 76)(  7, 18)(  8, 91)(  9, 64)
    ( 11, 41)( 12, 25)( 13, 85)( 14, 87)( 15, 57)( 19, 31)( 20, 38)( 21, 90)
    ( 22, 99)( 23, 97)( 24, 71)( 28, 79)( 29, 88)( 32, 83)( 33, 63)( 35, 65)
    ( 36, 61)( 37, 72)( 39, 52)( 43, 49)( 44, 78)( 45, 58)( 51, 80)( 54, 73)
    ( 55, 92)( 59, 94)( 60, 82)( 62, 74)( 67,100)( 68, 69)( 70, 93)( 89, 98),
  (  3, 22, 67)(  4, 65, 71)(  5, 78, 54)(  6, 40, 77)(  7, 49, 74)(  8, 55, 21)
    (  9, 11, 52)( 10, 48, 56)( 12, 20, 61)( 13, 70, 60)( 14, 59, 51)( 15, 58, 68)
    ( 16, 46, 66)( 17, 50, 53)( 18, 43, 62)( 19, 72, 63)( 24, 26, 35)( 25, 38, 36)
    ( 27, 34, 30)( 28, 32, 29)( 31, 37, 33)( 39, 64, 41)( 44, 73, 76)( 45, 69, 57)
    ( 79, 83, 88)( 80, 87, 94)( 81, 86, 84)( 82, 85, 93)( 90, 91, 92)( 96, 99,100),
  (  2,  4, 63)(  3, 14, 60)(  5, 37, 52)(  6, 75, 77)(  7, 78, 74)(  8, 29, 51)
    ( 10, 50, 66)( 11, 35, 54)( 12, 15, 68)( 13, 28, 55)( 16, 34, 53)( 17, 27, 56)
    ( 18, 44, 62)( 19, 43, 71)( 20, 45, 61)( 21, 22, 67)( 23, 59, 70)( 24, 31, 49)
    ( 25, 57, 69)( 26, 33, 47)( 30, 48, 46)( 36, 38, 58)( 39, 76, 72)( 41, 65, 73)
    ( 79, 92, 85)( 80, 91, 88)( 82, 96, 87)( 84, 95, 86)( 90, 99,100)( 93, 97, 94),
  (  2,  5)(  3, 13)(  4,  9)(  7, 18)( 10, 16)( 11, 19)( 12, 20)( 14, 22)( 23, 55)
    ( 24, 62)( 25, 38)( 26, 64)( 27, 50)( 29, 60)( 30, 53)( 31, 41)( 32, 59)
    ( 33, 65)( 35, 63)( 36, 61)( 37, 43)( 39, 52)( 40, 77)( 44, 78)( 45, 68)
    ( 46, 56)( 47, 76)( 49, 72)( 51, 70)( 54, 73)( 58, 69)( 71, 74)( 80, 93)
    ( 81, 86)( 82, 88)( 83, 94)( 85, 96)( 87, 99)( 89, 98)( 92, 97),
  (  1,  2)(  4, 13)(  5, 15)(  6, 17)(  7, 18)(  8, 20)( 12, 21)( 14, 19)( 23, 89)
    ( 24, 82)( 25, 90)( 26, 85)( 27, 86)( 28, 79)( 29, 88)( 30, 81)( 31, 87)
    ( 32, 83)( 33, 80)( 34, 84)( 35, 93)( 36, 92)( 37, 94)( 38, 91)( 40, 50)
    ( 42, 47)( 43, 49)( 44, 45)( 51, 63)( 53, 77)( 54, 68)( 55, 61)( 57, 76)
    ( 58, 78)( 59, 72)( 60, 71)( 62, 74)( 65, 70)( 69, 73)( 97, 98).
 The representatives of conjugacy classes of   $G^{s_4}$ are:\\
 (1),
 (1,42)(2,47)(3,96)(4,26)(5,76)(7,18)(8,91)(9,64)(11,41)(12,25)(13,85)(14,
    87)(15,57)(19,31)(20,38)(21,90)(22,\\99)(23,97)(24,71)(28,79)(29,88)(32,83)(33,
    63)(35,65)(36,61)(37,72)(39,52)(43,49)(44,78)(45,58)(51,80)(54,73)\\(55,92)(59,
    94)(60,82)(62,74)(67,100)(68,69)(70,93)(89,98),

  (2,76)(3,60)(4,62)(5,47)(7,52)(8,21)(9,24)(10,46)(11,37)(12,25)(13,29)(14,32)(16,
    56)(18,39)(19,43)(20,38)(22,\\59)(23,55)(26,74)(27,53)(28,67)(30,50)(31,49)(33,
    35)(41,72)(44,73)(45,58)(51,70)(54,78)(63,65)(64,71)(68,69)\\(79,100)(80,93)(82,
    96)(83,87)(85,88)(90,91)(92,97)(94,99),

    (1,15,42,57)(2,73,47,54)(3,43,96,
    49)(4,99,26,22)(5,78,76,44)(7,65,18,35)(8,97,91,23)(9,83,64,32)(10,46)(11,85,\\41,
    13)(12,20,25,38)(14,24,87,71)(17,34)(19,82,31,60)(21,92,90,55)(28,93,79,70)(29,
    37,88,72)(30,50)(33,52,63,39)\\(36,98,61,89)(40,81)(45,69,58,68)(48,66)(51,67,80,
    100)(59,62,94,74)(77,86),

     (1,98,42,89)(2,23,47,97)(3,49,96,43)(4,59,26,94)(5,
    92,76,55)(6,84)(7,100,18,67)(8,54,91,73)(9,83,64,32)(10,46)\\(11,88,41,29)(12,45,
    25,58)(13,37,85,72)(14,24,87,71)(15,61,57,36)(16,56)(17,34)(19,60,31,82)(20,69,
    38,68)(21,\\78,90,44)(22,74,99,62)(27,53)(28,52,79,39)(30,50)(33,93,63,70)(35,80,
    65,51)(40,81)(48,66)(75,95)(77,86),

    (2,54)(3,51)(4,37)(5,76)(6,40)(7,19)(8,
    23)(9,63)(11,41)(12,61)(13,59)(14,28)(15,45)(16,48)(17,30)(18,31)(20,\\38)(24,
    39)(25,36)(26,72)(32,70)(33,64)(34,50)(35,43)(44,78)(47,73)(49,65)(52,71)(56,
    66)(57,58)(60,67)(62,74)\\(79,87)(80,96)(81,84)(82,100)(83,93)(85,94)(89,98)(91,
    97),

    (2,78,73,5)(3,70,14,67)(4,11,72,74)(6,40)(7,43,33,71)(8,55,23,21)(9,65,
    31,39)(10,46)(12,61,25,36)(13,22,59,29)\\(15,58,57,45)(16,48,56,66)(17,50,34,
    30)(18,49,63,24)(19,52,64,35)(26,41,37,62)(27,53)(28,32,51,60)(44,54,76,47)\\(68,
    69)(79,83,80,82)(81,84)(85,99,94,88)(87,100,96,93)(89,98)(90,91,92,97),

  (1,57,69,58)(2,44,73,5)(3,51,64,33)(4,94)(6,81,84,40)(7,87,28,31)(8,55,23,90)(9,
    63,96,80)(11,22,41,99)(12,89,\\36,20)(13,72)(14,79,19,18)(15,68,45,42)(16,48)(17,
    50,34,30)(21,91,92,97)(24,52,60,100)(25,98,61,38)(26,59)(27,\\53)(29,74,88,62)(32,
    70,49,65)(35,83,93,43)(37,85)(39,82,67,71)(47,78,54,76)(56,66)(77,86),

  (1,20)(2,23)(3,70)(4,13)(5,21)(6,40)(7,71)(8,73)(9,35)(10,46)(11,29)(12,45)(14,
    67)(15,36)(16,66)(17,30)(18,\\24)(19,39)(22,74)(25,58)(26,85)(27,53)(28,82)(31,
    52)(32,80)(33,43)(34,50)(37,94)(38,42)(41,88)(44,92)(47,97)\\(48,56)(49,63)(51,
    83)(54,91)(55,78)(57,61)(59,72)(60,79)(62,99)(64,65)(68,98)(69,89)(75,95)(76,
    90)(77,86)(81,\\84)(87,100)(93,96),

    (2,41,73,62)(3,51,14,28)(4,76,72,44)(5,37,
    78,26)(6,40)(7,31,33,9)(8,22,23,29)(11,54,74,47)(12,45,25,58)(13,21,\\59,55)(15,
    61,57,36)(16,30,56,50)(17,48,34,66)(18,19,63,64)(20,38)(24,52,49,35)(32,67,60,
    70)(39,43,65,71)(68,69)\\(75,77)(79,96,80,87)(81,84)(82,93,83,100)(85,90,94,
    92)(86,95)(88,91,99,97),

     (2,72,73,4)(3,70,14,67)(5,11,78,74)(6,40)(7,49,33,
    24)(8,59,23,13)(9,52,31,35)(10,46)(12,58,25,45)(15,61,57,36)\\(16,50,56,30)(17,48,
    34,66)(18,43,63,71)(19,65,64,39)(21,22,55,29)(26,47,37,54)(27,53)(28,60,51,
    32)(41,44,62,76)\\(75,77)(79,82,80,83)(81,84)(85,91,94,97)(86,95)(87,100,96,
    93)(88,90,99,92),

    (1,42)(2,37,73,26)(3,93,14,100)(4,47,72,54)(5,41,78,62)(6,
    40)(7,43,33,71)(8,94,23,85)(9,39,31,65)(10,46)(11,44,\\74,76)(12,45,25,58)(13,91,
    59,97)(15,36,57,61)(16,50,56,30)(17,48,34,66)(18,49,63,24)(19,35,64,52)(20,
    38)(21,99,\\55,88)(22,92,29,90)(27,53)(28,82,51,83)(32,79,60,80)(67,96,70,87)(68,
    69)(75,77)(81,84)(86,95)(89,98),

     (1,57,98,36,42,15,89,61)(2,22,55,26,47,99,92,
    4)(3,51,31,33,96,80,19,63)(5,94,97,62,76,59,23,74)(6,81,84,40)(7,\\14,79,9,18,87,
    28,64)(8,41,78,13,91,11,44,85)(12,68,58,38,25,69,45,20)(16,50)(17,48,34,66)(21,
    72,54,29,90,37,73,\\88)(24,39,83,67,71,52,32,100)(27,53)(30,56)(35,82,93,43,65,60,
    70,49)(75,86,95,77),

     (1,69,42,68)(2,11,47,41)(3,33,96,63)(4,78,26,44)(5,37,76,
    72)(6,81)(7,87,18,14)(8,99,91,22)(9,28,64,79)(12,15,25,\\57)(13,92,85,55)(16,
    50)(17,48)(19,80,31,51)(20,98,38,89)(21,59,90,94)(23,88,97,29)(24,100,71,67)(30,
    56)(32,39,\\83,52)(34,66)(35,60,65,82)(36,45,61,58)(40,84)(43,70,49,93)(54,62,73,
    74)(75,77)(86,95),

    (1,68)(2,37)(3,35)(4,54)(5,11)(6,81)(7,83)(8,94)(9,67)(10,
    46)(12,15)(13,97)(14,52)(16,30)(17,48)(18,32)(19,93)\\(20,98)(21,22)(23,85)(24,
    79)(25,57)(26,73)(27,53)(28,71)(29,55)(31,70)(33,82)(34,66)(36,58)(38,89)(39,
    87)(40,84)\\(41,76)(42,69)(43,51)(44,62)(45,61)(47,72)(49,80)(50,56)(59,91)(60,
    63)(64,100)(65,96)(74,78)(75,77)(86,95)(88,\\92)(90,99),

  (2,11,65,24,76,37,63,9)(3,8,59,28,60,21,22,67)(4,52,19,78,62,7,43,54)(5,72,33,64,
    47,41,35,71)(6,40,75,77)(10,50,\\53,16,46,30,27,56)(12,20,58,69,25,38,45,68)(13,
    70,32,23,29,51,14,55)(15,61)(17,48,34,66)(18,49,73,26,39,31,44,74)\\(36,57)(79,82,
    90,99,100,96,91,94)(80,87,92,85,93,83,97,88)(81,95,86,84)(89,98),

  (2,26,63,49,76,74,65,31)(3,55,22,51,60,23,59,70)(4,33,43,5,62,35,19,47)(6,40,75,
    77)(7,64,78,72,52,71,54,41)(8,29,\\67,32,21,13,28,14)(9,44,37,39,24,73,11,18)(10,
    50,27,56,46,30,53,16)(12,38,45,69,25,20,58,68)(15,36)(17,66,34,48)\\(57,61)(79,87,
    91,88,100,83,90,85)(80,82,97,94,93,96,92,99)(81,95,86,84)(89,98),

  (1,58,57,38,42,45,15,20)(2,22,52,31,47,99,39,19)(3,55,26,100,96,92,4,67)(5,59,18,
    49,76,94,7,43)(6,40,75,86)(8,72,\\51,87,91,37,80,14)(9,78,29,33,64,44,88,63)(10,
    50,27,56)(11,70,32,90,41,93,83,21)(12,61,68,98,25,36,69,89)(13,65,71,\\54,85,35,
    24,73)(16,46,30,53)(17,48)(23,62,79,60,97,74,28,82)(34,66)(77,84,81,95),

  (1,58,15,38,42,45,57,20)(2,13,18,9,47,85,7,64)(3,8,62,93,96,91,74,70)(4,51,82,21,
    26,80,60,90)(5,29,52,24,76,88,39,\\71)(6,40,75,86)(10,50,53,16)(11,28,87,92,41,79,
    14,55)(12,36,68,98,25,61,69,89)(17,66)(19,44,94,35,31,78,59,65)(22,\\33,43,54,99,
    63,49,73)(23,72,67,32,97,37,100,83)(27,56,46,30)(34,48)(77,84,81,95),

  (2,65,19)(3,8,67)(4,72,62)(5,52,24)(6,40,75)(7,49,54)(9,78,33)(10,50,34)(12,61,
    69)(13,22,59)(14,55,28)(15,20,58)\\(16,48,53)(17,46,30)(18,43,73)(21,70,32)(23,51,
    60)(25,36,68)(26,37,74)(27,56,66)(31,47,35)(38,45,57)(39,71,76)(44,\\63,64)(79,87,
    92)(80,82,97)(81,95,84)(83,90,93)(85,99,94)(91,100,96),

  (1,42)(2,35,19,47,65,31)(3,91,67,96,8,100)(4,37,62,26,72,74)(5,39,24,76,52,71)(6,
    40,75)(7,43,54,18,49,73)(9,44,\\33,64,78,63)(10,50,34)(11,41)(12,36,69,25,61,
    68)(13,99,59,85,22,94)(14,92,28,87,55,79)(15,38,58,57,20,45)(16,48,\\53)(17,46,
    30)(21,93,32,90,70,83)(23,80,60,97,51,82)(27,56,66)(29,88)(81,95,84)(89,98),

  (1,57,12,89,61,58,42,15,25,98,36,45)(2,65,87,97,80,24,47,35,14,23,51,71)(3,90,93,
    43,78,33,96,21,70,49,44,63)(4,\\59,26,94)(5,39,60,55,79,19,76,52,82,92,28,31)(6,
    81,95,84,40,75)(7,32,91,67,64,54,18,83,8,100,9,73)(10,30,17,46,50,\\34)(11,99,37,
    29,74,13,41,22,72,88,62,85)(16,48,27,56,66,53)(20,69,38,68)(77,86),

  (1,2,41,69,54,62)(3,63,64,7,87,28)(4,57,92,26,15,55)(5,94,36)(6,50,16,81,34,66)(8,
    22,89,23,88,38)(9,18,14,79,96,\\33)(11,68,73,74,42,47)(12,44,13)(17,48,84,30,56,
    40)(19,51)(20,91,99,98,97,29)(21,72,45,90,37,58)(24,52,43,70,83,\\100)(25,78,
    85)(27,86,95)(31,80)(32,67,71,39,49,93)(35,82)(53,77,75)(59,61,76)(60,65),

  (1,76,22,89,55,62)(2,72,58,8,59,61)(3,65)(4,57,97,13,25,54)(5,99,98,92,74,42)(6,
    30,16,81,34,66)(7,83,79,82,33,24)\\(9,39,31,93,14,100)(10,46)(11,69,78,88,20,
    90)(12,73,26,15,23,85)(17,48,84,50,56,40)(18,32,28,60,63,71)(19,70,87,67,\\64,
    52)(21,41,68,44,29,38)(27,86,95,53,77,75)(35,96)(36,47,37,45,91,94)(43,51)(49,
    80),

    (1,2,94,32,65,42,47,59,83,35)(3,79,36,8,88,96,28,61,91,29)(4,9,51,89,97,
    26,64,80,98,23)(5,13,19,93,38,76,85,31,70,\\20)(6,75,53,48,46)(7,45,21,22,71,18,
    58,90,99,24)(10,84,95,27,66)(11,49,39,15,54,41,43,52,57,73)(12,44,62,87,100,25,
    78,74,14,67)(16,40,17,86,30)(33,69,92,37,60,63,68,55,72,82)(34,77,50,56,81),

  (1,47,94,83,65)(2,59,32,35,42)(3,28,36,91,88)(4,64,51,98,97)(5,85,19,70,38)(6,75,
    53,48,46)(7,58,21,99,71)(8,29,\\96,79,61)(9,80,89,23,26)(10,84,95,27,66)(11,43,39,
    57,54)(12,78,62,14,100)(13,31,93,20,76)(15,73,41,49,52)(16,40,\\17,86,30)(18,45,
    90,22,24)(25,44,74,87,67)(33,68,92,72,60)(34,77,50,56,81)(37,82,63,69,55),

  (1,44,72,87,80,98,21,13,71,65,42,78,37,14,51,89,90,85,24,35)(2,88,83,7,45,23,41,
    64,100,25,47,29,32,18,58,97,11,\\9,67,12)(3,93,20,54,4,49,63,69,91,59,96,70,38,73,
    26,43,33,68,8,94)(5,22,19,79,61,92,74,60,39,57,76,99,31,28,36,55,\\62,82,52,15)(6,
    75,27,48,10,84,95,53,66,46)(16,40,34,77,50,56,81,17,86,30),

  (1,78,72,14,80,89,21,85,71,35,42,44,37,87,51,98,90,13,24,65)(2,29,83,18,45,97,41,
    9,100,12,47,88,32,7,58,23,11,\\64,67,25)(3,70,20,73,4,43,63,68,91,94,96,93,38,54,
    26,49,33,69,8,59)(5,99,19,28,61,55,74,82,39,15,76,22,31,79,36,92,\\62,60,52,57)(6,
    75,27,48,10,84,95,53,66,46)(16,40,34,77,50,56,81,17,86,30).

The character table of $G^{s_4}$:\\
\begin{tabular}{c|cccccccccccccccccccccccc}
  & & & & & & & & & & 10& & & & & & & & & &20 & & & &\\\hline
$\chi_4^{(1)}$&1&1&1&1&1&1&1&1&1&1&1&1&1&1&1&1&1&1&1&1&1&1&1&1
\\$\chi_4^{(2)}$&1&1&1&1&1&-1&-1&-1&-1&1&1&1&1&1&1&-1&-1&-1&-1&1&1&1&-1&-1
\\$\chi_4^{(3)}$&4&4&4&4&4&-2&-2&-2&-2&.&.&.&.&.&.&.&.&.&.&1&1&1&1&1
\\$\chi_4^{(4)}$&4&4&4&4&4&2&2&2&2&.&.&.&.&.&.&.&.&.&.&1&1&1&-1&-1
\\$\chi_4^{(5)}$&5&5&5&5&5&-1&-1&-1&-1&1&1&1&1&1&1&1&1&1&1&-1&-1&-1&-1&-1
\\$\chi_4^{(6)}$&5&5&5&5&5&1&1&1&1&1&1&1&1&1&1&-1&-1&-1&-1&-1&-1&-1&1&1
\\$\chi_4^{(7)}$&6&6&6&6&6&.&.&.&.&-2&-2&-2&-2&-2&-2&.&.&.&.&.&.&.&.&.
\\$\chi_4^{(8)}$&6&6&-2&2&-6&.&.&.&.&-2&2&2&.&2&-2&.&.&-2&2&.&.&.&.&.
\\$\chi_4^{(9)}$&6&6&-2&2&-6&.&.&.&.&-2&2&2&.&2&-2&.&.&2&-2&.&.&.&.&.
\\$\chi_4^{(10)}$&6&6&-2&2&-6&.&.&.&.&2&-2&-2&.&-2&2&A&-A&.&.&.&.&.&.&.
\\$\chi_4^{(11)}$&6&6&-2&2&-6&.&.&.&.&2&-2&-2&.&-2&2&-A&A&.&.&.&.&.&.&.
\\$\chi_4^{(12)}$&8&-8&.&.&.&.&.&.&.&.&4&-4&.&.&.&.&.&.&.&2&-2&.&.&.
\\$\chi_4^{(13)}$&10&10&2&-2&-10&-4&.&.&4&2&-2&-2&.&2&-2&.&.&.&.&1&1&-1&-1&1
\\$\chi_4^{(14)}$&10&10&2&-2&-10&-2&-2&2&2&-2&2&2&.&-2&2&.&.&.&.&1&1&-1&1&-1
\\$\chi_4^{(15)}$&10&10&2&-2&-10&2&2&-2&-2&-2&2&2&.&-2&2&.&.&.&.&1&1&-1&-1&1
\\$\chi_4^{(16)}$&10&10&2&-2&-10&4&.&.&-4&2&-2&-2&.&2&-2&.&.&.&.&1&1&-1&1&-1
\\$\chi_4^{(17)}$&15&15&-1&-1&15&-3&1&1&-3&3&3&3&-1&-1&-1&-1&-1&1&1&.&.&.&.&.
\\$\chi_4^{(18)}$&15&15&-1&-1&15&3&-1&-1&3&3&3&3&-1&-1&-1&1&1&-1&-1&.&.&.&.&.
\\$\chi_4^{(19)}$&15&15&-1&-1&15&-3&1&1&-3&-1&-1&-1&-1&3&3&1&1&-1&-1&.&.&.&.&.
\\$\chi_4^{(20)}$&15&15&-1&-1&15&3&-1&-1&3&-1&-1&-1&-1&3&3&-1&-1&1&1&.&.&.&.&.
\\$\chi_4^{(21)}$&20&20&4&-4&-20&-2&2&-2&2&.&.&.&.&.&.&.&.&.&.&-1&-1&1&1&-1
\\$\chi_4^{(22)}$&20&20&4&-4&-20&2&-2&2&-2&.&.&.&.&.&.&.&.&.&.&-1&-1&1&-1&1
\\$\chi_4^{(23)}$&24&24&-8&8&-24&.&.&.&.&.&.&.&.&.&.&.&.&.&.&.&.&.&.&.
\\$\chi_4^{(24)}$&24&-24&.&.&.&.&.&.&.&.&-4&4&.&.&.&.&.&.&.&.&.&.&.&.
\\$\chi_4^{(25)}$&24&-24&.&.&.&.&.&.&.&.&-4&4&.&.&.&.&.&.&.&.&.&.&.&.
\\$\chi_4^{(26)}$&30&30&-2&-2&30&.&.&.&.&-2&-2&-2&2&-2&-2&.&.&.&.&.&.&.&.&.
\\$\chi_4^{(27)}$&32&-32&.&.&.&.&.&.&.&.&.&.&.&.&.&.&.&.&.&2&-2&.&.&.
\\$\chi_4^{(28)}$&40&-40&.&.&.&.&.&.&.&.&4&-4&.&.&.&.&.&.&.&-2&2&.&.&.
\end{tabular}

\begin{tabular}{c|cccc}
  & & & &\\\hline
$\chi_4^{(1)}$&1&1&1&1
\\$\chi_4^{(2)}$&1&1&1&1
\\$\chi_4^{(3)}$&-1&-1&-1&-1
\\$\chi_4^{(4)}$&-1&-1&-1&-1
\\$\chi_4^{(5)}$&.&.&.&.
\\$\chi_4^{(6)}$&.&.&.&.
\\$\chi_4^{(7)}$&1&1&1&1
\\$\chi_4^{(8)}$&1&1&-1&-1
\\$\chi_4^{(9)}$&1&1&-1&-1
\\$\chi_4^{(10)}$&1&1&-1&-1
\\$\chi_4^{(11)}$&1&1&-1&-1
\\$\chi_4^{(12)}$&2&-2&.&.
\\$\chi_4^{(13)}$&.&.&.&.
\\$\chi_4^{(14)}$&.&.&.&.
\\$\chi_4^{(15)}$&.&.&.&.
\\$\chi_4^{(16)}$&.&.&.&.
\\$\chi_4^{(17)}$&.&.&.&.
\\$\chi_4^{(18)}$&.&.&.&.
\\$\chi_4^{(19)}$&.&.&.&.
\\$\chi_4^{(20)}$&.&.&.&.
\\$\chi_4^{(21)}$&.&.&.&.
\\$\chi_4^{(22)}$&.&.&.&.
\\$\chi_4^{(23)}$&-1&-1&1&1
\\$\chi_4^{(24)}$&1&-1&B&-B
\\$\chi_4^{(25)}$&1&-1&-B&B
\\$\chi_4^{(26)}$&.&.&.&.
\\$\chi_4^{(27)}$&-2&2&.&.
\\$\chi_4^{(28)}$&.&.&.&.
\end{tabular}

 \noindent where   A = -2*E(4)
  = -2*ER(-1) = -2i;
B = -E(20)-E(20)$^9$+E(20)$^{13}$+E(20)$^{17}$
  = -ER(-5) = -i5.

The generators of $G^{s_5}$ are:\\
(  1, 70, 64, 90, 42, 93,  9, 21)(  2, 36,  7, 71, 47, 61, 18, 24)
    (  3, 73, 69, 79, 96, 54, 68, 28)(  4, 13, 74, 41, 26, 85, 62, 11)
    (  5, 45, 51, 31, 76, 58, 80, 19)(  6, 46, 50, 95)(  8, 20, 39, 49, 91, 38, 52,
      43)( 10, 30, 75, 84)( 12, 35, 60, 55, 25, 65, 82, 92)( 14, 23, 57,100, 87, 97,
     15, 67)( 16, 66, 53, 77)( 17, 40)( 22, 29, 59, 72, 99, 88, 94, 37)
    ( 27, 86, 56, 48)( 32, 78, 89, 33, 83, 44, 98, 63)( 34, 81),
  (  1, 14)(  2, 63)(  3, 45)(  4, 22)(  5, 28)(  6, 50)(  7, 78)(  8, 35)(  9, 15)
    ( 10, 75)( 11, 37)( 12, 43)( 13, 29)( 16, 27)( 17, 34)( 18, 44)( 19, 68)
    ( 20, 60)( 21, 67)( 23, 70)( 24, 98)( 25, 49)( 26, 99)( 30, 84)( 31, 69)
    ( 32, 36)( 33, 47)( 38, 82)( 39, 55)( 40, 81)( 41, 72)( 42, 87)( 46, 95)
    ( 48, 77)( 51, 73)( 52, 92)( 53, 56)( 54, 80)( 57, 64)( 58, 96)( 59, 74)
    ( 61, 83)( 62, 94)( 65, 91)( 66, 86)( 71, 89)( 76, 79)( 85, 88)( 90,100)
    ( 93, 97).

The representatives of conjugacy classes of   $G^{s_5}$ are:\\
(1),
(1,9,42,64)(2,18,47,7)(3,68,96,69)(4,62,26,74)(5,80,76,51)(6,50)(8,52,91,
    39)(10,75)(11,85,41,13)(12,82,25,60)\\(14,15,87,57)(16,53)(19,58,31,45)(20,43,38,
    49)(21,93,90,70)(22,94,99,59)(23,67,97,100)(24,61,71,36)(27,56)(28,\\54,79,73)(29,
    37,88,72)(30,84)(32,98,83,89)(33,78,63,44)(35,92,65,55)(46,95)(48,86)(66,77),

  (1,14)(2,63)(3,45)(4,22)(5,28)(6,50)(7,78)(8,35)(9,15)(10,75)(11,37)(12,43)(13,
    29)(16,27)(17,34)(18,44)(19,\\68)(20,60)(21,67)(23,70)(24,98)(25,49)(26,99)(30,
    84)(31,69)(32,36)(33,47)(38,82)(39,55)(40,81)(41,72)(42,87)\\(46,95)(48,77)(51,
    73)(52,92)(53,56)(54,80)(57,64)(58,96)(59,74)(61,83)(62,94)(65,91)(66,86)(71,
    89)(76,79)(85,\\88)(90,100)(93,97),

    (1,15,42,57)(2,44,47,78)(3,19,96,31)(4,94,
    26,59)(5,54,76,73)(7,63,18,33)(8,92,91,55)(9,87,64,14)(11,88,41,29)\\(12,38,25,
    20)(13,37,85,72)(16,56)(17,34)(21,97,90,23)(22,62,99,74)(24,83,71,32)(27,53)(28,
    80,79,51)(35,52,65,\\39)(36,98,61,89)(40,81)(43,82,49,60)(45,68,58,69)(48,66)(67,
    93,100,70)(77,86),

    (1,21,9,93,42,90,64,70)(2,24,18,61,47,71,7,36)(3,28,68,54,
    96,79,69,73)(4,11,62,85,26,41,74,13)(5,19,80,58,76,\\31,51,45)(6,95,50,46)(8,43,
    52,38,91,49,39,20)(10,84,75,30)(12,92,82,65,25,55,60,35)(14,67,15,97,87,100,57,
    23)\\(16,77,53,66)(17,40)(22,37,94,88,99,72,59,29)(27,48,56,86)(32,63,98,44,83,33,
    89,78)(34,81),

     (1,23,64,100,42,97,9,67)(2,32,7,89,47,83,18,98)(3,51,69,76,96,
    80,68,5)(4,29,74,72,26,88,62,37)(6,95,50,46)(8,\\60,39,25,91,82,52,12)(10,84,75,
    30)(11,22,13,59,41,99,85,94)(14,70,57,90,87,93,15,21)(16,86,53,48)(17,81)(19,28,
    45,73,31,79,58,54)(20,55,49,65,38,92,43,35)(24,63,36,78,71,33,61,44)(27,66,56,
    77)(34,40),

     (1,42)(2,47)(3,96)(4,26)(5,76)(7,18)(8,91)(9,64)(11,41)(12,25)(13,
    85)(14,87)(15,57)(19,31)(20,38)(21,90)(22,\\99)(23,97)(24,71)(28,79)(29,88)(32,
    83)(33,63)(35,65)(36,61)(37,72)(39,52)(43,49)(44,78)(45,58)(51,80)(54,73)\\(55,
    92)(59,94)(60,82)(62,74)(67,100)(68,69)(70,93)(89,98),

  (1,57,42,15)(2,78,47,44)(3,31,96,19)(4,59,26,94)(5,73,76,54)(7,33,18,63)(8,55,91,
    92)(9,14,64,87)(11,29,41,88)\\(12,20,25,38)(13,72,85,37)(16,56)(17,34)(21,23,90,
    97)(22,74,99,62)(24,32,71,83)(27,53)(28,51,79,80)(35,39,65,52)\\(36,89,61,98)(40,
    81)(43,60,49,82)(45,69,58,68)(48,66)(67,70,100,93)(77,86),

  (1,64,42,9)(2,7,47,18)(3,69,96,68)(4,74,26,62)(5,51,76,80)(6,50)(8,39,91,52)(10,
    75)(11,13,41,85)(12,60,25,82)\\(14,57,87,15)(16,53)(19,45,31,58)(20,49,38,43)(21,
    70,90,93)(22,59,99,94)(23,100,97,67)(24,36,71,61)(27,56)(28,\\73,79,54)(29,72,88,
    37)(30,84)(32,89,83,98)(33,44,63,78)(35,55,65,92)(46,95)(48,86)(66,77),

  (1,67,9,97,42,100,64,23)(2,98,18,83,47,89,7,32)(3,5,68,80,96,76,69,51)(4,37,62,88,
    26,72,74,29)(6,46,50,95)(8,\\12,52,82,91,25,39,60)(10,30,75,84)(11,94,85,99,41,59,
    13,22)(14,21,15,93,87,90,57,70)(16,48,53,86)(17,81)(19,54,\\58,79,31,73,45,28)(20,
    35,43,92,38,65,49,55)(24,44,61,33,71,78,36,63)(27,77,56,66)(34,40),

  (1,70,64,90,42,93,9,21)(2,36,7,71,47,61,18,24)(3,73,69,79,96,54,68,28)(4,13,74,41,
    26,85,62,11)(5,45,51,31,76,\\58,80,19)(6,46,50,95)(8,20,39,49,91,38,52,43)(10,30,
    75,84)(12,35,60,55,25,65,82,92)(14,23,57,100,87,97,15,67)\\(16,66,53,77)(17,
    40)(22,29,59,72,99,88,94,37)(27,86,56,48)(32,78,89,33,83,44,98,63)(34,81),

  (1,87)(2,33)(3,58)(4,99)(5,79)(6,50)(7,44)(8,65)(9,57)(10,75)(11,72)(12,49)(13,
    88)(14,42)(15,64)(16,27)(17,\\34)(18,78)(19,69)(20,82)(21,100)(22,26)(23,93)(24,
    89)(25,43)(28,76)(29,85)(30,84)(31,68)(32,61)(35,91)(36,83)\\(37,41)(38,60)(39,
    92)(40,81)(45,96)(46,95)(47,63)(48,77)(51,54)(52,55)(53,56)(59,62)(66,86)(67,
    90)(70,97)(71,\\98)(73,80)(74,94),

     (1,90,9,70,42,21,64,93)(2,71,18,36,47,24,7,
    61)(3,79,68,73,96,28,69,54)(4,41,62,13,26,11,74,85)(5,31,80,45,76,\\19,51,58)(6,
    95,50,46)(8,49,52,20,91,43,39,38)(10,84,75,30)(12,55,82,35,25,92,60,65)(14,100,
    15,23,87,67,57,97)\\(16,77,53,66)(17,40)(22,72,94,29,99,37,59,88)(27,48,56,86)(32,
    33,98,78,83,63,89,44)(34,81),

    (1,93,64,21,42,70,9,90)(2,61,7,24,47,36,18,
    71)(3,54,69,28,96,73,68,79)(4,85,74,11,26,13,62,41)(5,58,51,19,76,\\45,80,31)(6,
    46,50,95)(8,38,39,43,91,20,52,49)(10,30,75,84)(12,65,60,92,25,35,82,55)(14,97,
    57,67,87,23,15,100)\\(16,66,53,77)(17,40)(22,88,59,37,99,29,94,72)(27,86,56,
    48)(32,44,89,63,83,78,98,33)(34,81),

    (1,97,64,67,42,23,9,100)(2,83,7,98,47,32,
    18,89)(3,80,69,5,96,51,68,76)(4,88,74,37,26,29,62,72)(6,95,50,46)(8,\\82,39,12,91,
    60,52,25)(10,84,75,30)(11,99,13,94,41,22,85,59)(14,93,57,21,87,70,15,90)(16,86,
    53,48)(17,81)(19,79,\\45,54,31,28,58,73)(20,92,49,35,38,55,43,65)(24,33,36,44,71,
    63,61,78)(27,66,56,77)(34,40),

    (1,100,9,23,42,67,64,97)(2,89,18,32,47,98,7,
    83)(3,76,68,51,96,5,69,80)(4,72,62,29,26,37,74,88)(6,46,50,95)(8,\\25,52,60,91,12,
    39,82)(10,30,75,84)(11,59,85,22,41,94,13,99)(14,90,15,70,87,21,57,93)(16,48,53,
    86)(17,81)(19,73,\\58,28,31,54,45,79)(20,65,43,55,38,35,49,92)(24,78,61,63,71,44,
    36,33)(27,77,56,66)(34,40).

The character table of $G^{s_5}$:\\
\begin{tabular}{c|cccccccccccccccc}
  & & & & & & & & & &10 & & & & & &\\\hline

$\chi_5^{(1)}$&1&1&1&1&1&1&1&1&1&1&1&1&1&1&1&1
\\$\chi_5^{(2)}$&1&1&-1&-1&-1&1&1&-1&1&1&-1&-1&-1&-1&1&1
\\$\chi_5^{(3)}$&1&1&-1&-1&1&-1&1&-1&1&-1&1&-1&1&1&-1&-1
\\$\chi_5^{(4)}$&1&1&1&1&-1&-1&1&1&1&-1&-1&1&-1&-1&-1&-1
\\$\chi_5^{(5)}$&1&-1&-1&1&A&A&1&1&-1&-A&-A&-1&A&-A&A&-A
\\$\chi_5^{(6)}$&1&-1&-1&1&-A&-A&1&1&-1&A&A&-1&-A&A&-A&A
\\$\chi_5^{(7)}$&1&-1&1&-1&A&-A&1&-1&-1&A&-A&1&A&-A&-A&A
\\$\chi_5^{(8)}$&1&-1&1&-1&-A&A&1&-1&-1&-A&A&1&-A&A&A&-A
\\$\chi_5^{(9)}$&1&A&-1&-A&B&-/B&-1&A&-A&-B&/B&1&-B&-/B&/B&B
\\$\chi_5^{(10)}$&1&A&-1&-A&-B&/B&-1&A&-A&B&-/B&1&B&/B&-/B&-B
\\$\chi_5^{(11)}$&1&-A&-1&A&-/B&B&-1&-A&A&/B&-B&1&/B&B&-B&-/B
\\$\chi_5^{(12)}$&1&-A&-1&A&/B&-B&-1&-A&A&-/B&B&1&-/B&-B&B&/B
\\$\chi_5^{(13)}$&1&A&1&A&B&/B&-1&-A&-A&B&/B&-1&-B&-/B&-/B&-B
\\$\chi_5^{(14)}$&1&A&1&A&-B&-/B&-1&-A&-A&-B&-/B&-1&B&/B&/B&B
\\$\chi_5^{(15)}$&1&-A&1&-A&-/B&-B&-1&A&A&-/B&-B&-1&/B&B&B&/B
\\$\chi_5^{(16)}$&1&-A&1&-A&/B&B&-1&A&A&/B&B&-1&-/B&-B&-B&-/B

\end{tabular}

\noindent \noindent where   A = -E(4)
  = -ER(-1) = -i;
B = -E(8)$^3$.

The generators of $G^{s_6}$ are:\\
(  1, 57, 42, 15)(  2, 78, 47, 44)(  3, 31, 96, 19)(  4, 59, 26, 94)
    (  5, 73, 76, 54)(  7, 33, 18, 63)(  8, 55, 91, 92)(  9, 14, 64, 87)
    ( 11, 29, 41, 88)( 12, 20, 25, 38)( 13, 72, 85, 37)( 16, 56)( 17, 34)
    ( 21, 23, 90, 97)( 22, 74, 99, 62)( 24, 32, 71, 83)( 27, 53)( 28, 51, 79, 80)
    ( 35, 39, 65, 52)( 36, 89, 61, 98)( 40, 81)( 43, 60, 49, 82)( 45, 69, 58, 68)
    ( 48, 66)( 67, 70,100, 93)( 77, 86), (  2, 33, 76, 35)(  3, 32, 60, 14)
    (  4, 11, 62, 37)(  5, 65, 47, 63)(  6, 75)(  7, 73, 52, 44)(  8, 28, 21, 67)
    (  9, 19, 24, 43)( 10, 30, 46, 50)( 12, 69, 25, 68)( 13, 59, 29, 22)
    ( 16, 53, 56, 27)( 17, 34)( 18, 54, 39, 78)( 20, 58, 38, 45)( 23, 70, 55, 51)
    ( 26, 41, 74, 72)( 31, 71, 49, 64)( 36, 61)( 48, 66)( 79, 90,100, 91)
    ( 80, 97, 93, 92)( 82, 87, 96, 83)( 84, 95)( 85, 94, 88, 99)( 89, 98),
  (  1,  2, 15, 44, 42, 47, 57, 78)(  3, 67, 19, 93, 96,100, 31, 70)
    (  4, 62, 94, 99, 26, 74, 59, 22)(  5, 69, 54, 45, 76, 68, 73, 58)(  6, 50)
    (  7, 60, 63, 43, 18, 82, 33, 49)(  8, 12, 92, 38, 91, 25, 55, 20)
    (  9, 39, 87, 35, 64, 52, 14, 65)( 10, 46)( 11, 13, 88, 37, 41, 85, 29, 72)
    ( 16, 66, 56, 48)( 17, 81, 34, 40)( 21, 89, 97, 36, 90, 98, 23, 61)
    ( 24, 51, 83, 28, 71, 80, 32, 79)( 27, 77, 53, 86)( 30, 84).

The representatives of conjugacy classes of   $G^{s_6}$ are:\\
 (1),
 (2,33,76,35)(3,32,60,14)(4,11,62,37)(5,65,47,63)(6,75)(7,73,52,44)(8,28,21,
    67)(9,19,24,43)(10,30,46,50)(12,\\69,25,68)(13,59,29,22)(16,53,56,27)(17,34)(18,
    54,39,78)(20,58,38,45)(23,70,55,51)(26,41,74,72)(31,71,49,64)\\(36,61)(48,66)(79,
    90,100,91)(80,97,93,92)(82,87,96,83)(84,95)(85,94,88,99)(89,98),

  (2,35,76,33)(3,14,60,32)(4,37,62,11)(5,63,47,65)(6,75)(7,44,52,73)(8,67,21,28)(9,
    43,24,19)(10,50,46,30)(12,\\68,25,69)(13,22,29,59)(16,27,56,53)(17,34)(18,78,39,
    54)(20,45,38,58)(23,51,55,70)(26,72,74,41)(31,64,49,71)\\(36,61)(48,66)(79,91,100,
    90)(80,92,93,97)(82,83,96,87)(84,95)(85,99,88,94)(89,98),

  (2,76)(3,60)(4,62)(5,47)(7,52)(8,21)(9,24)(10,46)(11,37)(12,25)(13,29)(14,32)(16,
    56)(18,39)(19,43)(20,38)\\(22,59)(23,55)(26,74)(27,53)(28,67)(30,50)(31,49)(33,
    35)(41,72)(44,73)(45,58)(51,70)(54,78)(63,65)(64,71)(68,\\69)(79,100)(80,93)(82,
    96)(83,87)(85,88)(90,91)(92,97)(94,99),

     (1,2,15,44,42,47,57,78)(3,67,19,93,96,
    100,31,70)(4,62,94,99,26,74,59,22)(5,69,54,45,76,68,73,58)(6,50)(7,60,\\63,43,18,
    82,33,49)(8,12,92,38,91,25,55,20)(9,39,87,35,64,52,14,65)(10,46)(11,13,88,37,41,
    85,29,72)(16,66,56,\\48)(17,81,34,40)(21,89,97,36,90,98,23,61)(24,51,83,28,71,80,
    32,79)(27,77,53,86)(30,84),

    (1,2,49,52,42,47,43,39)(3,79,98,97,96,28,89,23)(4,
    13,22,88,26,85,99,29)(5,9,93,38,76,64,70,20)(6,75,50,46)\\(7,58,91,24,18,45,8,
    71)(10,84,95,30)(11,94,37,62,41,59,72,74)(12,54,87,100,25,73,14,67)(15,44,60,65,
    57,78,82,\\35)(16,86,27,66)(17,40)(19,51,61,90,31,80,36,21)(32,63,69,55,83,33,68,
    92)(34,81)(48,56,77,53),

     (1,2,64,7,42,47,9,18)(3,65,69,92,96,35,68,55)(4,41,
    74,85,26,11,62,13)(5,43,51,20,76,49,80,38)(6,75,50,10)(8,\\19,39,45,91,31,52,
    58)(12,73,60,79,25,54,82,28)(14,63,57,78,87,33,15,44)(16,77,53,66)(17,40)(21,71,
    70,61,90,24,\\93,36)(22,72,59,88,99,37,94,29)(23,83,100,98,97,32,67,89)(27,48,56,
    86)(30,46,84,95)(34,81),

    (1,2,68,54)(3,63,9,51)(4,94,26,59)(5,57,78,45)(6,50,
    84,30)(7,14,79,31)(8,89,97,38)(11,41)(12,55,61,21)(13,\\72,85,37)(15,44,58,76)(16,
    48)(17,81,34,40)(18,87,28,19)(20,91,98,23)(24,39,82,100)(25,92,36,90)(27,86)(29,
    88)\\(32,65,43,93)(33,64,80,96)(35,49,70,83)(42,47,69,73)(52,60,67,71)(53,77)(56,
    66),

    (1,3)(2,28)(4,62)(5,33)(6,30)(7,54)(8,67)(9,69)(10,75)(11,13)(12,24)(14,
    58)(15,19)(16,53)(17,34)(18,73)(20,\\32)(21,65)(22,59)(23,52)(25,71)(26,74)(27,
    56)(29,72)(31,57)(35,90)(36,82)(37,88)(38,83)(39,97)(40,81)(41,85)\\(42,96)(43,
    89)(44,80)(45,87)(46,95)(47,79)(48,86)(49,98)(50,84)(51,78)(55,70)(60,61)(63,
    76)(64,68)(66,77)(91,\\100)(92,93)(94,99),

    (1,3,61,60)(2,63,21,67)(5,79,91,
    35)(6,30,84,50)(7,23,70,78)(8,65,76,28)(9,12,71,68)(10,95,46,75)(11,88,41,
    29)\\(13,72,85,37)(14,20,83,45)(15,19,89,43)(16,27)(17,34)(18,97,93,44)(24,69,64,
    25)(31,98,49,57)(32,58,87,38)(33,\\90,100,47)(36,82,42,96)(39,73,80,92)(40,81)(48,
    86)(51,55,52,54)(53,56)(66,77),

    (1,5,69,44)(2,57,73,58)(3,79,64,18)(4,26)(6,
    30,84,50)(7,96,28,9)(8,12,23,36)(11,88,41,29)(14,33,19,51)(15,54,\\45,47)(16,
    66)(17,81,34,40)(20,90,89,55)(21,98,92,38)(22,74,99,62)(24,93,60,35)(25,97,61,
    91)(27,77)\\(31,80,87,63)(32,67,49,39)(42,76,68,78)(43,52,83,100)(48,56)(53,
    86)(59,94)(65,71,70,82),

    (1,7,57,33,42,18,15,63)(2,3,78,31,47,96,44,19)(4,22,
    59,74,26,99,94,62)(5,87,73,9,76,14,54,64)(8,83,55,24,91,\\32,92,71)(10,75)(11,85,
    29,37,41,13,88,72)(12,65,20,52,25,35,38,39)(16,86,56,77)(17,81,34,40)(21,60,23,
    49,90,\\82,97,43)(27,48,53,66)(28,45,51,69,79,58,80,68)(30,50)(36,93,89,67,61,70,
    98,100)(46,95),

    (1,7,25,65)(2,14,92,43)(3,23,24,76)(5,96,97,71)(8,60,78,64)(9,
    91,82,44)(10,95,46,75)(11,41)(12,35,42,18)(13,\\37,85,72)(15,63,20,39)(16,77)(17,
    81,34,40)(19,21,83,73)(22,74,99,62)(27,66)(28,61,70,69)(29,88)(31,90,32,54)\\(33,
    38,52,57)(36,93,68,79)(45,80,89,67)(47,87,55,49)(48,53)(51,98,100,58)(56,86),

  (1,7,19,76,42,18,31,5)(2,38,80,71,47,20,51,24)(3,54,15,63,96,73,57,33)(4,41,99,72,
    26,11,22,37)(6,10,30,75)\\(8,61,93,43,91,36,70,49)(9,21,45,39,64,90,58,52)(12,28,
    83,44,25,79,32,78)(13,59,88,62,85,94,29,74)(14,23,69,65,\\87,97,68,35)(16,48,27,
    77)(17,40)(34,81)(46,50,95,84)(53,86,56,66)(55,98,67,60,92,89,100,82),

  (1,7,9,2,42,18,64,47)(3,92,68,65,96,55,69,35)(4,85,62,41,26,13,74,11)(5,20,80,43,
    76,38,51,49)(6,10,50,75)(8,\\45,52,19,91,58,39,31)(12,79,82,73,25,28,60,54)(14,78,
    15,63,87,44,57,33)(16,66,53,77)(17,40)(21,61,93,71,90,36,\\70,24)(22,88,94,72,99,
    29,59,37)(23,98,67,83,97,89,100,32)(27,86,56,48)(30,95,84,46)(34,81),

  (1,8,69,23)(2,89,73,20)(3,65,64,70)(5,12,44,36)(6,50,84,30)(7,60,28,24)(9,93,96,
    35)(10,46)(11,88,41,29)(14,\\67,19,39)(15,92,45,21)(16,66)(17,81,34,40)(18,82,79,
    71)(22,62,99,74)(25,78,61,76)(27,77)(31,52,87,100)(32,33,\\49,51)(38,47,98,54)(42,
    91,68,97)(43,80,83,63)(48,56)(53,86)(55,58,90,57)(75,95),

  (1,9,42,64)(2,18,47,7)(3,68,96,69)(4,62,26,74)(5,80,76,51)(6,50)(8,52,91,39)(10,
    75)(11,85,41,13)(12,82,25,60)\\(14,15,87,57)(16,53)(19,58,31,45)(20,43,38,49)(21,
    93,90,70)(22,94,99,59)(23,67,97,100)(24,61,71,36)(27,56)(28,\\54,79,73)(29,37,88,
    72)(30,84)(32,98,83,89)(33,78,63,44)(35,92,65,55)(46,95)(48,86)(66,77),

  (1,9,61,71)(2,51,21,52)(3,12,60,68)(4,26)(5,7,91,70)(6,50,84,30)(8,93,76,18)(10,
    95,46,75)(11,88,41,29)(13,37,\\85,72)(14,98,83,57)(15,87,89,32)(16,27)(19,38,43,
    58)(20,49,45,31)(23,35,78,79)(24,42,64,36)(25,82,69,96)(28,97,\\65,44)(33,92,100,
    73)(39,47,80,90)(48,86)(53,56)(54,63,55,67)(59,94)(66,77),

  (1,9,58,49)(2,78,47,44)(3,98,32,12)(4,85,99,41)(5,55)(6,10,84,46)(7,28,93,65)(8,
    54)(11,26,13,22)(14,68,82,57)\\(15,87,69,60)(17,34)(18,79,70,35)(19,61,24,38)(20,
    31,36,71)(21,97,90,23)(25,96,89,83)(27,53)(29,94,72,74)(30,95)\\(33,51,67,52)(37,
    62,88,59)(39,63,80,100)(42,64,45,43)(48,77,66,86)(50,75)(73,91)(76,92),

  (1,9,38,31)(2,92,90,73)(3,15,87,25)(4,88,99,72)(5,44,91,23)(6,10)(7,33,18,63)(8,
    97,76,78)(11,62,85,59)(12,96,\\57,14)(13,94,41,74)(16,56)(17,34)(19,42,64,20)(21,
    54,47,55)(22,37,26,29)(24,58,43,61)(28,52)(30,75,50,95)(32,68,\\60,98)(35,51)(36,
    71,45,49)(39,79)(46,84)(48,77,66,86)(65,80)(67,93,100,70)(69,82,89,83),

  (1,12)(2,55)(3,71)(4,26)(5,23)(7,35)(8,44)(9,60)(10,46)(11,41)(14,49)(15,38)(17,
    34)(18,65)(19,32)(20,57)(21,\\54)(24,96)(25,42)(28,93)(29,88)(31,83)(33,39)(36,
    69)(40,81)(43,87)(45,98)(47,92)(51,67)(52,63)(58,89)(59,94)\\(61,68)(64,82)(70,
    79)(73,90)(75,95)(76,97)(78,91)(80,100),

    (1,12,42,25)(2,97,47,23)(3,9,96,
    64)(4,62,26,74)(5,92,76,55)(7,63,18,33)(8,54,91,73)(11,37,41,72)(13,88,85,
    29)\\(14,19,87,31)(15,38,57,20)(16,56)(17,34)(21,44,90,78)(22,94,99,59)(24,60,71,
    82)(27,53)(28,51,79,80)(30,50)(32,\\49,83,43)(35,39,65,52)(36,69,61,68)(40,81)(45,
    89,58,98)(67,93,100,70)(75,95),

     (1,15,42,57)(2,44,47,78)(3,19,96,31)(4,94,26,
    59)(5,54,76,73)(7,63,18,33)(8,92,91,55)(9,87,64,14)(11,88,41,29)\\(12,38,25,
    20)(13,37,85,72)(16,56)(17,34)(21,97,90,23)(22,62,99,74)(24,83,71,32)(27,53)(28,
    80,79,51)(35,52,65,\\39)(36,98,61,89)(40,81)(43,82,49,60)(45,68,58,69)(48,66)(67,
    93,100,70)(77,86),

    (1,18,25,35)(2,87,92,49)(3,97,24,5)(4,26)(7,12,65,42)(8,82,
    78,9)(10,95,46,75)(13,72,85,37)(14,55,43,47)(15,33,\\20,52)(16,77)(17,81,34,
    40)(19,90,83,54)(21,32,73,31)(22,62,99,74)(23,71,76,96)(27,66)(28,36,70,68)(38,
    39,57,63)\\(44,64,91,60)(45,51,89,100)(48,53)(56,86)(58,80,98,67)(59,94)(61,93,69,
    79),

     (1,20,42,38)(2,91,47,8)(3,83,96,32)(4,94,26,59)(5,90,76,21)(7,39,18,
    52)(9,49,64,43)(10,46)(11,88,41,29)(12,57,\\25,15)(13,72,85,37)(14,82,87,60)(16,
    56)(19,71,31,24)(22,74,99,62)(23,73,97,54)(27,53)(28,67,79,100)(33,65,63,35)\\(36,
    58,61,45)(44,55,78,92)(48,66)(51,70,80,93)(68,98,69,89)(75,95)(77,86),

  (1,23,69,8)(2,20,73,89)(3,70,64,65)(5,36,44,12)(6,30,84,50)(7,24,28,60)(9,35,96,
    93)(10,46)(11,29,41,88)(14,39,\\19,67)(15,21,45,92)(16,66)(17,40,34,81)(18,71,79,
    82)(22,74,99,62)(25,76,61,78)(27,77)(31,100,87,52)(32,51,49,33)\\(38,54,98,47)(42,
    97,68,91)(43,63,83,80)(48,56)(53,86)(55,57,90,58)(75,95),

  (1,31,36,43)(2,7,90,93)(3,98,82,15)(4,59,26,94)(5,80,8,52)(6,30,84,50)(9,20,24,
    58)(10,95,46,75)(12,83,69,87)\\(13,85)(14,25,32,68)(16,53)(18,21,70,47)(19,61,49,
    42)(22,74,99,62)(23,100,44,63)(27,56)(28,55,35,73)(33,97,67,\\78)(37,72)(38,71,45,
    64)(39,76,51,91)(48,77)(54,79,92,65)(57,96,89,60)(66,86),

  (1,36)(2,90)(3,82)(4,26)(5,8)(6,84)(7,93)(9,24)(10,46)(12,69)(14,32)(15,98)(18,
    70)(19,49)(20,58)(21,47)(22,\\99)(23,44)(25,68)(28,35)(30,50)(31,43)(33,67)(38,
    45)(39,51)(42,61)(52,80)(54,92)(55,73)(57,89)(59,94)(60,96)\\(62,74)(63,100)(64,
    71)(65,79)(75,95)(76,91)(78,97)(83,87),

     (1,42)(2,47)(3,96)(4,26)(5,76)(7,
    18)(8,91)(9,64)(11,41)(12,25)(13,85)(14,87)(15,57)(19,31)(20,38)(21,90)(22,\\
    99)(23,97)(24,71)(28,79)(29,88)(32,83)(33,63)(35,65)(36,61)(37,72)(39,52)(43,
    49)(44,78)(45,58)(51,80)(54,73)\\(55,92)(59,94)(60,82)(62,74)(67,100)(68,69)(70,
    93)(89,98),

    (1,57,42,15)(2,78,47,44)(3,31,96,19)(4,59,26,94)(5,73,76,54)(7,33,
    18,63)(8,55,91,92)(9,14,64,87)(11,29,41,88)\\(12,20,25,38)(13,72,85,37)(16,56)(17,
    34)(21,23,90,97)(22,74,99,62)(24,32,71,83)(27,53)(28,51,79,80)(35,39,65,52)\\(36,
    89,61,98)(40,81)(43,60,49,82)(45,69,58,68)(48,66)(67,70,100,93)(77,86),

  (1,79,25,93)(2,31,92,32)(3,91,24,44)(4,26)(5,64,97,60)(6,84)(7,36,65,68)(8,71,78,
    96)(9,23,82,76)(10,95,46,75)\\(11,41)(12,70,42,28)(13,72,85,37)(14,90,43,54)(15,
    51,20,100)(16,77)(17,40,34,81)(18,61,35,69)(19,55,83,47)(21,\\49,73,87)(22,74,99,
    62)(27,66)(29,88)(30,50)(33,89,52,45)(38,67,57,80)(39,58,63,98)(48,53)(56,
    86)(59,94),

    (1,80,12,100)(2,96,55,24)(3,92,71,47)(4,94,26,59)(5,87,23,43)(6,
    84)(7,89,35,58)(8,83,44,31)(9,90,60,73)(10,95,\\46,75)(11,88,41,29)(13,85)(14,97,
    49,76)(15,79,38,70)(16,86)(17,81,34,40)(18,98,65,45)(19,91,32,78)(20,93,57,28)\\(21,82,54,64)(22,99)(25,67,42,51)(27,48)(30,50)(33,61,39,68)(36,52,69,63)(37,
    72)(53,66)(56,77)(62,74),

    (1,82,61,96)(2,100,21,33)(3,42,60,36)(4,26)(5,65,91,
    28)(6,50,84,30)(7,44,70,97)(8,79,76,35)(9,69,71,25)(10,75,\\46,95)(11,88,41,
    29)(12,64,68,24)(13,72,85,37)(14,58,83,38)(15,49,89,31)(16,27)(17,34)(18,78,93,
    23)(19,57,43,98)\\(20,87,45,32)(22,99)(39,55,80,54)(40,81)(47,67,90,63)(48,86)(51,
    73,52,92)(53,56)(59,94)(62,74)(66,77),

    (1,89,42,98)(2,97,47,23)(3,43,96,49)(4,
    94,26,59)(5,55,76,92)(6,84)(7,67,18,100)(8,73,91,54)(9,32,64,83)(10,\\46)(11,29,
    41,88)(12,58,25,45)(13,72,85,37)(14,71,87,24)(15,36,57,61)(16,56)(17,34)(19,82,
    31,60)(20,68,38,69)\\(21,44,90,78)(22,62,99,74)(27,53)(28,39,79,52)(30,50)(33,70,
    63,93)(35,51,65,80)(40,81)(48,66)(75,95)(77,86).

The character table of $G^{s_6}$:\\
\begin{tabular}{c|cccccccccccccccccccccccc}
  & & & & & & & & & &10 & & & & & & & & & &20 & & & &\\\hline
$\chi_6^{(1)}$&1&1&1&1&1&1&1&1&1&1&1&1&1&1&1&1&1&1&1&1&1&1&1&1
\\$\chi_6^{(2)}$&1&-1&-1&1&-1&1&1&-1&1&1&-1&-1&-1&1&1&-1&1&1&-1&-1&1&1&1&-1
\\$\chi_6^{(3)}$&1&-1&-1&1&1&-1&-1&1&1&1&1&1&1&-1&-1&1&1&1&-1&-1&1&1&1&1
\\$\chi_6^{(4)}$&1&1&1&1&-1&-1&-1&-1&1&1&-1&-1&-1&-1&-1&-1&1&1&1&1&1&1&1&-1
\\$\chi_6^{(5)}$&1&-1&-1&1&A&-A&-A&A&1&1&A&-A&-A&A&A&A&-1&-1&1&1&-1&-1&-1&-A
\\$\chi_6^{(6)}$&1&-1&-1&1&-A&A&A&-A&1&1&-A&A&A&-A&-A&-A&-1&-1&1&1&-1&-1&-1&A
\\$\chi_6^{(7)}$&1&1&1&1&A&A&A&A&1&1&A&-A&-A&-A&-A&A&-1&-1&-1&-1&-1&-1&-1&-A
\\$\chi_6^{(8)}$&1&1&1&1&-A&-A&-A&-A&1&1&-A&A&A&A&A&-A&-1&-1&-1&-1&-1&-1&-1&A
\\$\chi_6^{(9)}$&1&A&-A&-1&-1&-A&A&1&-1&1&1&-1&1&A&-A&1&-1&1&-A&A&1&-1&1&1
\\$\chi_6^{(10)}$&1&-A&A&-1&-1&A&-A&1&-1&1&1&-1&1&-A&A&1&-1&1&A&-A&1&-1&1&1
\\$\chi_6^{(11)}$&1&A&-A&-1&1&A&-A&-1&-1&1&-1&1&-1&-A&A&-1&-1&1&-A&A&1&-1&1&-1
\\$\chi_6^{(12)}$&1&-A&A&-1&1&-A&A&-1&-1&1&-1&1&-1&A&-A&-1&-1&1&A&-A&1&-1&1&-1
\\$\chi_6^{(13)}$&1&A&-A&-1&A&-1&1&-A&-1&1&-A&-A&A&-1&1&-A&1&-1&A&-A&-1&1&-1&A
\\$\chi_6^{(14)}$&1&-A&A&-1&-A&-1&1&A&-1&1&A&A&-A&-1&1&A&1&-1&-A&A&-1&1&-1&-A
\\$\chi_6^{(15)}$&1&A&-A&-1&-A&1&-1&A&-1&1&A&A&-A&1&-1&A&1&-1&A&-A&-1&1&-1&-A
\\$\chi_6^{(16)}$&1&-A&A&-1&A&1&-1&-A&-1&1&-A&-A&A&1&-1&-A&1&-1&-A&A&-1&1&-1&A
\\$\chi_6^{(17)}$&2&.&.&-2&.&.&.&.&2&-2&.&.&.&.&.&.&-2&2&.&.&-2&2&-2&.
\\$\chi_6^{(18)}$&2&.&.&-2&.&.&.&.&2&-2&.&.&.&.&.&.&2&-2&.&.&2&-2&2&.
\\$\chi_6^{(19)}$&2&.&.&2&.&.&.&.&-2&-2&.&.&.&.&.&.&-2&-2&.&.&2&2&2&.
\\$\chi_6^{(20)}$&2&.&.&2&.&.&.&.&-2&-2&.&.&.&.&.&.&2&2&.&.&-2&-2&-2&.
\\$\chi_6^{(21)}$&4&.&.&.&.&.&.&-2&.&.&-2&.&-2&.&.&2&.&.&.&.&.&.&4&-2
\\$\chi_6^{(22)}$&4&.&.&.&.&.&.&2&.&.&2&.&2&.&.&-2&.&.&.&.&.&.&4&2
\\$\chi_6^{(23)}$&4&.&.&.&.&.&.&.&.&.&.&.&.&.&.&.&.&.&.&.&-4&.&4&.
\\$\chi_6^{(24)}$&4&.&.&.&.&.&.&.&.&.&.&.&.&.&.&.&.&.&.&.&4&.&-4&.
\\$\chi_6^{(25)}$&4&.&.&.&.&.&.&B&.&.&B&.&-B&.&.&-B&.&.&.&.&.&.&-4&-B
\\$\chi_6^{(26)}$&4&.&.&.&.&.&.&-B&.&.&-B&.&B&.&.&B&.&.&.&.&.&.&-4&B
\\$\chi_6^{(27)}$&4&.&.&.&.&.&.&.&.&2&.&.&.&.&.&C&.&-B&.&.&.&.&D&.
\\$\chi_6^{(28)}$&4&.&.&.&.&.&.&.&.&2&.&.&.&.&.&/C&.&B&.&.&.&.&-D&.
\\$\chi_6^{(29)}$&4&.&.&.&.&.&.&.&.&2&.&.&.&.&.&-/C&.&B&.&.&.&.&-D&.
\\$\chi_6^{(30)}$&4&.&.&.&.&.&.&.&.&2&.&.&.&.&.&-C&.&-B&.&.&.&.&D&.
\\$\chi_6^{(31)}$&4&.&.&.&.&.&.&C&.&-2&-C&.&-/C&.&.&.&.&B&.&.&.&.&D&/C
\\$\chi_6^{(32)}$&4&.&.&.&.&.&.&/C&.&-2&-/C&.&-C&.&.&.&.&-B&.&.&.&.&-D&C
\\$\chi_6^{(33)}$&4&.&.&.&.&.&.&-/C&.&-2&/C&.&C&.&.&.&.&-B&.&.&.&.&-D&-C
\\$\chi_6^{(34)}$&4&.&.&.&.&.&.&-C&.&-2&C&.&/C&.&.&.&.&B&.&.&.&.&D&-/C
\end{tabular}

\begin{tabular}{c|cccccccccc}
  & & & & & &30 & & & &\\\hline
$\chi_6^{(1)}$&1&1&1&1&1&1&1&1&1&1
\\$\chi_6^{(2)}$&1&-1&1&1&1&1&-1&-1&1&1
\\$\chi_6^{(3)}$&1&1&1&1&1&1&1&1&1&1
\\$\chi_6^{(4)}$&1&-1&1&1&1&1&-1&-1&1&1
\\$\chi_6^{(5)}$&1&-A&-1&1&1&-1&A&-A&1&-1
\\$\chi_6^{(6)}$&1&A&-1&1&1&-1&-A&A&1&-1
\\$\chi_6^{(7)}$&1&-A&-1&1&1&-1&A&-A&1&-1
\\$\chi_6^{(8)}$&1&A&-1&1&1&-1&-A&A&1&-1
\\$\chi_6^{(9)}$&1&1&1&1&1&1&1&1&1&1
\\$\chi_6^{(10)}$&1&1&1&1&1&1&1&1&1&1
\\$\chi_6^{(11)}$&1&-1&1&1&1&1&-1&-1&1&1
\\$\chi_6^{(12)}$&1&-1&1&1&1&1&-1&-1&1&1
\\$\chi_6^{(13)}$&1&A&-1&1&1&-1&-A&A&1&-1
\\$\chi_6^{(14)}$&1&-A&-1&1&1&-1&A&-A&1&-1
\\$\chi_6^{(15)}$&1&-A&-1&1&1&-1&A&-A&1&-1
\\$\chi_6^{(16)}$&1&A&-1&1&1&-1&-A&A&1&-1
\\$\chi_6^{(17)}$&2&.&2&2&2&-2&.&.&-2&-2
\\$\chi_6^{(18)}$&2&.&-2&2&2&2&.&.&-2&2
\\$\chi_6^{(19)}$&2&.&-2&2&2&2&.&.&-2&2
\\$\chi_6^{(20)}$&2&.&2&2&2&-2&.&.&-2&-2
\\$\chi_6^{(21)}$&.&2&.&-4&4&4&2&2&.&-4
\\$\chi_6^{(22)}$&.&-2&.&-4&4&4&-2&-2&.&-4
\\$\chi_6^{(23)}$&-4&.&.&4&4&4&.&.&.&4
\\$\chi_6^{(24)}$&-4&.&.&4&4&-4&.&.&.&-4
\\$\chi_6^{(25)}$&.&B&.&-4&4&-4&-B&B&.&4
\\$\chi_6^{(26)}$&.&-B&.&-4&4&-4&B&-B&.&4
\\$\chi_6^{(27)}$&.&/C&B&.&-4&-D&-C&-/C&-2&.
\\$\chi_6^{(28)}$&.&C&-B&.&-4&D&-/C&-C&-2&.
\\$\chi_6^{(29)}$&.&-C&-B&.&-4&D&/C&C&-2&.
\\$\chi_6^{(30)}$&.&-/C&B&.&-4&-D&C&/C&-2&.
\\$\chi_6^{(31)}$&.&.&-B&.&-4&-D&.&.&2&.
\\$\chi_6^{(32)}$&.&.&B&.&-4&D&.&.&2&.
\\$\chi_6^{(33)}$&.&.&B&.&-4&D&.&.&2&.
\\$\chi_6^{(34)}$&.&.&-B&.&-4&-D&.&.&2&.
\end{tabular}

\noindent \noindent where   A = -E(4)
  = -ER(-1) = -i;
B = -2*E(4)
  = -2*ER(-1) = -2i;
C = -2-2*E(4)
  = -2-2*ER(-1) = -2-2i;
D = 4*E(4)
  = 4*ER(-1) = 4i.

The generators of $G^{s_7}$ are:\\
(  1, 14)(  2, 63)(  3, 45)(  4, 22)(  5, 28)(  6, 50)(  7, 78)(  8,
35)
    (  9, 15)( 10, 75)( 11, 37)( 12, 43)( 13, 29)( 16, 27)( 17, 34)( 18, 44)
    ( 19, 68)( 20, 60)( 21, 67)( 23, 70)( 24, 98)( 25, 49)( 26, 99)( 30, 84)
    ( 31, 69)( 32, 36)( 33, 47)( 38, 82)( 39, 55)( 40, 81)( 41, 72)( 42, 87)
    ( 46, 95)( 48, 77)( 51, 73)( 52, 92)( 53, 56)( 54, 80)( 57, 64)( 58, 96)
    ( 59, 74)( 61, 83)( 62, 94)( 65, 91)( 66, 86)( 71, 89)( 76, 79)( 85, 88)
    ( 90,100)( 93, 97), (  3,  6)(  4,  9)(  5,  8)( 10, 20)( 11, 16)( 12, 19)
    ( 15, 22)( 18, 21)( 23, 33)( 24, 26)( 25, 31)( 27, 37)( 28, 35)( 30, 38)
    ( 39, 51)( 40, 66)( 41, 53)( 42, 61)( 43, 68)( 44, 67)( 45, 50)( 46, 58)
    ( 47, 70)( 49, 69)( 52, 76)( 54, 65)( 55, 73)( 56, 72)( 57, 59)( 60, 75)
    ( 62, 71)( 64, 74)( 79, 92)( 80, 91)( 81, 86)( 82, 84)( 83, 87)( 89, 94)
    ( 95, 96)( 98, 99), (  3, 60)(  4, 64)(  5, 49)(  6, 75)(  8, 69)(  9, 74)
    ( 10, 50)( 11, 47)( 12, 51)( 13, 17)( 15, 59)( 16, 70)( 18, 41)( 19, 39)
    ( 20, 45)( 21, 53)( 22, 57)( 23, 27)( 25, 28)( 29, 34)( 30, 38)( 31, 35)
    ( 33, 37)( 40, 42)( 43, 73)( 44, 72)( 46, 58)( 52, 54)( 55, 68)( 56, 67)
    ( 61, 66)( 65, 76)( 79, 91)( 80, 92)( 81, 87)( 82, 84)( 83, 86)( 90, 93)
    ( 95, 96)( 97,100), (  2,  7)(  4,  9)(  5, 19)(  8, 12)( 11, 18)( 13, 17)
    ( 15, 22)( 16, 21)( 23, 56)( 24, 71)( 25, 55)( 26, 62)( 27, 67)( 28, 68)
    ( 29, 34)( 30, 58)( 31, 73)( 33, 72)( 35, 43)( 37, 44)( 38, 46)( 39, 49)
    ( 40, 42)( 41, 47)( 51, 69)( 52, 65)( 53, 70)( 54, 76)( 57, 59)( 61, 66)
    ( 63, 78)( 64, 74)( 79, 80)( 81, 87)( 82, 95)( 83, 86)( 84, 96)( 89, 98)
    ( 91, 92)( 94, 99), (  2, 11)(  3,  5)(  6, 12)(  7, 21)(  8, 20)( 10, 19)
    ( 13, 16)( 17, 18)( 24, 66)( 25, 39)( 26, 71)( 27, 29)( 28, 45)( 30, 76)
    ( 31, 51)( 32, 57)( 33, 56)( 34, 44)( 35, 60)( 36, 64)( 37, 63)( 40, 61)
    ( 42, 62)( 43, 50)( 46, 65)( 47, 53)( 48, 59)( 49, 55)( 67, 78)( 68, 75)
    ( 69, 73)( 74, 77)( 79, 84)( 81, 83)( 85, 93)( 86, 98)( 87, 94)( 88, 97)
    ( 89, 99)( 91, 95), (  1,  2)(  3, 12)(  4, 53)(  5, 52)(  6, 51)(  7, 77)
    (  8, 54)(  9, 21)( 10, 55)( 11, 66)( 13, 62)( 14, 63)( 15, 67)( 16, 61)
    ( 17, 71)( 18, 74)( 19, 60)( 20, 68)( 22, 56)( 23, 87)( 24, 93)( 25, 79)
    ( 26, 90)( 27, 83)( 28, 92)( 29, 94)( 30, 82)( 31, 91)( 32, 88)( 33, 81)
    ( 34, 89)( 35, 80)( 36, 85)( 37, 86)( 38, 84)( 39, 75)( 40, 47)( 41, 64)
    ( 42, 70)( 43, 45)( 44, 59)( 46, 58)( 48, 78)( 49, 76)( 50, 73)( 57, 72)
    ( 65, 69)( 95, 96)( 97, 98)( 99,100).

The representatives of conjugacy classes of   $G^{s_7}$ are:\\
 (1),
 (1,14)(2,63)(3,45)(4,22)(5,28)(6,50)(7,78)(8,35)(9,15)(10,75)(11,37)(12,
    43)(13,29)(16,27)(17,34)(18,44)(19,\\68)(20,60)(21,67)(23,70)(24,98)(25,49)(26,
    99)(30,84)(31,69)(32,36)(33,47)(38,82)(39,55)(40,81)(41,72)(42,87)\\(46,95)(48,
    77)(51,73)(52,92)(53,56)(54,80)(57,64)(58,96)(59,74)(61,83)(62,94)(65,91)(66,
    86)(71,89)(76,79)(85,\\88)(90,100)(93,97),

     (1,7,61,41)(2,4)(3,79,68,75)(5,80,84,
    73)(6,20,91,43)(8,92,46,39)(9,85,71,18)(10,45,76,19)(11,74,70,77)(12,50,\\60,65)(13,66,93,62)(14,78,83,72)(15,88,89,44)(16,42,17,64)(21,40)(22,63)(23,48,37,
    59)(24,90,36,47)(25,49)(26,\\53)(27,87,34,57)(28,54,30,51)(29,86,97,94)(32,33,98,
    100)(35,52,95,55)(56,99)(58,96)(67,81),

     (1,78,61,72)(2,22)(3,76,68,10)(4,
    63)(5,54,84,51)(6,60,91,12)(7,83,41,14)(8,52,46,55)(9,88,71,44)(11,59,70,48)\\(13,
    86,93,94)(15,85,89,18)(16,87,17,57)(19,75,45,79)(20,65,43,50)(21,81)(23,77,37,
    74)(24,100,36,33)(26,56)(27,\\42,34,64)(28,80,30,73)(29,66,97,62)(31,69)(32,47,98,
    90)(35,92,95,39)(38,82)(40,67)(53,99),

    (1,61)(3,68)(5,84)(6,91)(7,41)(8,46)(9,
    71)(10,76)(11,70)(12,60)(13,93)(14,83)(15,89)(16,17)(18,85)(19,45)(20,\\43)(23,
    37)(24,36)(27,34)(28,30)(29,97)(32,98)(33,100)(35,95)(39,92)(42,64)(44,88)(47,
    90)(48,59)(50,65)(51,54)\\(52,55)(57,87)(62,66)(72,78)(73,80)(74,77)(75,79)(86,
    94),

     (1,83)(2,63)(3,19)(4,22)(5,30)(6,65)(7,72)(8,95)(9,89)(10,79)(11,23)(12,
    20)(13,97)(14,61)(15,71)(16,34)(17,\\27)(18,88)(21,67)(24,32)(25,49)(26,99)(28,
    84)(29,93)(31,69)(33,90)(35,46)(36,98)(37,70)(38,82)(39,52)(40,81)\\(41,78)(42,
    57)(43,60)(44,85)(45,68)(47,100)(48,74)(50,91)(51,80)(53,56)(54,73)(55,92)(58,
    96)(59,77)(62,86)(64,\\87)(66,94)(75,76),

    (1,62,61,66)(2,21)(3,55,68,52)(4,
    40)(5,6,84,91)(7,13,41,93)(8,10,46,76)(9,74,71,77)(11,85,70,18)(12,54,60,51)\\(14,
    94,83,86)(15,59,89,48)(16,47,17,90)(19,92,45,39)(20,73,43,80)(22,81)(23,44,37,
    88)(24,64,36,42)(25,58)(27,\\33,34,100)(28,50,30,65)(29,72,97,78)(31,82)(32,87,98,
    57)(35,75,95,79)(38,69)(49,96)(63,67),

    (1,94,61,86)(2,67)(3,39,68,92)(4,81)(5,
    50,84,65)(6,30,91,28)(7,29,41,97)(8,75,46,79)(9,59,71,48)(10,95,76,35)\\(11,88,70,
    44)(12,80,60,73)(13,72,93,78)(14,62,83,66)(15,74,89,77)(16,33,17,100)(18,37,85,
    23)(19,52,45,55)(20,51,\\43,54)(21,63)(22,40)(24,57,36,87)(25,96)(26,99)(27,47,34,
    90)(31,38)(32,42,98,64)(49,58)(53,56)(69,82),

  (1,74,61,77)(2,21)(3,54,68,51)(4,40)(5,10,84,76)(6,46,91,8)(7,11,41,70)(9,62,71,
    66)(12,55,60,52)(13,85,93,18)\\(14,59,83,48)(15,94,89,86)(16,47,17,90)(19,73,45,
    80)(20,92,43,39)(22,81)(23,78,37,72)(24,42,36,64)(25,49)(27,33,\\34,100)(28,75,30,
    79)(29,88,97,44)(31,69)(32,57,98,87)(35,50,95,65)(38,82)(58,96)(63,67),

  (1,59,61,48)(2,67)(3,80,68,73)(4,81)(5,75,84,79)(6,95,91,35)(7,37,41,23)(8,50,46,
    65)(9,94,71,86)(10,30,76,28)\\(11,72,70,78)(12,39,60,92)(13,88,93,44)(14,74,83,
    77)(15,62,89,66)(16,33,17,100)(18,29,85,97)(19,51,45,54)(20,52,\\43,55)(21,63)(22,
    40)(24,87,36,57)(26,99)(27,47,34,90)(32,64,98,42)(53,56),

  (1,9)(3,12)(5,8)(6,10)(7,18)(11,13)(14,15)(19,20)(23,97)(24,36)(25,96)(28,35)(29,
    37)(30,95)(31,38)(32,98)(39,\\80)(41,85)(42,64)(43,45)(44,78)(46,84)(48,86)(49,
    58)(50,75)(51,52)(54,55)(57,87)(59,94)(60,68)(61,71)(62,74)\\(65,79)(66,77)(69,
    82)(70,93)(72,88)(73,92)(76,91)(83,89),

    (1,15)(2,63)(3,43)(4,22)(5,35)(6,
    75)(7,44)(8,28)(9,14)(10,50)(11,29)(12,45)(13,37)(16,27)(17,34)(18,78)(19,\\
    60)(20,68)(21,67)(23,93)(24,32)(25,58)(26,99)(30,46)(31,82)(33,47)(36,98)(38,
    69)(39,54)(40,81)(41,88)(42,57)\\(48,66)(49,96)(51,92)(52,73)(53,56)(55,80)(59,
    62)(61,89)(64,87)(65,76)(70,97)(71,83)(72,85)(74,94)(77,86)(79,\\91)(84,95)(90,
    100),

     (1,42,24)(2,93,13)(3,68,96)(4,62,66)(5,75,54)(6,52,35)(7,11,47)(8,50,
    92)(10,80,28)(12,82,60)(14,87,98)(16,17,\\53)(18,85,21)(19,58,45)(20,43,38)(22,94,
    86)(23,72,100)(26,74,77)(27,34,56)(29,63,97)(30,73,76)(32,57,83)(33,78,\\37)(36,
    64,61)(39,65,46)(41,90,70)(44,88,67)(48,99,59)(51,79,84)(55,91,95),

  (1,87,24,14,42,98)(2,97,13,63,93,29)(3,19,96,45,68,58)(4,94,66,22,62,86)(5,10,54,
    28,75,80)(6,92,35,50,52,8)(7,\\37,47,78,11,33)(9,15)(12,38,60,43,82,20)(16,34,53,
    27,17,56)(18,88,21,44,85,67)(23,41,100,70,72,90)(25,49)(26,59,\\77,99,74,48)(30,
    51,76,84,73,79)(31,69)(32,64,83,36,57,61)(39,91,46,55,65,95)(40,81)(71,89),

  (1,42,24)(2,85,13,18,93,21)(3,50,96,8,68,92)(4,62,66)(5,43,54,20,75,38)(6,58,35,
    19,52,45)(7,17,47,16,11,53)(10,\\82,28,12,80,60)(14,87,98)(22,94,86)(23,72,
    100)(25,69)(26,36,77,61,74,64)(27,37,56,78,34,33)(29,44,97,67,63,88)\\(30,39,76,
    46,73,65)(31,49)(32,48,83,59,57,99)(40,71)(41,90,70)(51,91,84,55,79,95)(81,89),

(1,87,24,14,42,98)(2,88,13,44,93,67)(3,6,96,35,68,52)(4,94,66,22,62,86)(5,12,
    54,60,75,82)(7,34,47,27,11,56)(8,\\19,92,45,50,58)(9,15)(10,38,28,43,80,20)(16,37,
    53,78,17,33)(18,97,21,63,85,29)(23,41,100,70,72,90)(25,31)(26,32,\\77,83,74,
    57)(30,55,76,95,73,91)(36,48,61,59,64,99)(39,79,46,51,65,84)(40,89)(49,69)(71,
    81),

    (1,13)(2,40)(3,35)(4,21)(5,20)(6,73)(7,66)(8,45)(9,70)(10,39)(11,71)(12,
    30)(14,29)(15,23)(16,24)(17,36)(18,\\74)(19,46)(22,67)(25,96)(26,53)(27,98)(28,
    60)(31,82)(32,34)(33,57)(37,89)(38,69)(41,62)(42,90)(43,84)(44,59)\\(47,64)(48,
    88)(49,58)(50,51)(52,79)(54,65)(55,75)(56,99)(61,93)(63,81)(68,95)(72,94)(76,
    92)(77,85)(78,86)(80,\\91)(83,97)(87,100),

    (1,29)(2,81)(3,8)(4,67)(5,60)(6,
    51)(7,86)(9,23)(10,55)(11,89)(12,84)(13,14)(15,70)(16,98)(17,32)(18,59)(19,\\
    95)(20,28)(21,22)(24,27)(25,58)(26,56)(30,43)(31,38)(33,64)(34,36)(35,45)(37,
    71)(39,75)(40,63)(41,94)(42,100)\\(44,74)(46,68)(47,57)(48,85)(49,96)(50,73)(52,
    76)(53,99)(54,91)(61,97)(62,72)(65,80)(66,78)(69,82)(77,88)(79,\\92)(83,93)(87,
    90),

    (1,2,40,13,64,41,62,47)(3,91,30,82,31,12,80,35)(4,85,26,17,42,11,9,7)(5,
    20,79,95,96,25,68,52)(6,10,50,75)(8,\\45,65,84,38,69,43,54)(14,63,81,29,57,72,94,
    33)(15,78,22,88,99,34,87,37)(16,61,93,24)(18,74)(19,92,28,60,76,46,\\58,49)(21,66,
    70,71,90,36,53,77)(23,89,100,32,56,48,67,86)(27,83,97,98)(39,73,55,51)(44,59),

  (1,63,40,29,64,72,62,33)(2,81,13,57,41,94,47,14)(3,65,30,38,31,43,80,8)(4,88,26,
    34,42,37,9,78)(5,60,79,46,96,\\49,68,92)(6,75,50,10)(7,22,85,99,17,87,11,15)(12,
    54,35,45,91,84,82,69)(16,83,93,98)(18,59)(19,52,28,20,76,95,58,\\25)(21,86,70,89,
    90,32,53,48)(23,71,100,36,56,77,67,66)(24,27,61,97)(39,51,55,73)(44,74),

  (1,93,64,11,66,41,24,2,36,85)(3,12,96,6,5,91,31,52,51,95)(4,16,42,47,40,70,61,53,
    62,13)(7,77,90,74,21,71,17,9,\\18,26)(8,79,25,54,39,84,20,68,38,75)(10,35,76,49,
    80,55,30,60,19,82)(14,97,57,37,86,72,98,63,32,88)(15,44,99,78,\\48,100,59,67,89,
    34)(22,27,87,33,81,23,83,56,94,29)(28,65,69,92,73,46,45,43,58,50),

  (1,97,64,37,66,72,24,63,36,88)(2,32,85,14,93,57,11,86,41,98)(3,43,96,50,5,65,31,
    92,51,46)(4,27,42,33,40,23,61,\\56,62,29)(6,28,91,69,52,73,95,45,12,58)(7,48,90,
    59,21,89,17,15,18,99)(8,76,25,80,39,30,20,19,38,10)(9,44,26,78,\\77,100,74,67,71,
    34)(13,22,16,87,47,81,70,83,53,94)(35,79,49,54,55,84,60,68,82,75),

  (1,64,66,24,36)(2,85,93,11,41)(3,96,5,31,51)(4,42,40,61,62)(6,91,52,95,12)(7,90,
    21,17,18)(8,25,39,20,38)(9,26,\\77,74,71)(10,76,80,30,19)(13,16,47,70,53)(14,57,
    86,98,32)(15,99,48,59,89)(22,87,81,83,94)(23,56,29,27,33)(28,69,\\73,45,58)(34,44,
    78,100,67)(35,49,55,60,82)(37,72,63,88,97)(43,50,65,92,46)(54,84,68,75,79),

  (1,57,66,98,36,14,64,86,24,32)(2,88,93,37,41,63,85,97,11,72)(3,58,5,69,51,45,96,
    28,31,73)(4,87,40,83,62,22,\\42,81,61,94)(6,65,52,46,12,50,91,92,95,43)(7,100,21,
    34,18,78,90,67,17,44)(8,49,39,60,38,35,25,55,20,82)(9,99,\\77,59,71,15,26,48,74,
    89)(10,79,80,84,19,75,76,54,30,68)(13,27,47,23,53,29,16,33,70,56),

  (1,85,77,13,61,18,74,93)(2,40,21,4)(3,79,51,30,68,75,54,28)(5,45,76,73,84,19,10,
    80)(6,92,8,20,91,39,46,43)\\(7,66,70,71,41,62,11,9)(12,50,52,35,60,65,55,95)(14,
    88,48,29,83,44,59,97)(15,78,86,23,89,72,94,37)(16,64,90,\\36,17,42,47,24)(22,63,
    81,67)(25,82,49,38)(26,53)(27,57,100,32,34,87,33,98)(31,96,69,58)(56,99),

  (1,88,77,29,61,44,74,97)(2,81,21,22)(3,76,51,84,68,10,54,5)(4,63,40,67)(6,52,8,60,
    91,55,46,12)(7,86,70,89,41,\\94,11,15)(9,78,66,23,71,72,62,37)(13,83,18,59,93,14,
    85,48)(16,57,90,32,17,87,47,98)(19,75,80,28,45,79,73,30)\\(20,65,39,95,43,50,92,
    35)(24,27,64,100,36,34,42,33)(25,38,49,82)(26,56)(31,58,69,96)(53,99).

The character table of $G^{s_7}$:\\
\begin{tabular}{c|cccccccccccccccccccccccc}
  & & & & & & & & & & 10& & & & & & & & & & 20& & & &\\\hline

$\chi_7^{(1)}$&1&1&1&1&1&1&1&1&1&1&1&1&1&1&1&1&1&1&1&1&1&1&1&1
\\$\chi_7^{(2)}$&1&-1&-1&1&1&-1&-1&1&1&-1&-1&1&1&-1&-1&1&1&-1&1&-1&1&-1&1&-1
\\$\chi_7^{(3)}$&1&-1&-1&1&1&-1&1&-1&1&-1&1&-1&1&-1&1&-1&-1&1&-1&1&-1&1&1&-1
\\$\chi_7^{(4)}$&1&-1&1&-1&1&-1&-1&1&1&-1&-1&1&1&-1&-1&1&-1&1&-1&1&-1&1&1&-1
\\$\chi_7^{(5)}$&1&-1&1&-1&1&-1&1&-1&1&-1&1&-1&1&-1&1&-1&1&-1&1&-1&1&-1&1&-1
\\$\chi_7^{(6)}$&1&1&-1&-1&1&1&-1&-1&1&1&-1&-1&1&1&-1&-1&1&1&1&1&1&1&1&1
\\$\chi_7^{(7)}$&1&1&-1&-1&1&1&1&1&1&1&1&1&1&1&1&1&-1&-1&-1&-1&-1&-1&1&1
\\$\chi_7^{(8)}$&1&1&1&1&1&1&-1&-1&1&1&-1&-1&1&1&-1&-1&-1&-1&-1&-1&-1&-1&1&1
\\$\chi_7^{(9)}$&9&-9&-1&1&1&-1&-1&1&1&-1&3&-3&.&.&.&.&1&-1&-1&1&1&-1&-1&1
\\$\chi_7^{(10)}$&9&-9&-1&1&1&-1&1&-1&1&-1&-3&3&.&.&.&.&-1&1&1&-1&-1&1&-1&1
\\$\chi_7^{(11)}$&9&-9&1&-1&1&-1&-1&1&1&-1&3&-3&.&.&.&.&-1&1&1&-1&-1&1&-1&1
\\$\chi_7^{(12)}$&9&-9&1&-1&1&-1&1&-1&1&-1&-3&3&.&.&.&.&1&-1&-1&1&1&-1&-1&1
\\$\chi_7^{(13)}$&9&9&-1&-1&1&1&-1&-1&1&1&3&3&.&.&.&.&1&1&-1&-1&1&1&-1&-1
\\$\chi_7^{(14)}$&9&9&-1&-1&1&1&1&1&1&1&-3&-3&.&.&.&.&-1&-1&1&1&-1&-1&-1&-1
\\$\chi_7^{(15)}$&9&9&1&1&1&1&-1&-1&1&1&3&3&.&.&.&.&-1&-1&1&1&-1&-1&-1&-1
\\$\chi_7^{(16)}$&9&9&1&1&1&1&1&1&1&1&-3&-3&.&.&.&.&1&1&-1&-1&1&1&-1&-1
\\$\chi_7^{(17)}$&10&10&.&.&2&2&2&2&-2&-2&2&2&1&1&-1&-1&.&.&.&.&.&.&.&.
\\$\chi_7^{(18)}$&10&10&.&.&2&2&-2&-2&-2&-2&-2&-2&1&1&1&1&.&.&.&.&.&.&.&.
\\$\chi_7^{(19)}$&10&-10&.&.&2&-2&2&-2&-2&2&2&-2&1&-1&-1&1&.&.&.&.&.&.&.&.
\\$\chi_7^{(20)}$&10&-10&.&.&2&-2&-2&2&-2&2&-2&2&1&-1&1&-1&.&.&.&.&.&.&.&.
\\$\chi_7^{(21)}$&16&16&.&.&.&.&.&.&.&.&.&.&-2&-2&.&.&-4&-4&.&.&1&1&1&1
\\$\chi_7^{(22)}$&16&16&.&.&.&.&.&.&.&.&.&.&-2&-2&.&.&4&4&.&.&-1&-1&1&1
\\$\chi_7^{(23)}$&16&-16&.&.&.&.&.&.&.&.&.&.&-2&2&.&.&-4&4&.&.&1&-1&1&-1
\\$\chi_7^{(24)}$&16&-16&.&.&.&.&.&.&.&.&.&.&-2&2&.&.&4&-4&.&.&-1&1&1&-1
\\$\chi_7^{(25)}$&20&20&.&.&-4&-4&.&.&.&.&.&.&2&2&.&.&.&.&.&.&.&.&.&.
\\$\chi_7^{(26)}$&20&-20&.&.&-4&4&.&.&.&.&.&.&2&-2&.&.&.&.&.&.&.&.&.&.
\end{tabular}

\begin{tabular}{c|cc}
  & &\\\hline
$\chi_7^{(1)}$&1&1
\\$\chi_7^{(2)}$&-1&1
\\$\chi_7^{(3)}$&-1&1
\\$\chi_7^{(4)}$&1&-1
\\$\chi_7^{(5)}$&1&-1
\\$\chi_7^{(6)}$&-1&-1
\\$\chi_7^{(7)}$&-1&-1
\\$\chi_7^{(8)}$&1&1
\\$\chi_7^{(9)}$&1&-1
\\$\chi_7^{(10)}$&1&-1
\\$\chi_7^{(11)}$&-1&1
\\$\chi_7^{(12)}$&-1&1
\\$\chi_7^{(13)}$&1&1
\\$\chi_7^{(14)}$&1&1
\\$\chi_7^{(15)}$&-1&-1
\\$\chi_7^{(16)}$&-1&-1
\\$\chi_7^{(17)}$&.&.
\\$\chi_7^{(18)}$&.&.
\\$\chi_7^{(19)}$&.&.
\\$\chi_7^{(20)}$&.&.
\\$\chi_7^{(21)}$&.&.
\\$\chi_7^{(22)}$&.&.
\\$\chi_7^{(23)}$&.&.
\\$\chi_7^{(24)}$&.&.
\\$\chi_7^{(25)}$&.&.
\\$\chi_7^{(26)}$&.&.
\end{tabular}

The generators of $G^{s_8}$ are:\\
 (  1, 16, 61, 41, 26, 85, 62, 47, 66, 21)(  2, 36,  7, 71, 11, 64, 93, 74, 18, 24)
    (  3, 43, 69, 91, 30,  6, 58, 92, 51, 35)(  4, 17, 42, 90, 40, 13,  9, 53, 77,
      70)(  5, 20, 19, 49, 79, 95, 75, 82, 80, 39)(  8, 45, 12, 31, 65, 84, 50, 96,
      52, 73)( 10, 38, 54, 55, 28, 60, 68, 25, 76, 46)( 14, 27, 83, 72, 99, 88, 94,
      33, 86, 67)( 15, 56, 48, 23, 22, 34, 87,100, 81, 29)( 32, 78, 89, 37, 57, 97,
      59, 44, 98, 63), (  1, 14)(  2, 63)(  3, 45)(  4, 22)(  5, 28)(  6, 50)
    (  7, 78)(  8, 35)(  9, 15)( 10, 75)( 11, 37)( 12, 43)( 13, 29)( 16, 27)
    ( 17, 34)( 18, 44)( 19, 68)( 20, 60)( 21, 67)( 23, 70)( 24, 98)( 25, 49)
    ( 26, 99)( 30, 84)( 31, 69)( 32, 36)( 33, 47)( 38, 82)( 39, 55)( 40, 81)
    ( 41, 72)( 42, 87)( 46, 95)( 48, 77)( 51, 73)( 52, 92)( 53, 56)( 54, 80)
    ( 57, 64)( 58, 96)( 59, 74)( 61, 83)( 62, 94)( 65, 91)( 66, 86)( 71, 89)
    ( 76, 79)( 85, 88)( 90,100)( 93, 97)

The representatives of conjugacy classes of   $G^{s_8}$ are:\\
 (1),
  (1,14)(2,63)(3,45)(4,22)(5,28)(6,50)(7,78)(8,35)(9,15)(10,75)(11,37)(12,
    43)(13,29)(16,27)(17,34)(18,44)(19,\\68)(20,60)(21,67)(23,70)(24,98)(25,49)(26,
    99)(30,84)(31,69)(32,36)(33,47)(38,82)(39,55)(40,81)(41,72)(42,87)\\(46,95)(48,
    77)(51,73)(52,92)(53,56)(54,80)(57,64)(58,96)(59,74)(61,83)(62,94)(65,91)(66,
    86)(71,89)(76,79)(85,\\88)(90,100)(93,97),

    (1,16,61,41,26,85,62,47,66,21)(2,36,
    7,71,11,64,93,74,18,24)(3,43,69,91,30,6,58,92,51,35)(4,17,42,90,40,13,9,\\53,77,
    70)(5,20,19,49,79,95,75,82,80,39)(8,45,12,31,65,84,50,96,52,73)(10,38,54,55,28,
    60,68,25,76,46)(14,27,83,\\72,99,88,94,33,86,67)(15,56,48,23,22,34,87,100,81,
    29)(32,78,89,37,57,97,59,44,98,63),

    (1,21,66,47,62,85,26,41,61,16)(2,24,18,74,
    93,64,11,71,7,36)(3,35,51,92,58,6,30,91,69,43)(4,70,77,53,9,13,40,\\90,42,17)(5,
    39,80,82,75,95,79,49,19,20)(8,73,52,96,50,84,65,31,12,45)(10,46,76,25,68,60,28,
    55,54,38)(14,67,86,\\33,94,88,99,72,83,27)(15,29,81,100,87,34,22,23,48,56)(32,63,
    98,44,59,97,57,37,89,78),

    (1,26,66,61,62)(2,11,18,7,93)(3,30,51,69,58)(4,40,
    77,42,9)(5,79,80,19,75)(6,35,91,92,43)(8,65,52,12,50)(10,28,\\76,54,68)(13,70,90,
    53,17)(14,99,86,83,94)(15,22,81,48,87)(16,85,21,41,47)(20,95,39,49,82)(23,100,
    56,34,29)(24,\\71,74,36,64)(25,38,60,46,55)(27,88,67,72,33)(31,96,45,84,73)(32,57,
    98,89,59)(37,44,78,97,63),

     (1,27,61,72,26,88,62,33,66,67)(2,32,7,89,11,57,93,
    59,18,98)(3,12,69,65,30,50,58,52,51,8)(4,34,42,100,40,29,9,\\56,77,23)(5,60,19,25,
    79,46,75,38,80,55)(6,96,92,73,35,45,43,31,91,84)(10,82,54,39,28,20,68,49,76,
    95)(13,15,53,\\48,70,22,17,87,90,81)(14,16,83,41,99,85,94,47,86,21)(24,63,36,78,
    71,37,64,97,74,44),

    (1,33,26,27,66,88,61,67,62,72)(2,59,11,32,18,57,7,98,93,
    89)(3,52,30,12,51,50,69,8,58,65)(4,56,40,34,77,29,42,\\23,9,100)(5,38,79,60,80,46,
    19,55,75,25)(6,31,35,96,91,45,92,84,43,73)(10,49,28,82,76,20,54,95,68,39)(13,87,
    70,\\15,90,22,53,81,17,48)(14,47,99,16,86,85,83,21,94,41)(24,97,71,63,74,37,36,44,
    64,78),

     (1,41,62,21,61,85,66,16,26,47)(2,71,93,24,7,64,18,36,11,74)(3,91,58,
    35,69,6,51,43,30,92)(4,90,9,70,42,13,77,\\17,40,53)(5,49,75,39,19,95,80,20,79,
    82)(8,31,50,73,12,84,52,45,65,96)(10,55,68,46,54,60,76,38,28,25)(14,72,94,\\67,83,
    88,86,27,99,33)(15,23,87,29,48,34,81,56,22,100)(32,37,59,63,89,97,98,78,57,44),

    (1,47,26,16,66,85,61,21,62,41)(2,74,11,36,18,64,7,24,93,71)(3,92,30,43,51,6,
    69,35,58,91)(4,53,40,17,77,13,42,\\70,9,90)(5,82,79,20,80,95,19,39,75,49)(8,96,65,
    45,52,84,12,73,50,31)(10,25,28,38,76,60,54,46,68,55)(14,33,99,27,\\86,88,83,67,94,
    72)(15,100,22,56,81,34,48,29,87,23)(32,44,57,78,98,97,89,63,59,37),

  (1,61,26,62,66)(2,7,11,93,18)(3,69,30,58,51)(4,42,40,9,77)(5,19,79,75,80)(6,92,35,
    43,91)(8,12,65,50,52)(10,54,\\28,68,76)(13,53,70,17,90)(14,83,99,94,86)(15,48,22,
    87,81)(16,41,85,47,21)(20,49,95,82,39)(23,34,100,29,56)(24,36,\\71,64,74)(25,46,
    38,55,60)(27,72,88,33,67)(31,84,96,73,45)(32,89,57,59,98)(37,97,44,63,78),

  (1,62,61,66,26)(2,93,7,18,11)(3,58,69,51,30)(4,9,42,77,40)(5,75,19,80,79)(6,43,92,
    91,35)(8,50,12,52,65)(10,68,\\54,76,28)(13,17,53,90,70)(14,94,83,86,99)(15,87,48,
    81,22)(16,47,41,21,85)(20,82,49,39,95)(23,29,34,56,100)(24,\\64,36,74,71)(25,55,
    46,60,38)(27,33,72,67,88)(31,73,84,45,96)(32,59,89,98,57)(37,63,97,78,44),

  (1,66,62,26,61)(2,18,93,11,7)(3,51,58,30,69)(4,77,9,40,42)(5,80,75,79,19)(6,91,43,
    35,92)(8,52,50,65,12)(10,76,\\68,28,54)(13,90,17,70,53)(14,86,94,99,83)(15,81,87,
    22,48)(16,21,47,85,41)(20,39,82,95,49)(23,56,29,100,34)(24,\\74,64,71,36)(25,60,
    55,38,46)(27,67,33,88,72)(31,45,73,96,84)(32,98,59,57,89)(37,78,63,44,97),

  (1,67,66,33,62,88,26,72,61,27)(2,98,18,59,93,57,11,89,7,32)(3,8,51,52,58,50,30,65,
    69,12)(4,23,77,56,9,29,40,\\100,42,34)(5,55,80,38,75,46,79,25,19,60)(6,84,91,31,
    43,45,35,73,92,96)(10,95,76,49,68,20,28,39,54,82)(13,81,90,\\87,17,22,70,48,53,
    15)(14,21,86,47,94,85,99,41,83,16)(24,44,74,97,64,37,71,78,36,63),

  (1,72,62,67,61,88,66,27,26,33)(2,89,93,98,7,57,18,32,11,59)(3,65,58,8,69,50,51,12,
    30,52)(4,100,9,23,42,29,77,\\34,40,56)(5,25,75,55,19,46,80,60,79,38)(6,73,43,84,
    92,45,91,96,35,31)(10,39,68,95,54,20,76,82,28,49)(13,48,17,\\81,53,22,90,15,70,
    87)(14,41,94,21,83,85,86,16,99,47)(24,78,64,44,36,37,74,63,71,97),

  (1,83,26,94,66,14,61,99,62,86)(2,78,11,97,18,63,7,37,93,44)(3,31,30,96,51,45,69,
    84,58,73)(4,87,40,15,77,22,42,\\81,9,48)(5,68,79,10,80,28,19,76,75,54)(6,52,35,12,
    91,50,92,8,43,65)(13,56,70,34,90,29,53,23,17,100)(16,72,85,33,\\21,27,41,88,47,
    67)(20,25,95,38,39,60,49,46,82,55)(24,32,71,57,74,98,36,89,64,59),

  (1,85)(2,64)(3,6)(4,13)(5,95)(7,74)(8,84)(9,17)(10,60)(11,24)(12,96)(14,88)(15,
    34)(16,62)(18,71)(19,82)(20,75)\\(21,26)(22,29)(23,81)(25,54)(27,94)(28,46)(30,
    35)(31,52)(32,97)(33,83)(36,93)(37,98)(38,68)(39,79)(40,70)(41,66)\\(42,53)(43,
    58)(44,89)(45,50)(47,61)(48,100)(49,80)(51,91)(55,76)(56,87)(57,63)(59,78)(65,
    73)(67,99)(69,92)(72,\\86)(77,90),

    (1,86,62,99,61,14,66,94,26,83)(2,44,93,37,7,
    63,18,97,11,78)(3,73,58,84,69,45,51,96,30,31)(4,48,9,81,42,22,77,15,\\40,87)(5,54,
    75,76,19,28,80,10,79,68)(6,65,43,8,92,50,91,12,35,52)(13,100,17,23,53,29,90,34,
    70,56)(16,67,47,88,41,\\27,21,33,85,72)(20,55,82,46,49,60,39,38,95,25)(24,59,64,
    89,36,98,74,57,71,32),

     (1,88)(2,57)(3,50)(4,29)(5,46)(6,45)(7,59)(8,30)(9,
    34)(10,20)(11,98)(12,58)(13,22)(14,85)(15,17)(16,94)(18,89)\\(19,38)(21,99)(23,
    40)(24,37)(25,80)(26,67)(27,62)(28,95)(31,92)(32,93)(33,61)(35,84)(36,97)(39,
    76)(41,86)(42,56)\\(43,96)(44,71)(47,83)(48,90)(49,54)(51,65)(52,69)(53,87)(55,
    79)(60,75)(63,64)(66,72)(68,82)(70,81)(73,91)(74,78)\\(77,100),

  (1,94,61,86,26,14,62,83,66,99)(2,97,7,44,11,63,93,78,18,37)(3,96,69,73,30,45,58,
    31,51,84)(4,15,42,48,40,22,9,87,\\77,81)(5,10,19,54,79,28,75,68,80,76)(6,12,92,65,
    35,50,43,52,91,8)(13,34,53,100,70,29,17,56,90,23)(16,33,41,67,85,\\27,47,72,21,
    88)(20,38,49,55,95,60,82,25,39,46)(24,57,36,59,71,98,64,32,74,89),

  (1,99,66,83,62,14,26,86,61,94)(2,37,18,78,93,63,11,44,7,97)(3,84,51,31,58,45,30,
    73,69,96)(4,81,77,87,9,22,40,48,\\42,15)(5,76,80,68,75,28,79,54,19,10)(6,8,91,52,
    43,50,35,65,92,12)(13,23,90,56,17,29,70,100,53,34)(16,88,21,72,47,\\27,85,67,41,
    33)(20,46,39,25,82,60,95,55,49,38)(24,89,74,32,64,98,71,59,36,57).

The character table of $G^{s_8}$:\\
\begin{tabular}{c|cccccccccccccccccccc}
  & & & & & & & & & & 10& & & & & & & & & &20\\\hline
$\chi_8^{(1)}$&1&1&1&1&1&1&1&1&1&1&1&1&1&1&1&1&1&1&1&1
\\$\chi_8^{(2)}$&1&-1&-1&-1&1&1&1&-1&-1&1&1&1&1&1&-1&-1&-1&1&-1&-1
\\$\chi_8^{(3)}$&1&-1&1&1&1&-1&-1&1&1&1&1&1&-1&-1&-1&1&-1&-1&-1&-1
\\$\chi_8^{(4)}$&1&1&-1&-1&1&-1&-1&-1&-1&1&1&1&-1&-1&1&-1&1&-1&1&1
\\$\chi_8^{(5)}$&1&-1&A&/A&-/A&-A&-B&/B&B&-B&-A&-/B&-/A&-/B&B&-1&/B&1&A&/A
\\$\chi_8^{(6)}$&1&-1&B&/B&-/B&-B&-/A&A&/A&-/A&-B&-A&-/B&-A&/A&-1&A&1&B&/B
\\$\chi_8^{(7)}$&1&-1&/B&B&-B&-/B&-A&/A&A&-A&-/B&-/A&-B&-/A&A&-1&/A&1&/B&B
\\$\chi_8^{(8)}$&1&-1&/A&A&-A&-/A&-/B&B&/B&-/B&-/A&-B&-A&-B&/B&-1&B&1&/A&A
\\$\chi_8^{(9)}$&1&-1&-/A&-A&-A&/A&/B&-B&-/B&-/B&-/A&-B&A&B&/B&1&B&-1&/A&A
\\$\chi_8^{(10)}$&1&-1&-/B&-B&-B&/B&A&-/A&-A&-A&-/B&-/A&B&/A&A&1&/A&-1&/B&B
\\$\chi_8^{(11)}$&1&-1&-B&-/B&-/B&B&/A&-A&-/A&-/A&-B&-A&/B&A&/A&1&A&-1&B&/B
\\$\chi_8^{(12)}$&1&-1&-A&-/A&-/A&A&B&-/B&-B&-B&-A&-/B&/A&/B&B&1&/B&-1&A&/A
\\$\chi_8^{(13)}$&1&1&A&/A&-/A&A&B&/B&B&-B&-A&-/B&/A&/B&-B&-1&-/B&-1&-A&-/A
\\$\chi_8^{(14)}$&1&1&B&/B&-/B&B&/A&A&/A&-/A&-B&-A&/B&A&-/A&-1&-A&-1&-B&-/B
\\$\chi_8^{(15)}$&1&1&/B&B&-B&/B&A&/A&A&-A&-/B&-/A&B&/A&-A&-1&-/A&-1&-/B&-B
\\$\chi_8^{(16)}$&1&1&/A&A&-A&/A&/B&B&/B&-/B&-/A&-B&A&B&-/B&-1&-B&-1&-/A&-A
\\$\chi_8^{(17)}$&1&1&-/A&-A&-A&-/A&-/B&-B&-/B&-/B&-/A&-B&-A&-B&-/B&1&-B&1&-/A&-A
\\$\chi_8^{(18)}$&1&1&-/B&-B&-B&-/B&-A&-/A&-A&-A&-/B&-/A&-B&-/A&-A&1&-/A&1&-/B&-B
\\$\chi_8^{(19)}$&1&1&-B&-/B&-/B&-B&-/A&-A&-/A&-/A&-B&-A&-/B&-A&-/A&1&-A&1&-B&-/B
\\$\chi_8^{(20)}$&1&1&-A&-/A&-/A&-A&-B&-/B&-B&-B&-A&-/B&-/A&-/B&-B&1&-/B&1&-A&-/A
\end{tabular}

\noindent \noindent where   A = -E(5); B = -E(5)$^2$.

The generators of $G^{s_9}$ are:\\

 (  1, 61, 26, 62, 66)(  2,  7, 11, 93, 18)(  3, 69, 30, 58, 51)
    (  4, 42, 40,  9, 77)(  5, 19, 79, 75, 80)(  6, 92, 35, 43, 91)(  8, 12, 65, 50,
     52)( 10, 54, 28, 68, 76)( 13, 53, 70, 17, 90)( 14, 83, 99, 94, 86)
    ( 15, 48, 22, 87, 81)( 16, 41, 85, 47, 21)( 20, 49, 95, 82, 39)( 23, 34,100, 29,
     56)( 24, 36, 71, 64, 74)( 25, 46, 38, 55, 60)( 27, 72, 88, 33, 67)
    ( 31, 84, 96, 73, 45)( 32, 89, 57, 59, 98)( 37, 97, 44, 63, 78),
  (  2, 21, 53)(  3, 77, 60)(  4, 25, 69)(  5, 35, 57)(  6, 32, 75)(  7, 16, 70)
    (  8, 28, 64)(  9, 55, 51)( 10, 36, 50)( 11, 41, 17)( 12, 68, 74)( 13, 18, 47)
    ( 15, 73, 39)( 19, 43, 59)( 20, 48, 45)( 22, 31, 49)( 23, 27, 78)( 24, 65, 76)
    ( 29, 33, 44)( 30, 42, 46)( 34, 72, 37)( 38, 58, 40)( 52, 54, 71)( 56, 67, 63)
    ( 79, 91, 98)( 80, 92, 89)( 81, 96, 82)( 84, 95, 87)( 85, 90, 93)( 88, 97,100),
  (  1,  2, 51, 52, 70)(  3,  8, 17, 61,  7)(  4, 19, 88, 95, 24)(  5, 72, 49, 74,
      77)(  6, 84, 37, 83, 23)(  9, 80, 27, 20, 64)( 10, 87, 98, 25, 16)
    ( 11, 69, 12, 90, 26)( 13, 62, 93, 30, 65)( 14, 56, 91, 31, 78)( 15, 89, 38, 85,
     28)( 18, 58, 50, 53, 66)( 21, 76, 22, 59, 60)( 29, 43, 45, 63, 86)
    ( 32, 46, 41, 54, 81)( 33, 82, 36, 42, 79)( 34, 92, 96, 97, 99)( 35, 73, 44, 94,
     100)( 39, 71, 40, 75, 67)( 47, 68, 48, 57, 55).

The representatives of conjugacy classes of   $G^{s_9}$ are:\\
 (1),
  (2,21,53)(3,77,60)(4,25,69)(5,35,57)(6,32,75)(7,16,70)(8,28,64)(9,55,51)(10,
    36,50)(11,41,17)(12,68,74)(13,18,\\47)(15,73,39)(19,43,59)(20,48,45)(22,31,49)(23,
    27,78)(24,65,76)(29,33,44)(30,42,46)(34,72,37)(38,58,40)(52,54,\\71)(56,67,63)(79,
    91,98)(80,92,89)(81,96,82)(84,95,87)(85,90,93)(88,97,100),

  (1,2,51,52,70)(3,8,17,61,7)(4,19,88,95,24)(5,72,49,74,77)(6,84,37,83,23)(9,80,27,
    20,64)(10,87,98,25,16)(11,69,\\12,90,26)(13,62,93,30,65)(14,56,91,31,78)(15,89,38,
    85,28)(18,58,50,53,66)(21,76,22,59,60)(29,43,45,63,86)(32,46,\\41,54,81)(33,82,36,
    42,79)(34,92,96,97,99)(35,73,44,94,100)(39,71,40,75,67)(47,68,48,57,55),

  (1,2,66,18,62,93,26,11,61,7)(3,21,51,47,58,85,30,41,69,16)(4,12,77,8,9,52,40,50,
    42,65)(5,55,80,38,75,46,79,25,\\19,60)(6,67,91,33,43,88,35,72,92,27)(10,53,76,13,
    68,90,28,17,54,70)(14,56,86,29,94,100,99,34,83,23)(15,71,81,36,\\87,24,22,74,48,
    64)(20,63,39,44,82,97,95,37,49,78)(31,59,45,57,73,89,96,32,84,98),

  (1,3,90,93,50)(2,8,26,30,53)(4,33,71,80,49)(5,95,42,67,64)(6,97,100,45,14)(7,12,
    62,58,70)(9,72,24,79,39)(10,32,\\85,48,60)(11,65,66,51,17)(13,18,52,61,69)(15,55,
    76,98,41)(16,81,38,68,59)(19,82,40,27,74)(20,77,88,36,75)(21,87,\\46,28,57)(22,25,
    54,89,47)(23,96,94,43,78)(29,31,83,92,44)(34,73,86,91,37)(35,63,56,84,99),

  (1,3,88,100,36)(2,87,99,35,21)(4,54,80,44,31)(5,63,84,42,28)(6,85,93,48,14)(7,81,
    94,43,16)(8,57,95,46,53)(9,76,\\79,37,73)(10,75,97,45,77)(11,15,86,91,41)(12,59,
    82,38,70)(13,52,89,49,25)(17,65,98,39,55)(18,22,83,92,47)(19,78,\\96,40,68)(20,60,
    90,50,32)(23,74,62,58,27)(24,66,51,72,34)(26,30,67,56,64)(29,71,61,69,33),

  (1,4,2,16,76)(3,71,5,70,73)(6,22,46,78,65)(7,41,10,61,42)(8,43,15,60,44)(9,93,47,
    28,62)(11,85,54,26,40)(12,91,\\48,25,63)(13,84,58,24,75)(14,89,56,27,39)(17,45,69,
    64,19)(18,21,68,66,77)(20,83,57,23,72)(29,67,82,86,32)(30,74,\\79,90,31)(33,95,94,
    98,100)(34,88,49,99,59)(35,81,55,97,52)(36,80,53,96,51)(37,50,92,87,38),

  (1,4,63,27,65)(2,68,79,48,69)(3,18,28,19,15)(5,81,51,93,54)(6,29,95,38,24)(7,76,
    75,22,30)(8,62,9,97,33)(10,80,\\87,58,11)(12,66,77,44,67)(13,21,96,86,32)(14,89,
    53,16,73)(17,85,31,99,59)(20,25,64,43,34)(23,39,60,71,35)(26,40,\\37,88,52)(36,92,
    56,82,55)(41,45,83,57,70)(42,78,72,50,61)(46,74,91,100,49)(47,84,94,98,90),

  (1,4,64,62,9,36,61,42,74,66,77,71,26,40,24)(2,96,55,93,31,46,7,73,60,18,84,38,11,
    45,25)(3,44,95,58,37,20,69,63,\\82,51,97,49,30,78,39)(5,23,65,75,29,8,19,34,50,80,
    56,12,79,100,52)(6,13,28,43,17,10,92,53,68,91,90,54,35,70,76)\\(14,89,87,94,98,48,
    83,57,81,86,32,22,99,59,15)(16,41,85,47,21)(27,72,88,33,67),

  (1,5,98,61,19,32,26,79,89,62,75,57,66,80,59)(2,15,52,7,48,8,11,22,12,93,87,65,18,
    81,50)(3,30,51,69,58)(4,94,60,\\42,86,25,40,14,46,9,83,38,77,99,55)(6,35,91,92,
    43)(10,67,17,54,27,90,28,72,13,68,88,53,76,33,70)(16,97,95,41,44,\\82,85,63,39,47,
    78,20,21,37,49)(23,36,84,34,71,96,100,64,73,29,74,45,56,24,31),

  (1,6)(2,56)(3,14)(4,25)(5,57)(7,23)(8,63)(9,55)(10,20)(11,34)(12,78)(13,31)(15,
    72)(16,74)(17,73)(18,29)(19,\\59)(21,64)(22,33)(24,41)(26,35)(27,81)(28,95)(30,
    99)(32,75)(36,85)(37,65)(38,40)(39,76)(42,46)(43,62)(44,52)\\(45,90)(47,71)(48,
    88)(49,54)(50,97)(51,86)(53,84)(58,94)(60,77)(61,92)(66,91)(67,87)(68,82)(69,
    83)(70,96)(79,\\98)(80,89)(93,100),

    (1,7,61,11,26,93,62,18,66,2)(3,16,69,41,30,
    85,58,47,51,21)(4,65,42,50,40,52,9,8,77,12)(5,60,19,25,79,46,75,38,\\80,55)(6,27,
    92,72,35,88,43,33,91,67)(10,70,54,17,28,90,68,13,76,53)(14,23,83,34,99,100,94,
    29,86,56)(15,64,48,74,\\22,24,87,36,81,71)(20,78,49,37,95,97,82,44,39,63)(31,98,
    84,32,96,89,73,57,45,59),

     (1,7,69,65,53)(2,3,12,13,66)(4,79,67,20,74)(5,88,82,
    71,9)(6,96,44,86,56)(8,90,62,18,51)(10,81,89,55,21)(11,30,\\50,70,61)(14,23,92,73,
    63)(15,57,60,16,54)(17,26,93,58,52)(19,33,39,64,77)(22,98,46,85,68)(24,42,75,27,
    49)(25,41,\\28,48,59)(29,91,84,97,94)(31,37,99,100,43)(32,38,47,76,87)(34,35,45,
    78,83)(36,40,80,72,95),

     (1,8,11,13,58)(2,17,69,62,50)(3,26,65,18,70)(4,82,75,
    64,72)(5,24,33,40,20)(6,99,73,29,78)(7,90,30,66,52)(9,49,\\19,36,67)(10,46,15,47,
    59)(12,93,53,51,61)(14,84,34,44,43)(16,32,28,55,22)(21,98,54,38,48)(23,97,35,86,
    31)(25,81,\\85,57,76)(27,77,95,79,71)(37,92,94,45,56)(39,80,74,88,42)(41,89,68,60,
    87)(63,91,83,96,100),

    (1,9,74,26,4,36,66,40,64,61,77,24,62,42,71)(2,31,60,11,
    96,46,18,45,55,7,84,25,93,73,38)(3,37,82,30,44,20,51,78,\\95,69,97,39,58,63,49)(5,
    29,50,79,23,8,80,100,65,19,56,52,75,34,12)(6,17,68,35,13,10,91,70,28,92,90,76,
    43,53,54)\\(14,98,81,99,89,48,86,59,87,83,32,15,94,57,22)(16,21,47,85,41)(27,67,
    33,88,72),

    (1,9,7,21,54)(2,47,10,66,40)(3,24,19,53,31)(4,93,41,68,26)(5,13,45,
    51,74)(6,15,38,63,52)(8,92,48,55,78)(11,16,\\28,61,77)(12,35,22,60,37)(14,98,23,
    67,49)(17,96,30,71,75)(18,85,76,62,42)(20,86,59,56,33)(25,97,65,43,87)(27,95,\\83,
    32,34)(29,88,39,94,57)(36,79,70,84,69)(44,50,91,81,46)(58,64,80,90,73)(72,82,99,
    89,100),

     (1,9,78,67,52)(2,54,75,15,58)(3,11,68,5,22)(4,97,72,12,26)(6,34,82,
    46,71)(7,28,80,48,51)(8,61,77,37,27)(10,79,\\81,30,18)(13,85,73,94,57)(14,98,70,
    21,31)(16,84,83,32,17)(19,87,69,93,76)(20,55,74,35,29)(23,95,25,36,91)(24,43,\\56,49,60)(33,50,66,40,63)(38,64,92,100,39)(41,96,99,89,90)(42,44,88,65,62)(45,86,
    59,53,47),

    (1,11,62,2,61,93,66,7,26,18)(3,41,58,21,69,85,51,16,30,47)(4,50,9,
    12,42,52,77,65,40,8)(5,25,75,55,19,46,80,60,\\79,38)(6,72,43,67,92,88,91,27,35,
    33)(10,17,68,53,54,90,76,70,28,13)(14,34,94,56,83,100,86,23,99,29)(15,74,87,71,\\48,24,81,64,22,36)(20,37,82,63,49,97,39,78,95,44)(31,32,73,59,84,89,45,98,96,
    57),

    (1,12,25,66,8,60,62,52,55,26,50,38,61,65,46)(2,32,74,18,98,64,93,59,71,
    11,57,36,7,89,24)(3,82,47,51,95,85,58,\\49,41,30,20,16,69,39,21)(4,15,56,77,81,29,
    9,87,100,40,22,34,42,48,23)(5,14,96,80,86,84,75,94,31,79,99,45,19,83,\\73)(6,68,
    33,91,28,88,43,54,72,35,10,27,92,76,67)(13,17,53,90,70)(37,63,97,78,44),

  (1,12,18,17,30)(2,90,58,61,65)(3,62,52,11,53)(4,39,5,36,27)(6,94,31,34,63)(7,13,
    51,26,50)(8,93,70,69,66)(9,95,\\75,74,33)(10,38,22,41,57)(14,96,29,37,35)(15,21,
    32,68,25)(16,89,76,46,48)(19,71,72,42,20)(23,44,91,99,45)(24,67,\\77,82,80)(28,60,
    81,47,98)(40,49,79,64,88)(43,83,73,56,97)(54,55,87,85,59)(78,92,86,84,100),

  (1,13,26,70,66,90,61,53,62,17)(2,38,11,60,18,46,7,55,93,25)(3,71,30,74,51,36,69,
    64,58,24)(4,57,40,98,77,89,42,\\59,9,32)(5,23,79,100,80,56,19,34,75,29)(6,22,35,
    81,91,48,92,87,43,15)(8,16,65,85,52,21,12,41,50,47)(10,49,28,82,\\76,20,54,95,68,
    39)(14,44,99,78,86,97,83,63,94,37)(27,73,88,31,67,96,72,45,33,84),

  (1,26,66,61,62)(2,11,18,7,93)(3,30,51,69,58)(4,40,77,42,9)(5,79,80,19,75)(6,35,91,
    92,43)(8,65,52,12,50)(10,28,\\76,54,68)(13,70,90,53,17)(14,99,86,83,94)(15,22,81,
    48,87)(16,85,21,41,47)(20,95,39,49,82)(23,100,56,34,29)(24,\\71,74,36,64)(25,38,
    60,46,55)(27,88,67,72,33)(31,96,45,84,73)(32,57,98,89,59)(37,44,78,97,63),

  (1,61,26,62,66)(2,7,11,93,18)(3,69,30,58,51)(4,42,40,9,77)(5,19,79,75,80)(6,92,35,
    43,91)(8,12,65,50,52)(10,54,\\28,68,76)(13,53,70,17,90)(14,83,99,94,86)(15,48,22,
    87,81)(16,41,85,47,21)(20,49,95,82,39)(23,34,100,29,56)(24,\\36,71,64,74)(25,46,
    38,55,60)(27,72,88,33,67)(31,84,96,73,45)(32,89,57,59,98)(37,97,44,63,78),

  (1,62,61,66,26)(2,93,7,18,11)(3,58,69,51,30)(4,9,42,77,40)(5,75,19,80,79)(6,43,92,
    91,35)(8,50,12,52,65)(10,68,\\54,76,28)(13,17,53,90,70)(14,94,83,86,99)(15,87,48,
    81,22)(16,47,41,21,85)(20,82,49,39,95)(23,29,34,56,100)(24,\\64,36,74,71)(25,55,
    46,60,38)(27,33,72,67,88)(31,73,84,45,96)(32,59,89,98,57)(37,63,97,78,44),

  (1,66,62,26,61)(2,18,93,11,7)(3,51,58,30,69)(4,77,9,40,42)(5,80,75,79,19)(6,91,43,
    35,92)(8,52,50,65,12)(10,76,\\68,28,54)(13,90,17,70,53)(14,86,94,99,83)(15,81,87,
    22,48)(16,21,47,85,41)(20,39,82,95,49)(23,56,29,100,34)(24,\\74,64,71,36)(25,60,
    55,38,46)(27,67,33,88,72)(31,45,73,96,84)(32,98,59,57,89)(37,78,63,44,97).

The character table of $G^{s_9}$:\\
\begin{tabular}{c|cccccccccccccccccccc}
  & & & & & & & & & & 10& & & & & & & & & &20\\\hline

$\chi_9^{(1)}$&1&1&1&1&1&1&1&1&1&1&1&1&1&1&1&1&1&1&1&1
\\$\chi_9^{(2)}$&1&1&A&A&1&1&/B&/B&/B&A&1&/A&/A&A&B&B&B&B&/A&/A
\\$\chi_9^{(3)}$&1&1&B&B&1&1&A&A&A&B&1&/B&/B&B&/A&/A&/A&/A&/B&/B
\\$\chi_9^{(4)}$&1&1&/B&/B&1&1&/A&/A&/A&/B&1&B&B&/B&A&A&A&A&B&B
\\$\chi_9^{(5)}$&1&1&/A&/A&1&1&B&B&B&/A&1&A&A&/A&/B&/B&/B&/B&A&A
\\$\chi_9^{(6)}$&3&.&C&-A&E&*E&D&/F&.&.&-1&-/A&/C&/G&.&/D&F&-B&.&G
\\$\chi_9^{(7)}$&3&.&D&-/B&*E&E&/C&G&.&.&-1&-B&/D&/F&.&C&/G&-A&.&F
\\$\chi_9^{(8)}$&3&.&/D&-B&*E&E&C&/G&.&.&-1&-/B&D&F&.&/C&G&-/A&.&/F
\\$\chi_9^{(9)}$&3&.&/C&-/A&E&*E&/D&F&.&.&-1&-A&C&G&.&D&/F&-/B&.&/G
\\$\chi_9^{(10)}$&3&.&E&-1&*E&E&*E&E&.&.&-1&-1&E&*E&.&*E&E&-1&.&*E
\\$\chi_9^{(11)}$&3&.&*E&-1&E&*E&E&*E&.&.&-1&-1&*E&E&.&E&*E&-1&.&E
\\$\chi_9^{(12)}$&3&.&F&-B&E&*E&/G&C&.&.&-1&-/B&/F&/D&.&G&/C&-/A&.&D
\\$\chi_9^{(13)}$&3&.&G&-/A&*E&E&F&/D&.&.&-1&-A&/G&/C&.&/F&D&-/B&.&C
\\$\chi_9^{(14)}$&3&.&/G&-A&*E&E&/F&D&.&.&-1&-/A&G&C&.&F&/D&-B&.&/C
\\$\chi_9^{(15)}$&3&.&/F&-/B&E&*E&G&/C&.&.&-1&-B&F&D&.&/G&C&-A&.&/D
\\$\chi_9^{(16)}$&4&1&-1&.&-1&-1&-1&-1&1&1&.&.&-1&-1&1&-1&-1&.&1&-1
\\$\chi_9^{(17)}$&4&1&-A&.&-1&-1&-/B&-/B&/B&A&.&.&-/A&-A&B&-B&-B&.&/A&-/A
\\$\chi_9^{(18)}$&4&1&-B&.&-1&-1&-A&-A&A&B&.&.&-/B&-B&/A&-/A&-/A&.&/B&-/B
\\$\chi_9^{(19)}$&4&1&-/B&.&-1&-1&-/A&-/A&/A&/B&.&.&-B&-/B&A&-A&-A&.&B&-B
\\$\chi_9^{(20)}$&4&1&-/A&.&-1&-1&-B&-B&B&/A&.&.&-A&-/A&/B&-/B&-/B&.&A&-A
\\$\chi_9^{(21)}$&5&-1&.&1&.&.&.&.&-1&-1&1&1&.&.&-1&.&.&1&-1&.
\\$\chi_9^{(22)}$&5&-1&.&A&.&.&.&.&-/B&-A&1&/A&.&.&-B&.&.&B&-/A&.
\\$\chi_9^{(23)}$&5&-1&.&B&.&.&.&.&-A&-B&1&/B&.&.&-/A&.&.&/A&-/B&.
\\$\chi_9^{(24)}$&5&-1&.&/B&.&.&.&.&-/A&-/B&1&B&.&.&-A&.&.&A&-B&.
\\$\chi_9^{(25)}$&5&-1&.&/A&.&.&.&.&-B&-/A&1&A&.&.&-/B&.&.&/B&-A&.
\end{tabular}

\begin{tabular}{c|ccccc}
  & & & & &\\\hline
$\chi_9^{(1)}$&1&1&1&1&1
\\$\chi_9^{(2)}$&/B&A&/B&/A&B
\\$\chi_9^{(3)}$&A&B&A&/B&/A
\\$\chi_9^{(4)}$&/A&/B&/A&B&A
\\$\chi_9^{(5)}$&B&/A&B&A&/B
\\$\chi_9^{(6)}$&-/B&H&I&/H&/I
\\$\chi_9^{(7)}$&-/A&I&/H&/I&H
\\$\chi_9^{(8)}$&-A&/I&H&I&/H
\\$\chi_9^{(9)}$&-B&/H&/I&H&I
\\$\chi_9^{(10)}$&-1&3&3&3&3
\\$\chi_9^{(11)}$&-1&3&3&3&3
\\$\chi_9^{(12)}$&-A&/I&H&I&/H
\\$\chi_9^{(13)}$&-B&/H&/I&H&I
\\$\chi_9^{(14)}$&-/B&H&I&/H&/I
\\$\chi_9^{(15)}$&-/A&I&/H&/I&H
\\$\chi_9^{(16)}$&.&4&4&4&4
\\$\chi_9^{(17)}$&.&J&/K&/J&K
\\$\chi_9^{(18)}$&.&K&J&/K&/J
\\$\chi_9^{(19)}$&.&/K&/J&K&J
\\$\chi_9^{(20)}$&.&/J&K&J&/K
\\$\chi_9^{(21)}$&1&5&5&5&5
\\$\chi_9^{(22)}$&/B&L&/M&/L&M
\\$\chi_9^{(23)}$&A&M&L&/M&/L
\\$\chi_9^{(24)}$&/A&/M&/L&M&L
\\$\chi_9^{(25)}$&B&/L&M&L&/M
\end{tabular}

\noindent \noindent where   A = E(5)$^4$; B = E(5)$^3$; C =
-E(5)-E(5)$^2$; D = -E(5)-E(5)$^3$; E = -E(5)-E(5)$^4$
  = (1-ER(5))/2 = -b5;
F = E(5)$^2$+E(5)$^3$+E(5)$^4$; G = E(5)+E(5)$^3$+E(5)$^4$; H =
3*E(5)$^4$; I = 3*E(5)$^2$; J = 4*E(5)$^4$; K = 4*E(5)$^3$; L =
5*E(5)$^4$; M = 5*E(5)$^3$.

The generators of $G^{s_{10}}$ are:\\
(  1, 43, 80, 66, 35, 75, 62, 92, 79, 26,  6, 19, 61, 91,  5)(  2,
36, 78, 18, 24,
     63, 93, 74, 44, 11, 64, 97,  7, 71, 37)(  3, 38, 83, 51, 46, 14, 58, 25, 86,
      30, 60, 94, 69, 55, 99)(  4,  9, 42, 77, 40)(  8, 20, 27, 52, 39, 67, 50, 82,
      33, 65, 95, 88, 12, 49, 72)( 10, 16, 31, 76, 21, 45, 68, 47, 73, 28, 85, 96,
      54, 41, 84)( 13, 34, 87, 90, 23, 22, 17, 56, 48, 70, 29, 15, 53,100, 81)
    ( 32, 59, 89, 98, 57)

The representatives of conjugacy classes of   $G^{s_{10}}$ are:\\
 (1),
 (1,5,91,61,19,6,26,79,92,62,75,35,66,80,43)(2,37,71,7,97,64,11,44,74,93,63,
    24,18,78,36)(3,99,55,69,94,60,30,86,\\25,58,14,46,51,83,38)(4,40,77,42,9)(8,72,49,
    12,88,95,65,33,82,50,67,39,52,27,20)(10,84,41,54,96,85,28,73,47,68,45,\\21,76,31,
    16)(13,81,100,53,15,29,70,48,56,17,22,23,90,87,34)(32,57,98,89,59),

  (1,6,75)(2,64,63)(3,60,14)(5,26,35)(7,74,78)(8,95,67)(10,85,45)(11,24,37)(12,82,
    27)(13,29,22)(15,17,34)(16,96,\\68)(18,71,44)(19,62,43)(20,88,50)(21,84,28)(23,81,
    70)(25,83,69)(30,46,99)(31,54,47)(33,52,49)(36,97,93)(38,94,58)\\(39,72,65)(41,73,
    76)(48,90,100)(51,55,86)(53,56,87)(61,92,80)(66,91,79),

  (1,19,92,66,5,6,62,80,91,26,75,43,61,79,35)(2,97,74,18,37,64,93,78,71,11,63,36,7,
    44,24)(3,94,25,51,99,60,58,83,\\55,30,14,38,69,86,46)(4,9,42,77,40)(8,88,82,52,72,
    95,50,27,49,65,67,20,12,33,39)(10,96,47,76,84,85,68,31,41,28,45,\\16,54,73,21)(13,
    15,56,90,81,29,17,87,100,70,22,34,53,48,23)(32,59,89,98,57),

  (1,26,66,61,62)(2,11,18,7,93)(3,30,51,69,58)(4,40,77,42,9)(5,79,80,19,75)(6,35,91,
    92,43)(8,65,52,12,50)(10,28,\\76,54,68)(13,70,90,53,17)(14,99,86,83,94)(15,22,81,
    48,87)(16,85,21,41,47)(20,95,39,49,82)(23,100,56,34,29)(24,71,\\74,36,64)(25,38,
    60,46,55)(27,88,67,72,33)(31,96,45,84,73)(32,57,98,89,59)(37,44,78,97,63),

  (1,35,79,61,43,75,26,91,80,62,6,5,66,92,19)(2,24,44,7,36,63,11,71,78,93,64,37,18,
    74,97)(3,46,86,69,38,14,30,55,\\83,58,60,99,51,25,94)(4,40,77,42,9)(8,39,33,12,20,
    67,65,49,27,50,95,72,52,82,88)(10,21,73,54,16,45,28,41,31,68,85,\\84,76,47,96)(13,
    23,48,53,34,22,70,100,87,17,29,81,90,56,15)(32,57,98,89,59),

  (1,43,80,66,35,75,62,92,79,26,6,19,61,91,5)(2,36,78,18,24,63,93,74,44,11,64,97,7,
    71,37)(3,38,83,51,46,14,58,25,\\86,30,60,94,69,55,99)(4,9,42,77,40)(8,20,27,52,39,
    67,50,82,33,65,95,88,12,49,72)(10,16,31,76,21,45,68,47,73,28,85,\\96,54,41,84)(13,
    34,87,90,23,22,17,56,48,70,29,15,53,100,81)(32,59,89,98,57),

  (1,61,26,62,66)(2,7,11,93,18)(3,69,30,58,51)(4,42,40,9,77)(5,19,79,75,80)(6,92,35,
    43,91)(8,12,65,50,52)(10,54,\\28,68,76)(13,53,70,17,90)(14,83,99,94,86)(15,48,22,
    87,81)(16,41,85,47,21)(20,49,95,82,39)(23,34,100,29,56)(24,36,\\71,64,74)(25,46,
    38,55,60)(27,72,88,33,67)(31,84,96,73,45)(32,89,57,59,98)(37,97,44,63,78),

  (1,62,61,66,26)(2,93,7,18,11)(3,58,69,51,30)(4,9,42,77,40)(5,75,19,80,79)(6,43,92,
    91,35)(8,50,12,52,65)(10,68,\\54,76,28)(13,17,53,90,70)(14,94,83,86,99)(15,87,48,
    81,22)(16,47,41,21,85)(20,82,49,39,95)(23,29,34,56,100)(24,64,\\36,74,71)(25,55,
    46,60,38)(27,33,72,67,88)(31,73,84,45,96)(32,59,89,98,57)(37,63,97,78,44),

  (1,66,62,26,61)(2,18,93,11,7)(3,51,58,30,69)(4,77,9,40,42)(5,80,75,79,19)(6,91,43,
    35,92)(8,52,50,65,12)(10,76,\\68,28,54)(13,90,17,70,53)(14,86,94,99,83)(15,81,87,
    22,48)(16,21,47,85,41)(20,39,82,95,49)(23,56,29,100,34)(24,74,\\64,71,36)(25,60,
    55,38,46)(27,67,33,88,72)(31,45,73,96,84)(32,98,59,57,89)(37,78,63,44,97),

  (1,75,6)(2,63,64)(3,14,60)(5,35,26)(7,78,74)(8,67,95)(10,45,85)(11,37,24)(12,27,
    82)(13,22,29)(15,34,17)(16,68,\\96)(18,44,71)(19,43,62)(20,50,88)(21,28,84)(23,70,
    81)(25,69,83)(30,99,46)(31,47,54)(33,49,52)(36,93,97)(38,58,94)\\(39,65,72)(41,76,
    73)(48,100,90)(51,86,55)(53,87,56)(61,80,92)(66,79,91),

  (1,79,43,26,80,6,66,19,35,61,75,91,62,5,92)(2,44,36,11,78,64,18,97,24,7,63,71,93,
    37,74)(3,86,38,30,83,60,51,94,\\46,69,14,55,58,99,25)(4,77,9,40,42)(8,33,20,65,27,
    95,52,88,39,12,67,49,50,72,82)(10,73,16,28,31,85,76,96,21,54,45,\\41,68,84,47)(13,
    48,34,70,87,29,90,15,23,53,22,100,17,81,56)(32,98,59,57,89),

  (1,80,35,62,79,6,61,5,43,66,75,92,26,19,91)(2,78,24,93,44,64,7,37,36,18,63,74,11,
    97,71)(3,83,46,58,86,60,69,99,\\38,51,14,25,30,94,55)(4,42,40,9,77)(8,27,39,50,33,
    95,12,72,20,52,67,82,65,88,49)(10,31,21,68,73,85,54,84,16,76,45,\\47,28,96,41)(13,
    87,23,17,48,29,53,81,34,90,22,56,70,15,100)(32,89,57,59,98),

  (1,91,19,26,92,75,66,43,5,61,6,79,62,35,80)(2,71,97,11,74,63,18,36,37,7,64,44,93,
    24,78)(3,55,94,30,25,14,51,38,\\99,69,60,86,58,46,83)(4,77,9,40,42)(8,49,88,65,82,
    67,52,20,72,12,95,33,50,39,27)(10,41,96,28,47,45,76,16,84,54,85,\\73,68,21,31)(13,
    100,15,70,56,22,90,34,81,53,29,48,17,23,87)(32,98,59,57,89),

  (1,92,5,62,91,75,61,35,19,66,6,80,26,43,79)(2,74,37,93,71,63,7,24,97,18,64,78,11,
    36,44)(3,25,99,58,55,14,69,46,\\94,51,60,83,30,38,86)(4,42,40,9,77)(8,82,72,50,49,
    67,12,39,88,52,95,27,65,20,33)(10,47,84,68,41,45,54,21,96,76,85,\\31,28,16,73)(13,
    56,81,17,100,22,53,23,15,90,29,87,70,34,48)(32,89,57,59,98).

The character table of $G^{s_{10}}$:\\
\begin{tabular}{c|ccccccccccccccc}
  & & & & & & & & & &10 & & & & &\\\hline
$\chi_{10}^{(1)}$&1&1&1&1&1&1&1&1&1&1&1&1&1&1&1
\\$\chi_{10}^{(2)}$&1&A&/A&A&1&/A&/A&1&1&1&A&A&A&/A&/A
\\$\chi_{10}^{(3)}$&1&/A&A&/A&1&A&A&1&1&1&/A&/A&/A&A&A
\\$\chi_{10}^{(4)}$&1&B&1&/B&B&B&/B&/E&/B&E&1&E&/E&E&/E
\\$\chi_{10}^{(5)}$&1&C&/A&/D&B&D&/C&/E&/B&E&A&F&/G&G&/F
\\$\chi_{10}^{(6)}$&1&D&A&/C&B&C&/D&/E&/B&E&/A&G&/F&F&/G
\\$\chi_{10}^{(7)}$&1&E&1&/E&E&E&/E&B&/E&/B&1&/B&B&/B&B
\\$\chi_{10}^{(8)}$&1&F&/A&/G&E&G&/F&B&/E&/B&A&/D&C&/C&D
\\$\chi_{10}^{(9)}$&1&G&A&/F&E&F&/G&B&/E&/B&/A&/C&D&/D&C
\\$\chi_{10}^{(10)}$&1&/E&1&E&/E&/E&E&/B&E&B&1&B&/B&B&/B
\\$\chi_{10}^{(11)}$&1&/G&/A&F&/E&/F&G&/B&E&B&A&C&/D&D&/C
\\$\chi_{10}^{(12)}$&1&/F&A&G&/E&/G&F&/B&E&B&/A&D&/C&C&/D
\\$\chi_{10}^{(13)}$&1&/B&1&B&/B&/B&B&E&B&/E&1&/E&E&/E&E
\\$\chi_{10}^{(14)}$&1&/D&/A&C&/B&/C&D&E&B&/E&A&/G&F&/F&G
\\$\chi_{10}^{(15)}$&1&/C&A&D&/B&/D&C&E&B&/E&/A&/F&G&/G&F
\end{tabular}

\noindent \noindent where   A = E(3)$^2$
  = (-1-ER(-3))/2 = -1-b3;
B = E(5)$^3$; C = E(15)$^4$; D = E(15)$^{14}$; E = E(5); F =
E(15)$^{13}$; G = E(15)$^8$.

The generators of $G^{s_{11}}$ are:\\
 (  1, 75,  6)(  2, 63, 64)(  3, 14, 60)(  5, 35, 26)(  7, 78, 74)
    (  8, 67, 95)( 10, 45, 85)( 11, 37, 24)( 12, 27, 82)( 13, 22, 29)( 15, 34, 17)
    ( 16, 68, 96)( 18, 44, 71)( 19, 43, 62)( 20, 50, 88)( 21, 28, 84)( 23, 70, 81)
    ( 25, 69, 83)( 30, 99, 46)( 31, 47, 54)( 33, 49, 52)( 36, 93, 97)( 38, 58, 94)
    ( 39, 65, 72)( 41, 76, 73)( 48,100, 90)( 51, 86, 55)( 53, 87, 56)( 61, 80, 92)
    ( 66, 79, 91), (  2,  5)(  3, 13)(  4,  9)(  7, 18)( 10, 16)( 11, 19)( 12, 20)
    ( 14, 22)( 23, 55)( 24, 62)( 25, 38)( 26, 64)( 27, 50)( 29, 60)( 30, 53)
    ( 31, 41)( 32, 59)( 33, 65)( 35, 63)( 36, 61)( 37, 43)( 39, 52)( 40, 77)
    ( 44, 78)( 45, 68)( 46, 56)( 47, 76)( 49, 72)( 51, 70)( 54, 73)( 58, 69)
    ( 71, 74)( 80, 93)( 81, 86)( 82, 88)( 83, 94)( 85, 96)( 87, 99)( 89, 98)
    ( 92, 97), (  2,  7)(  4,  9)(  5, 19)(  8, 12)( 11, 18)( 13, 17)( 15, 22)
    ( 16, 21)( 23, 56)( 24, 71)( 25, 55)( 26, 62)( 27, 67)( 28, 68)( 29, 34)
    ( 30, 58)( 31, 73)( 33, 72)( 35, 43)( 37, 44)( 38, 46)( 39, 49)( 40, 42)
    ( 41, 47)( 51, 69)( 52, 65)( 53, 70)( 54, 76)( 57, 59)( 61, 66)( 63, 78)
    ( 64, 74)( 79, 80)( 81, 87)( 82, 95)( 83, 86)( 84, 96)( 89, 98)( 91, 92)
    ( 94, 99), (  1,  2,  6, 64, 75, 63)(  3, 14, 60)(  4, 77)(  5, 92, 26, 80, 35,
      61)(  7, 78, 74)(  8, 82, 95, 27, 67, 12)( 10, 47, 85, 31, 45, 54)
    ( 11, 91, 24, 79, 37, 66)( 13, 87, 29, 53, 22, 56)( 15, 55, 17, 86, 34, 51)
    ( 16, 84, 96, 28, 68, 21)( 18, 43, 71, 19, 44, 62)( 20, 49, 88, 33, 50, 52)
    ( 23, 69, 81, 25, 70, 83)( 30, 48, 46, 90, 99,100)( 32, 89)( 36, 93, 97)
    ( 38, 58, 94)( 39, 65, 72)( 40, 42)( 41, 76, 73)( 57, 59),
  (  1,  3, 75, 14,  6, 60)(  2, 15, 78, 29, 64, 17,  7, 22, 63, 34, 74, 13)
    (  4, 57,  9, 59)(  5, 86, 43, 25, 26, 51, 19, 83, 35, 55, 62, 69)
    (  8, 76, 27, 31, 95, 41, 12, 54, 67, 73, 82, 47)( 10, 88, 45, 20, 85, 50)
    ( 11, 81, 44, 56, 24, 70, 18, 87, 37, 23, 71, 53)( 16, 52, 28, 72, 96, 49, 21,
      65, 68, 33, 84, 39)( 30, 79, 94, 92, 46, 66, 58, 80, 99, 91, 38, 61)( 32, 77)
    ( 36, 90, 93, 48, 97,100)( 40, 89, 42, 98).

The representatives of conjugacy classes of   $G^{s_{11}}$ are:\\
(1),
(2,5)(3,13)(4,9)(7,18)(10,16)(11,19)(12,20)(14,22)(23,55)(24,62)(25,38)(26,
    64)(27,50)(29,60)(30,53)(31,41)\\(32,59)(33,65)(35,63)(36,61)(37,43)(39,52)(40,
    77)(44,78)(45,68)(46,56)(47,76)(49,72)(51,70)(54,73)(58,69)(71,\\74)(80,93)(81,
    86)(82,88)(83,94)(85,96)(87,99)(89,98)(92,97),

  (2,18,19)(3,13,17)(5,11,7)(8,20,12)(10,16,21)(14,22,15)(23,46,25)(24,74,26)(27,67,
    50)(28,45,68)(29,34,60)(30,\\69,70)(31,54,47)(32,59,57)(33,49,52)(35,37,78)(36,61,
    66)(38,56,55)(39,72,65)(40,42,77)(41,76,73)(43,63,44)(51,\\58,53)(62,64,71)(79,93,
    80)(81,99,83)(82,95,88)(84,85,96)(86,94,87)(91,97,92),

  (1,2,6,64,75,63)(3,14,60)(4,77)(5,92,26,80,35,61)(7,78,74)(8,82,95,27,67,12)(10,
    47,85,31,45,54)(11,91,24,79,37,\\66)(13,87,29,53,22,56)(15,55,17,86,34,51)(16,84,
    96,28,68,21)(18,43,71,19,44,62)(20,49,88,33,50,52)(23,69,81,25,70,\\83)(30,48,46,
    90,99,100)(32,89)(36,93,97)(38,58,94)(39,65,72)(40,42)(41,76,73)(57,59),

  (1,2,43)(3,87,34)(4,77,42)(5,91,36)(6,64,19)(7,92,24)(8,49,20)(10,84,31)(11,78,
    61)(12,82,27)(13,86,38)(14,56,\\17)(15,60,53)(18,44,71)(21,47,45)(22,55,58)(23,90,
    99)(25,69,83)(26,79,97)(28,54,85)(29,51,94)(30,81,100)(32,57,\\89)(33,88,95)(35,
    66,93)(37,74,80)(41,73,76)(46,70,48)(50,67,52)(62,75,63),

  (1,2,91,36,11,6,64,79,97,24,75,63,66,93,37)(3,86,100,30,15,60,51,48,46,17,14,55,
    90,99,34)(4,9,77,40,42)(5,43,\\74,18,92,26,19,78,71,80,35,62,7,44,61)(8,49,72,20,
    82,95,33,65,88,27,67,52,39,50,12)(10,16,84,31,73,85,96,28,54,\\76,45,68,21,47,
    41)(13,87,23,58,83,29,53,81,38,69,22,56,70,94,25)(32,59,57,89,98),

  (1,3,75,14,6,60)(2,15,78,29,64,17,7,22,63,34,74,13)(4,57,9,59)(5,86,43,25,26,51,
    19,83,35,55,62,69)(8,76,27,31,\\95,41,12,54,67,73,82,47)(10,88,45,20,85,50)(11,81,
    44,56,24,70,18,87,37,23,71,53)(16,52,28,72,96,49,21,65,68,33,\\84,39)(30,79,94,92,
    46,66,58,80,99,91,38,61)(32,77)(36,90,93,48,97,100)(40,89,42,98),

  (1,3,7,83,37,34)(2,81,92,46,62,13)(4,59)(5,87,35,56,26,53)(6,60,74,69,11,15)(8,85,
    39,84,50,73)(9,57,77,98,42,\\32)(10,65,21,88,41,67)(12,54)(14,78,25,24,17,75)(16,
    52,96,49,68,33)(18,86,91,100,36,58)(19,22,63,23,61,30)(20,\\76,95,45,72,28)(27,
    31)(29,64,70,80,99,43)(38,71,51,79,48,97)(40,89)(44,55,66,90,93,94)(47,82),

  (1,6,75)(2,24,63,11,64,37)(3,34,14,17,60,15)(4,9)(5,71,35,18,26,44)(7,62,78,19,74,
    43)(8,88,67,20,95,50)(10,84,\\45,21,85,28)(12,82,27)(13,29,22)(16,96,68)(23,94,70,
    38,81,58)(25,87,69,56,83,53)(30,55,99,51,46,86)(31,41,47,76,\\54,73)(32,57)(33,72,
    49,39,52,65)(36,91,93,66,97,79)(42,77)(48,90,100)(61,92,80)(89,98),

  (1,6,75)(2,62,44)(3,34,22)(5,74,37)(7,24,35)(8,82,50)(10,84,68)(11,26,78)(12,88,
    67)(13,60,15)(14,17,29)(16,85,\\28)(18,64,43)(19,71,63)(20,95,27)(21,96,45)(23,83,
    30)(25,99,70)(32,57,59)(33,49,52)(36,91,80)(38,86,53)(40,77,42)\\(41,76,73)(46,81,
    69)(48,90,100)(51,56,94)(55,87,58)(61,97,79)(66,92,93),

  (1,6,75)(2,64,63)(3,60,14)(5,26,35)(7,74,78)(8,95,67)(10,85,45)(11,24,37)(12,82,
    27)(13,29,22)(15,17,34)(16,96,\\68)(18,71,44)(19,62,43)(20,88,50)(21,84,28)(23,81,
    70)(25,83,69)(30,46,99)(31,54,47)(33,52,49)(36,97,93)(38,94,58)\\(39,72,65)(41,73,
    76)(48,90,100)(51,55,86)(53,56,87)(61,92,80)(66,91,79),

  (1,14,74,83,24,15)(2,23,80,46,19,29)(3,78,69,37,17,6)(4,59)(5,56)(7,25,11,34,75,
    60)(8,10,72,84,88,76)(9,57,77,\\98,42,32)(12,31,82,54,27,47)(13,63,70,92,30,
    43)(16,33,68,49,96,52)(18,55,79,100,93,38)(20,73,67,45,39,21)(22,64,\\81,61,99,
    62)(26,87)(28,50,41,95,85,65)(35,53)(36,94,71,86,66,48)(40,89)(44,51,91,90,97,
    58),

    (1,14)(2,29,7,34)(3,6)(4,59,9,57)(5,25,19,55)(8,31,12,73)(10,20)(11,56,
    18,23)(13,78,17,63)(15,64,22,74)(16,72,\\21,33)(24,87,71,81)(26,83,62,86)(27,41,
    67,47)(28,49,68,39)(30,92,58,91)(32,77)(35,69,43,51)(36,48)(37,53,44,70)\\(38,79,
    46,80)(40,98,42,89)(45,50)(52,96,65,84)(54,82,76,95)(60,75)(61,94,66,99)(85,
    88)(90,97)(93,100),

    (1,23,37,48,19,51)(2,3,74,25,92,87)(4,98,42,89,9,57)(5,30,
    36,38,44,22)(6,81,11,90,62,55)(7,69,61,56,63,14)(8,54,\\72,21,33,41)(10,12)(13,26,
    46,97,94,18)(15,79,17,66,34,91)(16,88,68,20,96,50)(24,100,43,86,75,70)(27,
    45)(28,49,76,\\67,31,39)(29,35,99,93,58,71)(32,40)(47,65,84,52,73,95)(53,64,60,78,
    83,80)(59,77)(82,85),

     (1,23,97,87,75,70,36,56,6,81,93,53)(2,13,66,38,63,22,79,
    58,64,29,91,94)(3,26,100,43,14,5,90,62,60,35,48,19)(4,\\98,42,59)(7,30,61,34,78,
    99,80,17,74,46,92,15)(8,96,52,41,67,16,33,76,95,68,49,73)(9,57,40,89)(10,88,85,
    50,45,20)\\(11,69,24,25,37,83)(12,54,65,28,27,31,72,84,82,47,39,21)(18,51,71,55,
    44,86)(32,77),

     (1,23,6,81,75,70)(2,30,64,46,63,99)(3,26,60,35,14,5)(4,57)(7,
    13,74,29,78,22)(8,54,95,47,67,31)(9,98)(10,88,85,50,\\45,20)(11,51,24,55,37,
    86)(12,96,82,68,27,16)(15,79,17,66,34,91)(18,69,71,25,44,83)(19,90,62,100,43,
    48)(21,33,84,52,\\28,49)(32,77)(36,56,97,87,93,53)(38,92,94,80,58,61)(39,73,72,76,
    65,41)(40,59)(42,89),

     (1,24,61,66,71)(2,5,19,7,93)(3,30,53,70,58)(4,42,77,40,
    9)(6,37,92,91,44)(8,88,12,39,49)(10,96,41,47,84)(11,80,79,\\18,75)(13,17,69,90, 51)(14,99,87,81,94)(15,83,48,86,22)(16,76,54,21,45)(20,27,65,52,67)(23,38,60,46,
    56)(25,100,55,\\29,34)(26,62,74,36,64)(28,85,68,73,31)(32,59,98,89,57)(33,95,50,82,72)(35,43,78,97,63),

     (1,30,75,99,6,46)(2,48,63,100,64,90)(3,93,14,97,60,
    36)(4,32)(5,15,35,34,26,17)(7,94,78,38,74,58)(8,21,67,28,95,84)\\(9,98)(10,52,45,
    33,85,49)(11,22,37,29,24,13)(12,68,27,96,82,16)(18,83,44,25,71,69)(19,81,43,23,
    62,70)(20,54,50,31,\\88,47)(39,73,65,41,72,76)(40,59)(42,57)(51,80,86,92,55,
    61)(53,79,87,91,56,66)(77,89),

    (1,35,93,61,63,75,26,97,80,64,6,5,36,92,2)(3,
    29,99,90,56,14,13,46,48,53,60,22,30,100,87)(4,40,42,77,9)(7,71,37,\\79,62,78,18,
    24,91,19,74,44,11,66,43)(8,20,39,49,12,67,50,65,52,27,95,88,72,33,82)(10,73,47,
    16,21,45,41,54,68,28,\\85,76,31,96,84)(15,70,38,83,51,34,81,58,25,86,17,23,94,69,
    55)(32,98,89,59,57),

     (1,48)(2,23)(3,91)(4,89)(5,56)(6,90)(7,55)(8,10)(9,
    98)(11,38)(12,16)(13,92)(14,66)(15,36)(17,97)(18,25)(19,46)\\(20,21)(22,61)(24,
    94)(26,87)(27,68)(28,50)(29,80)(30,43)(31,52)(32,42)(33,47)(34,93)(35,53)(37,
    58)(39,73)(40,59)\\(41,65)(44,69)(45,67)(49,54)(51,78)(57,77)(60,79)(62,99)(63,
    70)(64,81)(71,83)(72,76)(74,86)(75,100)(82,96)(84,\\88)(85,95),

  (1,75,6)(2,63,64)(3,14,60)(5,35,26)(7,78,74)(8,67,95)(10,45,85)(11,37,24)(12,27,
    82)(13,22,29)(15,34,17)(16,68,\\96)(18,44,71)(19,43,62)(20,50,88)(21,28,84)(23,70,
    81)(25,69,83)(30,99,46)(31,47,54)(33,49,52)(36,93,97)(38,58,94)\\(39,65,72)(41,76,
    73)(48,100,90)(51,86,55)(53,87,56)(61,80,92)(66,79,91).

The character table of $G^{s_{11}}$:\\
\begin{tabular}{c|ccccccccccccccccccccc}
  & & & & & & & & & &10 & & & & & & & & & & 20&\\\hline

$\chi_{11}^{(1)}$&1&1&1&1&1&1&1&1&1&1&1&1&1&1&1&1&1&1&1&1&1
\\$\chi_{11}^{(2)}$&1&1&1&1&1&1&-1&-1&1&1&1&-1&-1&-1&-1&-1&1&-1&1&-1&1
\\$\chi_{11}^{(3)}$&1&1&1&A&A&A&-/A&-/A&/A&/A&/A&-1&-1&-A&-A&-A&1&-/A&/A&-1&A
\\$\chi_{11}^{(4)}$&1&1&1&/A&/A&/A&-A&-A&A&A&A&-1&-1&-/A&-/A&-/A&1&-A&A&-1&/A
\\$\chi_{11}^{(5)}$&1&1&1&A&A&A&/A&/A&/A&/A&/A&1&1&A&A&A&1&/A&/A&1&A
\\$\chi_{11}^{(6)}$&1&1&1&/A&/A&/A&A&A&A&A&A&1&1&/A&/A&/A&1&A&A&1&/A
\\$\chi_{11}^{(7)}$&4&.&1&.&1&-1&.&1&.&1&4&1&.&1&.&-2&-1&-2&-1&-2&4
\\$\chi_{11}^{(8)}$&4&.&1&.&1&-1&.&-1&.&1&4&-1&.&-1&.&2&-1&2&-1&2&4
\\$\chi_{11}^{(9)}$&4&.&1&.&/A&-/A&.&A&.&A&C&1&.&/A&.&B&-1&/B&-A&-2&/C
\\$\chi_{11}^{(10)}$&4&.&1&.&A&-A&.&/A&.&/A&/C&1&.&A&.&/B&-1&B&-/A&-2&C
\\$\chi_{11}^{(11)}$&4&.&1&.&/A&-/A&.&-A&.&A&C&-1&.&-/A&.&-B&-1&-/B&-A&2&/C
\\$\chi_{11}^{(12)}$&4&.&1&.&A&-A&.&-/A&.&/A&/C&-1&.&-A&.&-/B&-1&-B&-/A&2&C
\\$\chi_{11}^{(13)}$&5&1&-1&1&-1&.&-1&1&1&-1&5&1&-1&1&-1&1&.&1&.&1&5
\\$\chi_{11}^{(14)}$&5&1&-1&1&-1&.&1&-1&1&-1&5&-1&1&-1&1&-1&.&-1&.&-1&5
\\$\chi_{11}^{(15)}$&5&1&-1&/A&-/A&.&-A&A&A&-A&D&1&-1&/A&-/A&/A&.&A&.&1&/D
\\$\chi_{11}^{(16)}$&5&1&-1&A&-A&.&-/A&/A&/A&-/A&/D&1&-1&A&-A&A&.&/A&.&1&D
\\$\chi_{11}^{(17)}$&5&1&-1&/A&-/A&.&A&-A&A&-A&D&-1&1&-/A&/A&-/A&.&-A&.&-1&/D
\\$\chi_{11}^{(18)}$&5&1&-1&A&-A&.&/A&-/A&/A&-/A&/D&-1&1&-A&A&-A&.&-/A&.&-1&D
\\$\chi_{11}^{(19)}$&6&-2&.&-2&.&1&.&.&-2&.&6&.&.&.&.&.&1&.&1&.&6
\\$\chi_{11}^{(20)}$&6&-2&.&B&.&/A&.&.&/B&.&E&.&.&.&.&.&1&.&A&.&/E
\\$\chi_{11}^{(21)}$&6&-2&.&/B&.&A&.&.&B&.&/E&.&.&.&.&.&1&.&/A&.&E
\end{tabular}

\noindent \noindent where   A = E(3)$^2$
  = (-1-ER(-3))/2 = -1-b3;
B = -2*E(3)
  = 1-ER(-3) = 1-i3;
C = 4*E(3)$^2$
  = -2-2*ER(-3) = -2-2i3;
D = 5*E(3)$^2$
  = (-5-5*ER(-3))/2 = -5-5b3;
E = 6*E(3)$^2$
  = -3-3*ER(-3) = -3-3i3.

The generators of $G^{s_{12}}$ are:\\
 (  1, 19, 23, 84, 81, 82, 52, 53, 38, 45, 10,  4, 72, 77, 32, 86, 98, 85, 50, 42)
    (  2, 12, 71, 51, 96, 29, 80, 25, 43, 48, 16, 22, 79, 26, 75, 33, 24, 59, 92, 46
     )(  3, 95,100, 47, 15, 20, 61, 41, 44,  5)(  6, 67, 90, 60, 94, 66, 17, 63, 36,
     76, 14, 97, 70, 87, 35, 58,  7, 34, 40, 69)(  8, 83, 55, 93, 57,  9, 30, 64,
      37, 73)( 11, 65, 49, 54, 13, 88, 74, 62, 39, 56, 21, 89, 91, 68, 18, 27, 31,
      28, 99, 78)

The representatives of conjugacy classes of   $G^{s_{12}}$ are:\\
 (1),
 (1,4,23,77,81,86,52,85,38,42,10,19,72,84,32,82,98,53,50,45)(2,22,71,26,96,
    33,80,59,43,46,16,12,79,51,75,29,24,\\25,92,48)(3,95,100,47,15,20,61,41,44,5)(6,
    97,90,87,94,58,17,34,36,69,14,67,70,60,35,66,7,63,40,76)(8,83,55,93,57,\\9,30,64,
    37,73)(11,89,49,68,13,27,74,28,39,78,21,65,91,54,18,88,31,62,99,56),

  (1,10)(2,16)(4,19)(6,14)(7,17)(11,21)(12,22)(13,18)(23,72)(24,80)(25,59)(26,
    51)(27,88)(28,62)(29,33)(31,74)\\(32,81)(34,63)(35,94)(36,40)(38,50)(39,99)(42,
    45)(43,92)(46,48)(49,91)(52,98)(53,85)(54,68)(56,78)(58,66)(60,\\87)(65,89)(67,
    97)(69,76)(70,90)(71,79)(75,96)(77,84)(82,86),

  (1,19,23,84,81,82,52,53,38,45,10,4,72,77,32,86,98,85,50,42)(2,12,71,51,96,29,80,
    25,43,48,16,22,79,26,75,33,24,\\59,92,46)(3,95,100,47,15,20,61,41,44,5)(6,67,90,
    60,94,66,17,63,36,76,14,97,70,87,35,58,7,34,40,69)(8,83,55,93,57,\\9,30,64,37,
    73)(11,65,49,54,13,88,74,62,39,56,21,89,91,68,18,27,31,28,99,78),

  (1,23,81,52,38,10,72,32,98,50)(2,71,96,80,43,16,79,75,24,92)(3,100,15,61,44)(4,77,
    86,85,42,19,84,82,53,45)(5,\\95,47,20,41)(6,90,94,17,36,14,70,35,7,40)(8,55,57,30,
    37)(9,64,73,83,93)(11,49,13,74,39,21,91,18,31,99)(12,51,29,\\25,48,22,26,33,59,
    46)(27,28,78,65,54,88,62,56,89,68)(34,69,67,60,66,63,76,97,87,58),

  (1,32,38,23,98,10,81,50,72,52)(2,75,43,71,24,16,96,92,79,80)(3,15,44,100,61)(4,82,
    42,77,53,19,86,45,84,85)(5,\\47,41,95,20)(6,35,36,90,7,14,94,40,70,17)(8,57,37,55,
    30)(9,73,93,64,83)(11,18,39,49,31,21,13,99,91,74)(12,33,48,\\51,59,22,29,46,26,
    25)(27,56,54,28,89,88,78,68,62,65)(34,97,66,69,87,63,67,58,76,60),

  (1,38,98,81,72)(2,43,24,96,79)(3,44,61,15,100)(4,42,53,86,84)(5,41,20,47,95)(6,36,
    7,94,70)(8,37,30,57,55)(9,\\93,83,73,64)(10,50,52,32,23)(11,39,31,13,91)(12,48,59,
    29,26)(14,40,17,35,90)(16,92,80,75,71)(18,49,21,99,74)(19,\\45,85,82,77)(22,46,25,
    33,51)(27,54,89,78,62)(28,88,68,65,56)(34,66,87,67,76)(58,60,97,69,63),

  (1,42,50,85,98,86,32,77,72,4,10,45,38,53,52,82,81,84,23,19)(2,46,92,59,24,33,75,
    26,79,22,16,48,43,25,80,29,96,\\51,71,12)(3,5,44,41,61,20,15,47,100,95)(6,69,40,
    34,7,58,35,87,70,97,14,76,36,63,17,66,94,60,90,67)(8,73,37,64,30,\\9,57,93,55,
    83)(11,78,99,28,31,27,18,68,91,89,21,56,39,62,74,88,13,54,49,65),

  (1,45,50,53,98,82,32,84,72,19,10,42,38,85,52,86,81,77,23,4)(2,48,92,25,24,29,75,
    51,79,12,16,46,43,59,80,33,96,\\26,71,22)(3,5,44,41,61,20,15,47,100,95)(6,76,40,
    63,7,66,35,60,70,67,14,69,36,34,17,58,94,87,90,97)(8,73,37,64,30,\\9,57,93,55,
    83)(11,56,99,62,31,88,18,54,91,65,21,78,39,28,74,27,13,68,49,89),

  (1,50,98,32,72,10,38,52,81,23)(2,92,24,75,79,16,43,80,96,71)(3,44,61,15,100)(4,45,
    53,82,84,19,42,85,86,77)(5,\\41,20,47,95)(6,40,7,35,70,14,36,17,94,90)(8,37,30,57,
    55)(9,93,83,73,64)(11,99,31,18,91,21,39,74,13,49)(12,46,59,\\33,26,22,48,25,29,
    51)(27,68,89,56,62,88,54,65,78,28)(34,58,87,97,76,63,66,60,67,69),

  (1,52,72,50,81,10,98,23,38,32)(2,80,79,92,96,16,24,71,43,75)(3,61,100,44,15)(4,85,
    84,45,86,19,53,77,42,82)(5,\\20,95,41,47)(6,17,70,40,94,14,7,90,36,35)(8,30,55,37,
    57)(9,83,64,93,73)(11,74,91,99,13,21,31,49,39,18)(12,25,26,\\46,29,22,59,51,48,
    33)(27,65,62,68,78,88,89,28,54,56)(34,60,76,58,67,63,87,69,66,97),

  (1,53,32,19,38,86,23,45,98,84,10,85,81,4,50,82,72,42,52,77)(2,25,75,12,43,33,71,
    48,24,51,16,59,96,22,92,29,79,\\46,80,26)(3,41,15,95,44,20,100,5,61,47)(6,63,35,
    67,36,58,90,76,7,60,14,34,94,97,40,66,70,69,17,87)(8,64,57,83,37,\\9,55,73,30,
    93)(11,62,18,65,39,27,49,56,31,54,21,28,13,89,99,88,91,78,74,68),

  (1,72,81,98,38)(2,79,96,24,43)(3,100,15,61,44)(4,84,86,53,42)(5,95,47,20,41)(6,70,
    94,7,36)(8,55,57,30,37)(9,64,\\73,83,93)(10,23,32,52,50)(11,91,13,31,39)(12,26,29,
    59,48)(14,90,35,17,40)(16,71,75,80,92)(18,74,99,21,49)(19,77,\\82,85,45)(22,51,33,
    25,46)(27,62,78,89,54)(28,56,65,68,88)(34,76,67,87,66)(58,63,69,97,60),

  (1,77,52,42,72,82,50,4,81,85,10,84,98,45,23,86,38,19,32,53)(2,26,80,46,79,29,92,
    22,96,59,16,51,24,48,71,33,43,\\12,75,25)(3,47,61,5,100,20,44,95,15,41)(6,87,17,
    69,70,66,40,97,94,34,14,60,7,76,90,58,36,67,35,63)(8,93,30,73,55,\\9,37,83,57,
    64)(11,68,74,78,91,88,99,89,13,28,21,54,31,56,49,27,39,65,18,62),

  (1,81,38,72,98)(2,96,43,79,24)(3,15,44,100,61)(4,86,42,84,53)(5,47,41,95,20)(6,94,
    36,70,7)(8,57,37,55,30)(9,73,\\93,64,83)(10,32,50,23,52)(11,13,39,91,31)(12,29,48,
    26,59)(14,35,40,90,17)(16,75,92,71,80)(18,99,49,74,21)(19,82,\\45,77,85)(22,33,46,
    51,25)(27,78,54,62,89)(28,65,88,56,68)(34,67,66,76,87)(58,69,60,63,97),

  (1,82,10,86)(2,29,16,33)(3,20)(4,98,19,52)(5,15)(6,66,14,58)(7,67,17,97)(8,9)(11,
    88,21,27)(12,80,22,24)(13,56,\\18,78)(23,53,72,85)(25,79,59,71)(26,92,51,43)(28,
    49,62,91)(30,83)(31,65,74,89)(32,42,81,45)(34,90,63,70)(35,69,\\94,76)(36,87,40,
    60)(37,93)(38,77,50,84)(39,68,99,54)(41,100)(44,47)(46,96,48,75)(55,64)(57,
    73)(61,95),

    (1,84,52,45,72,86,50,19,81,53,10,77,98,42,23,82,38,4,32,85)(2,51,
    80,48,79,33,92,12,96,25,16,26,24,46,71,29,43,\\22,75,59)(3,47,61,5,100,20,44,95,
    15,41)(6,60,17,76,70,58,40,67,94,63,14,87,7,69,90,66,36,97,35,34)(8,93,30,73,55,\\
    9,37,83,57,64)(11,54,74,56,91,27,99,65,13,62,21,68,31,78,49,88,39,89,18,28),

  (1,85,32,4,38,82,23,42,98,77,10,53,81,19,50,86,72,45,52,84)(2,59,75,22,43,29,71,
    46,24,26,16,25,96,12,92,33,79,\\48,80,51)(3,41,15,95,44,20,100,5,61,47)(6,34,35,
    97,36,66,90,69,7,87,14,63,94,67,40,58,70,76,17,60)(8,64,57,83,37,\\9,55,73,30,
    93)(11,28,18,89,39,88,49,78,31,68,21,62,13,65,99,27,91,56,74,54),

  (1,86,10,82)(2,33,16,29)(3,20)(4,52,19,98)(5,15)(6,58,14,66)(7,97,17,67)(8,9)(11,
    27,21,88)(12,24,22,80)(13,78,\\18,56)(23,85,72,53)(25,71,59,79)(26,43,51,92)(28,
    91,62,49)(30,83)(31,89,74,65)(32,45,81,42)(34,70,63,90)(35,76,\\94,69)(36,60,40,
    87)(37,93)(38,84,50,77)(39,54,99,68)(41,100)(44,47)(46,75,48,96)(55,64)(57,
    73)(61,95),

    (1,98,72,38,81)(2,24,79,43,96)(3,61,100,44,15)(4,53,84,42,86)(5,
    20,95,41,47)(6,7,70,36,94)(8,30,55,37,57)(9,83,\\64,93,73)(10,52,23,50,32)(11,31,
    91,39,13)(12,59,26,48,29)(14,17,90,40,35)(16,80,71,92,75)(18,21,74,49,99)(19,85,\\
    77,45,82)(22,25,51,46,33)(27,89,62,54,78)(28,68,56,88,65)(34,87,76,66,67)(58,97,
    63,60,69).

The character table of $G^{s_{12}}$:\\
\begin{tabular}{c|cccccccccccccccccccc}
  & & & & & & & & & &10 & & & & & & & & & &20\\\hline

$\chi_{12}^{(1)}$&1&1&1&1&1&1&1&1&1&1&1&1&1&1&1&1&1&1&1&1
\\$\chi_{12}^{(2)}$&1&-1&1&-1&1&1&1&-1&-1&1&1&-1&1&-1&1&-1&-1&-1&-1&1
\\$\chi_{12}^{(3)}$&1&A&1&A&B&/A&/B&/A&/A&/B&A&B&B&/B&/A&1&/B&B&1&A
\\$\chi_{12}^{(4)}$&1&-A&1&-A&B&/A&/B&-/A&-/A&/B&A&-B&B&-/B&/A&-1&-/B&-B&-1&A
\\$\chi_{12}^{(5)}$&1&B&1&B&/A&/B&A&/B&/B&A&B&/A&/A&A&/B&1&A&/A&1&B
\\$\chi_{12}^{(6)}$&1&-B&1&-B&/A&/B&A&-/B&-/B&A&B&-/A&/A&-A&/B&-1&-A&-/A&-1&B
\\$\chi_{12}^{(7)}$&1&/B&1&/B&A&B&/A&B&B&/A&/B&A&A&/A&B&1&/A&A&1&/B
\\$\chi_{12}^{(8)}$&1&-/B&1&-/B&A&B&/A&-B&-B&/A&/B&-A&A&-/A&B&-1&-/A&-A&-1&/B
\\$\chi_{12}^{(9)}$&1&/A&1&/A&/B&A&B&A&A&B&/A&/B&/B&B&A&1&B&/B&1&/A
\\$\chi_{12}^{(10)}$&1&-/A&1&-/A&/B&A&B&-A&-A&B&/A&-/B&/B&-B&A&-1&-B&-/B&-1&/A
\\$\chi_{12}^{(11)}$&1&C&-1&-C&-1&-1&1&C&-C&-1&-1&C&1&-C&1&-C&C&-C&C&1
\\$\chi_{12}^{(12)}$&1&-C&-1&C&-1&-1&1&-C&C&-1&-1&-C&1&C&1&C&-C&C&-C&1
\\$\chi_{12}^{(13)}$&1&D&-1&-D&-B&-/A&/B&-/D&/D&-/B&-A&E&B&/E&/A&-C&-/E&-E&C&A
\\$\chi_{12}^{(14)}$&1&-D&-1&D&-B&-/A&/B&/D&-/D&-/B&-A&-E&B&-/E&/A&C&/E&E&-C&A
\\$\chi_{12}^{(15)}$&1&E&-1&-E&-/A&-/B&A&-/E&/E&-A&-B&-/D&/A&-D&/B&-C&D&/D&C&B
\\$\chi_{12}^{(16)}$&1&-E&-1&E&-/A&-/B&A&/E&-/E&-A&-B&/D&/A&D&/B&C&-D&-/D&-C&B
\\$\chi_{12}^{(17)}$&1&-/E&-1&/E&-A&-B&/A&E&-E&-/A&-/B&D&A&/D&B&-C&-/D&-D&C&/B
\\$\chi_{12}^{(18)}$&1&/E&-1&-/E&-A&-B&/A&-E&E&-/A&-/B&-D&A&-/D&B&C&/D&D&-C&/B
\\$\chi_{12}^{(19)}$&1&-/D&-1&/D&-/B&-A&B&D&-D&-B&-/A&-/E&/B&-E&A&-C&E&/E&C&/A
\\$\chi_{12}^{(20)}$&1&/D&-1&-/D&-/B&-A&B&-D&D&-B&-/A&/E&/B&E&A&C&-E&-/E&-C&/A

\end{tabular}

\noindent \noindent where   A = E(5)$^4$; B = E(5)$^3$; C = -E(4)
  = -ER(-1) = -i;
D = -E(20); E = -E(20)$^{17}$.

The generators of $G^{s_{13}}$ are:\\
 (  1,  4, 23, 77, 81, 86, 52, 85, 38, 42, 10, 19, 72, 84, 32, 82, 98, 53, 50, 45)
    (  2, 22, 71, 26, 96, 33, 80, 59, 43, 46, 16, 12, 79, 51, 75, 29, 24, 25, 92, 48
     )(  3, 95,100, 47, 15, 20, 61, 41, 44,  5)(  6, 97, 90, 87, 94, 58, 17, 34, 36,
     69, 14, 67, 70, 60, 35, 66,  7, 63, 40, 76)(  8, 83, 55, 93, 57,  9, 30, 64,
      37, 73)( 11, 89, 49, 68, 13, 27, 74, 28, 39, 78, 21, 65, 91, 54, 18, 88, 31,
      62, 99, 56)

The representatives of conjugacy classes of   $G^{s_{13}}$ are:\\
(1),
(1,4,23,77,81,86,52,85,38,42,10,19,72,84,32,82,98,53,50,45)(2,22,71,26,96,
    33,80,59,43,46,16,12,79,51,75,29,24,\\25,92,48)(3,95,100,47,15,20,61,41,44,5)(6,
    97,90,87,94,58,17,34,36,69,14,67,70,60,35,66,7,63,40,76)(8,83,55,93,57,\\9,30,64,
    37,73)(11,89,49,68,13,27,74,28,39,78,21,65,91,54,18,88,31,62,99,56),

  (1,10)(2,16)(4,19)(6,14)(7,17)(11,21)(12,22)(13,18)(23,72)(24,80)(25,59)(26,
    51)(27,88)(28,62)(29,33)(31,74)\\(32,81)(34,63)(35,94)(36,40)(38,50)(39,99)(42,
    45)(43,92)(46,48)(49,91)(52,98)(53,85)(54,68)(56,78)(58,66)(60,\\87)(65,89)(67,
    97)(69,76)(70,90)(71,79)(75,96)(77,84)(82,86),

  (1,19,23,84,81,82,52,53,38,45,10,4,72,77,32,86,98,85,50,42)(2,12,71,51,96,29,80,
    25,43,48,16,22,79,26,75,33,24,\\59,92,46)(3,95,100,47,15,20,61,41,44,5)(6,67,90,
    60,94,66,17,63,36,76,14,97,70,87,35,58,7,34,40,69)(8,83,55,93,57,\\9,30,64,37,
    73)(11,65,49,54,13,88,74,62,39,56,21,89,91,68,18,27,31,28,99,78),

  (1,23,81,52,38,10,72,32,98,50)(2,71,96,80,43,16,79,75,24,92)(3,100,15,61,44)(4,77,
    86,85,42,19,84,82,53,45)(5,\\95,47,20,41)(6,90,94,17,36,14,70,35,7,40)(8,55,57,30,
    37)(9,64,73,83,93)(11,49,13,74,39,21,91,18,31,99)(12,51,29,\\25,48,22,26,33,59,
    46)(27,28,78,65,54,88,62,56,89,68)(34,69,67,60,66,63,76,97,87,58),

  (1,32,38,23,98,10,81,50,72,52)(2,75,43,71,24,16,96,92,79,80)(3,15,44,100,61)(4,82,
    42,77,53,19,86,45,84,85)(5,\\47,41,95,20)(6,35,36,90,7,14,94,40,70,17)(8,57,37,55,
    30)(9,73,93,64,83)(11,18,39,49,31,21,13,99,91,74)(12,33,48,\\51,59,22,29,46,26,
    25)(27,56,54,28,89,88,78,68,62,65)(34,97,66,69,87,63,67,58,76,60),

  (1,38,98,81,72)(2,43,24,96,79)(3,44,61,15,100)(4,42,53,86,84)(5,41,20,47,95)(6,36,
    7,94,70)(8,37,30,57,55)(9,\\93,83,73,64)(10,50,52,32,23)(11,39,31,13,91)(12,48,59,
    29,26)(14,40,17,35,90)(16,92,80,75,71)(18,49,21,99,74)\\(19,45,85,82,77)(22,46,25,
    33,51)(27,54,89,78,62)(28,88,68,65,56)(34,66,87,67,76)(58,60,97,69,63),

  (1,42,50,85,98,86,32,77,72,4,10,45,38,53,52,82,81,84,23,19)(2,46,92,59,24,33,75,
    26,79,22,16,48,43,25,80,29,96,\\51,71,12)(3,5,44,41,61,20,15,47,100,95)(6,69,40,
    34,7,58,35,87,70,97,14,76,36,63,17,66,94,60,90,67)(8,73,37,64,30,\\9,57,93,55,
    83)(11,78,99,28,31,27,18,68,91,89,21,56,39,62,74,88,13,54,49,65),

  (1,45,50,53,98,82,32,84,72,19,10,42,38,85,52,86,81,77,23,4)(2,48,92,25,24,29,75,
    51,79,12,16,46,43,59,80,33,96,\\26,71,22)(3,5,44,41,61,20,15,47,100,95)(6,76,40,
    63,7,66,35,60,70,67,14,69,36,34,17,58,94,87,90,97)(8,73,37,64,30,\\9,57,93,55,
    83)(11,56,99,62,31,88,18,54,91,65,21,78,39,28,74,27,13,68,49,89),

  (1,50,98,32,72,10,38,52,81,23)(2,92,24,75,79,16,43,80,96,71)(3,44,61,15,100)(4,45,
    53,82,84,19,42,85,86,77)(5,\\41,20,47,95)(6,40,7,35,70,14,36,17,94,90)(8,37,30,57,
    55)(9,93,83,73,64)(11,99,31,18,91,21,39,74,13,49)(12,46,59,\\33,26,22,48,25,29,
    51)(27,68,89,56,62,88,54,65,78,28)(34,58,87,97,76,63,66,60,67,69),

  (1,52,72,50,81,10,98,23,38,32)(2,80,79,92,96,16,24,71,43,75)(3,61,100,44,15)(4,85,
    84,45,86,19,53,77,42,82)(5,\\20,95,41,47)(6,17,70,40,94,14,7,90,36,35)(8,30,55,37,
    57)(9,83,64,93,73)(11,74,91,99,13,21,31,49,39,18)(12,25,26,\\46,29,22,59,51,48,
    33)(27,65,62,68,78,88,89,28,54,56)(34,60,76,58,67,63,87,69,66,97),

  (1,53,32,19,38,86,23,45,98,84,10,85,81,4,50,82,72,42,52,77)(2,25,75,12,43,33,71,
    48,24,51,16,59,96,22,92,29,79,\\46,80,26)(3,41,15,95,44,20,100,5,61,47)(6,63,35,
    67,36,58,90,76,7,60,14,34,94,97,40,66,70,69,17,87)(8,64,57,83,37,\\9,55,73,30,
    93)(11,62,18,65,39,27,49,56,31,54,21,28,13,89,99,88,91,78,74,68),

  (1,72,81,98,38)(2,79,96,24,43)(3,100,15,61,44)(4,84,86,53,42)(5,95,47,20,41)(6,70,
    94,7,36)(8,55,57,30,37)(9,\\64,73,83,93)(10,23,32,52,50)(11,91,13,31,39)(12,26,29,
    59,48)(14,90,35,17,40)(16,71,75,80,92)(18,74,99,21,49)(19,\\77,82,85,45)(22,51,33,
    25,46)(27,62,78,89,54)(28,56,65,68,88)(34,76,67,87,66)(58,63,69,97,60),

  (1,77,52,42,72,82,50,4,81,85,10,84,98,45,23,86,38,19,32,53)(2,26,80,46,79,29,92,
    22,96,59,16,51,24,48,71,33,43,\\12,75,25)(3,47,61,5,100,20,44,95,15,41)(6,87,17,
    69,70,66,40,97,94,34,14,60,7,76,90,58,36,67,35,63)(8,93,30,73,55,\\9,37,83,57,
    64)(11,68,74,78,91,88,99,89,13,28,21,54,31,56,49,27,39,65,18,62),

  (1,81,38,72,98)(2,96,43,79,24)(3,15,44,100,61)(4,86,42,84,53)(5,47,41,95,20)(6,94,
    36,70,7)(8,57,37,55,30)(9,73,\\93,64,83)(10,32,50,23,52)(11,13,39,91,31)(12,29,48,
    26,59)(14,35,40,90,17)(16,75,92,71,80)(18,99,49,74,21)(19,82,\\45,77,85)(22,33,46,
    51,25)(27,78,54,62,89)(28,65,88,56,68)(34,67,66,76,87)(58,69,60,63,97),

  (1,82,10,86)(2,29,16,33)(3,20)(4,98,19,52)(5,15)(6,66,14,58)(7,67,17,97)(8,9)(11,
    88,21,27)(12,80,22,24)(13,56,\\18,78)(23,53,72,85)(25,79,59,71)(26,92,51,43)(28,
    49,62,91)(30,83)(31,65,74,89)(32,42,81,45)(34,90,63,70)(35,69,\\94,76)(36,87,40,
    60)(37,93)(38,77,50,84)(39,68,99,54)(41,100)(44,47)(46,96,48,75)(55,64)(57,
    73)(61,95),

    (1,84,52,45,72,86,50,19,81,53,10,77,98,42,23,82,38,4,32,85)(2,51,
    80,48,79,33,92,12,96,25,16,26,24,46,71,29,43,\\22,75,59)(3,47,61,5,100,20,44,95,
    15,41)(6,60,17,76,70,58,40,67,94,63,14,87,7,69,90,66,36,97,35,34)(8,93,30,73,55,\\
    9,37,83,57,64)(11,54,74,56,91,27,99,65,13,62,21,68,31,78,49,88,39,89,18,28),

  (1,85,32,4,38,82,23,42,98,77,10,53,81,19,50,86,72,45,52,84)(2,59,75,22,43,29,71,
    46,24,26,16,25,96,12,92,33,79,\\48,80,51)(3,41,15,95,44,20,100,5,61,47)(6,34,35,
    97,36,66,90,69,7,87,14,63,94,67,40,58,70,76,17,60)(8,64,57,83,37,\\9,55,73,30,
    93)(11,28,18,89,39,88,49,78,31,68,21,62,13,65,99,27,91,56,74,54),

  (1,86,10,82)(2,33,16,29)(3,20)(4,52,19,98)(5,15)(6,58,14,66)(7,97,17,67)(8,9)(11,
    27,21,88)(12,24,22,80)(13,78,\\18,56)(23,85,72,53)(25,71,59,79)(26,43,51,92)(28,
    91,62,49)(30,83)(31,89,74,65)(32,45,81,42)(34,70,63,90)(35,76,\\94,69)(36,60,40,
    87)(37,93)(38,84,50,77)(39,54,99,68)(41,100)(44,47)(46,75,48,96)(55,64)(57,
    73)(61,95),

     (1,98,72,38,81)(2,24,79,43,96)(3,61,100,44,15)(4,53,84,42,86)(5,
    20,95,41,47)(6,7,70,36,94)(8,30,55,37,57)(9,83,\\64,93,73)(10,52,23,50,32)(11,31,
    91,39,13)(12,59,26,48,29)(14,17,90,40,35)(16,80,71,92,75)(18,21,74,49,99)(19,85,\\
    77,45,82)(22,25,51,46,33)(27,89,62,54,78)(28,68,56,88,65)(34,87,76,66,67)(58,97,
    63,60,69).

The character table of $G^{s_{13}}$:\\
\begin{tabular}{c|cccccccccccccccccccc}
  & & & & & & & & & &10 & & & & & & & & & &20\\\hline
$\chi_{13}^{(1)}$&1&1&1&1&1&1&1&1&1&1&1&1&1&1&1&1&1&1&1&1
\\$\chi_{13}^{(2)}$&1&-1&1&-1&1&1&1&-1&-1&1&1&-1&1&-1&1&-1&-1&-1&-1&1
\\$\chi_{13}^{(3)}$&1&A&1&A&B&/A&/B&/A&/A&/B&A&B&B&/B&/A&1&/B&B&1&A
\\$\chi_{13}^{(4)}$&1&-A&1&-A&B&/A&/B&-/A&-/A&/B&A&-B&B&-/B&/A&-1&-/B&-B&-1&A
\\$\chi_{13}^{(5)}$&1&B&1&B&/A&/B&A&/B&/B&A&B&/A&/A&A&/B&1&A&/A&1&B
\\$\chi_{13}^{(6)}$&1&-B&1&-B&/A&/B&A&-/B&-/B&A&B&-/A&/A&-A&/B&-1&-A&-/A&-1&B
\\$\chi_{13}^{(7)}$&1&/B&1&/B&A&B&/A&B&B&/A&/B&A&A&/A&B&1&/A&A&1&/B
\\$\chi_{13}^{(8)}$&1&-/B&1&-/B&A&B&/A&-B&-B&/A&/B&-A&A&-/A&B&-1&-/A&-A&-1&/B
\\$\chi_{13}^{(9)}$&1&/A&1&/A&/B&A&B&A&A&B&/A&/B&/B&B&A&1&B&/B&1&/A
\\$\chi_{13}^{(10)}$&1&-/A&1&-/A&/B&A&B&-A&-A&B&/A&-/B&/B&-B&A&-1&-B&-/B&-1&/A
\\$\chi_{13}^{(11)}$&1&C&-1&-C&-1&-1&1&C&-C&-1&-1&C&1&-C&1&-C&C&-C&C&1
\\$\chi_{13}^{(12)}$&1&-C&-1&C&-1&-1&1&-C&C&-1&-1&-C&1&C&1&C&-C&C&-C&1
\\$\chi_{13}^{(13)}$&1&D&-1&-D&-B&-/A&/B&-/D&/D&-/B&-A&E&B&/E&/A&-C&-/E&-E&C&A
\\$\chi_{13}^{(14)}$&1&-D&-1&D&-B&-/A&/B&/D&-/D&-/B&-A&-E&B&-/E&/A&C&/E&E&-C&A
\\$\chi_{13}^{(15)}$&1&E&-1&-E&-/A&-/B&A&-/E&/E&-A&-B&-/D&/A&-D&/B&-C&D&/D&C&B
\\$\chi_{13}^{(16)}$&1&-E&-1&E&-/A&-/B&A&/E&-/E&-A&-B&/D&/A&D&/B&C&-D&-/D&-C&B
\\$\chi_{13}^{(17)}$&1&-/E&-1&/E&-A&-B&/A&E&-E&-/A&-/B&D&A&/D&B&-C&-/D&-D&C&/B
\\$\chi_{13}^{(18)}$&1&/E&-1&-/E&-A&-B&/A&-E&E&-/A&-/B&-D&A&-/D&B&C&/D&D&-C&/B
\\$\chi_{13}^{(19)}$&1&-/D&-1&/D&-/B&-A&B&D&-D&-B&-/A&-/E&/B&-E&A&-C&E&/E&C&/A
\\$\chi_{13}^{(20)}$&1&/D&-1&-/D&-/B&-A&B&-D&D&-B&-/A&/E&/B&E&A&C&-E&-/E&-C&/A

\end{tabular}

\noindent \noindent where   A = E(5)$^4$; B = E(5)$^3$; C = E(4)
  = ER(-1) = i;
D = E(20); E = E(20)$^{17}$.

The generators of $G^{s_{14}}$ are:\\
 (  1, 23, 81, 52, 38, 10, 72, 32, 98, 50)(  2, 71, 96, 80, 43, 16, 79, 75, 24, 92)
    (  3,100, 15, 61, 44)(  4, 77, 86, 85, 42, 19, 84, 82, 53, 45)(  5, 95, 47, 20,
      41)(  6, 90, 94, 17, 36, 14, 70, 35,  7, 40)(  8, 55, 57, 30, 37)
    (  9, 64, 73, 83, 93)( 11, 49, 13, 74, 39, 21, 91, 18, 31, 99)( 12, 51, 29, 25,
      48, 22, 26, 33, 59, 46)( 27, 28, 78, 65, 54, 88, 62, 56, 89, 68)
    ( 34, 69, 67, 60, 66, 63, 76, 97, 87, 58),
  (  1,  4, 23, 77, 81, 86, 52, 85, 38, 42, 10, 19, 72, 84, 32, 82, 98, 53, 50, 45)
    (  2, 22, 71, 26, 96, 33, 80, 59, 43, 46, 16, 12, 79, 51, 75, 29, 24, 25, 92, 48
     )(  3, 95,100, 47, 15, 20, 61, 41, 44,  5)(  6, 97, 90, 87, 94, 58, 17, 34, 36,
     69, 14, 67, 70, 60, 35, 66,  7, 63, 40, 76)(  8, 83, 55, 93, 57,  9, 30, 64,
      37, 73)( 11, 89, 49, 68, 13, 27, 74, 28, 39, 78, 21, 65, 91, 54, 18, 88, 31,
      62, 99, 56).

The representatives of conjugacy classes of   $G^{s_{14}}$ are:\\
 (1),
  (1,4,23,77,81,86,52,85,38,42,10,19,72,84,32,82,98,53,50,45)(2,22,71,26,96,
    33,80,59,43,46,16,12,79,51,75,29,24,25,\\92,48)(3,95,100,47,15,20,61,41,44,5)(6,
    97,90,87,94,58,17,34,36,69,14,67,70,60,35,66,7,63,40,76)(8,83,55,93,57,9,30,\\64,
    37,73)(11,89,49,68,13,27,74,28,39,78,21,65,91,54,18,88,31,62,99,56),

  (1,10)(2,16)(4,19)(6,14)(7,17)(11,21)(12,22)(13,18)(23,72)(24,80)(25,59)(26,
    51)(27,88)(28,62)(29,33)(31,74)(32,\\81)(34,63)(35,94)(36,40)(38,50)(39,99)(42,
    45)(43,92)(46,48)(49,91)(52,98)(53,85)(54,68)(56,78)(58,66)(60,87)(65,\\89)(67,
    97)(69,76)(70,90)(71,79)(75,96)(77,84)(82,86),

  (1,19,23,84,81,82,52,53,38,45,10,4,72,77,32,86,98,85,50,42)(2,12,71,51,96,29,80,
    25,43,48,16,22,79,26,75,33,24,59,\\92,46)(3,95,100,47,15,20,61,41,44,5)(6,67,90,
    60,94,66,17,63,36,76,14,97,70,87,35,58,7,34,40,69)(8,83,55,93,57,9,30,\\64,37,
    73)(11,65,49,54,13,88,74,62,39,56,21,89,91,68,18,27,31,28,99,78),

  (1,23,81,52,38,10,72,32,98,50)(2,71,96,80,43,16,79,75,24,92)(3,100,15,61,44)(4,77,
    86,85,42,19,84,82,53,45)(5,95,\\47,20,41)(6,90,94,17,36,14,70,35,7,40)(8,55,57,30,
    37)(9,64,73,83,93)(11,49,13,74,39,21,91,18,31,99)(12,51,29,25,48,\\22,26,33,59,
    46)(27,28,78,65,54,88,62,56,89,68)(34,69,67,60,66,63,76,97,87,58),

  (1,32,38,23,98,10,81,50,72,52)(2,75,43,71,24,16,96,92,79,80)(3,15,44,100,61)(4,82,
    42,77,53,19,86,45,84,85)(5,47,\\41,95,20)(6,35,36,90,7,14,94,40,70,17)(8,57,37,55,
    30)(9,73,93,64,83)(11,18,39,49,31,21,13,99,91,74)(12,33,48,51,59,\\22,29,46,26,
    25)(27,56,54,28,89,88,78,68,62,65)(34,97,66,69,87,63,67,58,76,60),

  (1,38,98,81,72)(2,43,24,96,79)(3,44,61,15,100)(4,42,53,86,84)(5,41,20,47,95)(6,36,
    7,94,70)(8,37,30,57,55)(9,93,\\83,73,64)(10,50,52,32,23)(11,39,31,13,91)(12,48,59,
    29,26)(14,40,17,35,90)(16,92,80,75,71)(18,49,21,99,74)(19,45,\\85,82,77)(22,46,25,
    33,51)(27,54,89,78,62)(28,88,68,65,56)(34,66,87,67,76)(58,60,97,69,63),

  (1,42,50,85,98,86,32,77,72,4,10,45,38,53,52,82,81,84,23,19)(2,46,92,59,24,33,75,
    26,79,22,16,48,43,25,80,29,96,\\51,71,12)(3,5,44,41,61,20,15,47,100,95)(6,69,40,
    34,7,58,35,87,70,97,14,76,36,63,17,66,94,60,90,67)(8,73,37,64,30,\\9,57,93,55,
    83)(11,78,99,28,31,27,18,68,91,89,21,56,39,62,74,88,13,54,49,65),

  (1,45,50,53,98,82,32,84,72,19,10,42,38,85,52,86,81,77,23,4)(2,48,92,25,24,29,75,
    51,79,12,16,46,43,59,80,33,96,\\26,71,22)(3,5,44,41,61,20,15,47,100,95)(6,76,40,
    63,7,66,35,60,70,67,14,69,36,34,17,58,94,87,90,97)(8,73,37,64,30,\\9,57,93,55,
    83)(11,56,99,62,31,88,18,54,91,65,21,78,39,28,74,27,13,68,49,89),

  (1,50,98,32,72,10,38,52,81,23)(2,92,24,75,79,16,43,80,96,71)(3,44,61,15,100)(4,45,
    53,82,84,19,42,85,86,77)(5,\\41,20,47,95)(6,40,7,35,70,14,36,17,94,90)(8,37,30,57,
    55)(9,93,83,73,64)(11,99,31,18,91,21,39,74,13,49)(12,46,59,\\33,26,22,48,25,29,
    51)(27,68,89,56,62,88,54,65,78,28)(34,58,87,97,76,63,66,60,67,69),

  (1,52,72,50,81,10,98,23,38,32)(2,80,79,92,96,16,24,71,43,75)(3,61,100,44,15)(4,85,
    84,45,86,19,53,77,42,82)(5,\\20,95,41,47)(6,17,70,40,94,14,7,90,36,35)(8,30,55,37,
    57)(9,83,64,93,73)(11,74,91,99,13,21,31,49,39,18)(12,25,26,\\46,29,22,59,51,48,
    33)(27,65,62,68,78,88,89,28,54,56)(34,60,76,58,67,63,87,69,66,97),

  (1,53,32,19,38,86,23,45,98,84,10,85,81,4,50,82,72,42,52,77)(2,25,75,12,43,33,71,
    48,24,51,16,59,96,22,92,29,79,\\46,80,26)(3,41,15,95,44,20,100,5,61,47)(6,63,35,
    67,36,58,90,76,7,60,14,34,94,97,40,66,70,69,17,87)(8,64,57,83,37,\\9,55,73,30,
    93)(11,62,18,65,39,27,49,56,31,54,21,28,13,89,99,88,91,78,74,68),

  (1,72,81,98,38)(2,79,96,24,43)(3,100,15,61,44)(4,84,86,53,42)(5,95,47,20,41)(6,70,
    94,7,36)(8,55,57,30,37)(9,64,\\73,83,93)(10,23,32,52,50)(11,91,13,31,39)(12,26,29,
    59,48)(14,90,35,17,40)(16,71,75,80,92)(18,74,99,21,49)(19,77,\\82,85,45)(22,51,33,
    25,46)(27,62,78,89,54)(28,56,65,68,88)(34,76,67,87,66)(58,63,69,97,60),

  (1,77,52,42,72,82,50,4,81,85,10,84,98,45,23,86,38,19,32,53)(2,26,80,46,79,29,92,
    22,96,59,16,51,24,48,71,33,43,\\12,75,25)(3,47,61,5,100,20,44,95,15,41)(6,87,17,
    69,70,66,40,97,94,34,14,60,7,76,90,58,36,67,35,63)(8,93,30,73,55,\\9,37,83,57,
    64)(11,68,74,78,91,88,99,89,13,28,21,54,31,56,49,27,39,65,18,62),

  (1,81,38,72,98)(2,96,43,79,24)(3,15,44,100,61)(4,86,42,84,53)(5,47,41,95,20)(6,94,
    36,70,7)(8,57,37,55,30)(9,73,\\93,64,83)(10,32,50,23,52)(11,13,39,91,31)(12,29,48,
    26,59)(14,35,40,90,17)(16,75,92,71,80)(18,99,49,74,21)(19,82,\\45,77,85)(22,33,46,
    51,25)(27,78,54,62,89)(28,65,88,56,68)(34,67,66,76,87)(58,69,60,63,97),

  (1,82,10,86)(2,29,16,33)(3,20)(4,98,19,52)(5,15)(6,66,14,58)(7,67,17,97)(8,9)(11,
    88,21,27)(12,80,22,24)(13,56,\\18,78)(23,53,72,85)(25,79,59,71)(26,92,51,43)(28,
    49,62,91)(30,83)(31,65,74,89)(32,42,81,45)(34,90,63,70)(35,69,\\94,76)(36,87,40,
    60)(37,93)(38,77,50,84)(39,68,99,54)(41,100)(44,47)(46,96,48,75)(55,64)(57,
    73)(61,95),

     (1,84,52,45,72,86,50,19,81,53,10,77,98,42,23,82,38,4,32,85)(2,51,
    80,48,79,33,92,12,96,25,16,26,24,46,71,29,43,\\22,75,59)(3,47,61,5,100,20,44,95,
    15,41)(6,60,17,76,70,58,40,67,94,63,14,87,7,69,90,66,36,97,35,34)(8,93,30,73,55,\\
    9,37,83,57,64)(11,54,74,56,91,27,99,65,13,62,21,68,31,78,49,88,39,89,18,28),

  (1,85,32,4,38,82,23,42,98,77,10,53,81,19,50,86,72,45,52,84)(2,59,75,22,43,29,71,
    46,24,26,16,25,96,12,92,33,79,\\48,80,51)(3,41,15,95,44,20,100,5,61,47)(6,34,35,
    97,36,66,90,69,7,87,14,63,94,67,40,58,70,76,17,60)(8,64,57,83,37,\\9,55,73,30,
    93)(11,28,18,89,39,88,49,78,31,68,21,62,13,65,99,27,91,56,74,54),

  (1,86,10,82)(2,33,16,29)(3,20)(4,52,19,98)(5,15)(6,58,14,66)(7,97,17,67)(8,9)(11,
    27,21,88)(12,24,22,80)(13,78,\\18,56)(23,85,72,53)(25,71,59,79)(26,43,51,92)(28,
    91,62,49)(30,83)(31,89,74,65)(32,45,81,42)(34,70,63,90)(35,76,\\94,69)(36,60,40,
    87)(37,93)(38,84,50,77)(39,54,99,68)(41,100)(44,47)(46,75,48,96)(55,64)(57,
    73)(61,95),

     (1,98,72,38,81)(2,24,79,43,96)(3,61,100,44,15)(4,53,84,42,86)(5,
    20,95,41,47)(6,7,70,36,94)(8,30,55,37,57)(9,83,\\64,93,73)(10,52,23,50,32)(11,31,
    91,39,13)(12,59,26,48,29)(14,17,90,40,35)(16,80,71,92,75)(18,21,74,49,99)(19,85,\\
    77,45,82)(22,25,51,46,33)(27,89,62,54,78)(28,68,56,88,65)(34,87,76,66,67)(58,97,
    63,60,69).

The character table of $G^{s_{14}}$:\\
\begin{tabular}{c|cccccccccccccccccccc}
  & & & & & & & & & &10 & & & & & & & & & &20\\\hline

$\chi_{14}^{(1)}$&1&1&1&1&1&1&1&1&1&1&1&1&1&1&1&1&1&1&1&1
\\$\chi_{14}^{(2)}$&1&-1&1&-1&1&1&1&-1&-1&1&1&-1&1&-1&1&-1&-1&-1&-1&1
\\$\chi_{14}^{(3)}$&1&A&-1&-A&-1&-1&1&A&-A&-1&-1&A&1&-A&1&-A&A&-A&A&1
\\$\chi_{14}^{(4)}$&1&-A&-1&A&-1&-1&1&-A&A&-1&-1&-A&1&A&1&A&-A&A&-A&1
\\$\chi_{14}^{(5)}$&1&B&1&B&-C&-/B&-/C&/B&/B&-/C&-B&C&-C&/C&-/B&-1&/C&C&-1&-B
\\$\chi_{14}^{(6)}$&1&C&1&C&-/B&-/C&-B&/C&/C&-B&-C&/B&-/B&B&-/C&-1&B&/B&-1&-C
\\$\chi_{14}^{(7)}$&1&/C&1&/C&-B&-C&-/B&C&C&-/B&-/C&B&-B&/B&-C&-1&/B&B&-1&-/C
\\$\chi_{14}^{(8)}$&1&/B&1&/B&-/C&-B&-C&B&B&-C&-/B&/C&-/C&C&-B&-1&C&/C&-1&-/B
\\$\chi_{14}^{(9)}$&1&-/B&1&-/B&-/C&-B&-C&-B&-B&-C&-/B&-/C&-/C&-C&-B&1&-C&-/C&1&-/B
\\$\chi_{14}^{(10)}$&1&-/C&1&-/C&-B&-C&-/B&-C&-C&-/B&-/C&-B&-B&-/B&-C&1&-/B&-B&1&-/C
\\$\chi_{14}^{(11)}$&1&-C&1&-C&-/B&-/C&-B&-/C&-/C&-B&-C&-/B&-/B&-B&-/C&1&-B&-/B&1&-C
\\$\chi_{14}^{(12)}$&1&-B&1&-B&-C&-/B&-/C&-/B&-/B&-/C&-B&-C&-C&-/C&-/B&1&-/C&-C&1&-B
\\$\chi_{14}^{(13)}$&1&D&-1&-D&/C&B&-C&-/D&/D&C&/B&-/E&-/C&-E&-B&-A&E&/E&A&-/B
\\$\chi_{14}^{(14)}$&1&-/D&-1&/D&C&/B&-/C&D&-D&/C&B&E&-C&/E&-/B&-A&-/E&-E&A&-B
\\$\chi_{14}^{(15)}$&1&E&-1&-E&/B&/C&-B&-/E&/E&B&C&D&-/B&/D&-/C&-A&-/D&-D&A&-C
\\$\chi_{14}^{(16)}$&1&-/E&-1&/E&B&C&-/B&E&-E&/B&/C&-/D&-B&-D&-C&-A&D&/D&A&-/C
\\$\chi_{14}^{(17)}$&1&/E&-1&-/E&B&C&-/B&-E&E&/B&/C&/D&-B&D&-C&A&-D&-/D&-A&-/C
\\$\chi_{14}^{(18)}$&1&-E&-1&E&/B&/C&-B&/E&-/E&B&C&-D&-/B&-/D&-/C&A&/D&D&-A&-C
\\$\chi_{14}^{(19)}$&1&/D&-1&-/D&C&/B&-/C&-D&D&/C&B&-E&-C&-/E&-/B&A&/E&E&-A&-B
\\$\chi_{14}^{(20)}$&1&-D&-1&D&/C&B&-C&/D&-/D&C&/B&/E&-/C&E&-B&A&-E&-/E&-A&-/B

\end{tabular}

\noindent where   A = -E(4)
  = -ER(-1) = -i;
B = -E(5); C = -E(5)$^2$; D = -E(20); E = -E(20)$^{13}$.

The generators of $G^{s_{15}}$ are:\\
(  1, 81, 38, 72, 98)(  2, 96, 43, 79, 24)(  3, 15, 44,100, 61)
    (  4, 86, 42, 84, 53)(  5, 47, 41, 95, 20)(  6, 94, 36, 70,  7)(  8, 57, 37, 55,
     30)(  9, 73, 93, 64, 83)( 10, 32, 50, 23, 52)( 11, 13, 39, 91, 31)
    ( 12, 29, 48, 26, 59)( 14, 35, 40, 90, 17)( 16, 75, 92, 71, 80)( 18, 99, 49, 74,
     21)( 19, 82, 45, 77, 85)( 22, 33, 46, 51, 25)( 27, 78, 54, 62, 89)
    ( 28, 65, 88, 56, 68)( 34, 67, 66, 76, 87)( 58, 69, 60, 63, 97),
  (  3, 69)(  4, 75)(  5, 74)(  6, 39)(  7, 13)(  8, 59)(  9, 45)( 10, 51)( 11, 70)
    ( 12, 57)( 15, 60)( 16, 53)( 18, 41)( 19, 64)( 20, 49)( 21, 47)( 22, 50)
    ( 23, 33)( 25, 32)( 26, 30)( 27, 34)( 29, 37)( 31, 36)( 42, 71)( 44, 63)
    ( 46, 52)( 48, 55)( 54, 66)( 58, 61)( 62, 76)( 67, 78)( 73, 77)( 80, 84)
    ( 82, 83)( 85, 93)( 86, 92)( 87, 89)( 91, 94)( 95, 99)( 97,100),
  (  2,  4)(  5, 14)(  6, 18)(  7, 21)(  8, 12)(  9, 16)( 10, 22)( 17, 20)( 23, 51)
    ( 24, 53)( 25, 52)( 26, 55)( 27, 56)( 28, 54)( 29, 57)( 30, 59)( 32, 33)
    ( 34, 58)( 35, 47)( 36, 49)( 37, 48)( 40, 41)( 42, 43)( 46, 50)( 60, 66)
    ( 62, 65)( 63, 76)( 64, 71)( 67, 69)( 68, 78)( 70, 74)( 73, 75)( 79, 84)
    ( 80, 83)( 86, 96)( 87, 97)( 88, 89)( 90, 95)( 92, 93)( 94, 99),
  (  1,  2, 10, 73, 81, 96, 32, 93, 38, 43, 50, 64, 72, 79, 23, 83, 98, 24, 52,  9)
    (  3,  5, 54, 21, 15, 47, 62, 18, 44, 41, 89, 99,100, 95, 27, 49, 61, 20, 78, 74
     )(  4, 57, 84, 30, 86, 37, 53,  8, 42, 55)(  6, 68, 13, 76, 94, 28, 39, 87, 36,
     65, 91, 34, 70, 88, 31, 67,  7, 56, 11, 66)( 12, 16, 46, 85, 29, 75, 51, 19,
      48, 92, 25, 82, 26, 71, 22, 45, 59, 80, 33, 77)( 14, 69, 90, 97, 35, 60, 17,
      58, 40, 63).

The representatives of conjugacy classes of   $G^{s_{15}}$ are:\\
 (1),
 (3,69)(4,75)(5,74)(6,39)(7,13)(8,59)(9,45)(10,51)(11,70)(12,57)(15,60)(16,
    53)(18,41)(19,64)(20,49)(21,47)(22,\\50)(23,33)(25,32)(26,30)(27,34)(29,37)(31,
    36)(42,71)(44,63)(46,52)(48,55)(54,66)(58,61)(62,76)(67,78)(73,77)(80,\\84)(82,
    83)(85,93)(86,92)(87,89)(91,94)(95,99)(97,100),

  (2,4,73,77,75)(3,67,68,78,69)(5,70,11,74,14)(6,39,18,40,41)(7,13,21,35,47)(8,30,
    55,37,57)(9,45,16,24,53)(10,23,\\32,52,50)(12,29,48,26,59)(15,66,28,54,60)(17,20,
    36,31,49)(19,71,43,42,64)(22,46,25,33,51)(27,58,61,34,56)(44,76,\\65,62,63)(79,84,
    83,82,80)(85,92,96,86,93)(87,88,89,97,100)(90,95,94,91,99),

  (2,73,75,4,77)(3,68,69,67,78)(5,11,14,70,74)(6,18,41,39,40)(7,21,47,13,35)(8,55,
    57,30,37)(9,16,53,45,24)(10,32,\\50,23,52)(12,48,59,29,26)(15,28,60,66,54)(17,36,
    49,20,31)(19,43,64,71,42)(22,25,51,46,33)(27,61,56,58,34)(44,65,\\63,76,62)(79,83,
    80,84,82)(85,96,93,92,86)(87,89,100,88,97)(90,94,99,95,91),

  (1,2,10,73,81,96,32,93,38,43,50,64,72,79,23,83,98,24,52,9)(3,5,54,21,15,47,62,18,
    44,41,89,99,100,95,27,49,61,20,\\78,74)(4,57,84,30,86,37,53,8,42,55)(6,68,13,76,
    94,28,39,87,36,65,91,34,70,88,31,67,7,56,11,66)(12,16,46,85,29,75,\\51,19,48,92,
    25,82,26,71,22,45,59,80,33,77)(14,69,90,97,35,60,17,58,40,63),

  (1,2,10,19,72,79,23,77,81,96,32,82,98,24,52,85,38,43,50,45)(3,90,97,95,100,35,60,
    47,15,17,58,20,61,40,63,41,44,\\14,69,5)(4,51,73,12,84,33,83,26,86,25,93,29,53,46,
    9,59,42,22,64,48)(6,87,99,27,70,66,21,62,94,34,49,78,7,76,18,89,\\36,67,74,54)(8,
    80,30,71,55,92,37,75,57,16)(11,88,31,65,91,28,39,68,13,56),

  (1,2,12,4,81,96,29,86,38,43,48,42,72,79,26,84,98,24,59,53)(3,13,60,21,15,39,63,18,
    44,91,97,99,100,31,58,49,61,\\11,69,74)(5,88,95,28,47,56,20,65,41,68)(6,89,35,76,
    94,27,40,87,36,78,90,34,70,54,17,67,7,62,14,66)(8,92,25,73,57,\\71,22,93,37,80,33,
    64,55,16,46,83,30,75,51,9)(10,82,23,85,32,45,52,19,50,77),

  (1,2,12,71,72,79,26,75,81,96,29,80,98,24,59,92,38,43,48,16)(3,74,88,95,100,99,28,
    47,15,21,56,20,61,49,65,41,44,\\18,68,5)(4,51,82,30,84,33,85,37,86,25,45,8,53,46,
    19,55,42,22,77,57)(6,63,91,27,70,69,13,62,94,97,31,78,7,60,39,89,\\36,58,11,54)(9,
    52,83,23,64,50,93,32,73,10)(14,66,17,67,90,34,40,87,35,76),

  (1,2,81,96,38,43,72,79,98,24)(3,7,60,39,44,94,97,31,61,70,69,13,15,6,63,91,100,36,
    58,11)(4,57,92,29,42,55,80,\\26,53,8,75,12,86,37,71,48,84,30,16,59)(5,54,21,76,41,
    89,99,34,20,78,74,66,47,62,18,87,95,27,49,67)(9,52,77,51,93,\\32,19,22,83,23,45,
    46,73,10,85,25,64,50,82,33)(14,88,35,56,40,68,90,28,17,65),

  (1,2,81,96,38,43,72,79,98,24)(3,13,60,6,44,91,97,36,61,11,69,7,15,39,63,94,100,31,
    58,70)(4,12,92,37,42,48,80,\\30,53,59,75,57,86,29,71,55,84,26,16,8)(5,66,21,62,41,
    87,99,27,20,67,74,54,47,76,18,89,95,34,49,78)(9,46,77,10,93,\\25,19,50,83,33,45,
    52,73,51,85,32,64,22,82,23)(14,88,35,56,40,68,90,28,17,65),

  (1,8,10,22,12)(2,16,19,9,4)(5,18,6,14,13)(7,17,11,20,21)(23,51,26,72,55)(24,80,85,
    83,53)(25,59,98,30,52)(27,\\62,78,89,54)(28,88,68,65,56)(29,81,57,32,33)(31,95,74,
    70,90)(34,67,66,76,87)(35,39,47,99,94)(36,40,91,41,49)\\(37,50,46,48,38)(42,43,92,
    45,93)(58,97,63,60,69)(64,84,79,71,77)(73,86,96,75,82),

  (1,8,98,30,72,55,38,37,81,57)(2,16,24,80,79,71,43,92,96,75)(3,58,61,97,100,63,44,
    60,15,69)(4,82,53,19,84,85,\\42,77,86,45)(5,70,20,36,95,94,41,6,47,7)(9,93,83,73,
    64)(10,26,52,48,23,29,50,12,32,59)(11,90,31,40,91,35,39,14,\\13,17)(18,49,21,99,
    74)(22,46,25,33,51)(27,67,89,34,62,87,54,76,78,66)(28,88,68,65,56),

  (1,8,38,37,98,30,81,57,72,55)(2,64,43,9,24,93,96,83,79,73)(3,65,44,56,61,28,15,88,
    100,68)(4,86,42,84,53)(5,\\90,41,14,20,40,47,17,95,35)(6,49,36,21,7,99,94,74,70,
    18)(10,59,50,29,52,26,32,12,23,48)(11,13,39,91,31)(16,19,\\92,45,80,85,75,82,71,
    77)(22,33,46,51,25)(27,97,54,69,89,63,78,58,62,60)(34,67,66,76,87),

  (1,8,50,33,48)(2,64,92,42,77)(3,89,34,97,68)(4,82,79,73,16)(5,70,13,99,35)(6,91,
    74,90,41)(7,39,49,40,47)(9,\\80,53,19,43)(10,25,12,72,55)(11,18,14,20,36)(15,27,
    67,58,28)(17,95,94,31,21)(22,29,98,30,32)(23,46,26,81,57)\\(24,93,75,86,45)(37,52,
    51,59,38)(44,78,66,69,65)(54,76,60,88,100)(56,61,62,87,63)(71,84,85,96,83),

  (1,8,32,51,29)(2,82,53,80,73)(3,58,62,56,67)(4,16,93,96,45)(5,13,17,36,21)(6,49,
    95,31,40)(7,99,41,91,35)(9,\\24,19,84,71)(10,46,12,98,30)(11,90,94,74,20)(14,70,
    18,47,39)(15,69,89,68,66)(22,26,38,37,23)(25,48,81,57,50)\\(27,28,76,44,60)(33,59,
    72,55,52)(34,61,97,54,88)(42,92,83,79,85)(43,77,86,75,64)(63,78,65,87,100),

  (1,8,72,55,81,57,98,30,38,37)(2,82,79,85,96,45,24,19,43,77)(3,89,100,54,15,27,61,
    62,44,78)(4,64,84,73,86,83,\\53,93,42,9)(5,20,95,41,47)(6,14,70,90,94,35,7,17,36,
    40)(10,29,23,59,32,48,52,12,50,26)(11,18,91,74,13,99,31,21,\\39,49)(16,80,71,92,
    75)(22,25,51,46,33)(28,76,56,67,65,87,68,66,88,34)(58,97,63,60,69),

  (1,8,81,57,38,37,72,55,98,30)(2,86,96,42,43,84,79,53,24,4)(3,66,15,76,44,87,100,
    34,61,67)(5,18,47,99,41,49,\\95,74,20,21)(6,91,94,31,36,11,70,13,7,39)(9,80,73,16,
    93,75,64,92,83,71)(10,48,32,26,50,59,23,12,52,29)(14,90,35,\\17,40)(19,77,82,85,
    45)(22,51,33,25,46)(27,62,78,89,54)(28,69,65,60,88,63,56,97,68,58),

  (1,9,52,24,98,83,23,79,72,64,50,43,38,93,32,96,81,73,10,2)(3,74,78,20,61,49,27,95,
    100,99,89,41,44,18,62,47,\\15,21,54,5)(4,55,42,8,53,37,86,30,84,57)(6,66,11,56,7,
    67,31,88,70,34,91,65,36,87,39,28,94,76,13,68)(12,77,33,80,\\59,45,22,71,26,82,25,
    92,48,19,51,75,29,85,46,16)(14,63,40,58,17,60,35,97,90,69),

  (1,9,37,45,98,83,57,82,72,64,8,19,38,93,30,85,81,73,55,77)(2,33,43,51,24,22,96,46,
    79,25)(3,6,44,36,61,7,15,\\94,100,70)(4,10,92,48,53,52,75,29,84,23,16,12,42,50,80,
    59,86,32,71,26)(5,34,91,69,20,87,39,58,95,76,13,97,41,66,\\11,63,47,67,31,60)(14,
    56,49,27,17,88,99,89,90,65,18,62,40,28,21,54,35,68,74,78),

  (1,10,12,8,22)(2,19,4,16,9)(5,6,13,18,14)(7,11,21,17,20)(23,26,55,51,72)(24,85,53,
    80,83)(25,98,52,59,30)(27,\\78,54,62,89)(28,68,56,88,65)(29,57,33,81,32)(31,74,90,
    95,70)(34,66,87,67,76)(35,47,94,39,99)(36,91,49,40,41)\\(37,46,38,50,48)(42,92,93,
    43,45)(58,63,69,97,60)(64,79,77,84,71)(73,96,82,86,75),

  (1,19)(2,32,4,33)(3,11)(5,62,14,65)(6,69,18,67)(7,58,21,34)(8,73,12,75)(9,59,16,
    30)(10,53,22,24)(13,15)(17,\\28,20,54)(23,42,51,43)(25,79,52,84)(26,80,55,83)(27,
    40,56,41)(29,92,57,93)(31,61)(35,88,47,89)(36,63,49,76)\\(37,64,48,71)(38,45)(39,
    44)(46,96,50,86)(60,99,66,94)(68,95,78,90)(70,97,74,87)(72,77)(81,82)(85,98)(91,
    100),

     (1,19)(2,33,4,32)(3,11)(5,65,14,62)(6,67,18,69)(7,34,21,58)(8,75,12,
    73)(9,30,16,59)(10,24,22,53)(13,15)(17,\\54,20,28)(23,43,51,42)(25,84,52,79)(26,
    83,55,80)(27,41,56,40)(29,93,57,92)(31,61)(35,89,47,88)(36,76,49,63)\\(37,71,48,
    64)(38,45)(39,44)(46,86,50,96)(60,94,66,99)(68,90,78,95)(70,87,74,97)(72,77)(81,
    82)(85,98)(91,100),

    (1,38,98,81,72)(2,43,24,96,79)(3,44,61,15,100)(4,42,53,86,
    84)(5,41,20,47,95)(6,36,7,94,70)(8,37,30,57,55)(9,\\93,83,73,64)(10,50,52,32,
    23)(11,39,31,13,91)(12,48,59,29,26)(14,40,17,35,90)(16,92,80,75,71)(18,49,21,99,
    74)\\(19,45,85,82,77)(22,46,25,33,51)(27,54,89,78,62)(28,88,68,65,56)(34,66,87,67,
    76)(58,60,97,69,63),

     (1,72,81,98,38)(2,79,96,24,43)(3,100,15,61,44)(4,84,86,
    53,42)(5,95,47,20,41)(6,70,94,7,36)(8,55,57,30,37)(9,\\64,73,83,93)(10,23,32,52,
    50)(11,91,13,31,39)(12,26,29,59,48)(14,90,35,17,40)(16,71,75,80,92)(18,74,99,21,
    49)\\(19,77,82,85,45)(22,51,33,25,46)(27,62,78,89,54)(28,56,65,68,88)(34,76,67,87,
    66)(58,63,69,97,60),

    (1,81,38,72,98)(2,96,43,79,24)(3,15,44,100,61)(4,86,42,
    84,53)(5,47,41,95,20)(6,94,36,70,7)(8,57,37,55,30)(9,\\73,93,64,83)(10,32,50,23,
    52)(11,13,39,91,31)(12,29,48,26,59)(14,35,40,90,17)(16,75,92,71,80)(18,99,49,74,
    21)\\(19,82,45,77,85)(22,33,46,51,25)(27,78,54,62,89)(28,65,88,56,68)(34,67,66,76,
    87)(58,69,60,63,97),

    (1,98,72,38,81)(2,24,79,43,96)(3,61,100,44,15)(4,53,84,
    42,86)(5,20,95,41,47)(6,7,70,36,94)(8,30,55,37,57)(9,\\83,64,93,73)(10,52,23,50,
    32)(11,31,91,39,13)(12,59,26,48,29)(14,17,90,40,35)(16,80,71,92,75)(18,21,74,49,
    99)\\(19,85,77,45,82)(22,25,51,46,33)(27,89,62,54,78)(28,68,56,88,65)(34,87,76,66,
    67)(58,97,63,60,69).

The character table of $G^{s_{15}}$:\\
\begin{tabular}{c|ccccccccccccccccccccc}
  & & & & & & & & & &10 & & & & & & & & & &20 &\\\hline
$\chi_{15}^{(1)}$&1&1&1&1&1&1&1&1&1&1&1&1&1&1&1&1&1&1&1&1&1
\\$\chi_{15}^{(2)}$&1&1&1&1&-1&-1&-1&-1&-1&-1&1&1&1&1&1&1&1&-1&-1&1&-1
\\$\chi_{15}^{(3)}$&1&-1&1&1&C&-C&-C&C&-C&C&1&-1&-1&1&1&-1&-1&-C&C&1&-C
\\$\chi_{15}^{(4)}$&1&-1&1&1&-C&C&C&-C&C&-C&1&-1&-1&1&1&-1&-1&C&-C&1&C
\\$\chi_{15}^{(5)}$&4&.&-1&-1&.&.&.&.&.&.&-1&.&.&4&-1&.&.&.&.&-1&.
\\$\chi_{15}^{(6)}$&4&.&-1&-1&.&.&.&.&.&.&-1&.&.&-1&4&.&.&.&.&-1&.
\\$\chi_{15}^{(7)}$&4&.&A&*A&.&.&.&.&.&.&*B&.&.&-1&-1&.&.&.&.&B&.
\\$\chi_{15}^{(8)}$&4&.&*A&A&.&.&.&.&.&.&B&.&.&-1&-1&.&.&.&.&*B&.
\\$\chi_{15}^{(9)}$&4&.&B&*B&.&.&.&.&.&.&A&.&.&-1&-1&.&.&.&.&*A&.
\\$\chi_{15}^{(10)}$&4&.&*B&B&.&.&.&.&.&.&*A&.&.&-1&-1&.&.&.&.&A&.
\\$\chi_{15}^{(11)}$&5&-1&.&.&D&-E&-D&E&/E&-/E&.&G&F&.&.&/F&/G&/D&-/D&.&-C
\\$\chi_{15}^{(12)}$&5&-1&.&.&-/D&/E&/D&-/E&-E&E&.&/G&/F&.&.&F&G&-D&D&.&-C
\\$\chi_{15}^{(13)}$&5&-1&.&.&E&/D&-E&-/D&-D&D&.&F&/G&.&.&G&/F&/E&-/E&.&-C
\\$\chi_{15}^{(14)}$&5&-1&.&.&-/E&-D&/E&D&/D&-/D&.&/F&G&.&.&/G&F&-E&E&.&-C
\\$\chi_{15}^{(15)}$&5&-1&.&.&/E&D&-/E&-D&-/D&/D&.&/F&G&.&.&/G&F&E&-E&.&C
\\$\chi_{15}^{(16)}$&5&-1&.&.&-E&-/D&E&/D&D&-D&.&F&/G&.&.&G&/F&-/E&/E&.&C
\\$\chi_{15}^{(17)}$&5&-1&.&.&/D&-/E&-/D&/E&E&-E&.&/G&/F&.&.&F&G&D&-D&.&C
\\$\chi_{15}^{(18)}$&5&-1&.&.&-D&E&D&-E&-/E&/E&.&G&F&.&.&/F&/G&-/D&/D&.&C
\\$\chi_{15}^{(19)}$&5&1&.&.&F&/G&F&/G&G&G&.&-/G&-/F&.&.&-F&-G&/F&/F&.&-1
\\$\chi_{15}^{(20)}$&5&1&.&.&G&F&G&F&/F&/F&.&-F&-/G&.&.&-G&-/F&/G&/G&.&-1
\\$\chi_{15}^{(21)}$&5&1&.&.&/G&/F&/G&/F&F&F&.&-/F&-G&.&.&-/G&-F&G&G&.&-1
\\$\chi_{15}^{(22)}$&5&1&.&.&/F&G&/F&G&/G&/G&.&-G&-F&.&.&-/F&-/G&F&F&.&-1
\\$\chi_{15}^{(23)}$&5&1&.&.&-/F&-G&-/F&-G&-/G&-/G&.&-G&-F&.&.&-/F&-/G&-F&-F&.&1
\\$\chi_{15}^{(24)}$&5&1&.&.&-/G&-/F&-/G&-/F&-F&-F&.&-/F&-G&.&.&-/G&-F&-G&-G&.&1
\\$\chi_{15}^{(25)}$&5&1&.&.&-G&-F&-G&-F&-/F&-/F&.&-F&-/G&.&.&-G&-/F&-/G&-/G&.&1
\\$\chi_{15}^{(26)}$&5&1&.&.&-F&-/G&-F&-/G&-G&-G&.&-/G&-/F&.&.&-F&-G&-/F&-/F&.&1
\end{tabular}

\begin{tabular}{c|ccccc}
  & & & & &\\\hline
$\chi_{15}^{(1)}$&1&1&1&1&1
\\$\chi_{15}^{(2)}$&-1&1&1&1&1
\\$\chi_{15}^{(3)}$&C&1&1&1&1
\\$\chi_{15}^{(4)}$&-C&1&1&1&1
\\$\chi_{15}^{(5)}$&.&4&4&4&4
\\$\chi_{15}^{(6)}$&.&4&4&4&4
\\$\chi_{15}^{(7)}$&.&4&4&4&4
\\$\chi_{15}^{(8)}$&.&4&4&4&4
\\$\chi_{15}^{(9)}$&.&4&4&4&4
\\$\chi_{15}^{(10)}$&.&4&4&4&4
\\$\chi_{15}^{(11)}$&C&H&/H&I&/I
\\$\chi_{15}^{(12)}$&C&/H&H&/I&I
\\$\chi_{15}^{(13)}$&C&I&/I&/H&H
\\$\chi_{15}^{(14)}$&C&/I&I&H&/H
\\$\chi_{15}^{(15)}$&-C&/I&I&H&/H
\\$\chi_{15}^{(16)}$&-C&I&/I&/H&H
\\$\chi_{15}^{(17)}$&-C&/H&H&/I&I
\\$\chi_{15}^{(18)}$&-C&H&/H&I&/I
\\$\chi_{15}^{(19)}$&-1&/H&H&/I&I
\\$\chi_{15}^{(20)}$&-1&I&/I&/H&H
\\$\chi_{15}^{(21)}$&-1&/I&I&H&/H
\\$\chi_{15}^{(22)}$&-1&H&/H&I&/I
\\$\chi_{15}^{(23)}$&1&H&/H&I&/I
\\$\chi_{15}^{(24)}$&1&/I&I&H&/H
\\$\chi_{15}^{(25)}$&1&I&/I&/H&H
\\$\chi_{15}^{(26)}$&1&/H&H&/I&I

\end{tabular}

\noindent \noindent where   A = -E(5)-2*E(5)$^2$-2*E(5)$^3$-E(5)$^4$
  = (3+ER(5))/2 = 2+b5;
B = 2*E(5)$^2$+2*E(5)$^3$
  = -1-ER(5) = -1-r5;
C = -E(4)
  = -ER(-1) = -i;
D = -E(20); E = -E(20)$^{13}$; F = -E(5); G = -E(5)$^2$; H =
5*E(5)$^2$; I = 5*E(5).

The generators of $G^{s_{16}}$ are:\\
 (  1, 82, 10, 86)(  2, 29, 16, 33)(  3, 20)(  4, 98, 19, 52)(  5, 15)
    (  6, 66, 14, 58)(  7, 67, 17, 97)(  8,  9)( 11, 88, 21, 27)( 12, 80, 22, 24)
    ( 13, 56, 18, 78)( 23, 53, 72, 85)( 25, 79, 59, 71)( 26, 92, 51, 43)
    ( 28, 49, 62, 91)( 30, 83)( 31, 65, 74, 89)( 32, 42, 81, 45)( 34, 90, 63, 70)
    ( 35, 69, 94, 76)( 36, 87, 40, 60)( 37, 93)( 38, 77, 50, 84)( 39, 68, 99, 54)
    ( 41,100)( 44, 47)( 46, 96, 48, 75)( 55, 64)( 57, 73)( 61, 95),
  (  3, 20)(  4, 17)(  6, 14)(  7, 19)(  8,  9)( 11, 12)( 13, 18)( 21, 22)( 23, 31)
    ( 24, 27)( 25, 34)( 26, 36)( 28, 38)( 32, 35)( 39, 75)( 40, 51)( 42, 69)
    ( 43, 60)( 44, 47)( 45, 76)( 46, 68)( 48, 54)( 49, 77)( 50, 62)( 52, 67)
    ( 53, 65)( 56, 78)( 57, 73)( 58, 66)( 59, 63)( 70, 71)( 72, 74)( 79, 90)
    ( 80, 88)( 81, 94)( 84, 91)( 85, 89)( 87, 92)( 96, 99)( 97, 98),
  (  3, 47, 20, 44)(  4, 68, 17, 46)(  5, 15)(  6, 66, 14, 58)(  7, 48, 19, 54)
    (  8, 73,  9, 57)( 11, 42, 12, 69)( 13, 78, 18, 56)( 21, 45, 22, 76)
    ( 23, 28, 31, 38)( 24, 35, 27, 32)( 25, 36, 34, 26)( 30, 37)( 39, 67, 75, 52)
    ( 40, 63, 51, 59)( 41, 61)( 43, 71, 60, 70)( 49, 65, 77, 53)( 50, 72, 62, 74)
    ( 55, 64)( 79, 87, 90, 92)( 80, 94, 88, 81)( 83, 93)( 84, 85, 91, 89)( 95,100)
    ( 96, 98, 99, 97), (  2,  6, 13)(  3,  8,  5)(  4, 12, 17)(  7, 19, 22)
    (  9, 15, 20)( 14, 18, 16)( 23, 51, 76)( 24, 67, 52)( 25, 36, 68)( 26, 69, 72)
    ( 28, 48, 31)( 29, 66, 56)( 30, 61, 37)( 32, 70, 77)( 33, 58, 78)( 34, 50, 42)
    ( 35, 53, 43)( 38, 45, 63)( 39, 71, 60)( 40, 54, 59)( 44, 57, 64)( 46, 74, 62)
    ( 47, 73, 55)( 49, 75, 65)( 79, 87, 99)( 80, 97, 98)( 81, 90, 84)( 83, 95, 93)
    ( 85, 92, 94)( 89, 91, 96), (  1,  2)(  3, 20)(  6, 18)(  8,  9)( 10, 16)
    ( 11, 21)( 12, 22)( 13, 14)( 23, 79)( 24, 80)( 25, 85)( 26, 91)( 27, 88)
    ( 28, 92)( 29, 82)( 30, 83)( 31, 90)( 32, 94)( 33, 86)( 34, 89)( 35, 81)
    ( 36, 84)( 37, 93)( 38, 87)( 39, 75)( 40, 77)( 42, 76)( 43, 62)( 45, 69)
    ( 46, 68)( 48, 54)( 49, 51)( 50, 60)( 53, 59)( 56, 58)( 63, 65)( 66, 78)
    ( 70, 74)( 71, 72)( 96, 99), (  1,  4, 13, 12, 10, 19, 18, 22)(  2,  7,  6, 11,
      16, 17, 14, 21)(  3,  9, 20,  8)(  5, 15)( 23, 62, 79, 38, 72, 28, 71, 50)
    ( 24, 82, 98, 56, 80, 86, 52, 78)( 25, 84, 53, 91, 59, 77, 85, 49)
    ( 26, 34, 36, 89, 51, 63, 40, 65)( 27, 29, 67, 66, 88, 33, 97, 58)( 30, 44)
    ( 31, 43, 70, 60, 74, 92, 90, 87)( 32, 68, 45, 39, 81, 54, 42, 99)
    ( 35, 46, 69, 96, 94, 48, 76, 75)( 37, 73, 93, 57)( 47, 83)( 55, 61, 64,
    95).

The representatives of conjugacy classes of   $G^{s_{16}}$ are:\\
 (1),
 (1,10)(2,16)(4,19)(6,14)(7,17)(11,21)(12,22)(13,18)(23,72)(24,80)(25,59)(26,
    51)(27,88)(28,62)(29,33)(31,74)\\(32,81)(34,63)(35,94)(36,40)(38,50)(39,99)(42,
    45)(43,92)(46,48)(49,91)(52,98)(53,85)(54,68)(56,78)(58,66)(60,\\87)(65,89)(67,
    97)(69,76)(70,90)(71,79)(75,96)(77,84)(82,86),

  (1,86,10,82)(2,33,16,29)(3,20)(4,52,19,98)(5,15)(6,58,14,66)(7,97,17,67)(8,9)(11,
    27,21,88)(12,24,22,80)(13,78,\\18,56)(23,85,72,53)(25,71,59,79)(26,43,51,92)(28,
    91,62,49)(30,83)(31,89,74,65)(32,45,81,42)(34,70,63,90)(35,76,\\94,69)(36,60,40,
    87)(37,93)(38,84,50,77)(39,54,99,68)(41,100)(44,47)(46,75,48,96)(55,64)(57,
    73)(61,95),

    (1,82,10,86)(2,29,16,33)(3,20)(4,98,19,52)(5,15)(6,66,14,58)(7,67,
    17,97)(8,9)(11,88,21,27)(12,80,22,24)(13,56,\\18,78)(23,53,72,85)(25,79,59,71)(26,
    92,51,43)(28,49,62,91)(30,83)(31,65,74,89)(32,42,81,45)(34,90,63,70)(35,69,\\94,
    76)(36,87,40,60)(37,93)(38,77,50,84)(39,68,99,54)(41,100)(44,47)(46,96,48,
    75)(55,64)(57,73)(61,95),

    (2,16)(4,12)(5,15)(7,21)(8,9)(11,17)(13,18)(19,
    22)(23,63)(24,52)(25,74)(26,62)(27,67)(28,51)(29,33)(31,59)(32,\\75)(34,72)(35,
    39)(36,50)(38,40)(42,46)(43,49)(45,48)(53,70)(54,76)(55,64)(56,78)(57,73)(60,
    77)(65,71)(68,69)(79,\\89)(80,98)(81,96)(84,87)(85,90)(88,97)(91,92)(94,99),

  (1,86,10,82)(2,29,16,33)(3,20)(4,24,19,80)(6,58,14,66)(7,88,17,27)(11,67,21,
    97)(12,52,22,98)(13,56,18,78)(23,90,\\72,70)(25,65,59,89)(26,49,51,91)(28,92,62,
    43)(30,83)(31,79,74,71)(32,48,81,46)(34,53,63,85)(35,54,94,68)(36,77,40,\\84)(37,
    93)(38,87,50,60)(39,76,99,69)(41,100)(42,75,45,96)(44,47)(61,95),

  (1,27,18,80,10,88,13,24)(2,98,6,67,16,52,14,97)(4,58,7,29,19,66,17,33)(5,9,15,
    8)(11,78,22,86,21,56,12,82)(23,71,\\53,25,72,79,85,59)(26,75,28,76,51,96,62,
    69)(30,64)(31,63,89,90,74,34,65,70)(32,77,68,87,81,84,54,60)(35,43,48,91,\\94,92,
    46,49)(36,42,50,99,40,45,38,39)(37,44,93,47)(41,100)(55,83)(57,95,73,61),

  (1,67,2,80,10,97,16,24)(4,58,21,78,19,66,11,56)(5,8,15,9)(6,27,13,52,14,88,18,
    98)(7,33,12,82,17,29,22,86)(23,65,\\53,74,72,89,85,31)(25,34,71,70,59,63,79,
    90)(26,32,60,75,51,81,87,96)(28,54,77,69,62,68,84,76)(30,55,83,64)(35,49,\\39,50,
    94,91,99,38)(36,46,43,45,40,48,92,42)(37,44,93,47)(41,100)(57,95)(61,73),

  (1,24)(2,97)(3,20)(4,66)(5,8)(6,52)(7,29)(9,15)(10,80)(11,56)(12,82)(13,27)(14,
    98)(16,67)(17,33)(18,88)(19,58)\\(21,78)(22,86)(23,70)(25,74)(26,39)(28,45)(30,
    64)(31,59)(32,49)(34,53)(35,60)(36,69)(37,47)(38,75)(40,76)(41,\\100)(42,62)(43,
    54)(44,93)(46,77)(48,84)(50,96)(51,99)(55,83)(57,95)(61,73)(63,85)(65,71)(68,
    92)(72,90)(79,89)\\(81,91)(87,94),

     (1,22,10,12)(2,17,16,7)(4,6,19,14)(5,9)(8,
    15)(11,13,21,18)(23,63,72,34)(24,86,80,82)(25,65,59,89)(26,54,51,68)\\(27,78,88,
    56)(28,81,62,32)(29,97,33,67)(30,55)(31,79,74,71)(35,40,94,36)(37,44)(38,48,50,
    46)(39,43,99,92)(42,49,\\45,91)(47,93)(52,58,98,66)(53,70,85,90)(57,61)(60,76,87,
    69)(64,83)(73,95)(75,84,96,77),

    (1,52,6,24)(2,88,13,67)(3,20)(4,66,12,82)(5,
    9)(7,33,11,78)(8,15)(10,98,14,80)(16,27,18,97)(17,29,21,56)(19,58,\\22,86)(23,53,
    72,85)(26,35,77,42)(28,48,87,99)(30,55,83,64)(32,43,76,38)(34,70,63,90)(36,68,
    91,96)(37,47)(39,62,\\46,60)(40,54,49,75)(41,100)(44,93)(45,51,94,84)(50,81,92,
    69)(57,95,73,61),

    (1,98,6,80)(2,27,13,97)(3,20)(4,58,12,86)(5,9)(7,29,11,
    56)(8,15)(10,52,14,24)(16,88,18,67)(17,33,21,78)(19,66,\\22,82)(23,85,72,53)(25,
    59)(26,94,77,45)(28,46,87,39)(30,55,83,64)(31,74)(32,92,76,50)(34,90,63,70)(35,
    84,42,51)\\(36,54,91,75)(37,47)(38,81,43,69)(40,68,49,96)(41,100)(44,93)(48,60,99,
    62)(57,95,73,61)(65,89)(71,79),

    (1,19,14,12)(2,11,18,17)(4,6,22,10)(5,8)(7,16,
    21,13)(9,15)(24,86,98,66)(25,71,59,79)(26,76,84,81)(27,56,67,33)\\(28,96,60,
    54)(29,88,78,97)(30,64,83,55)(31,89,74,65)(32,51,69,77)(34,63)(35,38,45,92)(36,
    39,49,48)(37,44)(40,99,\\91,46)(42,43,94,50)(47,93)(52,58,80,82)(57,61,73,95)(62,
    75,87,68)(70,90),

    (1,4,14,22)(2,21,18,7)(5,8)(6,12,10,19)(9,15)(11,13,17,
    16)(23,72)(24,82,98,58)(25,79,59,71)(26,69,84,32)(27,78,\\67,29)(28,75,60,68)(30,
    64,83,55)(31,65,74,89)(33,88,56,97)(35,50,45,43)(36,99,49,46)(37,44)(38,42,92,
    94)(39,91,\\48,40)(47,93)(51,76,77,81)(52,66,80,86)(53,85)(54,62,96,87)(57,61,73,
    95),

    (1,24,32)(2,88,94)(3,30,73)(4,96,18)(6,17,39)(7,99,14)(8,37,47)(9,93,
    44)(10,80,81)(11,69,33)(12,42,82)(13,19,75)\\(16,27,35)(20,83,57)(21,76,29)(22,45,
    86)(23,31,70)(26,36,77)(34,53,65)(38,43,60)(40,84,51)(41,64,95)(46,56,52)(48,\\78,
    98)(50,92,87)(54,58,67)(55,61,100)(63,85,89)(66,97,68)(72,74,90),

  (1,80,32,10,24,81)(2,27,94,16,88,35)(3,30,73)(4,75,18,19,96,13)(6,7,39,14,17,
    99)(8,37,47)(9,93,44)(11,76,33,21,\\69,29)(12,45,82,22,42,86)(20,83,57)(23,74,70,
    72,31,90)(25,59)(26,40,77,51,36,84)(28,62)(34,85,65,63,53,89)(38,92,\\60,50,43,
    87)(41,64,95)(46,78,52,48,56,98)(49,91)(54,66,67,68,58,97)(55,61,100)(71,79),

  (1,22,81,82,24,45,10,12,32,86,80,42)(2,11,35,29,88,69,16,21,94,33,27,76)(3,83,73,
    20,30,57)(4,46,13,98,96,56,19,\\48,18,52,75,78)(5,15)(6,67,99,66,17,54,14,97,39,
    58,7,68)(8,93,47,9,37,44)(23,89,90,53,31,63,72,65,70,85,74,34)(25,\\71,59,79)(26,
    60,84,92,36,38,51,87,77,43,40,50)(28,91,62,49)(41,55,95,100,64,61),

  (1,12,81,86,24,42,10,22,32,82,80,45)(2,21,35,33,88,76,16,11,94,29,27,69)(3,83,73,
    20,30,57)(4,48,13,52,96,78,19,\\46,18,98,75,56)(5,15)(6,97,99,58,17,68,14,67,39,
    66,7,54)(8,93,47,9,37,44)(23,65,90,85,31,34,72,89,70,53,74,63)(25,\\79,59,71)(26,
    87,84,43,36,50,51,60,77,92,40,38)(28,49,62,91)(41,55,95,100,64,61),

  (1,24,32,38,31,86,22,45,84,89,10,80,81,50,74,82,12,42,77,65)(2,4,69,62,25,33,52,
    35,49,71,16,19,76,28,59,29,98,\\94,91,79)(3,15,37,100,55,20,5,93,41,64)(6,27,54,
    43,23,58,21,99,51,85,14,88,68,92,72,66,11,39,26,53)(7,75,60,34,56,\\97,48,40,70,
    13,17,96,87,63,78,67,46,36,90,18)(8,30,47,95,57,9,83,44,61,73),

  (1,80,32,50,31,82,22,42,84,65,10,24,81,38,74,86,12,45,77,89)(2,19,69,28,25,29,52,
    94,49,79,16,4,76,62,59,33,98,\\35,91,71)(3,15,37,100,55,20,5,93,41,64)(6,88,54,92,
    23,66,21,39,51,53,14,27,68,43,72,58,11,99,26,85)(7,96,60,63,56,\\67,48,36,70,18,
    17,75,87,34,78,97,46,40,90,13)(8,30,47,95,57,9,83,44,61,73),

  (1,22,81,77,31,10,12,32,84,74)(2,52,76,91,25,16,98,69,49,59)(3,5,37,41,55)(4,35,
    28,79,33,19,94,62,71,29)(6,21,\\68,26,23,14,11,54,51,72)(7,48,87,90,56,17,46,60,
    70,78)(8,83,47,61,57)(9,30,44,95,73)(13,67,75,40,63,18,97,96,36,\\34)(15,93,100,
    64,20)(24,45,50,65,86,80,42,38,89,82)(27,99,92,53,58,88,39,43,85,66),

  (1,12,81,84,31)(2,98,76,49,25)(3,5,37,41,55)(4,94,28,71,33)(6,11,68,51,23)(7,46,
    87,70,56)(8,83,47,61,57)(9,30,\\44,95,73)(10,22,32,77,74)(13,97,75,36,63)(14,21,
    54,26,72)(15,93,100,64,20)(16,52,69,91,59)(17,48,60,90,78)(18,\\67,96,40,34)(19,
    35,62,79,29)(24,42,50,89,86)(27,39,92,85,58)(38,65,82,80,45)(43,53,66,88,99),

  (1,24,32,34,36,82,12,42,90,87,10,80,81,63,40,86,22,45,70,60)(2,7,76,59,28,29,67,
    35,71,49,16,17,69,25,62,33,97,\\94,79,91)(3,37,73,95,44,20,93,57,61,47)(4,39,65,
    38,58,98,68,74,77,6,19,99,89,50,66,52,54,31,84,14)(5,83,41,64,8,15,\\30,100,55,
    9)(11,96,72,92,18,88,48,85,51,78,21,75,23,43,13,27,46,53,26,56),

  (1,80,32,63,36,86,12,45,90,60,10,24,81,34,40,82,22,42,70,87)(2,17,76,25,28,33,67,
    94,71,91,16,7,69,59,62,29,97,\\35,79,49)(3,37,73,95,44,20,93,57,61,47)(4,99,65,50,
    58,52,68,31,77,14,19,39,89,38,66,98,54,74,84,6)(5,83,41,64,8,15,\\30,100,55,9)(11,
    75,72,43,18,27,48,53,51,56,21,96,23,92,13,88,46,85,26,78),

  (1,22,81,90,36)(2,97,69,71,28)(3,93,73,61,44)(4,54,89,77,58)(5,30,41,55,8)(6,98,
    39,31,50)(7,94,25,49,29)(9,15,\\83,100,64)(10,12,32,70,40)(11,46,23,51,18)(13,21,
    48,72,26)(14,52,99,74,38)(16,67,76,79,62)(17,35,59,91,33)(19,68,\\65,84,66)(20,37,
    57,95,47)(24,45,63,87,82)(27,75,85,92,56)(34,60,86,80,42)(43,78,88,96,53),

  (1,12,81,70,36,10,22,32,90,40)(2,67,69,79,28,16,97,76,71,62)(3,93,73,61,44)(4,68,
    89,84,58,19,54,65,77,66)(5,30,\\41,55,8)(6,52,39,74,50,14,98,99,31,38)(7,35,25,91,
    29,17,94,59,49,33)(9,15,83,100,64)(11,48,23,26,18,21,46,72,51,13)\\(20,37,57,95,
    47)(24,42,63,60,82,80,45,34,87,86)(27,96,85,43,56,88,75,53,92,78).

The character table of $G^{s_{16}}$:\\
\begin{tabular}{c|cccccccccccccccccccccc}
  & & & & & & & & & &10 & & & & & & & & & &20 & &\\\hline

$\chi_{16}^{(1)}$&1&1&1&1&1&1&1&1&1&1&1&1&1&1&1&1&1&1&1&1&1&1
\\$\chi_{16}^{(2)}$&3&3&3&3&3&3&-1&-1&-1&-1&-1&-1&-1&-1&.&.&.&.&G&G&G&G
\\$\chi_{16}^{(3)}$&3&3&3&3&3&3&-1&-1&-1&-1&-1&-1&-1&-1&.&.&.&.&*G&*G&*G&*G
\\$\chi_{16}^{(4)}$&4&4&4&4&4&4&.&.&.&.&.&.&.&.&1&1&1&1&-1&-1&-1&-1
\\$\chi_{16}^{(5)}$&4&-4&A&-A&.&.&.&.&.&.&2&-2&E&-E&1&-1&F&-F&-F&F&1&-1
\\$\chi_{16}^{(6)}$&4&-4&-A&A&.&.&.&.&.&.&2&-2&-E&E&1&-1&-F&F&F&-F&1&-1
\\$\chi_{16}^{(7)}$&5&5&5&5&5&5&1&1&1&1&1&1&1&1&-1&-1&-1&-1&.&.&.&.
\\$\chi_{16}^{(8)}$&6&6&-6&-6&-2&2&.&.&-2&2&2&2&-2&-2&.&.&.&.&-1&-1&1&1
\\$\chi_{16}^{(9)}$&6&6&-6&-6&-2&2&.&.&2&-2&-2&-2&2&2&.&.&.&.&-1&-1&1&1
\\$\chi_{16}^{(10)}$&10&10&-10&-10&2&-2&.&.&2&-2&2&2&-2&-2&1&1&-1&-1&.&.&.&.
\\$\chi_{16}^{(11)}$&10&10&-10&-10&2&-2&.&.&-2&2&-2&-2&2&2&1&1&-1&-1&.&.&.&.
\\$\chi_{16}^{(12)}$&12&-12&B&-B&.&.&.&.&.&.&-2&2&-E&E&.&.&.&.&H&-H&G&-G
\\$\chi_{16}^{(13)}$&12&-12&B&-B&.&.&.&.&.&.&-2&2&-E&E&.&.&.&.&I&-I&*G&-*G
\\$\chi_{16}^{(14)}$&12&-12&-B&B&.&.&.&.&.&.&-2&2&E&-E&.&.&.&.&-I&I&*G&-*G
\\$\chi_{16}^{(15)}$&12&-12&-B&B&.&.&.&.&.&.&-2&2&E&-E&.&.&.&.&-H&H&G&-G
\\$\chi_{16}^{(16)}$&12&12&-12&-12&-4&4&.&.&.&.&.&.&.&.&.&.&.&.&G&G&-G&-G
\\$\chi_{16}^{(17)}$&12&12&-12&-12&-4&4&.&.&.&.&.&.&.&.&.&.&.&.&*G&*G&-*G&-*G
\\$\chi_{16}^{(18)}$&15&15&15&15&-1&-1&-1&-1&-1&-1&3&3&3&3&.&.&.&.&.&.&.&.
\\$\chi_{16}^{(19)}$&15&15&15&15&-1&-1&-1&-1&3&3&-1&-1&-1&-1&.&.&.&.&.&.&.&.
\\$\chi_{16}^{(20)}$&15&15&15&15&-1&-1&-1&3&-1&-1&-1&-1&-1&-1&.&.&.&.&.&.&.&.
\\$\chi_{16}^{(21)}$&15&15&15&15&-1&-1&3&-1&-1&-1&-1&-1&-1&-1&.&.&.&.&.&.&.&.
\\$\chi_{16}^{(22)}$&16&-16&C&-C&.&.&.&.&.&.&.&.&.&.&1&-1&F&-F&F&-F&-1&1
\\$\chi_{16}^{(23)}$&16&-16&-C&C&.&.&.&.&.&.&.&.&.&.&1&-1&-F&F&-F&F&-1&1
\\$\chi_{16}^{(24)}$&20&20&-20&-20&4&-4&.&.&.&.&.&.&.&.&-1&-1&1&1&.&.&.&.
\\$\chi_{16}^{(25)}$&20&-20&D&-D&.&.&.&.&.&.&2&-2&E&-E&-1&1&-F&F&.&.&.&.
\\$\chi_{16}^{(26)}$&20&-20&-D&D&.&.&.&.&.&.&2&-2&-E&E&-1&1&F&-F&.&.&.&.
\end{tabular}

\begin{tabular}{c|cccc}
  & & & &\\\hline
$\chi_{16}^{(1)}$&1&1&1&1
\\$\chi_{16}^{(2)}$&*G&*G&*G&*G
\\$\chi_{16}^{(3)}$&G&G&G&G
\\$\chi_{16}^{(4)}$&-1&-1&-1&-1
\\$\chi_{16}^{(5)}$&F&-F&-1&1
\\$\chi_{16}^{(6)}$&-F&F&-1&1
\\$\chi_{16}^{(7)}$&.&.&.&.
\\$\chi_{16}^{(8)}$&-1&-1&1&1
\\$\chi_{16}^{(9)}$&-1&-1&1&1
\\$\chi_{16}^{(10)}$&.&.&.&.
\\$\chi_{16}^{(11)}$&.&.&.&.
\\$\chi_{16}^{(12)}$&-I&I&-*G&*G
\\$\chi_{16}^{(13)}$&-H&H&-G&G
\\$\chi_{16}^{(14)}$&H&-H&-G&G
\\$\chi_{16}^{(15)}$&I&-I&-*G&*G
\\$\chi_{16}^{(16)}$&*G&*G&-*G&-*G
\\$\chi_{16}^{(17)}$&G&G&-G&-G
\\$\chi_{16}^{(18)}$&.&.&.&.
\\$\chi_{16}^{(19)}$&.&.&.&.
\\$\chi_{16}^{(20)}$&.&.&.&.
\\$\chi_{16}^{(21)}$&.&.&.&.
\\$\chi_{16}^{(22)}$&-F&F&1&-1
\\$\chi_{16}^{(23)}$&F&-F&1&-1
\\$\chi_{16}^{(24)}$&.&.&.&.
\\$\chi_{16}^{(25)}$&.&.&.&.
\\$\chi_{16}^{(26)}$&.&.&.&.
\end{tabular}

\noindent \noindent where  A = -4*E(4)
  = -4*ER(-1) = -4i;
B = -12*E(4)
  = -12*ER(-1) = -12i;
C = -16*E(4)
  = -16*ER(-1) = -16i;
D = -20*E(4)
  = -20*ER(-1) = -20i;
E = -2*E(4)
  = -2*ER(-1) = -2i;
F = -E(4)
  = -ER(-1) = -i;
G = -E(5)-E(5)$^4$
  = (1-ER(5))/2 = -b5;
H = -E(20)-E(20)$^9$; I = -E(20)$^{13}$-E(20)$^{17}$.

The generators of $G^{s_{17}}$ are:\\
 (  1, 37, 46, 32, 43, 24, 76)(  2, 33, 91, 13, 87, 20, 96)(  3, 63, 74, 81, 55,
      95, 85)(  4, 54, 72, 35, 84, 65, 31)(  5, 44, 49, 99, 16, 19, 42)
    (  7, 10, 12, 94, 82, 86, 64)(  8, 38, 17, 15, 56, 39, 47)(  9, 88, 57, 14, 92,
      90, 18)( 11, 97, 71, 75, 78, 23, 41)( 21, 62, 52, 29, 60, 59, 25)
    ( 22, 34, 83, 73, 77, 40, 26)( 27, 45, 69, 48, 36, 93, 58)( 28, 61, 30, 68, 67,
     100, 70)( 50, 98, 89, 53, 80, 66, 51)

The representatives of conjugacy classes of   $G^{s_{17}}$ are:\\
 (1),
 (1,24,32,37,76,43,46)(2,20,13,33,96,87,91)(3,95,81,63,85,55,74)(4,65,35,54,
    31,84,72)(5,19,99,44,42,16,49)(7,\\86,94,10,64,82,12)(8,39,15,38,47,56,17)(9,90,
    14,88,18,92,57)(11,23,75,97,41,78,71)(21,59,29,62,25,60,52)(22,40,\\73,34,26,77,
    83)(27,93,48,45,58,36,69)(28,100,68,61,70,67,30)(50,66,53,98,51,80,89),

  (1,32,76,46,24,37,43)(2,13,96,91,20,33,87)(3,81,85,74,95,63,55)(4,35,31,72,65,54,
    84)(5,99,42,49,19,44,16)(7,\\94,64,12,86,10,82)(8,15,47,17,39,38,56)(9,14,18,57,
    90,88,92)(11,75,41,71,23,97,78)(21,29,25,52,59,62,60)(22,73,\\26,83,40,34,77)(27,
    48,58,69,93,45,36)(28,68,70,30,100,61,67)(50,53,51,89,66,98,80),

  (1,37,46,32,43,24,76)(2,33,91,13,87,20,96)(3,63,74,81,55,95,85)(4,54,72,35,84,65,
    31)(5,44,49,99,16,19,42)(7,\\10,12,94,82,86,64)(8,38,17,15,56,39,47)(9,88,57,14,
    92,90,18)(11,97,71,75,78,23,41)(21,62,52,29,60,59,25)(22,34,\\83,73,77,40,26)(27,
    45,69,48,36,93,58)(28,61,30,68,67,100,70)(50,98,89,53,80,66,51),

  (1,43,37,24,46,76,32)(2,87,33,20,91,96,13)(3,55,63,95,74,85,81)(4,84,54,65,72,31,
    35)(5,16,44,19,49,42,99)(7,\\82,10,86,12,64,94)(8,56,38,39,17,47,15)(9,92,88,90,
    57,18,14)(11,78,97,23,71,41,75)(21,60,62,59,52,25,29)(22,77,\\34,40,83,26,73)(27,
    36,45,93,69,58,48)(28,67,61,100,30,70,68)(50,80,98,66,89,51,53),

  (1,46,43,76,37,32,24)(2,91,87,96,33,13,20)(3,74,55,85,63,81,95)(4,72,84,31,54,35,
    65)(5,49,16,42,44,99,19)(7,\\12,82,64,10,94,86)(8,17,56,47,38,15,39)(9,57,92,18,
    88,14,90)(11,71,78,41,97,75,23)(21,52,60,25,62,29,59)(22,83,\\77,26,34,73,40)(27,
    69,36,58,45,48,93)(28,30,67,70,61,68,100)(50,89,80,51,98,53,66),

  (1,76,24,43,32,46,37)(2,96,20,87,13,91,33)(3,85,95,55,81,74,63)(4,31,65,84,35,72,
    54)(5,42,19,16,99,49,44)(7,\\64,86,82,94,12,10)(8,47,39,56,15,17,38)(9,18,90,92,
    14,57,88)(11,41,23,78,75,71,97)(21,25,59,60,29,52,62)(22,26,\\40,77,73,83,34)(27,
    58,93,36,48,69,45)(28,70,100,67,68,30,61)(50,51,66,80,53,89,98).

The character table of $G^{s_{17}}$:\\
\begin{tabular}{c|ccccccc}
  & & & & & & &\\\hline
$\chi_{17}^{(1)}$&1&1&1&1&1&1&1
\\$\chi_{17}^{(2)}$&1&A&B&C&/B&/A&/C
\\$\chi_{17}^{(3)}$&1&B&/C&/A&C&/B&A
\\$\chi_{17}^{(4)}$&1&C&/A&B&A&/C&/B
\\$\chi_{17}^{(5)}$&1&/C&A&/B&/A&C&B
\\$\chi_{17}^{(6)}$&1&/B&C&A&/C&B&/A
\\$\chi_{17}^{(7)}$&1&/A&/B&/C&B&A&C
\end{tabular}

\noindent \noindent where   A = E(7)$^5$; B = E(7)$^3$; C = E(7).

The generators of $G^{s_{18}}$ are:\\
(  1, 87, 46, 96, 20, 51, 78, 32, 42, 89, 74, 86)(  2, 60,  6, 38,
10, 98,  5, 27,
     72,100, 73, 53)(  3, 52, 76, 25, 90, 66)(  4, 95, 40, 28, 11, 93)
    (  7,  8, 26, 97, 30, 80, 70, 77, 65, 62, 85, 21)(  9, 22, 37, 99, 36, 79, 67,
      75, 43, 41, 92, 13)( 12, 88, 29, 31, 48, 69, 82, 71, 63, 59, 47, 16)
    ( 14, 49, 44, 34, 94, 45)( 15, 17, 91, 23, 33, 58, 68, 83, 39, 64, 55, 54)
    ( 18, 61, 81, 24)( 19, 50, 84, 35)( 56, 57)

The representatives of conjugacy classes of   $G^{s_{18}}$ are:\\
 (1),
 (1,20,42)(2,10,72)(3,90,76)(4,11,40)(5,73,6)(7,30,65)(8,80,62)(9,36,43)(12,
    48,63)(13,99,75)(14,94,44)(15,33,\\39)(16,31,71)(17,58,64)(21,97,77)(22,79,41)(23,
    83,54)(25,52,66)(26,70,85)(27,53,38)(28,95,93)(29,82,47)(32,86,\\96)(34,49,45)(37,
    67,92)(46,78,74)(51,89,87)(55,91,68)(59,88,69)(60,98,100),

  (1,32,46,89,20,86,78,87,42,96,74,51)(2,27,6,100,10,53,5,60,72,38,73,98)(3,52,76,
    25,90,66)(4,95,40,28,11,93)\\(7,77,26,62,30,21,70,8,65,97,85,80)(9,75,37,41,36,13,
    67,22,43,99,92,79)(12,71,29,59,48,16,82,88,63,31,47,69)(14,\\49,44,34,94,45)(15,
    83,91,64,33,54,68,17,39,23,55,58)(18,24,81,61)(19,35,84,50)(56,57),

  (1,42,20)(2,72,10)(3,76,90)(4,40,11)(5,6,73)(7,65,30)(8,62,80)(9,43,36)(12,63,
    48)(13,75,99)(14,44,94)(15,39,\\33)(16,71,31)(17,64,58)(21,77,97)(22,41,79)(23,54,
    83)(25,66,52)(26,85,70)(27,38,53)(28,93,95)(29,47,82)(32,96,\\86)(34,45,49)(37,92,
    67)(46,74,78)(51,87,89)(55,68,91)(59,69,88)(60,100,98),

  (1,46,20,78,42,74)(2,6,10,5,72,73)(3,76,90)(4,40,11)(7,26,30,70,65,85)(8,97,80,77,
    62,21)(9,37,36,67,43,92)(12,\\29,48,82,63,47)(13,22,99,79,75,41)(14,44,94)(15,91,
    33,68,39,55)(16,88,31,69,71,59)(17,23,58,83,64,54)(18,81)(19,\\84)(24,61)(25,66,
    52)(27,100,53,60,38,98)(28,93,95)(32,89,86,87,96,51)(34,45,49)(35,50),

  (1,51,74,96,42,87,78,86,20,89,46,32)(2,98,73,38,72,60,5,53,10,100,6,27)(3,66,90,
    25,76,52)(4,93,11,28,40,95)(7,\\80,85,97,65,8,70,21,30,62,26,77)(9,79,92,99,43,22,
    67,13,36,41,37,75)(12,69,47,31,63,88,82,16,48,59,29,71)(14,45,\\94,34,44,49)(15,
    58,55,23,39,17,68,54,33,64,91,83)(18,61,81,24)(19,50,84,35)(56,57),

  (1,74,42,78,20,46)(2,73,72,5,10,6)(3,90,76)(4,11,40)(7,85,65,70,30,26)(8,21,62,77,
    80,97)(9,92,43,67,36,37)(12,\\47,63,82,48,29)(13,41,75,79,99,22)(14,94,44)(15,55,
    39,68,33,91)(16,59,71,69,31,88)(17,54,64,83,58,23)(18,81)\\(19,84)(24,61)(25,52,
    66)(27,98,38,60,53,100)(28,95,93)(32,51,96,87,86,89)(34,49,45)(35,50),

  (1,78)(2,5)(6,72)(7,70)(8,77)(9,67)(10,73)(12,82)(13,79)(15,68)(16,69)(17,83)(18,
    81)(19,84)(20,74)(21,80)(22,\\75)(23,64)(24,61)(26,65)(27,60)(29,63)(30,85)(31,
    59)(32,87)(33,55)(35,50)(36,92)(37,43)(38,100)(39,91)(41,99)\\(42,46)(47,48)(51,
    86)(53,98)(54,58)(62,97)(71,88)(89,96),

    (1,86,74,89,42,32,78,51,20,96,46,
    87)(2,53,73,100,72,27,5,98,10,38,6,60)(3,66,90,25,76,52)(4,93,11,28,40,95)(7,\\21,
    85,62,65,77,70,80,30,97,26,8)(9,13,92,41,43,75,67,79,36,99,37,22)(12,16,47,59,
    63,71,82,69,48,31,29,88)(14,45,\\94,34,44,49)(15,54,55,64,39,83,68,58,33,23,91,
    17)(18,24,81,61)(19,35,84,50)(56,57),

    (1,87,46,96,20,51,78,32,42,89,74,86)(2,
    60,6,38,10,98,5,27,72,100,73,53)(3,52,76,25,90,66)(4,95,40,28,11,93)(7,\\8,26,97,
    30,80,70,77,65,62,85,21)(9,22,37,99,36,79,67,75,43,41,92,13)(12,88,29,31,48,69,
    82,71,63,59,47,16)(14,49,\\44,34,94,45)(15,17,91,23,33,58,68,83,39,64,55,54)(18,
    61,81,24)(19,50,84,35)(56,57),

     (1,89,78,96)(2,100,5,38)(3,25)(4,28)(6,53,72,
    98)(7,62,70,97)(8,85,77,30)(9,41,67,99)(10,60,73,27)(11,95)(12,\\59,82,31)(13,43,
    79,37)(14,34)(15,64,68,23)(16,63,69,29)(17,55,83,33)(18,61,81,24)(19,50,84,
    35)(20,87,74,32)(21,\\65,80,26)(22,92,75,36)(39,58,91,54)(40,93)(42,51,46,86)(44,
    45)(47,71,48,88)(49,94)(52,90)(56,57)(66,76),

  (1,96,78,89)(2,38,5,100)(3,25)(4,28)(6,98,72,53)(7,97,70,62)(8,30,77,85)(9,99,67,
    41)(10,27,73,60)(11,95)(12,\\31,82,59)(13,37,79,43)(14,34)(15,23,68,64)(16,29,69,
    63)(17,33,83,55)(18,24,81,61)(19,35,84,50)(20,32,74,87)(21,\\26,80,65)(22,36,75,
    92)(39,54,91,58)(40,93)(42,86,46,51)(44,45)(47,88,48,71)(49,94)(52,90)(56,
    57)(66,76).

The character table of $G^{s_{18}}$:\\
\begin{tabular}{c|cccccccccccc}
  & & & & & & & & & &10 & &\\\hline

$\chi_{18}^{(1)}$&1&1&1&1&1&1&1&1&1&1&1&1
\\$\chi_{18}^{(2)}$&1&1&-1&1&1&-1&1&1&-1&-1&-1&-1
\\$\chi_{18}^{(3)}$&1&A&A&/A&/A&/A&A&1&/A&A&1&1
\\$\chi_{18}^{(4)}$&1&A&-A&/A&/A&-/A&A&1&-/A&-A&-1&-1
\\$\chi_{18}^{(5)}$&1&/A&/A&A&A&A&/A&1&A&/A&1&1
\\$\chi_{18}^{(6)}$&1&/A&-/A&A&A&-A&/A&1&-A&-/A&-1&-1
\\$\chi_{18}^{(7)}$&1&1&B&1&-1&-B&-1&-1&B&-B&-B&B
\\$\chi_{18}^{(8)}$&1&1&-B&1&-1&B&-1&-1&-B&B&B&-B
\\$\chi_{18}^{(9)}$&1&A&C&/A&-/A&/C&-A&-1&-/C&-C&-B&B
\\$\chi_{18}^{(10)}$&1&A&-C&/A&-/A&-/C&-A&-1&/C&C&B&-B
\\$\chi_{18}^{(11)}$&1&/A&-/C&A&-A&-C&-/A&-1&C&/C&-B&B
\\$\chi_{18}^{(12)}$&1&/A&/C&A&-A&C&-/A&-1&-C&-/C&B&-B
\end{tabular}

\noindent \noindent where   A = E(3)
  = (-1+ER(-3))/2 = b3;
B = -E(4)
  = -ER(-1) = -i;
C = -E(12)$^7$.

The generators of $G^{s_{19}}$ are:\\
(  1, 46, 20, 78, 42, 74)(  2,  6, 10,  5, 72, 73)(  3, 76, 90)
    (  4, 40, 11)(  7, 26, 30, 70, 65, 85)(  8, 97, 80, 77, 62, 21)(  9, 37, 36, 67,
     43, 92)( 12, 29, 48, 82, 63, 47)( 13, 22, 99, 79, 75, 41)( 14, 44, 94)
    ( 15, 91, 33, 68, 39, 55)( 16, 88, 31, 69, 71, 59)( 17, 23, 58, 83, 64, 54)
    ( 18, 81)( 19, 84)( 24, 61)( 25, 66, 52)( 27,100, 53, 60, 38, 98)( 28, 93, 95)
    ( 32, 89, 86, 87, 96, 51)( 34, 45, 49)( 35, 50),
  (  1,  2, 20, 10, 42, 72)(  3, 45, 90, 34, 76, 49)(  4, 40, 11)(  5, 74, 73, 46,
       6, 78)(  7, 79, 30, 41, 65, 22)(  8,  9, 80, 36, 62, 43)( 12, 63, 48)
    ( 13, 85, 99, 26, 75, 70)( 14, 66, 94, 25, 44, 52)( 15, 91, 33, 68, 39, 55)
    ( 16, 88, 31, 69, 71, 59)( 17, 64, 58)( 18, 50)( 19, 61)( 21, 92, 97, 37, 77, 67
     )( 23, 54, 83)( 24, 84)( 27, 51, 53, 89, 38, 87)( 28, 93, 95)( 29, 47, 82)
    ( 32, 60, 86, 98, 96,100)( 35, 81), (  1, 27)(  2, 86)(  3, 14)(  4, 28)
    (  5, 51)(  6, 87)(  7,  9)(  8, 75)( 10, 96)( 11, 95)( 12, 31)( 13, 80)
    ( 15, 64)( 16, 63)( 17, 33)( 18, 19)( 20, 53)( 21, 79)( 22, 77)( 23, 68)
    ( 24, 50)( 25, 34)( 26, 37)( 29, 69)( 30, 36)( 32, 72)( 35, 61)( 38, 42)
    ( 39, 58)( 40, 93)( 41, 97)( 43, 65)( 44, 76)( 45, 66)( 46,100)( 47, 88)
    ( 48, 71)( 49, 52)( 54, 91)( 55, 83)( 56, 57)( 59, 82)( 60, 78)( 62, 99)
    ( 67, 70)( 73, 89)( 74, 98)( 81, 84)( 85, 92)( 90, 94).

The representatives of conjugacy classes of   $G^{s_{19}}$ are:\\
 (1),
 (1,2,20,10,42,72)(3,45,90,34,76,49)(4,40,11)(5,74,73,46,6,78)(7,79,30,41,65,
    22)(8,9,80,36,62,43)(12,63,48)(13,\\85,99,26,75,70)(14,66,94,25,44,52)(15,91,33,
    68,39,55)(16,88,31,69,71,59)(17,64,58)(18,50)(19,61)(21,92,97,37,77,\\67)(23,54,
    83)(24,84)(27,51,53,89,38,87)(28,93,95)(29,47,82)(32,60,86,98,96,100)(35,81),

  (1,6,42,73,20,5)(2,78,72,46,10,74)(3,49,76,34,90,45)(4,11,40)(7,75,65,99,30,13)(8,
    37,62,92,80,67)(9,77,43,97,\\36,21)(12,47,63,82,48,29)(14,52,44,25,94,66)(15,33,
    39)(16,31,71)(17,54,64,83,58,23)(18,35)(19,24)(22,26,41,85,\\79,70)(27,32,38,96,
    53,86)(28,95,93)(50,81)(51,60,87,100,89,98)(55,91,68)(59,88,69)(61,84),

  (1,10)(2,42)(3,34)(5,46)(6,74)(7,41)(8,36)(9,62)(13,26)(14,25)(15,68)(16,69)(18,
    50)(19,61)(20,72)(21,37)(22,\\30)(24,84)(27,89)(31,59)(32,98)(33,55)(35,81)(38,
    51)(39,91)(43,80)(44,66)(45,76)(49,90)(52,94)(53,87)(60,96)\\(65,79)(67,97)(70,
    99)(71,88)(73,78)(75,85)(77,92)(86,100),

    (1,20,42)(2,10,72)(3,90,76)(4,11,
    40)(5,73,6)(7,30,65)(8,80,62)(9,36,43)(12,48,63)(13,99,75)(14,94,44)(15,33,\\
    39)(16,31,71)(17,58,64)(21,97,77)(22,79,41)(23,83,54)(25,52,66)(26,70,85)(27,53,
    38)(28,95,93)(29,82,47)(32,86,\\96)(34,49,45)(37,67,92)(46,78,74)(51,89,87)(55,91,
    68)(59,88,69)(60,98,100),

    (1,27)(2,86)(3,14)(4,28)(5,51)(6,87)(7,9)(8,75)(10,
    96)(11,95)(12,31)(13,80)(15,64)(16,63)(17,33)(18,19)(20,\\53)(21,79)(22,77)(23,
    68)(24,50)(25,34)(26,37)(29,69)(30,36)(32,72)(35,61)(38,42)(39,58)(40,93)(41,
    97)(43,65)\\(44,76)(45,66)(46,100)(47,88)(48,71)(49,52)(54,91)(55,83)(56,57)(59,
    82)(60,78)(62,99)(67,70)(73,89)(74,98)(81,\\84)(85,92)(90,94),

  (1,32,46,89,20,86,78,87,42,96,74,51)(2,27,6,100,10,53,5,60,72,38,73,98)(3,52,76,
    25,90,66)(4,95,40,28,11,93)(7,\\77,26,62,30,21,70,8,65,97,85,80)(9,75,37,41,36,13,
    67,22,43,99,92,79)(12,71,29,59,48,16,82,88,63,31,47,69)(14,49,\\44,34,94,45)(15,
    83,91,64,33,54,68,17,39,23,55,58)(18,24,81,61)(19,35,84,50)(56,57),

  (1,38,20,27,42,53)(2,32,10,86,72,96)(3,44,90,14,76,94)(4,93,11,28,40,95)(5,87,73,
    51,6,89)(7,43,30,9,65,36)(8,99,\\80,75,62,13)(12,16,48,31,63,71)(15,58,33,64,39,
    17)(18,19)(21,22,97,79,77,41)(23,91,83,68,54,55)(24,50)(25,45,52,\\34,66,49)(26,
    92,70,37,85,67)(29,88,82,69,47,59)(35,61)(46,98,78,100,74,60)(56,57)(81,84),

  (1,42,20)(2,72,10)(3,76,90)(4,40,11)(5,6,73)(7,65,30)(8,62,80)(9,43,36)(12,63,
    48)(13,75,99)(14,44,94)(15,39,33)\\(16,71,31)(17,64,58)(21,77,97)(22,41,79)(23,54,
    83)(25,66,52)(26,85,70)(27,38,53)(28,93,95)(29,47,82)(32,96,86)\\(345,49)(37,92,
    67)(46,74,78)(51,87,89)(55,68,91)(59,69,88)(60,100,98),

  (1,46,20,78,42,74)(2,6,10,5,72,73)(3,76,90)(4,40,11)(7,26,30,70,65,85)(8,97,80,77,
    62,21)(9,37,36,67,43,92)(12,29,\\48,82,63,47)(13,22,99,79,75,41)(14,44,94)(15,91,
    33,68,39,55)(16,88,31,69,71,59)(17,23,58,83,64,54)(18,81)(19,84)\\(24,61)(25,66,
    52)(27,100,53,60,38,98)(28,93,95)(32,89,86,87,96,51)(34,45,49)(35,50),

  (1,51,74,96,42,87,78,86,20,89,46,32)(2,98,73,38,72,60,5,53,10,100,6,27)(3,66,90,
    25,76,52)(4,93,11,28,40,95)(7,80,\\85,97,65,8,70,21,30,62,26,77)(9,79,92,99,43,22,
    67,13,36,41,37,75)(12,69,47,31,63,88,82,16,48,59,29,71)(14,45,94,34,\\44,49)(15,
    58,55,23,39,17,68,54,33,64,91,83)(18,61,81,24)(19,50,84,35)(56,57),

  (1,53,42,27,20,38)(2,96,72,86,10,32)(3,94,76,14,90,44)(4,95,40,28,11,93)(5,89,6,
    51,73,87)(7,36,65,9,30,43)(8,13,\\62,75,80,99)(12,71,63,31,48,16)(15,17,39,64,33,
    58)(18,19)(21,41,77,79,97,22)(23,55,54,68,83,91)(24,50)(25,49,66,\\34,52,45)(26,
    67,85,37,70,92)(29,59,47,69,82,88)(35,61)(46,60,74,100,78,98)(56,57)(81,84),

  (1,74,42,78,20,46)(2,73,72,5,10,6)(3,90,76)(4,11,40)(7,85,65,70,30,26)(8,21,62,77,
    80,97)(9,92,43,67,36,37)(12,47,\\63,82,48,29)(13,41,75,79,99,22)(14,94,44)(15,55,
    39,68,33,91)(16,59,71,69,31,88)(17,54,64,83,58,23)(18,81)(19,84)\\(24,61)(25,52,
    66)(27,98,38,60,53,100)(28,95,93)(32,51,96,87,86,89)(34,49,45)(35,50),

  (1,78)(2,5)(6,72)(7,70)(8,77)(9,67)(10,73)(12,82)(13,79)(15,68)(16,69)(17,83)(18,
    81)(19,84)(20,74)(21,80)(22,75)\\(23,64)(24,61)(26,65)(27,60)(29,63)(30,85)(31,
    59)(32,87)(33,55)(35,50)(36,92)(37,43)(38,100)(39,91)(41,99)(42,\\46)(47,48)(51,
    86)(53,98)(54,58)(62,97)(71,88)(89,96),

    (1,89,78,96)(2,100,5,38)(3,25)(4,
    28)(6,53,72,98)(7,62,70,97)(8,85,77,30)(9,41,67,99)(10,60,73,27)(11,95)(12,59,\\
    82,31)(13,43,79,37)(14,34)(15,64,68,23)(16,63,69,29)(17,55,83,33)(18,61,81,
    24)(19,50,84,35)(20,87,74,32)(21,65,\\80,26)(22,92,75,36)(39,58,91,54)(40,93)(42,
    51,46,86)(44,45)(47,71,48,88)(49,94)(52,90)(56,57)(66,76).

The character table of $G^{s_{19}}$:\\
\begin{tabular}{c|ccccccccccccccc}
  & & & & & & & & & &10 & & & & &\\\hline

$\chi_{19}^{(1)}$&1&1&1&1&1&1&1&1&1&1&1&1&1&1&1
\\$\chi_{19}^{(2)}$&1&-1&-1&-1&1&-1&1&-1&1&1&1&-1&1&1&1
\\$\chi_{19}^{(3)}$&1&-1&-1&-1&1&1&-1&1&1&1&-1&1&1&1&-1
\\$\chi_{19}^{(4)}$&1&1&1&1&1&-1&-1&-1&1&1&-1&-1&1&1&-1
\\$\chi_{19}^{(5)}$&1&A&/A&-1&-/A&-1&-/A&A&-A&-A&-A&/A&-/A&1&1
\\$\chi_{19}^{(6)}$&1&/A&A&-1&-A&-1&-A&/A&-/A&-/A&-/A&A&-A&1&1
\\$\chi_{19}^{(7)}$&1&A&/A&-1&-/A&1&/A&-A&-A&-A&A&-/A&-/A&1&-1
\\$\chi_{19}^{(8)}$&1&/A&A&-1&-A&1&A&-/A&-/A&-/A&/A&-A&-A&1&-1
\\$\chi_{19}^{(9)}$&1&-/A&-A&1&-A&-1&A&/A&-/A&-/A&/A&A&-A&1&-1
\\$\chi_{19}^{(10)}$&1&-A&-/A&1&-/A&-1&/A&A&-A&-A&A&/A&-/A&1&-1
\\$\chi_{19}^{(11)}$&1&-/A&-A&1&-A&1&-A&-/A&-/A&-/A&-/A&-A&-A&1&1
\\$\chi_{19}^{(12)}$&1&-A&-/A&1&-/A&1&-/A&-A&-A&-A&-A&-/A&-/A&1&1
\\$\chi_{19}^{(13)}$&2&.&.&.&2&.&.&.&2&-2&.&.&-2&-2&.
\\$\chi_{19}^{(14)}$&2&.&.&.&B&.&.&.&/B&-/B&.&.&-B&-2&.
\\$\chi_{19}^{(15)}$&2&.&.&.&/B&.&.&.&B&-B&.&.&-/B&-2&.
\end{tabular}

\noindent \noindent where   A = -E(3)
  = (1-ER(-3))/2 = -b3;
B = 2*E(3)$^2$
  = -1-ER(-3) = -1-i3.

The generators of $G^{s_{20}}$ are:\\
(  1, 53, 42, 27, 20, 38)(  2, 96, 72, 86, 10, 32)(  3, 94, 76, 14,
90, 44)
    (  4, 95, 40, 28, 11, 93)(  5, 89,  6, 51, 73, 87)(  7, 36, 65,  9, 30, 43)
    (  8, 13, 62, 75, 80, 99)( 12, 71, 63, 31, 48, 16)( 15, 17, 39, 64, 33, 58)
    ( 18, 19)( 21, 41, 77, 79, 97, 22)( 23, 55, 54, 68, 83, 91)( 24, 50)
    ( 25, 49, 66, 34, 52, 45)( 26, 67, 85, 37, 70, 92)( 29, 59, 47, 69, 82, 88)
    ( 35, 61)( 46, 60, 74,100, 78, 98)( 56, 57)( 81, 84),
  (  2, 73)(  4, 78)(  5, 72)(  6, 10)(  7, 25)(  8, 23)(  9, 34)( 11, 74)( 12, 33)
    ( 13, 55)( 15, 63)( 16, 64)( 17, 31)( 18, 57)( 19, 56)( 21, 59)( 22, 29)
    ( 24, 61)( 28, 60)( 30, 52)( 32, 51)( 35, 50)( 36, 49)( 39, 48)( 40, 46)
    ( 41, 47)( 43, 45)( 54, 62)( 58, 71)( 65, 66)( 68, 75)( 69, 77)( 79, 82)
    ( 80, 83)( 86, 89)( 87, 96)( 88, 97)( 91, 99)( 93,100)( 95, 98),
  (  1,  4)(  2,  6)(  5, 10)(  8, 16)( 11, 20)( 12, 13)( 15, 22)( 17, 21)( 23, 69)
    ( 24, 61)( 25, 70)( 26, 66)( 27, 28)( 29, 68)( 31, 80)( 32, 89)( 33, 79)
    ( 34, 67)( 35, 50)( 37, 45)( 38, 93)( 39, 41)( 40, 42)( 47, 91)( 48, 99)
    ( 49, 92)( 51, 96)( 52, 85)( 53, 95)( 54, 88)( 55, 82)( 56, 84)( 57, 81)
    ( 58, 97)( 59, 83)( 62, 71)( 63, 75)( 64, 77)( 72, 73)( 86, 87).

The representatives of conjugacy classes of   $G^{s_{20}}$ are:\\
 (1),
 (2,73)(4,78)(5,72)(6,10)(7,25)(8,23)(9,34)(11,74)(12,33)(13,55)(15,63)(16,
    64)(17,31)(18,57)(19,56)(21,59)(22,\\29)(24,61)(28,60)(30,52)(32,51)(35,50)(36,
    49)(39,48)(40,46)(41,47)(43,45)(54,62)(58,71)(65,66)(68,75)(69,77)\\(79,82)(80,
    83)(86,89)(87,96)(88,97)(91,99)(93,100)(95,98),

  (1,4,78)(2,72,10)(5,73,6)(7,70,25)(8,69,64)(9,67,34)(11,74,20)(12,79,55)(13,82,
    33)(15,75,29)(16,77,23)(17,80,59)\\(18,81,57)(19,84,56)(21,83,31)(22,68,63)(26,66,
    65)(27,28,60)(30,85,52)(32,96,86)(36,92,49)(37,45,43)(38,93,100)\\(39,99,47)(40,
    46,42)(41,91,48)(51,89,87)(53,95,98)(54,71,97)(58,62,88),

  (1,11,46)(3,90,76)(4,74,42)(5,6,73)(7,85,66)(8,59,58)(9,92,45)(12,41,68)(13,47,
    15)(14,94,44)(16,21,54)(17,62,69)\\(18,81,57)(19,84,56)(20,40,78)(22,55,48)(23,31,
    97)(25,30,26)(27,95,100)(28,98,38)(29,33,99)(34,36,37)(39,75,82)\\(43,67,49)(51,
    87,89)(52,65,70)(53,93,60)(63,79,91)(64,80,88)(71,77,83),

  (1,11,42,4,20,40)(2,5,72,6,10,73)(3,90,76)(7,30,65)(8,31,62,16,80,71)(9,36,43)(12,
    99,63,13,48,75)(14,94,44)(15,79,\\39,22,33,41)(17,97,64,21,58,77)(23,59,54,69,83,
    88)(24,61)(25,85,66,70,52,26)(27,95,38,28,53,93)(29,55,47,68,82,91)\\(32,87,96,89,
    86,51)(34,92,45,67,49,37)(35,50)(46,78,74)(56,84)(57,81)(60,98,100),

  (1,20,42)(2,10,72)(3,90,76)(4,11,40)(5,73,6)(7,30,65)(8,80,62)(9,36,43)(12,48,
    63)(13,99,75)(14,94,44)(15,33,39)\\(16,31,71)(17,58,64)(21,97,77)(22,79,41)(23,83,
    54)(25,52,66)(26,70,85)(27,53,38)(28,95,93)(29,82,47)(32,86,96)\\(34,49,45)(37,67,
    92)(46,78,74)(51,89,87)(55,91,68)(59,88,69)(60,98,100),

  (1,27)(2,86)(3,14)(4,28)(5,51)(6,87)(7,9)(8,75)(10,96)(11,95)(12,31)(13,80)(15,
    64)(16,63)(17,33)(18,19)(20,53)\\(21,79)(22,77)(23,68)(24,50)(25,34)(26,37)(29,
    69)(30,36)(32,72)(35,61)(38,42)(39,58)(40,93)(41,97)(43,65)(44,76)\\(45,66)(46,
    100)(47,88)(48,71)(49,52)(54,91)(55,83)(56,57)(59,82)(60,78)(62,99)(67,70)(73,
    89)(74,98)(81,84)(85,\\92)(90,94),

    (1,27)(2,89)(3,14)(4,60)(5,32)(6,96)(7,34)(8,
    68)(9,25)(10,87)(11,98)(12,17)(13,83)(15,16)(18,56)(19,57)(20,53)\\(21,82)(22,
    69)(23,75)(24,35)(26,37)(28,78)(29,77)(30,49)(31,33)(36,52)(38,42)(39,71)(40,
    100)(41,88)(43,66)(44,\\76)(45,65)(46,93)(47,97)(48,58)(50,61)(51,72)(54,99)(55,
    80)(59,79)(62,91)(63,64)(67,70)(73,86)(74,95)(81,84)(85,\\92)(90,94),

  (1,28,78,27,4,60)(2,32,10,86,72,96)(3,14)(5,89,6,51,73,87)(7,67,25,9,70,34)(8,29,
    64,75,69,15)(11,98,20,95,74,53)\\(12,21,55,31,79,83)(13,59,33,80,82,17)(16,22,23,
    63,77,68)(18,84,57,19,81,56)(24,50)(26,45,65,37,66,43)(30,92,52,36,\\85,49)(35,
    61)(38,40,100,42,93,46)(39,62,47,58,99,88)(41,54,48,97,91,71)(44,76)(90,94),

  (1,38,20,27,42,53)(2,32,10,86,72,96)(3,44,90,14,76,94)(4,93,11,28,40,95)(5,87,73,
    51,6,89)(7,43,30,9,65,36)(8,99,\\80,75,62,13)(12,16,48,31,63,71)(15,58,33,64,39,
    17)(18,19)(21,22,97,79,77,41)(23,91,83,68,54,55)(24,50)(25,45,52,\\34,66,49)(26,
    92,70,37,85,67)(29,88,82,69,47,59)(35,61)(46,98,78,100,74,60)(56,57)(81,84),

  (1,38,20,27,42,53)(2,51,10,89,72,87)(3,44,90,14,76,94)(4,100,11,60,40,98)(5,96,73,
    32,6,86)(7,45,30,34,65,49)(8,\\91,80,68,62,55)(9,66,36,25,43,52)(12,64,48,17,63,
    58)(13,23,99,83,75,54)(15,71,33,16,39,31)(18,56)(19,57)(21,29,97,\\82,77,47)(22,
    88,79,69,41,59)(24,35)(26,92,70,37,85,67)(28,46,95,78,93,74)(50,61)(81,84),

  (1,40,74)(2,10,72)(3,76,90)(4,46,20)(7,26,52)(8,88,17)(9,37,49)(11,78,42)(12,22,
    91)(13,29,39)(14,44,94)(15,99,\\82)(16,97,83)(18,81,57)(19,84,56)(21,23,71)(25,65,
    85)(27,93,98)(28,100,53)(30,70,66)(31,77,54)(32,86,96)(33,75,\\47)(34,43,92)(36,
    67,45)(38,95,60)(41,55,63)(48,79,68)(58,80,69)(59,64,62),

  (1,40,20,4,42,11)(2,73,10,6,72,5)(3,76,90)(7,65,30)(8,71,80,16,62,31)(9,43,36)(12,
    75,48,13,63,99)(14,44,94)(15,\\41,33,22,39,79)(17,77,58,21,64,97)(23,88,83,69,54,
    59)(24,61)(25,26,52,70,66,85)(27,93,53,28,38,95)(29,91,82,68,47,\\55)(32,51,86,89,
    96,87)(34,37,49,67,45,92)(35,50)(46,74,78)(56,84)(57,81)(60,100,98),

  (1,42,20)(2,72,10)(3,76,90)(4,40,11)(5,6,73)(7,65,30)(8,62,80)(9,43,36)(12,63,
    48)(13,75,99)(14,44,94)(15,39,33)\\(16,71,31)(17,64,58)(21,77,97)(22,41,79)(23,54,
    83)(25,66,52)(26,85,70)(27,38,53)(28,93,95)(29,47,82)(32,96,86)(34,\\45,49)(37,92,
    67)(46,74,78)(51,87,89)(55,68,91)(59,69,88)(60,100,98),

  (1,53,42,27,20,38)(2,87,72,89,10,51)(3,94,76,14,90,44)(4,98,40,60,11,100)(5,86,6,
    32,73,96)(7,49,65,34,30,45)(8,\\55,62,68,80,91)(9,52,43,25,36,66)(12,58,63,17,48,
    64)(13,54,75,83,99,23)(15,31,39,16,33,71)(18,56)(19,57)(21,47,77,\\82,97,29)(22,
    59,41,69,79,88)(24,35)(26,67,85,37,70,92)(28,74,93,78,95,46)(50,61)(81,84),

  (1,53,42,27,20,38)(2,96,72,86,10,32)(3,94,76,14,90,44)(4,95,40,28,11,93)(5,89,6,
    51,73,87)(7,36,65,9,30,43)(8,13,\\62,75,80,99)(12,71,63,31,48,16)(15,17,39,64,33,
    58)(18,19)(21,41,77,79,97,22)(23,55,54,68,83,91)(24,50)(25,49,66,34,\\52,45)(26,
    67,85,37,70,92)(29,59,47,69,82,88)(35,61)(46,60,74,100,78,98)(56,57)(81,84),

  (1,93,74,27,40,98)(2,96,72,86,10,32)(3,44,90,14,76,94)(4,100,20,28,46,53)(5,51)(6,
    87)(7,37,52,9,26,49)(8,47,17,\\75,88,33)(11,60,42,95,78,38)(12,77,91,31,22,54)(13,
    69,39,80,29,58)(15,62,82,64,99,59)(16,41,83,63,97,55)(18,84,57,\\19,81,56)(21,68,
    71,79,23,48)(24,50)(25,43,85,34,65,92)(30,67,66,36,70,45)(35,61)(73,89),

  (1,95,46,27,11,100)(2,86)(3,94,76,14,90,44)(4,98,42,28,74,38)(5,87,73,51,6,89)(7,
    92,66,9,85,45)(8,82,58,75,59,39)\\(10,96)(12,97,68,31,41,23)(13,88,15,80,47,
    64)(16,79,54,63,21,91)(17,99,69,33,62,29)(18,84,57,19,81,56)(20,93,78,53,\\40,
    60)(22,83,48,77,55,71)(24,50)(25,36,26,34,30,37)(32,72)(35,61)(43,70,49,65,67,
    52).

The character table of $G^{s_{20}}$:\\
\begin{tabular}{c|cccccccccccccccccc}
  & & & & & & & & & &10 & & & & & & & &\\\hline
$\chi_{20}^{(1)}$&1&1&1&1&1&1&1&1&1&1&1&1&1&1&1&1&1&1
\\$\chi_{20}^{(2)}$&1&-1&1&1&-1&1&-1&1&-1&-1&1&1&-1&1&1&-1&-1&-1
\\$\chi_{20}^{(3)}$&1&-1&1&1&-1&1&1&-1&1&1&-1&1&-1&1&-1&1&1&1
\\$\chi_{20}^{(4)}$&1&1&1&1&1&1&-1&-1&-1&-1&-1&1&1&1&-1&-1&-1&-1
\\$\chi_{20}^{(5)}$&1&-1&1&A&-A&A&-1&1&-1&-/A&/A&/A&-/A&/A&A&-A&-/A&-A
\\$\chi_{20}^{(6)}$&1&-1&1&/A&-/A&/A&-1&1&-1&-A&A&A&-A&A&/A&-/A&-A&-/A
\\$\chi_{20}^{(7)}$&1&-1&1&A&-A&A&1&-1&1&/A&-/A&/A&-/A&/A&-A&A&/A&A
\\$\chi_{20}^{(8)}$&1&-1&1&/A&-/A&/A&1&-1&1&A&-A&A&-A&A&-/A&/A&A&/A
\\$\chi_{20}^{(9)}$&1&1&1&A&A&A&-1&-1&-1&-/A&-/A&/A&/A&/A&-A&-A&-/A&-A
\\$\chi_{20}^{(10)}$&1&1&1&/A&/A&/A&-1&-1&-1&-A&-A&A&A&A&-/A&-/A&-A&-/A
\\$\chi_{20}^{(11)}$&1&1&1&A&A&A&1&1&1&/A&/A&/A&/A&/A&A&A&/A&A
\\$\chi_{20}^{(12)}$&1&1&1&/A&/A&/A&1&1&1&A&A&A&A&A&/A&/A&A&/A
\\$\chi_{20}^{(13)}$&2&.&-1&-1&.&2&-2&.&1&-2&.&-1&.&2&.&-2&1&1
\\$\chi_{20}^{(14)}$&2&.&-1&-1&.&2&2&.&-1&2&.&-1&.&2&.&2&-1&-1
\\$\chi_{20}^{(15)}$&2&.&-1&-A&.&B&-2&.&1&-/B&.&-/A&.&/B&.&-B&/A&A
\\$\chi_{20}^{(16)}$&2&.&-1&-/A&.&/B&-2&.&1&-B&.&-A&.&B&.&-/B&A&/A
\\$\chi_{20}^{(17)}$&2&.&-1&-A&.&B&2&.&-1&/B&.&-/A&.&/B&.&B&-/A&-A
\\$\chi_{20}^{(18)}$&2&.&-1&-/A&.&/B&2&.&-1&B&.&-A&.&B&.&/B&-A&-/A

\end{tabular}

\noindent \noindent where   A = E(3)$^2$
  = (-1-ER(-3))/2 = -1-b3;
B = 2*E(3)$^2$
  = -1-ER(-3) = -1-i3.

The generators of $G^{s_{21}}$ are:\\
(  1,  5, 22, 79, 61, 90, 59, 33)(  2, 88, 65, 69,  8, 13, 67, 12)(
3, 49)
    (  4, 93, 89, 92, 62, 39, 15, 44)(  6, 75, 27, 16, 84, 95, 53, 56)
    (  7, 58, 23, 11, 51, 20, 54, 72)(  9, 71, 83, 87, 64, 24, 32, 14)
    ( 10, 40, 30, 48)( 17, 86, 34, 77)( 18, 45, 97, 41, 80, 38, 73, 37)
    ( 21, 94, 63, 42, 76, 99, 28, 36)( 25, 47, 29, 35, 68, 91, 85,100)
    ( 26, 70, 98, 55, 74, 52, 57, 78)( 43, 96)( 46, 81, 50, 66)( 60, 82),
  (  1, 21, 22, 63, 61, 76, 59, 28)(  2, 85, 65, 25,  8, 29, 67, 68)(  3, 43)
    (  4, 52, 89, 78, 62, 70, 15, 55)(  5, 94, 79, 42, 90, 99, 33, 36)
    (  6, 95, 27, 56, 84, 75, 53, 16)(  7, 38, 23, 37, 51, 45, 54, 41)
    (  9, 71, 83, 87, 64, 24, 32, 14)( 10, 40, 30, 48)( 11, 18, 20, 97, 72, 80, 58,
      73)( 12, 91, 88,100, 69, 47, 13, 35)( 17, 86, 34, 77)( 19, 31)
    ( 26, 39, 98, 44, 74, 93, 57, 92)( 46, 81, 50, 66)( 49, 96).

The representatives of conjugacy classes of   $G^{s_{21}}$ are:\\
 (1),
  (1,5,22,79,61,90,59,33)(2,88,65,69,8,13,67,12)(3,49)(4,93,89,92,62,39,15,
    44)(6,75,27,16,84,95,53,56)(7,58,23,\\11,51,20,54,72)(9,71,83,87,64,24,32,14)(10,
    40,30,48)(17,86,34,77)(18,45,97,41,80,38,73,37)(21,94,63,42,76,99,28,\\36)(25,47,
    29,35,68,91,85,100)(26,70,98,55,74,52,57,78)(43,96)(46,81,50,66)(60,82),

  (1,21,22,63,61,76,59,28)(2,85,65,25,8,29,67,68)(3,43)(4,52,89,78,62,70,15,55)(5,
    94,79,42,90,99,33,36)(6,95,27,\\56,84,75,53,16)(7,38,23,37,51,45,54,41)(9,71,83,
    87,64,24,32,14)(10,40,30,48)(11,18,20,97,72,80,58,73)(12,91,88,\\100,69,47,13,
    35)(17,86,34,77)(19,31)(26,39,98,44,74,93,57,92)(46,81,50,66)(49,96),

  (1,22,61,59)(2,65,8,67)(4,89,62,15)(5,79,90,33)(6,27,84,53)(7,23,51,54)(9,83,64,
    32)(10,30)(11,20,72,58)(12,88,\\69,13)(14,71,87,24)(16,95,56,75)(17,34)(18,97,80,
    73)(21,63,76,28)(25,29,68,85)(26,98,74,57)(35,91,100,47)(36,94,\\42,99)(37,45,41,
    38)(39,44,93,92)(40,48)(46,50)(52,78,70,55)(66,81)(77,86),

  (1,28,59,76,61,63,22,21)(2,68,67,29,8,25,65,85)(3,43)(4,55,15,70,62,78,89,52)(5,
    36,33,99,90,42,79,94)(6,16,53,\\75,84,56,27,95)(7,41,54,45,51,37,23,38)(9,14,32,
    24,64,87,83,71)(10,48,30,40)(11,73,58,80,72,97,20,18)(12,35,13,\\47,69,100,88,
    91)(17,77,34,86)(19,31)(26,92,57,93,74,44,98,39)(46,66,50,81)(49,96),

  (1,33,59,90,61,79,22,5)(2,12,67,13,8,69,65,88)(3,49)(4,44,15,39,62,92,89,93)(6,56,
    53,95,84,16,27,75)(7,72,54,\\20,51,11,23,58)(9,14,32,24,64,87,83,71)(10,48,30,
    40)(17,77,34,86)(18,37,73,38,80,41,97,45)(21,36,28,99,76,42,63,\\94)(25,100,85,91,
    68,35,29,47)(26,78,57,52,74,55,98,70)(43,96)(46,66,50,81)(60,82),

  (1,36)(2,91)(3,96)(4,74)(5,21)(6,84)(7,80)(8,47)(11,37)(12,68)(13,29)(15,98)(16,
    56)(18,51)(19,31)(20,45)(22,\\94)(23,73)(25,69)(26,62)(27,53)(28,33)(35,67)(38,
    58)(39,70)(41,72)(42,61)(43,49)(44,55)(52,93)(54,97)(57,89)\\(59,99)(60,82)(63,
    79)(65,100)(75,95)(76,90)(78,92)(85,88),

    (1,42)(2,47)(3,96)(4,26)(5,76)(7,
    18)(8,91)(9,64)(11,41)(12,25)(13,85)(14,87)(15,57)(19,31)(20,38)(21,90)(22,\\
    99)(23,97)(24,71)(28,79)(29,88)(32,83)(33,63)(35,65)(36,61)(37,72)(39,52)(43,
    49)(44,78)(45,58)(51,80)(54,73)\\(55,92)(59,94)(60,82)(62,74)(67,100)(68,69)(70,
    93)(89,98),

     (1,59,61,22)(2,67,8,65)(4,15,62,89)(5,33,90,79)(6,53,84,27)(7,54,
    51,23)(9,32,64,83)(10,30)(11,58,72,20)(12,13,\\69,88)(14,24,87,71)(16,75,56,
    95)(17,34)(18,73,80,97)(21,28,76,63)(25,85,68,29)(26,57,74,98)(35,47,100,91)(36,\\
    99,42,94)(37,38,41,45)(39,92,93,44)(40,48)(46,50)(52,55,70,78)(66,81)(77,86),

  (1,61)(2,8)(4,62)(5,90)(6,84)(7,51)(9,64)(11,72)(12,69)(13,88)(14,87)(15,89)(16,
    56)(18,80)(20,58)(21,76)(22,\\59)(23,54)(24,71)(25,68)(26,74)(27,53)(28,63)(29,
    85)(32,83)(33,79)(35,100)(36,42)(37,41)(38,45)(39,93)(44,92)\\(47,91)(52,70)(55,
    78)(57,98)(65,67)(73,97)(75,95)(94,99),

     (1,63,59,21,61,28,22,76)(2,25,67,85,8,
    68,65,29)(3,43)(4,78,15,52,62,55,89,70)(5,42,33,94,90,36,79,99)(6,56,53,\\95,84,
    16,27,75)(7,37,54,38,51,41,23,45)(9,87,32,71,64,14,83,24)(10,48,30,40)(11,97,58,
    18,72,73,20,80)(12,100,13,\\91,69,35,88,47)(17,77,34,86)(19,31)(26,44,57,39,74,92,
    98,93)(46,66,50,81)(49,96),

     (1,76,22,28,61,21,59,63)(2,29,65,68,8,85,67,25)(3,
    43)(4,70,89,55,62,52,15,78)(5,99,79,36,90,94,33,42)(6,75,27,\\16,84,95,53,56)(7,
    45,23,41,51,38,54,37)(9,24,83,14,64,71,32,87)(10,40,30,48)(11,80,20,73,72,18,58,
    97)(12,47,88,\\35,69,91,13,100)(17,86,34,77)(19,31)(26,93,98,92,74,39,57,44)(46,
    81,50,66)(49,96),

    (1,79,59,5,61,33,22,90)(2,69,67,88,8,12,65,13)(3,49)(4,92,
    15,93,62,44,89,39)(6,16,53,75,84,56,27,95)(7,11,54,\\58,51,72,23,20)(9,87,32,71,
    64,14,83,24)(10,48,30,40)(17,77,34,86)(18,41,73,45,80,37,97,38)(21,42,28,94,76,
    36,63,\\99)(25,35,85,47,68,100,29,91)(26,55,57,70,74,78,98,52)(43,96)(46,66,50,
    81)(60,82),

     (1,90,22,33,61,5,59,79)(2,13,65,12,8,88,67,69)(3,49)(4,39,89,44,
    62,93,15,92)(6,95,27,56,84,75,53,16)(7,20,23,\\72,51,58,54,11)(9,24,83,14,64,71,
    32,87)(10,40,30,48)(17,86,34,77)(18,38,97,37,80,45,73,41)(21,99,63,36,76,94,28,
    42)(25,91,29,100,68,47,85,35)(26,52,98,78,74,70,57,55)(43,96)(46,81,50,66)(60,
    82),

     (1,94,61,99)(2,100,8,35)(3,96)(4,57,62,98)(5,63,90,28)(6,53,84,27)(7,73,
    51,97)(9,83,64,32)(10,30)(11,45,72,38)\\(12,85,69,29)(13,68,88,25)(14,71,87,
    24)(15,74,89,26)(16,75,56,95)(17,34)(18,54,80,23)(19,31)(20,41,58,37)(21,79,\\76, 33)(22,42,59,36)(39,55,93,78)(40,48)(43,49)(44,52,92,70)(46,50)(47,67,91,65)(60,
    82)(66,81)(77,86),

    (1,99,61,94)(2,35,8,100)(3,96)(4,98,62,57)(5,28,90,63)(6,
    27,84,53)(7,97,51,73)(9,32,64,83)(10,30)(11,38,72,45)\\(12,29,69,85)(13,25,88,
    68)(14,24,87,71)(15,26,89,74)(16,95,56,75)(17,34)(18,23,80,54)(19,31)(20,37,58,
    41)(21,33,\\76,79)(22,36,59,42)(39,78,93,55)(40,48)(43,49)(44,70,92,52)(46,50)(47,
    65,91,67)(60,82)(66,81)(77,86).

The character table of $G^{s_{21}}$:\\
\begin{tabular}{c|cccccccccccccccc}
  & & & & & & & & & &10 & & & & & &\\\hline
$\chi_{21}^{(1)}$&1&1&1&1&1&1&1&1&1&1&1&1&1&1&1&1
\\$\chi_{21}^{(2)}$&1&-1&-1&1&-1&-1&1&1&1&1&-1&-1&-1&-1&1&1
\\$\chi_{21}^{(3)}$&1&-1&1&1&1&-1&-1&-1&1&1&1&1&-1&-1&-1&-1
\\$\chi_{21}^{(4)}$&1&1&-1&1&-1&1&-1&-1&1&1&-1&-1&1&1&-1&-1
\\$\chi_{21}^{(5)}$&1&A&A&-1&-A&-A&1&1&-1&1&-A&A&-A&A&-1&-1
\\$\chi_{21}^{(6)}$&1&-A&-A&-1&A&A&1&1&-1&1&A&-A&A&-A&-1&-1
\\$\chi_{21}^{(7)}$&1&A&-A&-1&A&-A&-1&-1&-1&1&A&-A&-A&A&1&1
\\$\chi_{21}^{(8)}$&1&-A&A&-1&-A&A&-1&-1&-1&1&-A&A&A&-A&1&1
\\$\chi_{21}^{(9)}$&1&B&B&-A&/B&/B&1&-1&A&-1&-/B&-B&-/B&-B&-A&A
\\$\chi_{21}^{(10)}$&1&-/B&-/B&A&-B&-B&1&-1&-A&-1&B&/B&B&/B&A&-A
\\$\chi_{21}^{(11)}$&1&/B&/B&A&B&B&1&-1&-A&-1&-B&-/B&-B&-/B&A&-A
\\$\chi_{21}^{(12)}$&1&-B&-B&-A&-/B&-/B&1&-1&A&-1&/B&B&/B&B&-A&A
\\$\chi_{21}^{(13)}$&1&B&-B&-A&-/B&/B&-1&1&A&-1&/B&B&-/B&-B&A&-A
\\$\chi_{21}^{(14)}$&1&-/B&/B&A&B&-B&-1&1&-A&-1&-B&-/B&B&/B&-A&A
\\$\chi_{21}^{(15)}$&1&/B&-/B&A&-B&B&-1&1&-A&-1&B&/B&-B&-/B&-A&A
\\$\chi_{21}^{(16)}$&1&-B&B&-A&/B&-/B&-1&1&A&-1&-/B&-B&/B&B&A&-A
\end{tabular}

\noindent \noindent where   A = -E(4)
  = -ER(-1) = -i;
B = -E(8).

The generators of $G^{s_{22}}$ are:\\
(  1, 42, 17, 71, 24)(  2, 96, 70, 54, 25)(  3, 66, 47, 41, 83)
    (  4,  5, 26, 46,  6)(  7, 68, 90, 60, 73)(  8, 44, 99, 22, 97)(  9, 36, 49, 34,
     20)( 10, 91, 84, 23, 78)( 11, 12, 82, 13, 35)( 14, 79, 98, 55, 92)
    ( 15, 18, 93, 32, 65)( 19, 27, 58, 38, 43)( 21, 77, 69, 57, 85)( 28, 48, 87, 89,
     86)( 29, 40, 37, 72, 88)( 30, 52, 62, 39, 75)( 31, 61, 56, 64, 45)
    ( 33, 50, 94, 95, 67)( 51,100, 59, 63, 80),
  (  1,  4, 35, 15, 75)(  2, 70, 25, 96, 54)(  3, 41, 66, 83, 47)(  5, 11, 18, 30,
      42)(  6, 13, 65, 39, 24)(  8, 95, 38, 88, 49)(  9, 22, 50, 27, 37)
    ( 10, 79, 48, 80, 31)( 12, 93, 52, 17, 26)( 14, 28, 63, 45, 78)( 16, 81, 53, 76,
     74)( 19, 40, 20, 99, 33)( 21, 85, 57, 69, 77)( 23, 92, 86, 59, 64)
    ( 29, 34, 44, 67, 43)( 32, 62, 71, 46, 82)( 36, 97, 94, 58, 72)( 51, 61, 91, 98,
     87)( 55, 89,100, 56, 84).

The representatives of conjugacy classes of   $G^{s_{22}}$ are:\\
 (1),
 (1,4,35,15,75)(2,70,25,96,54)(3,41,66,83,47)(5,11,18,30,42)(6,13,65,39,
    24)(8,95,38,88,49)(9,22,50,27,37)(10,79,\\48,80,31)(12,93,52,17,26)(14,28,63,45,
    78)(16,81,53,76,74)(19,40,20,99,33)(21,85,57,69,77)(23,92,86,59,64)(29,34,\\44,67,
    43)(32,62,71,46,82)(36,97,94,58,72)(51,61,91,98,87)(55,89,100,56,84),

  (1,5,12,32,39)(2,54,96,25,70)(3,83,41,47,66)(4,11,93,62,24)(6,35,18,52,71)(7,68,
    90,60,73)(8,67,19,37,36)(9,97,\\95,43,40)(10,98,89,59,45)(13,15,30,17,46)(14,48,
    51,56,23)(16,81,53,76,74)(20,22,94,38,29)(26,82,65,75,42)(27,72,\\49,44,33)(28,80,
    61,84,92)(31,91,55,86,63)(34,99,50,58,88)(64,78,79,87,100),

  (1,6,82,93,30)(2,96,70,54,25)(3,47,83,66,41)(4,13,32,52,42)(5,35,65,62,17)(7,73,
    60,90,68)(8,94,27,40,34)(9,99,\\67,38,72)(10,14,86,100,61)(11,15,39,71,26)(12,18,
    75,24,46)(16,81,53,76,74)(19,29,49,97,50)(20,44,95,58,37)(21,57,\\77,85,69)(22,33,
    43,88,36)(23,55,87,80,45)(28,59,56,91,79)(31,78,92,89,51)(48,63,64,84,98),

  (1,11,52,46,65)(3,47,83,66,41)(4,18,17,82,39)(5,93,71,13,75)(6,15,42,12,62)(7,68,
    90,60,73)(8,43,20,50,72)(9,94,\\88,44,19)(10,87,56,92,63)(14,80,91,89,64)(16,53,
    74,81,76)(21,85,57,69,77)(22,58,49,67,40)(23,28,31,98,100)(24,35,\\30,26,32)(27,
    36,95,29,99)(33,37,97,38,34)(45,79,51,84,86)(48,61,55,59,78),

  (1,12,39,5,32)(2,96,70,54,25)(3,41,66,83,47)(4,93,24,11,62)(6,18,71,35,52)(7,90,
    73,68,60)(8,19,36,67,37)(9,95,\\40,97,43)(10,89,45,98,59)(13,30,46,15,17)(14,51,
    23,48,56)(16,53,74,81,76)(20,94,29,22,38)(26,65,42,82,75)(27,49,\\33,72,44)(28,61,
    92,80,84)(31,55,63,91,86)(34,50,88,99,58)(64,79,100,78,87),

  (1,13,62,26,18)(2,54,96,25,70)(4,65,71,12,30)(5,15,24,82,52)(6,32,17,11,75)(7,73,
    60,90,68)(8,58,9,33,29)(10,28,\\64,55,51)(14,59,84,87,31)(16,53,74,81,76)(19,34,
    95,72,22)(20,67,88,97,27)(21,69,85,77,57)(23,89,61,79,63)(35,39,\\46,93,42)(36,50,
    40,44,38)(37,99,43,49,94)(45,92,100,91,48)(56,98,80,78,86),

  (1,15,4,75,35)(2,96,70,54,25)(3,83,41,47,66)(5,30,11,42,18)(6,39,13,24,65)(8,88,
    95,49,38)(9,27,22,37,50)(10,80,\\79,31,48)(12,17,93,26,52)(14,45,28,78,63)(16,76,
    81,74,53)(19,99,40,33,20)(21,69,85,77,57)(23,59,92,64,86)(29,67,\\34,43,44)(32,46,
    62,82,71)(36,58,97,72,94)(51,98,61,87,91)(55,56,89,84,100),

  (1,17,24,42,71)(2,70,25,96,54)(3,47,83,66,41)(4,26,6,5,46)(7,90,73,68,60)(8,99,97,
    44,22)(9,49,20,36,34)(10,84,\\78,91,23)(11,82,35,12,13)(14,98,92,79,55)(15,93,65,
    18,32)(19,58,43,27,38)(21,69,85,77,57)(28,87,86,48,89)(29,37,\\88,40,72)(30,62,75,
    52,39)(31,56,45,61,64)(33,94,67,50,95)(51,59,80,100,63),

  (1,18,26,62,13)(2,70,25,96,54)(4,30,12,71,65)(5,52,82,24,15)(6,75,11,17,32)(7,68,
    90,60,73)(8,29,33,9,58)(10,51,\\55,64,28)(14,31,87,84,59)(16,76,81,74,53)(19,22,
    72,95,34)(20,27,97,88,67)(21,57,77,85,69)(23,63,79,61,89)(35,42,\\93,46,39)(36,38,
    44,40,50)(37,94,49,43,99)(45,48,91,100,92)(56,86,78,80,98),

  (1,24,71,17,42)(2,25,54,70,96)(3,83,41,47,66)(4,6,46,26,5)(7,73,60,90,68)(8,97,22,
    99,44)(9,20,34,49,36)(10,78,\\23,84,91)(11,35,13,82,12)(14,92,55,98,79)(15,65,32,
    93,18)(19,43,38,58,27)(21,85,57,69,77)(28,86,89,87,48)(29,88,\\72,37,40)(30,75,39,
    62,52)(31,45,64,56,61)(33,67,95,94,50)(51,80,63,59,100),

  (1,26,13,18,62)(2,25,54,70,96)(4,12,65,30,71)(5,82,15,52,24)(6,11,32,75,17)(7,90,
    73,68,60)(8,33,58,29,9)(10,55,\\28,51,64)(14,87,59,31,84)(16,81,53,76,74)(19,72,
    34,22,95)(20,97,67,27,88)(21,77,69,57,85)(23,79,89,63,61)(35,93,\\39,42,46)(36,44,
    50,38,40)(37,49,99,94,43)(45,91,92,48,100)(56,78,98,86,80),

  (1,30,93,82,6)(2,25,54,70,96)(3,41,66,83,47)(4,42,52,32,13)(5,17,62,65,35)(7,68,
    90,60,73)(8,34,40,27,94)(9,72,\\38,67,99)(10,61,100,86,14)(11,26,71,39,15)(12,46,
    24,75,18)(16,74,76,53,81)(19,50,97,49,29)(20,37,58,95,44)(21,69,\\85,77,57)(22,36,
    88,43,33)(23,45,80,87,55)(28,79,91,56,59)(31,51,89,92,78)(48,98,84,64,63),

  (1,32,5,39,12)(2,25,54,70,96)(3,47,83,66,41)(4,62,11,24,93)(6,52,35,71,18)(7,60,
    68,73,90)(8,37,67,36,19)(9,43,\\97,40,95)(10,59,98,45,89)(13,17,15,46,30)(14,56,
    48,23,51)(16,76,81,74,53)(20,38,22,29,94)(26,75,82,42,65)(27,44,\\72,33,49)(28,84,
    80,92,61)(31,86,91,63,55)(34,58,99,88,50)(64,87,78,100,79),

  (1,35,75,4,15)(2,25,54,70,96)(3,66,47,41,83)(5,18,42,11,30)(6,65,24,13,39)(8,38,
    49,95,88)(9,50,37,22,27)(10,48,\\31,79,80)(12,52,26,93,17)(14,63,78,28,45)(16,53,
    74,81,76)(19,20,33,40,99)(21,57,77,85,69)(23,86,64,92,59)(29,44,\\43,34,67)(32,71,
    82,62,46)(36,94,72,97,58)(51,91,87,61,98)(55,100,84,89,56),

  (1,39,32,12,5)(2,70,25,96,54)(3,66,47,41,83)(4,24,62,93,11)(6,71,52,18,35)(7,73,
    60,90,68)(8,36,37,19,67)(9,40,\\43,95,97)(10,45,59,89,98)(13,46,17,30,15)(14,23,
    56,51,48)(16,74,76,53,81)(20,29,38,94,22)(26,42,75,65,82)(27,33,\\44,49,72)(28,92,
    84,61,80)(31,63,86,55,91)(34,88,58,50,99)(64,100,87,79,78),

  (1,42,17,71,24)(2,96,70,54,25)(3,66,47,41,83)(4,5,26,46,6)(7,68,90,60,73)(8,44,99,
    22,97)(9,36,49,34,20)(10,91,\\84,23,78)(11,12,82,13,35)(14,79,98,55,92)(15,18,93,
    32,65)(19,27,58,38,43)(21,77,69,57,85)(28,48,87,89,86)(29,40,\\37,72,88)(30,52,62,
    39,75)(31,61,56,64,45)(33,50,94,95,67)(51,100,59,63,80),

  (1,46,11,65,52)(3,66,47,41,83)(4,82,18,39,17)(5,13,93,75,71)(6,12,15,62,42)(7,60,
    68,73,90)(8,50,43,72,20)(9,44,\\94,19,88)(10,92,87,63,56)(14,89,80,64,91)(16,81,
    53,76,74)(21,69,85,77,57)(22,67,58,40,49)(23,98,28,100,31)(24,26,\\35,32,30)(27,
    29,36,99,95)(33,38,37,34,97)(45,84,79,86,51)(48,59,61,78,55),

  (1,52,65,11,46)(3,83,41,47,66)(4,17,39,18,82)(5,71,75,93,13)(6,42,62,15,12)(7,90,
    73,68,60)(8,20,72,43,50)(9,88,\\19,94,44)(10,56,63,87,92)(14,91,64,80,89)(16,74,
    76,53,81)(21,57,77,85,69)(22,49,40,58,67)(23,31,100,28,98)(24,30,\\32,35,26)(27,
    95,99,36,29)(33,97,34,37,38)(45,51,86,79,84)(48,55,78,61,59),

  (1,62,18,13,26)(2,96,70,54,25)(4,71,30,65,12)(5,24,52,15,82)(6,17,75,32,11)(7,60,
    68,73,90)(8,9,29,58,33)(10,64,\\51,28,55)(14,84,31,59,87)(16,74,76,53,81)(19,95,
    22,34,72)(20,88,27,67,97)(21,85,57,69,77)(23,61,63,89,79)(35,46,\\42,39,93)(36,40,
    38,50,44)(37,43,94,99,49)(45,100,48,92,91)(56,80,86,98,78),

  (1,65,46,52,11)(3,41,66,83,47)(4,39,82,17,18)(5,75,13,71,93)(6,62,12,42,15)(7,73,
    60,90,68)(8,72,50,20,43)(9,19,\\44,88,94)(10,63,92,56,87)(14,64,89,91,80)(16,76,
    81,74,53)(21,77,69,57,85)(22,40,67,49,58)(23,100,98,31,28)(24,32,\\26,30,35)(27,
    99,29,95,36)(33,34,38,97,37)(45,86,84,51,79)(48,78,59,55,61),

  (1,71,42,24,17)(2,54,96,25,70)(3,41,66,83,47)(4,46,5,6,26)(7,60,68,73,90)(8,22,44,
    97,99)(9,34,36,20,49)(10,23,91,\\78,84)(11,13,12,35,82)(14,55,79,92,98)(15,32,18,
    65,93)(19,38,27,43,58)(21,57,77,85,69)(28,89,48,86,87)(29,72,40,88,\\37)(30,39,52,
    75,62)(31,64,61,45,56)(33,95,50,67,94)(51,63,100,80,59),

  (1,75,15,35,4)(2,54,96,25,70)(3,47,83,66,41)(5,42,30,18,11)(6,24,39,65,13)(8,49,
    88,38,95)(9,37,27,50,22)(10,31,80,\\48,79)(12,26,17,52,93)(14,78,45,63,28)(16,74,
    76,53,81)(19,33,99,20,40)(21,77,69,57,85)(23,64,59,86,92)(29,43,67,44,\\34)(32,82,
    46,71,62)(36,72,58,94,97)(51,87,98,91,61)(55,84,56,100,89),

  (1,82,30,6,93)(2,70,25,96,54)(3,83,41,47,66)(4,32,42,13,52)(5,65,17,35,62)(7,60,
    68,73,90)(8,27,34,94,40)(9,67,72,\\99,38)(10,86,61,14,100)(11,39,26,15,71)(12,75,
    46,18,24)(16,53,74,81,76)(19,49,50,29,97)(20,95,37,44,58)(21,77,69,\\57,85)(22,43,
    36,33,88)(23,87,45,55,80)(28,56,79,59,91)(31,92,51,78,89)(48,64,98,63,84),

  (1,93,6,30,82)(2,54,96,25,70)(3,66,47,41,83)(4,52,13,42,32)(5,62,35,17,65)(7,90,
    73,68,60)(8,40,94,34,27)(9,38,99,\\72,67)(10,100,14,61,86)(11,71,15,26,39)(12,24,
    18,46,75)(16,76,81,74,53)(19,97,29,50,49)(20,58,44,37,95)(21,85,57,69,\\77)(22,88,
    33,36,43)(23,80,55,45,87)(28,91,59,79,56)(31,89,78,51,92)(48,84,63,98,64).

The character table of $G^{s_{22}}$:\\

${}_{\!\!\!\!\!\!\!
 \!\!\!\!\!\!\!\!\!\!_{
\!\!\!\!_{ \!\!\!\!\!\!\!\!\!\!\! _{ \small
\begin{tabular}{c|ccccccccccccccccccccccccc}
  & & & & & & & & & & 10& & & & & & & & & &20 & & & & &\\\hline
$\chi_{22}^{(1)}$&1&1&1&1&1&1&1&1&1&1&1&1&1&1&1&1&1&1&1&1&1&1&1&1&1
\\$\chi_{22}^{(2)}$&1&1&A&/A&A&B&/A&1&B&A&/A&B&A&/B&1&/A&A&/B&B&/B&/A&/B&1&/B&B
\\$\chi_{22}^{(3)}$&1&1&B&/B&B&/A&/B&1&/A&B&/B&/A&B&A&1&/B&B&A&/A&A&/B&A&1&A&/A
\\$\chi_{22}^{(4)}$&1&1&/B&B&/B&A&B&1&A&/B&B&A&/B&/A&1&B&/B&/A&A&/A&B&/A&1&/A&A
\\$\chi_{22}^{(5)}$&1&1&/A&A&/A&/B&A&1&/B&/A&A&/B&/A&B&1&A&/A&B&/B&B&A&B&1&B&/B
\\$\chi_{22}^{(6)}$&1&A&1&B&A&1&/B&/B&/B&B&A&/A&/B&1&B&1&/A&/B&B&A&/A&B&/A&/A&A
\\$\chi_{22}^{(7)}$&1&B&1&/A&B&1&A&A&A&/A&B&/B&A&1&/A&1&/B&A&/A&B&/B&/A&/B&/B&B
\\$\chi_{22}^{(8)}$&1&/B&1&A&/B&1&/A&/A&/A&A&/B&B&/A&1&A&1&B&/A&A&/B&B&A&B&B&/B
\\$\chi_{22}^{(9)}$&1&/A&1&/B&/A&1&B&B&B&/B&/A&A&B&1&/B&1&A&B&/B&/A&A&/B&A&A&/A
\\$\chi_{22}^{(10)}$&1&A&A&A&B&B&B&/B&1&/B&1&A&/A&/B&B&/A&1&A&/A&/A&/B&1&/A&B&/B
\\$\chi_{22}^{(11)}$&1&B&B&B&/A&/A&/A&A&1&A&1&B&/B&A&/A&/B&1&B&/B&/B&A&1&/B&/A&A
\\$\chi_{22}^{(12)}$&1&/B&/B&/B&A&A&A&/A&1&/A&1&/B&B&/A&A&B&1&/B&B&B&/A&1&B&A&/A
\\$\chi_{22}^{(13)}$&1&/A&/A&/A&/B&/B&/B&B&1&B&1&/A&A&B&/B&A&1&/A&A&A&B&1&A&/B&B
\\$\chi_{22}^{(14)}$&1&A&B&1&/B&/A&A&/B&B&/A&/A&/B&1&A&B&/B&A&/A&A&B&B&/B&/A&1&1
\\$\chi_{22}^{(15)}$&1&B&/A&1&A&/B&B&A&/A&/B&/B&A&1&B&/A&A&B&/B&B&/A&/A&A&/B&1&1
\\$\chi_{22}^{(16)}$&1&/B&A&1&/A&B&/B&/A&A&B&B&/A&1&/B&A&/A&/B&B&/B&A&A&/A&B&1&1
\\$\chi_{22}^{(17)}$&1&/A&/B&1&B&A&/A&B&/B&A&A&B&1&/A&/B&B&/A&A&/A&/B&/B&B&A&1&1
\\$\chi_{22}^{(18)}$&1&A&/B&/A&/A&A&1&/B&/A&1&/B&1&A&/A&B&B&B&B&/B&1&A&A&/A&/B&B
\\$\chi_{22}^{(19)}$&1&B&A&/B&/B&B&1&A&/B&1&A&1&B&/B&/A&/A&/A&/A&A&1&B&B&/B&A&/A
\\$\chi_{22}^{(20)}$&1&/B&/A&B&B&/B&1&/A&B&1&/A&1&/B&B&A&A&A&A&/A&1&/B&/B&B&/A&A
\\$\chi_{22}^{(21)}$&1&/A&B&A&A&/A&1&B&A&1&B&1&/A&A&/B&/B&/B&/B&B&1&/A&/A&A&B&/B
\\$\chi_{22}^{(22)}$&1&A&/A&/B&1&/B&/A&/B&A&A&B&B&B&B&B&A&/B&1&1&/B&1&/A&/A&A&/A
\\$\chi_{22}^{(23)}$&1&B&/B&A&1&A&/B&A&B&B&/A&/A&/A&/A&/A&B&A&1&1&A&1&/B&/B&B&/B
\\$\chi_{22}^{(24)}$&1&/B&B&/A&1&/A&B&/A&/B&/B&A&A&A&A&A&/B&/A&1&1&/A&1&B&B&/B&B
\\$\chi_{22}^{(25)}$&1&/A&A&B&1&B&A&B&/A&/A&/B&/B&/B&/B&/B&/A&B&1&1&B&1&A&A&/A&A
\end{tabular}
 }}}}$

\noindent \noindent where   A = E(5)$^4$; B = E(5)$^3$.

The generators of $G^{s_{23}}$ are:\\
 (  1, 44, 84, 95, 20, 57, 68, 28, 92, 90, 52)(  2, 26, 43, 23, 94,  3, 50, 16,100,
     21, 33)(  4, 99, 83, 30, 15, 69, 93, 24, 39,  5, 75)(  6, 72, 22, 34, 89, 98,
      53, 66, 47, 77, 40)(  7, 42, 76, 37, 19, 88, 17, 27, 45, 96, 80)
    (  8, 49, 71, 29, 14, 79, 85, 25, 48, 54, 62)(  9, 64, 91, 31, 63, 10, 87, 61,
      38, 56, 51)( 11, 65, 46, 32, 67, 13, 82, 70, 35, 73, 81)( 12, 58, 60, 36, 86,
      74, 18, 78, 41, 97, 59)

The representatives of conjugacy classes of   $G^{s_{23}}$ are:\\
 (1),
 (1,20,92,44,57,90,84,68,52,95,28)(2,94,100,26,3,21,43,50,33,23,16)(4,15,39,
    99,69,5,83,93,75,30,24)(6,89,47,72,\\98,77,22,53,40,34,66)(7,19,45,42,88,96,76,17,
    80,37,27)(8,14,48,49,79,54,71,85,62,29,25)(9,63,38,64,10,56,91,87,51,\\31,61)(11,
    67,35,65,13,73,46,82,81,32,70)(12,86,41,58,74,97,60,18,59,36,78),

  (1,28,95,52,68,84,90,57,44,92,20)(2,16,23,33,50,43,21,3,26,100,94)(4,24,30,75,93,
    83,5,69,99,39,15)(6,66,34,40,\\53,22,77,98,72,47,89)(7,27,37,80,17,76,96,88,42,45,
    19)(8,25,29,62,85,71,54,79,49,48,14)(9,61,31,51,87,91,56,10,64,\\38,63)(11,70,32,
    81,82,46,73,13,65,35,67)(12,78,36,59,18,60,97,74,58,41,86),

  (1,44,84,95,20,57,68,28,92,90,52)(2,26,43,23,94,3,50,16,100,21,33)(4,99,83,30,15,
    69,93,24,39,5,75)(6,72,22,34,\\89,98,53,66,47,77,40)(7,42,76,37,19,88,17,27,45,96,
    80)(8,49,71,29,14,79,85,25,48,54,62)(9,64,91,31,63,10,87,61,38,\\56,51)(11,65,46,
    32,67,13,82,70,35,73,81)(12,58,60,36,86,74,18,78,41,97,59),

  (1,52,90,92,28,68,57,20,95,84,44)(2,33,21,100,16,50,3,94,23,43,26)(4,75,5,39,24,
    93,69,15,30,83,99)(6,40,77,47,\\66,53,98,89,34,22,72)(7,80,96,45,27,17,88,19,37,
    76,42)(8,62,54,48,25,85,79,14,29,71,49)(9,51,56,38,61,87,10,63,31,\\91,64)(11,81,
    73,35,70,82,13,67,32,46,65)(12,59,97,41,78,18,74,86,36,60,58),

  (1,57,52,20,90,95,92,84,28,44,68)(2,3,33,94,21,23,100,43,16,26,50)(4,69,75,15,5,
    30,39,83,24,99,93)(6,98,40,89,77,\\34,47,22,66,72,53)(7,88,80,19,96,37,45,76,27,
    42,17)(8,79,62,14,54,29,48,71,25,49,85)(9,10,51,63,56,31,38,91,61,64,\\87)(11,13,
    81,67,73,32,35,46,70,65,82)(12,74,59,86,97,36,41,60,78,58,18),

  (1,68,44,28,84,92,95,90,20,52,57)(2,50,26,16,43,100,23,21,94,33,3)(4,93,99,24,83,
    39,30,5,15,75,69)(6,53,72,66,22,\\47,34,77,89,40,98)(7,17,42,27,76,45,37,96,19,80,
    88)(8,85,49,25,71,48,29,54,14,62,79)(9,87,64,61,91,38,31,56,63,51,\\10)(11,82,65,
    70,46,35,32,73,67,81,13)(12,18,58,78,60,41,36,97,86,59,74),

  (1,84,20,68,92,52,44,95,57,28,90)(2,43,94,50,100,33,26,23,3,16,21)(4,83,15,93,39,
    75,99,30,69,24,5)(6,22,89,53,47,\\40,72,34,98,66,77)(7,76,19,17,45,80,42,37,88,27,
    96)(8,71,14,85,48,62,49,29,79,25,54)(9,91,63,87,38,51,64,31,10,61,\\56)(11,46,67,
    82,35,81,65,32,13,70,73)(12,60,86,18,41,59,58,36,74,78,97),

  (1,90,28,57,95,44,52,92,68,20,84)(2,21,16,3,23,26,33,100,50,94,43)(4,5,24,69,30,
    99,75,39,93,15,83)(6,77,66,98,34,\\72,40,47,53,89,22)(7,96,27,88,37,42,80,45,17,
    19,76)(8,54,25,79,29,49,62,48,85,14,71)(9,56,61,10,31,64,51,38,87,63,\\91)(11,73,
    70,13,32,65,81,35,82,67,46)(12,97,78,74,36,58,59,41,18,86,60),

  (1,92,57,84,52,28,20,44,90,68,95)(2,100,3,43,33,16,94,26,21,50,23)(4,39,69,83,75,
    24,15,99,5,93,30)(6,47,98,22,40,\\66,89,72,77,53,34)(7,45,88,76,80,27,19,42,96,17,
    37)(8,48,79,71,62,25,14,49,54,85,29)(9,38,10,91,51,61,63,64,56,87,31)\\(11,35,13,
    46,81,70,67,65,73,82,32)(12,41,74,60,59,78,86,58,97,18,36),

  (1,95,68,90,44,20,28,52,84,57,92)(2,23,50,21,26,94,16,33,43,3,100)(4,30,93,5,99,
    15,24,75,83,69,39)(6,34,53,77,72,\\89,66,40,22,98,47)(7,37,17,96,42,19,27,80,76,
    88,45)(8,29,85,54,49,14,25,62,71,79,48)(9,31,87,56,64,63,61,51,91,10,38)\\(11,32,
    82,73,65,67,70,81,46,13,35)(12,36,18,97,58,86,78,59,60,74,41).

The character table of $G^{s_{23}}$:\\
\begin{tabular}{c|ccccccccccc}
  & & & & & & & & & &10 &\\\hline
$\chi_{23}^{(1)}$&1&1&1&1&1&1&1&1&1&1&1
\\$\chi_{23}^{(2)}$&1&A&/A&C&/C&D&/D&/E&E&B&/B
\\$\chi_{23}^{(3)}$&1&B&/B&/E&E&/C&C&A&/A&D&/D
\\$\chi_{23}^{(4)}$&1&C&/C&/B&B&A&/A&/D&D&/E&E
\\$\chi_{23}^{(5)}$&1&D&/D&A&/A&E&/E&B&/B&/C&C
\\$\chi_{23}^{(6)}$&1&E&/E&D&/D&/B&B&/C&C&/A&A
\\$\chi_{23}^{(7)}$&1&/E&E&/D&D&B&/B&C&/C&A&/A
\\$\chi_{23}^{(8)}$&1&/D&D&/A&A&/E&E&/B&B&C&/C
\\$\chi_{23}^{(9)}$&1&/C&C&B&/B&/A&A&D&/D&E&/E
\\$\chi_{23}^{(10)}$&1&/B&B&E&/E&C&/C&/A&A&/D&D
\\$\chi_{23}^{(11)}$&1&/A&A&/C&C&/D&D&E&/E&/B&B
\end{tabular}

\noindent \noindent where   A = E(11)$^4$; B = E(11)$^8$; C = E(11);
D = E(11)$^5$; E = E(11)$^9$.

The generators of $G^{s_{24}}$ are:\\
 (  1, 84, 20, 68, 92, 52, 44, 95, 57, 28, 90)(  2, 43, 94, 50,100, 33, 26, 23,  3,
     16, 21)(  4, 83, 15, 93, 39, 75, 99, 30, 69, 24,  5)(  6, 22, 89, 53, 47, 40,
      72, 34, 98, 66, 77)(  7, 76, 19, 17, 45, 80, 42, 37, 88, 27, 96)
    (  8, 71, 14, 85, 48, 62, 49, 29, 79, 25, 54)(  9, 91, 63, 87, 38, 51, 64, 31,
      10, 61, 56)( 11, 46, 67, 82, 35, 81, 65, 32, 13, 70, 73)( 12, 60, 86, 18, 41,
      59, 58, 36, 74, 78, 97)

The representatives of conjugacy classes of   $G^{s_{24}}$ are:\\
 (1),
  (1,20,92,44,57,90,84,68,52,95,28)(2,94,100,26,3,21,43,50,33,23,16)(4,15,39,
    99,69,5,83,93,75,30,24)(6,89,47,72,\\98,77,22,53,40,34,66)(7,19,45,42,88,96,76,17,
    80,37,27)(8,14,48,49,79,54,71,85,62,29,25)(9,63,38,64,10,56,91,87,51,\\31,61)(11,
    67,35,65,13,73,46,82,81,32,70)(12,86,41,58,74,97,60,18,59,36,78),

  (1,28,95,52,68,84,90,57,44,92,20)(2,16,23,33,50,43,21,3,26,100,94)(4,24,30,75,93,
    83,5,69,99,39,15)(6,66,34,40,\\53,22,77,98,72,47,89)(7,27,37,80,17,76,96,88,42,45,
    19)(8,25,29,62,85,71,54,79,49,48,14)(9,61,31,51,87,91,56,10,64,\\38,63)(11,70,32,
    81,82,46,73,13,65,35,67)(12,78,36,59,18,60,97,74,58,41,86),

  (1,44,84,95,20,57,68,28,92,90,52)(2,26,43,23,94,3,50,16,100,21,33)(4,99,83,30,15,
    69,93,24,39,5,75)(6,72,22,34,\\89,98,53,66,47,77,40)(7,42,76,37,19,88,17,27,45,96,
    80)(8,49,71,29,14,79,85,25,48,54,62)(9,64,91,31,63,10,87,61,38,\\56,51)(11,65,46,
    32,67,13,82,70,35,73,81)(12,58,60,36,86,74,18,78,41,97,59),

  (1,52,90,92,28,68,57,20,95,84,44)(2,33,21,100,16,50,3,94,23,43,26)(4,75,5,39,24,
    93,69,15,30,83,99)(6,40,77,47,\\66,53,98,89,34,22,72)(7,80,96,45,27,17,88,19,37,
    76,42)(8,62,54,48,25,85,79,14,29,71,49)(9,51,56,38,61,87,10,63,31,\\91,64)(11,81,
    73,35,70,82,13,67,32,46,65)(12,59,97,41,78,18,74,86,36,60,58),

  (1,57,52,20,90,95,92,84,28,44,68)(2,3,33,94,21,23,100,43,16,26,50)(4,69,75,15,5,
    30,39,83,24,99,93)(6,98,40,89,\\77,34,47,22,66,72,53)(7,88,80,19,96,37,45,76,27,
    42,17)(8,79,62,14,54,29,48,71,25,49,85)(9,10,51,63,56,31,38,91,61,\\64,87)(11,13,
    81,67,73,32,35,46,70,65,82)(12,74,59,86,97,36,41,60,78,58,18),

  (1,68,44,28,84,92,95,90,20,52,57)(2,50,26,16,43,100,23,21,94,33,3)(4,93,99,24,83,
    39,30,5,15,75,69)(6,53,72,66,\\22,47,34,77,89,40,98)(7,17,42,27,76,45,37,96,19,80,
    88)(8,85,49,25,71,48,29,54,14,62,79)(9,87,64,61,91,38,31,56,63,\\51,10)(11,82,65,
    70,46,35,32,73,67,81,13)(12,18,58,78,60,41,36,97,86,59,74),

  (1,84,20,68,92,52,44,95,57,28,90)(2,43,94,50,100,33,26,23,3,16,21)(4,83,15,93,39,
    75,99,30,69,24,5)(6,22,89,53,\\47,40,72,34,98,66,77)(7,76,19,17,45,80,42,37,88,27,
    96)(8,71,14,85,48,62,49,29,79,25,54)(9,91,63,87,38,51,64,31,10,\\61,56)(11,46,67,
    82,35,81,65,32,13,70,73)(12,60,86,18,41,59,58,36,74,78,97),

  (1,90,28,57,95,44,52,92,68,20,84)(2,21,16,3,23,26,33,100,50,94,43)(4,5,24,69,30,
    99,75,39,93,15,83)(6,77,66,98,\\34,72,40,47,53,89,22)(7,96,27,88,37,42,80,45,17,
    19,76)(8,54,25,79,29,49,62,48,85,14,71)(9,56,61,10,31,64,51,38,87,\\63,91)(11,73,
    70,13,32,65,81,35,82,67,46)(12,97,78,74,36,58,59,41,18,86,60),

  (1,92,57,84,52,28,20,44,90,68,95)(2,100,3,43,33,16,94,26,21,50,23)(4,39,69,83,75,
    24,15,99,5,93,30)(6,47,98,22,\\40,66,89,72,77,53,34)(7,45,88,76,80,27,19,42,96,17,
    37)(8,48,79,71,62,25,14,49,54,85,29)(9,38,10,91,51,61,63,64,56,\\87,31)(11,35,13,
    46,81,70,67,65,73,82,32)(12,41,74,60,59,78,86,58,97,18,36),

  (1,95,68,90,44,20,28,52,84,57,92)(2,23,50,21,26,94,16,33,43,3,100)(4,30,93,5,99,
    15,24,75,83,69,39)(6,34,53,77,\\72,89,66,40,22,98,47)(7,37,17,96,42,19,27,80,76,
    88,45)(8,29,85,54,49,14,25,62,71,79,48)(9,31,87,56,64,63,61,51,91,\\10,38)(11,32,
    82,73,65,67,70,81,46,13,35)(12,36,18,97,58,86,78,59,60,74,41).

The character table of $G^{s_{24}}$:\\
\begin{tabular}{c|ccccccccccc}
  & & & & & & & & & &10 &\\\hline
$\chi_{24}^{(1)}$&1&1&1&1&1&1&1&1&1&1&1
\\$\chi_{24}^{(2)}$&1&A&/A&C&/C&D&/D&/E&E&B&/B
\\$\chi_{24}^{(3)}$&1&B&/B&/E&E&/C&C&A&/A&D&/D
\\$\chi_{24}^{(4)}$&1&C&/C&/B&B&A&/A&/D&D&/E&E
\\$\chi_{24}^{(5)}$&1&D&/D&A&/A&E&/E&B&/B&/C&C
\\$\chi_{24}^{(6)}$&1&E&/E&D&/D&/B&B&/C&C&/A&A
\\$\chi_{24}^{(7)}$&1&/E&E&/D&D&B&/B&C&/C&A&/A
\\$\chi_{24}^{(8)}$&1&/D&D&/A&A&/E&E&/B&B&C&/C
\\$\chi_{24}^{(9)}$&1&/C&C&B&/B&/A&A&D&/D&E&/E
\\$\chi_{24}^{(10)}$&1&/B&B&E&/E&C&/C&/A&A&/D&D
\\$\chi_{24}^{(11)}$&1&/A&A&/C&C&/D&D&E&/E&/B&B
\end{tabular}
$s_9$=$\left(\begin{array}{cccc} Z(2^3)^5& 0*Z(2)& Z(2^3)^4& Z(2^3)
\\ Z(2^3)^4& Z(2^3)& Z(2^3)^3& 0*Z(2)\\Z(2^3)^5& Z(2^3)^4& Z(2^3)^3&
Z(2)^0\\Z(2^3)^2& Z(2)^0& Z(2^3)^6& Z(2^3)^4
\end{array}\right)$,
$s_{10}$= $\left(\begin{array}{cccc} Z(2^3)^2& 0*Z(2)& Z(2^3)^2&
Z(2^3)^2 \\ Z(2)^0& Z(2^3)^6& Z(2)^0& Z(2^3)^2\\Z(2^3)^5& Z(2^3)^5&
Z(2^3)^6& 0*Z(2) \\Z(2^3)^2& Z(2^3)^5& Z(2)^0& Z(2^3)^2
\end{array}\right)$,

$s_11$=$\left(\begin{array}{cccc} Z(2^3)^4& Z(2^3)& Z(2^3)^4&
Z(2^3)^6\\Z(2^3)^5& Z(2^3)^3& Z(2^3)^4& 0*Z(2)\\Z(2^3)^6& Z(2^3)^5&
Z(2^3)^5& Z(2)^0 \\Z(2^3)^6& Z(2)^0& Z(2^3)^6& Z(2^3)^3
\end{array}\right)$.

The character table of $G^{s_1}$:\\
\begin{tabular}{c|ccccccccccc}
  & & & & & & & & & & 10&\\\hline

$\chi_1^{(1)}$&1&1&1&1&1&1&1&1&1&1&1
\\$\chi_1^{(2)}$&14&.&.&.&1&1&1&G&-G&-2&-1
\\$\chi_1^{(3)}$&14&.&.&.&1&1&1&-G&G&-2&-1
\\$\chi_1^{(4)}$&35&.&.&.&D&E&F&-1&-1&3&.
\\$\chi_1^{(5)}$&35&.&.&.&E&F&D&-1&-1&3&.
\\$\chi_1^{(6)}$&35&.&.&.&F&D&E&-1&-1&3&.
\\$\chi_1^{(7)}$&64&1&1&1&-1&-1&-1&.&.&.&-1
\\$\chi_1^{(8)}$&65&A&C&B&.&.&.&1&1&1&.
\\$\chi_1^{(9)}$&65&B&A&C&.&.&.&1&1&1&.
\\$\chi_1^{(10)}$&65&C&B&A&.&.&.&1&1&1&.
\\$\chi_1^{(11)}$&91&.&.&.&.&.&.&-1&-1&-5&1
  \end{tabular}

   \noindent where  A = E(7)$^2$+E(7)$^5$; B = E(7)+E(7)$^6$;
C = E(7)$^3$+E(7)$^4$; D = -E(13)$^4$-E(13)$^6$-E(13)$^7$-E(13)$^9$;
E = -E(13)-E(13)$^5$-E(13)$^8$-E(13)$^{12}$; F =
-E(13)$^2$-E(13)$^3$-E(13)$^{10}$-E(13)$^{11}$; G = -2*E(4)
  = -2*ER(-1) = -2i.

 The generators,  the representatives of conjugacy classes and character table  of $G^{s_{1}}$ are the same as
 $G$ since $G= G^{s_1}=G$.

The generators of $G^{s_2}$ are:\\
 $\left(\begin{array}{cccc} Z(2^3)^4& Z(2^3)^6& Z(2^3)^4& Z(2^3)^5 \\ 0*Z(2)& Z(2^3)^5& Z(2^3)^6& Z(2^3
             )^6 \\ Z(2^3)^4& Z(2^3)^6& Z(2^3)& Z(2)^0 \\Z(2^3)^4& Z(2^3)& Z(2^3)^3& Z(2)^0 \end{array}\right)$

The representatives of conjugacy classes of   $G^{s_2}$ are:\\
 $\left(\begin{array}{cccc} Z(2)^0& 0*Z(2)& 0*Z(2)& 0*Z(2) \\ 0*Z(2)& Z(2)^0& 0*Z(2)& 0*Z(2)\\0*Z(2)& 0*Z(2)& Z(2)^0& 0*Z(2) \\ 0*Z(2)& 0*Z(2)& 0*Z(2)& Z(2)^0 \end{array}\right)$,
 $\left(\begin{array}{cccc} Z(2)^0& Z(2)^0& Z(2^3)^6& Z(2^3)^5 \\ Z(2^3)^3& Z(2^3)& Z(2^3)^6& Z(2^3)^4 \\Z(2^3)& Z(2^3)^6& Z(2^3)^5& Z(2^3)^6\\Z(2^3)^4& Z(2^3)^4& 0*Z(2)& Z(2^3)^4 \end{array}\right)$,

  $\left(\begin{array}{cccc} Z(2^3)& Z(2^3)^2& Z(2^3)^5& Z(2^3)^5 \\ Z(2^3)^3& Z(2^3)^6& 0*Z(2)& Z(2^3)\\Z(2^3)^4& Z(2)^0& Z(2^3)^4& Z(2^3)^3\\Z(2^3)^6& Z(2)^0& Z(2^3)^4& Z(2^3)^5 \end{array}\right)$,
$\left(\begin{array}{cccc} Z(2^3)^2& Z(2^3)^3& 0*Z(2)& Z(2^3) \\
0*Z(2)& Z(2^3)^6& Z(2^3)^3& Z(2^3)^4\\Z(2^3)^3& Z(2^3)^4& Z(2^3)&
0*Z(2) \\Z(2^3)^6& 0*Z(2)& Z(2^3)^4& Z(2^3)^5 \end{array}\right)$,

  $\left(\begin{array}{cccc} Z(2^3)^4& Z(2^3)^6& Z(2^3)^4& Z(2^3)^5 \\0*Z(2)& Z(2^3)^5& Z(2^3)^6& Z(2^3)^6\\Z(2^3)^4& Z(2^3)^6& Z(2^3)& Z(2)^0\\Z(2^3)^4& Z(2^3)& Z(2^3)^3& Z(2)^0 \end{array}\right)$,
$\left(\begin{array}{cccc} Z(2^3)^5& 0*Z(2)& Z(2^3)^4& Z(2^3) \\
Z(2^3)^4& Z(2^3)& Z(2^3)^3& 0*Z(2)\\0*Z(2)& Z(2^3)^4& Z(2^3)^6&
Z(2^3)^3 \\Z(2^3)^6& Z(2^3)^3& 0*Z(2)& Z(2^3)^2 \end{array}\right)$,

  $\left(\begin{array}{cccc} Z(2^3)^5& Z(2^3)^3& Z(2^3)& Z(2^3)^5\\Z(2^3)^4& Z(2^3)^4& 0*Z(2)& Z(2^3)^5\\Z(2)^0& Z(2)^0& Z(2^3)^6& Z(2^3)^2\\Z(2^3)^6& Z(2^3)^4& Z(2^3)^3& Z(2^3) \end{array}\right)$.

The character table of $G^{s_2}$:\\
\begin{tabular}{c|ccccccc}
  & & & & & & &\\\hline
$\chi_2^{(1)}$&1&1&1&1&1&1&1
\\$\chi_2^{(2)}$&1&A&/C&/B&/A&B&C
\\$\chi_2^{(3)}$&1&B&A&C&/B&/C&/A
\\$\chi_2^{(4)}$&1&C&/B&A&/C&/A&B
\\$\chi_2^{(5)}$&1&/C&B&/A&C&A&/B
\\$\chi_2^{(6)}$&1&/B&/A&/C&B&C&A
\\$\chi_2^{(7)}$&1&/A&C&B&A&/B&/C
\end{tabular}

   \noindent where  A = E(7)$^6$; B = E(7)$^5$; C = E(7)$^4$.

The generators of $G^{s_3}$ are:\\
 $\left(\begin{array}{cccc} Z(2^3)^2& Z(2^3)^3& 0*Z(2)& Z(2^3) \\ 0*Z(2)& Z(2^3)^6& Z(2^3)^3& Z(2^3)^4 \\Z(2^3)^3& Z(2^3)^4& Z(2^3)& 0*Z(2)\\Z(2^3)^6& 0*Z(2)& Z(2^3)^4& Z(2^3)^5 \end{array}\right)$

The representatives of conjugacy classes of   $G^{s_3}$ are:\\
 $\left(\begin{array}{cccc} Z(2)^0& 0*Z(2)& 0*Z(2)& 0*Z(2) \\ 0*Z(2)& Z(2)^0& 0*Z(2)& 0*Z(2) \\0*Z(2)& 0*Z(2)& Z(2)^0& 0*Z(2) \\ 0*Z(2)& 0*Z(2)& 0*Z(2)& Z(2)^0 \end{array}\right)$,
$\left(\begin{array}{cccc} Z(2)^0& Z(2)^0& Z(2^3)^6& Z(2^3)^5 \\
Z(2^3)^3& Z(2^3)& Z(2^3)^6& Z(2^3)^4 \\Z(2^3)& Z(2^3)^6& Z(2^3)^5&
Z(2^3)^6 \\Z(2^3)^4& Z(2^3)^4& 0*Z(2)& Z(2^3)^4 \end{array}\right)$,

  $\left(\begin{array}{cccc} Z(2^3)& Z(2^3)^2& Z(2^3)^5& Z(2^3)^5 \\ Z(2^3)^3& Z(2^3)^6& 0*Z(2)& Z(2^3) \\Z(2^3)^4& Z(2)^0& Z(2^3)^4& Z(2^3)^3 \\Z(2^3)^6& Z(2)^0& Z(2^3)^4& Z(2^3)^5 \end{array}\right)$,
$\left(\begin{array}{cccc} Z(2^3)^2& Z(2^3)^3& 0*Z(2)& Z(2^3) \\
0*Z(2)& Z(2^3)^6& Z(2^3)^3& Z(2^3)^4\\Z(2^3)^3& Z(2^3)^4& Z(2^3)&
0*Z(2) \\Z(2^3)^6& 0*Z(2)& Z(2^3)^4& Z(2^3)^5 \end{array}\right)$,

  $\left(\begin{array}{cccc} Z(2^3)^4& Z(2^3)^6& Z(2^3)^4& Z(2^3)^5\\0*Z(2)& Z(2^3)^5& Z(2^3)^6& Z(2^3)^6\\Z(2^3)^4& Z(2^3)^6& Z(2^3)& Z(2)^0 \\Z(2^3)^4& Z(2^3)& Z(2^3)^3& Z(2)^0 \end{array}\right)$,
 $\left(\begin{array}{cccc} Z(2^3)^5& 0*Z(2)& Z(2^3)^4& Z(2^3) \\ Z(2^3)^4& Z(2^3)& Z(2^3)^3& 0*Z(2)\\0*Z(2)& Z(2^3)^4& Z(2^3)^6& Z(2^3)^3 \\Z(2^3)^6& Z(2^3)^3& 0*Z(2)& Z(2^3)^2 \end{array}\right)$,

  $\left(\begin{array}{cccc} Z(2^3)^5& Z(2^3)^3& Z(2^3)& Z(2^3)^5 \\Z(2^3)^4& Z(2^3)^4& 0*Z(2)& Z(2^3)^5\\Z(2)^0& Z(2)^0& Z(2^3)^6& Z(2^3)^2 \\Z(2^3)^6& Z(2^3)^4& Z(2^3)^3& Z(2^3) \end{array}\right)$.

The character table of $G^{s_3}$:\\
\begin{tabular}{c|ccccccc}
  & & & & & & &\\\hline
$\chi_3^{(1)}$&1&1&1&1&1&1&1
\\$\chi_3^{(2)}$&1&A&/C&/B&/A&B&C
\\$\chi_3^{(3)}$&1&B&A&C&/B&/C&/A
\\$\chi_3^{(4)}$&1&C&/B&A&/C&/A&B
\\$\chi_3^{(5)}$&1&/C&B&/A&C&A&/B
\\$\chi_3^{(6)}$&1&/B&/A&/C&B&C&A
\\$\chi_3^{(7)}$&1&/A&C&B&A&/B&/C

  \end{tabular}

   \noindent where  A = E(7)$^3$; B = E(7)$^6$; C = E(7)$^2$.

The generators of $G^{s_4}$ are:\\
$\left(\begin{array}{cccc} Z(2^3)& Z(2^3)^2& Z(2^3)^5& Z(2^3)^5 \\
Z(2^3)^3& Z(2^3)^6& 0*Z(2)& Z(2^3) \\Z(2^3)^4& Z(2)^0& Z(2^3)^4&
Z(2^3)^3 \\Z(2^3)^6& Z(2)^0& Z(2^3)^4& Z(2^3)^5 \end{array}\right)$

The representatives of conjugacy classes of   $G^{s_4}$ are:\\
 $\left(\begin{array}{cccc} Z(2)^0& 0*Z(2)& 0*Z(2)& 0*Z(2) \\ 0*Z(2)& Z(2)^0& 0*Z(2)& 0*Z(2)\\0*Z(2)& 0*Z(2)& Z(2)^0& 0*Z(2) \\ 0*Z(2)& 0*Z(2)& 0*Z(2)& Z(2)^0 \end{array}\right)$,
$\left(\begin{array}{cccc} Z(2)^0& Z(2)^0& Z(2^3)^6& Z(2^3)^5 \\
Z(2^3)^3& Z(2^3)& Z(2^3)^6& Z(2^3)^4\\Z(2^3)& Z(2^3)^6& Z(2^3)^5&
Z(2^3)^6\\Z(2^3)^4& Z(2^3)^4& 0*Z(2)& Z(2^3)^4 \end{array}\right)$,

  $\left(\begin{array}{cccc} Z(2^3)& Z(2^3)^2& Z(2^3)^5& Z(2^3)^5 \\ Z(2^3)^3& Z(2^3)^6& 0*Z(2)& Z(2^3)\\Z(2^3)^4& Z(2)^0& Z(2^3)^4& Z(2^3)^3\\Z(2^3)^6& Z(2)^0& Z(2^3)^4& Z(2^3)^5 \end{array}\right)$,
$\left(\begin{array}{cccc} Z(2^3)^2& Z(2^3)^3& 0*Z(2)& Z(2^3) \\
0*Z(2)& Z(2^3)^6& Z(2^3)^3& Z(2^3)^4 \\Z(2^3)^3& Z(2^3)^4& Z(2^3)&
0*Z(2) \\Z(2^3)^6& 0*Z(2)& Z(2^3)^4& Z(2^3)^5 \end{array}\right)$,

  $\left(\begin{array}{cccc} Z(2^3)^4& Z(2^3)^6& Z(2^3)^4& Z(2^3)^5 \\0*Z(2)& Z(2^3)^5& Z(2^3)^6& Z(2^3)^6 \\Z(2^3)^4& Z(2^3)^6& Z(2^3)& Z(2)^0 \\Z(2^3)^4& Z(2^3)& Z(2^3)^3& Z(2)^0 \end{array}\right)$,
$\left(\begin{array}{cccc} Z(2^3)^5& 0*Z(2)& Z(2^3)^4& Z(2^3) \\
Z(2^3)^4& Z(2^3)& Z(2^3)^3& 0*Z(2) \\0*Z(2)& Z(2^3)^4& Z(2^3)^6&
Z(2^3)^3 \\Z(2^3)^6& Z(2^3)^3& 0*Z(2)& Z(2^3)^2 \end{array}\right)$,

  $\left(\begin{array}{cccc} Z(2^3)^5& Z(2^3)^3& Z(2^3)& Z(2^3)^5 \\Z(2^3)^4& Z(2^3)^4& 0*Z(2)& Z(2^3)^5 \\Z(2)^0& Z(2)^0& Z(2^3)^6& Z(2^3)^2 \\Z(2^3)^6& Z(2^3)^4& Z(2^3)^3& Z(2^3) \end{array}\right)$.

The character table of $G^{s_4}$:\\
\begin{tabular}{c|ccccccc}
  & & & & & & &\\\hline

$\chi_4^{(1)}$&1&1&1&1&1&1&1
\\$\chi_4^{(2)}$&1&A&/C&/B&/A&B&C
\\$\chi_4^{(3)}$&1&B&A&C&/B&/C&/A
\\$\chi_4^{(4)}$&1&C&/B&A&/C&/A&B
\\$\chi_4^{(5)}$&1&/C&B&/A&C&A&/B
\\$\chi_4^{(6)}$&1&/B&/A&/C&B&C&A
\\$\chi_4^{(7)}$&1&/A&C&B&A&/B&/C
  \end{tabular}

   \noindent where  A = E(7)$^2$; B = E(7)$^4$; C = E(7)$^6$.

The generators of $G^{s_5}$ are:\\
 $\left(\begin{array}{cccc} Z(2^3)& Z(2^3)& Z(2^3)^4& Z(2)^0 \\ Z(2^3)^4& Z(2^3)^6& Z(2^3)& Z(2^3)^6 \\Z(2^3)^3& Z(2)^0& Z(2)^0& Z(2^3)^4\\Z(2^3)^4& Z(2^3)^5& Z(2)^0& Z(2^3)^6 \end{array}\right)$

The representatives of conjugacy classes of   $G^{s_5}$ are:\\
$\left(\begin{array}{cccc} Z(2)^0& 0*Z(2)& 0*Z(2)& 0*Z(2) \\ 0*Z(2)&
Z(2)^0& 0*Z(2)& 0*Z(2) \\0*Z(2)& 0*Z(2)& Z(2)^0& 0*Z(2) \\ 0*Z(2)&
0*Z(2)& 0*Z(2)& Z(2)^0 \end{array}\right)$,
$\left(\begin{array}{cccc} Z(2)^0& 0*Z(2)& Z(2)^0& Z(2)^0 \\
Z(2^3)^5& Z(2)^0& Z(2^3)^5& Z(2^3)^4\\Z(2^3)^3& Z(2^3)^6& Z(2^3)&
Z(2^3)^4 \\Z(2)^0& Z(2^3)^6& Z(2^3)^4& Z(2^3)^6 \end{array}\right)$,

  $\left(\begin{array}{cccc} Z(2)^0& Z(2^3)^3& Z(2)^0& Z(2^3)^2 \\Z(2^3)^3& Z(2^3)^2& Z(2^3)^2& Z(2^3)^2\\Z(2^3)^4& Z(2^3)^4& Z(2^3)^2& Z(2^3)^2 \\Z(2^3)^3& Z(2)^0& Z(2)^0& Z(2^3)^2 \end{array}\right)$,
$\left(\begin{array}{cccc} Z(2^3)& Z(2^3)& Z(2^3)^4& Z(2)^0 \\
Z(2^3)^4& Z(2^3)^6& Z(2^3)& Z(2^3)^6\\Z(2^3)^3& Z(2)^0& Z(2)^0&
Z(2^3)^4 \\Z(2^3)^4& Z(2^3)^5& Z(2)^0& Z(2^3)^6 \end{array}\right)$,

  $\left(\begin{array}{cccc} Z(2^3)^2& Z(2^3)^2& Z(2^3)^2& Z(2^3)^2 \\Z(2)^0& Z(2^3)^2& Z(2^3)^2& Z(2)^0 \\Z(2)^0& Z(2^3)^4& Z(2^3)^2& Z(2^3)^3 \\Z(2^3)^3& Z(2^3)^4& Z(2^3)^3& Z(2)^0 \end{array}\right)$,
$\left(\begin{array}{cccc} Z(2^3)^2& Z(2^3)^2& Z(2^3)^4& Z(2)^0
\\Z(2^3)^4& Z(2^3)^3& 0*Z(2)& Z(2^3)^3 \\Z(2^3)^4& Z(2^3)^5&
Z(2^3)^4& Z(2^3)^3\\Z(2^3)^5& Z(2)^0& Z(2^3)^2& Z(2^3)^2
\end{array}\right)$,

  $\left(\begin{array}{cccc} Z(2^3)^2& Z(2^3)^3& Z(2^3)^3& Z(2)^0 \\Z(2^3)^2& Z(2^3)^4& 0*Z(2)& Z(2^3)^4\\Z(2)^0& Z(2^3)^5& Z(2^3)^3& Z(2^3)^2 \\Z(2^3)^5& Z(2^3)^4& Z(2^3)^4& Z(2^3)^2 \end{array}\right)$,
$\left(\begin{array}{cccc} Z(2^3)^3& 0*Z(2)& Z(2^3)^6&
Z(2^3)\\Z(2^3)^2& Z(2^3)^2& Z(2^3)^5& Z(2^3)^4 \\Z(2^3)^6& Z(2^3)^5&
Z(2^3)^3& Z(2^3)^2 \\Z(2^3)^3& Z(2)^0& Z(2^3)^3& Z(2^3)^5
\end{array}\right)$,

  $\left(\begin{array}{cccc} Z(2^3)^4& Z(2^3)^2& Z(2^3)^3& Z(2^3)^4\\Z(2^3)^4& Z(2^3)^3& Z(2^3)^6& Z(2^3)^2 \\Z(2)^0& 0*Z(2)& Z(2^3)^5& Z(2^3)\\Z(2^3)^3& Z(2^3)^5& Z(2^3)^5& Z(2^3)^6 \end{array}\right)$,
$\left(\begin{array}{cccc} Z(2^3)^5& Z(2^3)^2& Z(2^3)^4&
Z(2^3)\\Z(2^3)^3& Z(2^3)^3& Z(2^3)^5& Z(2^3)^6\\Z(2)^0& Z(2^3)^5&
Z(2^3)^2& 0*Z(2) \\Z(2^3)^3& Z(2^3)^6& Z(2^3)^2& Z(2^3)^3
\end{array}\right)$,

  $\left(\begin{array}{cccc} Z(2^3)^6& Z(2^3)& Z(2^3)^2& Z(2^3)^4 \\Z(2^3)^5& Z(2^3)^5& Z(2^3)^6& Z(2^3)^3\\Z(2^3)^5& 0*Z(2)& Z(2^3)^3& Z(2^3)^2 \\Z(2^3)^3& Z(2)^0& Z(2^3)^4& Z(2^3)^4 \end{array}\right)$,
$\left(\begin{array}{cccc} Z(2^3)^6& Z(2^3)^4& Z(2^3)^4& Z(2)^0 \\
Z(2^3)^4& Z(2^3)& Z(2^3)^5& Z(2)^0\\Z(2^3)^6& Z(2^3)^6& Z(2)^0&
0*Z(2)\\Z(2)^0& Z(2^3)^3& Z(2^3)^5& Z(2)^0 \end{array}\right)$,

  $\left(\begin{array}{cccc} Z(2^3)^6& Z(2^3)^4& Z(2^3)^6& Z(2)^0 \\ Z(2)^0& Z(2)^0& Z(2^3)& Z(2^3)^4 \\Z(2^3)^5& Z(2)^0& Z(2^3)^6& Z(2^3) \\Z(2^3)^4& Z(2^3)^3& Z(2^3)^4& Z(2^3) \end{array}\right)$.

The character table of $G^{s_5}$:\\
\begin{tabular}{c|ccccccccccccc}
  & & & & & & & & & & 10& & &\\\hline
$\chi_5^{(1)}$&1&1&1&1&1&1&1&1&1&1&1&1&1
\\$\chi_5^{(2)}$&1&A&F&/C&/F&D&/D&B&/E&/B&E&/A&C
\\$\chi_5^{(3)}$&1&B&/A&/F&A&/E&E&D&C&/D&/C&/B&F
\\$\chi_5^{(4)}$&1&C&E&D&/E&/A&A&F&/B&/F&B&/C&/D
\\$\chi_5^{(5)}$&1&D&/B&A&B&C&/C&/E&F&E&/F&/D&/A
\\$\chi_5^{(6)}$&1&E&D&/B&/D&/F&F&/C&A&C&/A&/E&B
\\$\chi_5^{(7)}$&1&F&/C&/E&C&/B&B&/A&/D&A&D&/F&E
\\$\chi_5^{(8)}$&1&/F&C&E&/C&B&/B&A&D&/A&/D&F&/E
\\$\chi_5^{(9)}$&1&/E&/D&B&D&F&/F&C&/A&/C&A&E&/B
\\$\chi_5^{(10)}$&1&/D&B&/A&/B&/C&C&E&/F&/E&F&D&A
\\$\chi_5^{(11)}$&1&/C&/E&/D&E&A&/A&/F&B&F&/B&C&D
\\$\chi_5^{(12)}$&1&/B&A&F&/A&E&/E&/D&/C&D&C&B&/F
\\$\chi_5^{(13)}$&1&/A&/F&C&F&/D&D&/B&E&B&/E&A&/C

  \end{tabular}

   \noindent where  A = E(13)$^4$; B = E(13)$^8$; C =
E(13)$^{12}$; D = E(13)$^3$; E = E(13)$^7$; F = E(13)$^{11}$.

The generators of $G^{s_6}$ are:\\
 $\left(\begin{array}{cccc} Z(2^3)^2& Z(2^3)^2& Z(2^3)^2& Z(2^3)^2 \\ Z(2)^0& Z(2^3)^2& Z(2^3)^2& Z(2)^
            0 \\ Z(2)^0& Z(2^3)^4& Z(2^3)^2& Z(2^3)^3 \\Z(2^3)^3& Z(2^3)^4& Z(2^3)^3& Z(2)^0 \end{array}\right)$

The representatives of conjugacy classes of   $G^{s_6}$ are:\\
 $\left(\begin{array}{cccc} Z(2)^0& 0*Z(2)& 0*Z(2)& 0*Z(2) \\ 0*Z(2)& Z(2)^0& 0*Z(2)& 0*Z(2) \\ 0*Z(2)& 0*Z(2)& Z(2)^0& 0*Z(2) \\ 0*Z(2)& 0*Z(2)& 0*Z(2)& Z(2)^0 \end{array}\right)$,
$\left(\begin{array}{cccc} Z(2)^0& 0*Z(2)& Z(2)^0& Z(2)^0 \\
Z(2^3)^5& Z(2)^0& Z(2^3)^5& Z(2^3)^4 \\Z(2^3)^3& Z(2^3)^6& Z(2^3)&
Z(2^3)^4 \\Z(2)^0& Z(2^3)^6& Z(2^3)^4& Z(2^3)^6 \end{array}\right)$,

  $\left(\begin{array}{cccc} Z(2)^0& Z(2^3)^3& Z(2)^0& Z(2^3)^2 \\Z(2^3)^3& Z(2^3)^2& Z(2^3)^2& Z(2^3)^2 \\Z(2^3)^4& Z(2^3)^4& Z(2^3)^2& Z(2^3)^2 \\Z(2^3)^3& Z(2)^0& Z(2)^0& Z(2^3)^2 \end{array}\right)$,
$\left(\begin{array}{cccc} Z(2^3)& Z(2^3)& Z(2^3)^4& Z(2)^0 \\
Z(2^3)^4& Z(2^3)^6& Z(2^3)& Z(2^3)^6\\Z(2^3)^3& Z(2)^0& Z(2)^0&
Z(2^3)^4 \\Z(2^3)^4& Z(2^3)^5& Z(2)^0& Z(2^3)^6 \end{array}\right)$,

  $\left(\begin{array}{cccc} Z(2^3)^2& Z(2^3)^2& Z(2^3)^2& Z(2^3)^2 \\Z(2)^0& Z(2^3)^2& Z(2^3)^2& Z(2)^0 \\Z(2)^0& Z(2^3)^4& Z(2^3)^2& Z(2^3)^3\\Z(2^3)^3& Z(2^3)^4& Z(2^3)^3& Z(2)^0 \end{array}\right)$,
$\left(\begin{array}{cccc} Z(2^3)^2& Z(2^3)^2& Z(2^3)^4& Z(2)^0
\\Z(2^3)^4& Z(2^3)^3& 0*Z(2)& Z(2^3)^3 \\Z(2^3)^4& Z(2^3)^5&
Z(2^3)^4& Z(2^3)^3 \\Z(2^3)^5& Z(2)^0& Z(2^3)^2& Z(2^3)^2
\end{array}\right)$,

  $\left(\begin{array}{cccc} Z(2^3)^2& Z(2^3)^3& Z(2^3)^3& Z(2)^0 \\Z(2^3)^2& Z(2^3)^4& 0*Z(2)& Z(2^3)^4 \\Z(2)^0& Z(2^3)^5& Z(2^3)^3& Z(2^3)^2 \\Z(2^3)^5& Z(2^3)^4& Z(2^3)^4& Z(2^3)^2 \end{array}\right)$,
$\left(\begin{array}{cccc} Z(2^3)^3& 0*Z(2)& Z(2^3)^6& Z(2^3)
\\Z(2^3)^2& Z(2^3)^2& Z(2^3)^5& Z(2^3)^4 \\Z(2^3)^6& Z(2^3)^5&
Z(2^3)^3& Z(2^3)^2 \\Z(2^3)^3& Z(2)^0& Z(2^3)^3& Z(2^3)^5
\end{array}\right)$,

  $\left(\begin{array}{cccc} Z(2^3)^4& Z(2^3)^2& Z(2^3)^3& Z(2^3)^4 \\Z(2^3)^4& Z(2^3)^3& Z(2^3)^6& Z(2^3)^2 \\Z(2)^0& 0*Z(2)& Z(2^3)^5& Z(2^3) \\Z(2^3)^3& Z(2^3)^5& Z(2^3)^5& Z(2^3)^6 \end{array}\right)$,
$\left(\begin{array}{cccc} Z(2^3)^5& Z(2^3)^2& Z(2^3)^4& Z(2^3)
\\Z(2^3)^3& Z(2^3)^3& Z(2^3)^5& Z(2^3)^6\\Z(2)^0& Z(2^3)^5&
Z(2^3)^2& 0*Z(2)\\Z(2^3)^3& Z(2^3)^6& Z(2^3)^2& Z(2^3)^3
\end{array}\right)$,

  $\left(\begin{array}{cccc} Z(2^3)^6& Z(2^3)& Z(2^3)^2& Z(2^3)^4 \\Z(2^3)^5& Z(2^3)^5& Z(2^3)^6& Z(2^3)^3\\Z(2^3)^5& 0*Z(2)& Z(2^3)^3& Z(2^3)^2 \\Z(2^3)^3& Z(2)^0& Z(2^3)^4& Z(2^3)^4 \end{array}\right)$,
  $\left(\begin{array}{cccc} Z(2^3)^6& Z(2^3)^4& Z(2^3)^4& Z(2)^0 \\
Z(2^3)^4& Z(2^3)& Z(2^3)^5& Z(2)^0\\Z(2^3)^6& Z(2^3)^6& Z(2)^0&
0*Z(2)\\Z(2)^0& Z(2^3)^3& Z(2^3)^5& Z(2)^0 \end{array}\right)$,

  $\left(\begin{array}{cccc} Z(2^3)^6& Z(2^3)^4& Z(2^3)^6& Z(2)^0 \\ Z(2)^0& Z(2)^0& Z(2^3)& Z(2^3)^4 \\Z(2^3)^5& Z(2)^0& Z(2^3)^6& Z(2^3) \\Z(2^3)^4& Z(2^3)^3& Z(2^3)^4& Z(2^3) \end{array}\right)$.

The character table of $G^{s_6}$:\\
 \begin{tabular}{c|ccccccccccccc}
  & & & & & & & & & & 10& & &\\\hline
$\chi_6^{(1)}$&1&1&1&1&1&1&1&1&1&1&1&1&1
\\$\chi_6^{(2)}$&1&A&F&/C&/F&D&/D&B&/E&/B&E&/A&C
\\$\chi_6^{(3)}$&1&B&/A&/F&A&/E&E&D&C&/D&/C&/B&F
\\$\chi_6^{(4)}$&1&C&E&D&/E&/A&A&F&/B&/F&B&/C&/D
\\$\chi_6^{(5)}$&1&D&/B&A&B&C&/C&/E&F&E&/F&/D&/A
\\$\chi_6^{(6)}$&1&E&D&/B&/D&/F&F&/C&A&C&/A&/E&B
\\$\chi_6^{(7)}$&1&F&/C&/E&C&/B&B&/A&/D&A&D&/F&E
\\$\chi_6^{(8)}$&1&/F&C&E&/C&B&/B&A&D&/A&/D&F&/E
\\$\chi_6^{(9)}$&1&/E&/D&B&D&F&/F&C&/A&/C&A&E&/B
\\$\chi_6^{(10)}$&1&/D&B&/A&/B&/C&C&E&/F&/E&F&D&A
\\$\chi_6^{(11)}$&1&/C&/E&/D&E&A&/A&/F&B&F&/B&C&D
\\$\chi_6^{(12)}$&1&/B&A&F&/A&E&/E&/D&/C&D&C&B&/F
\\$\chi_6^{(13)}$&1&/A&/F&C&F&/D&D&/B&E&B&/E&A&/C

  \end{tabular}

   \noindent where  A = E(13)$^2$; B = E(13)$^4$; C =
E(13)$^6$; D = E(13)$^8$; E = E(13)$^{10}$; F = E(13)$^{12}$.

The generators of $G^{s_7}$ are:\\
 $\left(\begin{array}{cccc} Z(2)^0& 0*Z(2)& Z(2)^0& Z(2)^0 \\ Z(2^3)^5& Z(2)^0& Z(2^3)^5& Z(2^3)^4 \\Z(2^3)^3& Z(2^3)^6& Z(2^3)& Z(2^3)^4\\Z(2)^0& Z(2^3)^6& Z(2^3)^4& Z(2^3)^6 \end{array}\right)$

The representatives of conjugacy classes of   $G^{s_7}$ are:\\
$\left(\begin{array}{cccc} Z(2)^0& 0*Z(2)& 0*Z(2)& 0*Z(2) \\ 0*Z(2)&
Z(2)^0& 0*Z(2)& 0*Z(2) \\0*Z(2)& 0*Z(2)& Z(2)^0& 0*Z(2) \\ 0*Z(2)&
0*Z(2)& 0*Z(2)& Z(2)^0 \end{array}\right)$,
$\left(\begin{array}{cccc} Z(2)^0& 0*Z(2)& Z(2)^0& Z(2)^0 \\
Z(2^3)^5& Z(2)^0& Z(2^3)^5& Z(2^3)^4\\Z(2^3)^3& Z(2^3)^6& Z(2^3)&
Z(2^3)^4 \\Z(2)^0& Z(2^3)^6& Z(2^3)^4& Z(2^3)^6 \end{array}\right)$,

  $\left(\begin{array}{cccc} Z(2)^0& Z(2^3)^3& Z(2)^0& Z(2^3)^2\\Z(2^3)^3& Z(2^3)^2& Z(2^3)^2& Z(2^3)^2 \\Z(2^3)^4& Z(2^3)^4& Z(2^3)^2& Z(2^3)^2 \\Z(2^3)^3& Z(2)^0& Z(2)^0& Z(2^3)^2 \end{array}\right)$,
$\left(\begin{array}{cccc} Z(2^3)& Z(2^3)& Z(2^3)^4& Z(2)^0 \\
Z(2^3)^4& Z(2^3)^6& Z(2^3)& Z(2^3)^6 \\Z(2^3)^3& Z(2)^0& Z(2)^0&
Z(2^3)^4 \\Z(2^3)^4& Z(2^3)^5& Z(2)^0& Z(2^3)^6 \end{array}\right)$,

  $\left(\begin{array}{cccc} Z(2^3)^2& Z(2^3)^2& Z(2^3)^2& Z(2^3)^2 \\Z(2)^0& Z(2^3)^2& Z(2^3)^2& Z(2)^0 \\Z(2)^0& Z(2^3)^4& Z(2^3)^2& Z(2^3)^3 \\Z(2^3)^3& Z(2^3)^4& Z(2^3)^3& Z(2)^0 \end{array}\right)$,
$\left(\begin{array}{cccc} Z(2^3)^2& Z(2^3)^2& Z(2^3)^4& Z(2)^0
\\Z(2^3)^4& Z(2^3)^3& 0*Z(2)& Z(2^3)^3\\Z(2^3)^4& Z(2^3)^5&
Z(2^3)^4& Z(2^3)^3 \\Z(2^3)^5& Z(2)^0& Z(2^3)^2& Z(2^3)^2
\end{array}\right)$,

  $\left(\begin{array}{cccc} Z(2^3)^2& Z(2^3)^3& Z(2^3)^3& Z(2)^0 \\Z(2^3)^2& Z(2^3)^4& 0*Z(2)& Z(2^3)^4 \\Z(2)^0& Z(2^3)^5& Z(2^3)^3& Z(2^3)^2\\Z(2^3)^5& Z(2^3)^4& Z(2^3)^4& Z(2^3)^2 \end{array}\right)$,
$\left(\begin{array}{cccc} Z(2^3)^3& 0*Z(2)& Z(2^3)^6& Z(2^3)
\\Z(2^3)^2& Z(2^3)^2& Z(2^3)^5& Z(2^3)^4 \\Z(2^3)^6& Z(2^3)^5&
Z(2^3)^3& Z(2^3)^2 \\Z(2^3)^3& Z(2)^0& Z(2^3)^3& Z(2^3)^5
\end{array}\right)$,

  $\left(\begin{array}{cccc} Z(2^3)^4& Z(2^3)^2& Z(2^3)^3& Z(2^3)^4 \\Z(2^3)^4& Z(2^3)^3& Z(2^3)^6& Z(2^3)^2\\Z(2)^0& 0*Z(2)& Z(2^3)^5& Z(2^3)\\Z(2^3)^3& Z(2^3)^5& Z(2^3)^5& Z(2^3)^6 \end{array}\right)$,
$\left(\begin{array}{cccc} Z(2^3)^5& Z(2^3)^2& Z(2^3)^4& Z(2^3)
\\Z(2^3)^3& Z(2^3)^3& Z(2^3)^5& Z(2^3)^6\\Z(2)^0& Z(2^3)^5&
Z(2^3)^2& 0*Z(2) \\Z(2^3)^3& Z(2^3)^6& Z(2^3)^2& Z(2^3)^3
\end{array}\right)$,

  $\left(\begin{array}{cccc} Z(2^3)^6& Z(2^3)& Z(2^3)^2& Z(2^3)^4\\Z(2^3)^5& Z(2^3)^5& Z(2^3)^6& Z(2^3)^3 \\Z(2^3)^5& 0*Z(2)& Z(2^3)^3& Z(2^3)^2 \\Z(2^3)^3& Z(2)^0& Z(2^3)^4& Z(2^3)^4 \end{array}\right)$,
$\left(\begin{array}{cccc} Z(2^3)^6& Z(2^3)^4& Z(2^3)^4& Z(2)^0 \\
Z(2^3)^4& Z(2^3)& Z(2^3)^5& Z(2)^0 \\Z(2^3)^6& Z(2^3)^6& Z(2)^0&
0*Z(2)\\Z(2)^0& Z(2^3)^3& Z(2^3)^5& Z(2)^0 \end{array}\right)$,

  $\left(\begin{array}{cccc} Z(2^3)^6& Z(2^3)^4& Z(2^3)^6& Z(2)^0 \\ Z(2)^0& Z(2)^0& Z(2^3)& Z(2^3)^4 \\Z(2^3)^5& Z(2)^0& Z(2^3)^6& Z(2^3) \\Z(2^3)^4& Z(2^3)^3& Z(2^3)^4& Z(2^3) \end{array}\right)$.

The character table of $G^{s_7}$:\\
 \begin{tabular}{c|ccccccccccccc}
  & & & & & & & & & & 10& & &\\\hline

$\chi_7^{(1)}$&1&1&1&1&1&1&1&1&1&1&1&1&1
\\$\chi_7^{(2)}$&1&A&F&/C&/F&D&/D&B&/E&/B&E&/A&C
\\$\chi_7^{(3)}$&1&B&/A&/F&A&/E&E&D&C&/D&/C&/B&F
\\$\chi_7^{(4)}$&1&C&E&D&/E&/A&A&F&/B&/F&B&/C&/D
\\$\chi_7^{(5)}$&1&D&/B&A&B&C&/C&/E&F&E&/F&/D&/A
\\$\chi_7^{(6)}$&1&E&D&/B&/D&/F&F&/C&A&C&/A&/E&B
\\$\chi_7^{(7)}$&1&F&/C&/E&C&/B&B&/A&/D&A&D&/F&E
\\$\chi_7^{(8)}$&1&/F&C&E&/C&B&/B&A&D&/A&/D&F&/E
\\$\chi_7^{(9)}$&1&/E&/D&B&D&F&/F&C&/A&/C&A&E&/B
\\$\chi_7^{(10)}$&1&/D&B&/A&/B&/C&C&E&/F&/E&F&D&A
\\$\chi_7^{(11)}$&1&/C&/E&/D&E&A&/A&/F&B&F&/B&C&D
\\$\chi_7^{(12)}$&1&/B&A&F&/A&E&/E&/D&/C&D&C&B&/F
\\$\chi_7^{(13)}$&1&/A&/F&C&F&/D&D&/B&E&B&/E&A&/C
  \end{tabular}

   \noindent where  A = E(13); B = E(13)$^2$; C = E(13)$^3$;
D = E(13)$^4$; E = E(13)$^5$; F = E(13)$^6$.

The generators of $G^{s_8}$ are:\\
 $\left(\begin{array}{cccc} Z(2^3)^4& Z(2)^0& 0*Z(2)& Z(2^3) \\ Z(2^3)^6& Z(2^3)^3& Z(2^3)^3& Z(2^3)^4\\Z(2)^0& Z(2^3)^4& Z(2^3)& 0*Z(2) \\Z(2^3)^2& Z(2^3)^5& Z(2^3)^4& Z(2^3)^5 \end{array}\right)$
$\left(\begin{array}{cccc} Z(2^3)^4& Z(2^3)& Z(2^3)^2&
Z(2^3)^5\\Z(2^3)^4& Z(2^3)^6& Z(2^3)^6& Z(2^3) \\Z(2^3)^2& Z(2^3)&
Z(2^3)^2& Z(2^3)^3 \\Z(2)^0& Z(2^3)& Z(2^3)^6& Z(2^3)^5
\end{array}\right)$

  $\left(\begin{array}{cccc} Z(2^3)^3& Z(2^3)^2& Z(2^3)^2& Z(2)^0 \\ 0*Z(2)& Z(2^3)^2& Z(2^3)^5& Z(2^3)^2 \\Z(2^3)^4& Z(2^3)^3& Z(2^3)^2& Z(2^3)^2\\Z(2)^0& Z(2^3)^4& 0*Z(2)& Z(2^3)^3 \end{array}\right)$

The representatives of conjugacy classes of   $G^{s_8}$ are:\\
$\left(\begin{array}{cccc} Z(2)^0& 0*Z(2)& 0*Z(2)& 0*Z(2) \\ 0*Z(2)&
Z(2)^0& 0*Z(2)& 0*Z(2)\\0*Z(2)& 0*Z(2)& Z(2)^0& 0*Z(2) \\ 0*Z(2)&
0*Z(2)& 0*Z(2)& Z(2)^0 \end{array}\right)$,
$\left(\begin{array}{cccc} Z(2)^0& Z(2^3)& Z(2)^0& Z(2^3)^4 \\
Z(2^3)& Z(2^3)^6& Z(2^3)^2& Z(2)^0\\Z(2^3)^4& Z(2)^0& Z(2^3)^6&
Z(2^3)\\Z(2^3)^4& Z(2^3)^4& Z(2^3)& Z(2)^0 \end{array}\right)$,

  $\left(\begin{array}{cccc} Z(2^3)& Z(2^3)^2& 0*Z(2)& Z(2)^0 \\ Z(2^3)^3& Z(2^3)^6& Z(2)^0& Z(2^3)^4 \\Z(2^3)^5& Z(2^3)^3& Z(2^3)^2& Z(2^3)^6\\Z(2^3)^5& Z(2)^0& 0*Z(2)& Z(2^3)^3 \end{array}\right)$,
$\left(\begin{array}{cccc} Z(2^3)& Z(2^3)^5& Z(2^3)^2& Z(2^3)^2 \\
0*Z(2)& Z(2^3)^3& Z(2^3)^5& Z(2^3) \\Z(2^3)& Z(2^3)^5& Z(2^3)&
Z(2^3)^4 \\Z(2^3)& Z(2^3)^2& Z(2^3)^3& Z(2^3)^3 \end{array}\right)$,

  $\left(\begin{array}{cccc} Z(2^3)^2& 0*Z(2)& Z(2^3)^2& Z(2^3)^2 \\ Z(2)^0& Z(2^3)^6& Z(2)^0& Z(2^3)^2 \\Z(2^3)^5& Z(2^3)^5& Z(2^3)^6& 0*Z(2) \\Z(2^3)^2& Z(2^3)^5& Z(2)^0& Z(2^3)^2 \end{array}\right)$,
$\left(\begin{array}{cccc} Z(2^3)^2& Z(2^3)& Z(2^3)^6& Z(2^3) \\
Z(2^3)^3& Z(2)^0& Z(2^3)^6& Z(2^3)^6\\Z(2)^0& Z(2^3)^4& Z(2)^0&
Z(2^3)\\Z(2^3)& Z(2)^0& Z(2^3)^3& Z(2^3)^2 \end{array}\right)$,

  $\left(\begin{array}{cccc} Z(2^3)^3& Z(2^3)^2& Z(2^3)^2& Z(2)^0 \\0*Z(2)& Z(2^3)^2& Z(2^3)^5& Z(2^3)^2 \\Z(2^3)^4& Z(2^3)^3& Z(2^3)^2& Z(2^3)^2\\Z(2)^0& Z(2^3)^4& 0*Z(2)& Z(2^3)^3 \end{array}\right)$,
$\left(\begin{array}{cccc} Z(2^3)^3& Z(2^3)^4& Z(2^3)& Z(2^3)^2
\\Z(2^3)^3& Z(2^3)& Z(2^3)^5& Z(2^3)^2 \\Z(2^3)^2& Z(2^3)^5&
Z(2^3)^3& Z(2^3)^5\\Z(2^3)& Z(2^3)& 0*Z(2)& Z(2^3)
\end{array}\right)$,

  $\left(\begin{array}{cccc} Z(2^3)^3& Z(2^3)^4& Z(2^3)^6& Z(2^3)^5\\Z(2^3)& 0*Z(2)& Z(2^3)^3& Z(2^3)^6 \\ 0*Z(2)& Z(2^3)& 0*Z(2)& Z(2^3)^4 \\Z(2^3)^5& 0*Z(2)& Z(2^3)& Z(2^3)^3 \end{array}\right)$,
$\left(\begin{array}{cccc} Z(2^3)^3& Z(2^3)^6& Z(2^3)^4& Z(2)^0 \\
0*Z(2)& Z(2^3)^2& Z(2)^0& 0*Z(2) \\Z(2)^0& Z(2^3)^3& Z(2^3)^6&
Z(2^3)^2 \\Z(2^3)^5& Z(2^3)^5& Z(2^3)^3& Z(2^3) \end{array}\right)$,

  $\left(\begin{array}{cccc} Z(2^3)^4& Z(2)^0& 0*Z(2)& Z(2^3) \\ Z(2^3)^6& Z(2^3)^3& Z(2^3)^3& Z(2^3)^4 \\Z(2)^0& Z(2^3)^4& Z(2^3)& 0*Z(2)\\Z(2^3)^2& Z(2^3)^5& Z(2^3)^4& Z(2^3)^5 \end{array}\right)$,
$\left(\begin{array}{cccc} Z(2^3)^4& Z(2^3)& Z(2^3)^2& Z(2^3)^5
\\Z(2^3)^4& Z(2^3)^6& Z(2^3)^6& Z(2^3) \\Z(2^3)^2& Z(2^3)& Z(2^3)^2&
Z(2^3)^3 \\Z(2)^0& Z(2^3)& Z(2^3)^6& Z(2^3)^5 \end{array}\right)$,

  $\left(\begin{array}{cccc} Z(2^3)^4& Z(2^3)^2& 0*Z(2)& Z(2^3)^6 \\ Z(2)^0& 0*Z(2)& Z(2^3)^4& 0*Z(2)\\Z(2)^0& Z(2^3)^2& 0*Z(2)& Z(2^3)^2\\Z(2^3)^6& Z(2)^0& Z(2)^0& Z(2^3)^4 \end{array}\right)$,
$\left(\begin{array}{cccc} Z(2^3)^4& Z(2^3)^4& Z(2)^0& Z(2^3)^3 \\
Z(2^3)^3& Z(2^3)^2& Z(2^3)& Z(2)^0\\Z(2^3)^5& Z(2^3)^6& Z(2^3)^2&
Z(2^3)^4\\Z(2^3)^3& Z(2^3)^5& Z(2^3)^3& Z(2^3)^4
\end{array}\right)$,

  $\left(\begin{array}{cccc} Z(2^3)^5& 0*Z(2)& Z(2^3)^4& Z(2^3) \\ Z(2^3)^4& Z(2^3)& Z(2^3)^3& 0*Z(2) \\Z(2^3)^5& Z(2^3)^4& Z(2^3)^3& Z(2)^0 \\Z(2^3)^2& Z(2)^0& Z(2^3)^6& Z(2^3)^4 \end{array}\right)$,
$\left(\begin{array}{cccc} Z(2^3)^5& Z(2^3)^3& Z(2^3)&
Z(2^3)^5\\Z(2^3)^6& Z(2^3)^2& Z(2^3)^6& Z(2^3)^2\\Z(2^3)& Z(2^3)&
Z(2^3)^6& Z(2^3)\\Z(2)^0& Z(2^3)^2& Z(2^3)^4& Z(2^3)^4
\end{array}\right)$.

The character table of $G^{s_8}$:\\
\begin{tabular}{c|cccccccccccccccc}
  & & & & & & & & & & 10& & & & & &\\\hline

$\chi_8^{(1)}$&1&1&1&1&1&1&1&1&1&1&1&1&1&1&1&1
\\$\chi_8^{(2)}$&1&-1&1&1&1&-1&-1&1&1&1&-1&-1&-1&1&-1&-1
\\$\chi_8^{(3)}$&1&1&1&-1&1&1&-1&-1&-1&1&1&-1&-1&-1&1&-1
\\$\chi_8^{(4)}$&1&-1&-1&1&1&-1&1&1&-1&-1&1&-1&1&-1&1&-1
\\$\chi_8^{(5)}$&1&1&-1&-1&1&1&1&-1&1&-1&-1&-1&1&1&-1&-1
\\$\chi_8^{(6)}$&1&-1&-1&-1&1&-1&-1&-1&1&-1&1&1&-1&1&1&1
\\$\chi_8^{(7)}$&1&1&-1&1&1&1&-1&1&-1&-1&-1&1&-1&-1&-1&1
\\$\chi_8^{(8)}$&1&-1&1&-1&1&-1&1&-1&-1&1&-1&1&1&-1&-1&1
\\$\chi_8^{(9)}$&1&1&A&-A&-1&-1&-1&A&-1&-A&-A&A&1&1&A&-A
\\$\chi_8^{(10)}$&1&1&-A&A&-1&-1&-1&-A&-1&A&A&-A&1&1&-A&A
\\$\chi_8^{(11)}$&1&-1&A&A&-1&1&-1&-A&1&-A&A&A&1&-1&-A&-A
\\$\chi_8^{(12)}$&1&-1&-A&-A&-1&1&-1&A&1&A&-A&-A&1&-1&A&A
\\$\chi_8^{(13)}$&1&1&-A&-A&-1&-1&1&A&1&A&A&A&-1&-1&-A&-A
\\$\chi_8^{(14)}$&1&1&A&A&-1&-1&1&-A&1&-A&-A&-A&-1&-1&A&A
\\$\chi_8^{(15)}$&1&-1&-A&A&-1&1&1&-A&-1&A&-A&A&-1&1&A&-A
\\$\chi_8^{(16)}$&1&-1&A&-A&-1&1&1&A&-1&-A&A&-A&-1&1&-A&A
  \end{tabular}

   \noindent where  A = -E(4)
  = -ER(-1) = -i.

The generators of $G^{s_9}$ are:\\
$\left(\begin{array}{cccc} Z(2^3)^5& 0*Z(2)& Z(2^3)^4& Z(2^3) \\
Z(2^3)^4& Z(2^3)& Z(2^3)^3& 0*Z(2) \\Z(2^3)^5& Z(2^3)^4& Z(2^3)^3&
Z(2)^0 \\Z(2^3)^2& Z(2)^0& Z(2^3)^6& Z(2^3)^4 \end{array}\right)$,
$\left(\begin{array}{cccc} Z(2^3)^4& Z(2^3)& Z(2^3)^2&
Z(2^3)^5\\Z(2^3)^4& Z(2^3)^6& Z(2^3)^6& Z(2^3) \\Z(2^3)^2& Z(2^3)&
Z(2^3)^2& Z(2^3)^3 \\Z(2)^0& Z(2^3)& Z(2^3)^6& Z(2^3)^5
  \end{array}\right)$,

  $\left(\begin{array}{cccc} Z(2^3)^3& Z(2^3)^2& Z(2^3)^2& Z(2)^0\\0*Z(2)& Z(2^3)^2& Z(2^3)^5& Z(2^3)^2 \\Z(2^3)^4& Z(2^3)^3& Z(2^3)^2& Z(2^3)^2 \\Z(2)^0& Z(2^3)^4& 0*Z(2)& Z(2^3)^3 \end{array}\right)$

The representatives of conjugacy classes of   $G^{s_9}$ are:\\
 $\left(\begin{array}{cccc} Z(2)^0& 0*Z(2)& 0*Z(2)& 0*Z(2) \\ 0*Z(2)& Z(2)^0& 0*Z(2)& 0*Z(2) \\0*Z(2)& 0*Z(2)& Z(2)^0& 0*Z(2) \\ 0*Z(2)& 0*Z(2)& 0*Z(2)& Z(2)^0 \end{array}\right)$,
$\left(\begin{array}{cccc} Z(2)^0& Z(2^3)& Z(2)^0& Z(2^3)^4 \\
Z(2^3)& Z(2^3)^6& Z(2^3)^2& Z(2)^0 \\Z(2^3)^4& Z(2)^0& Z(2^3)^6&
Z(2^3)\\Z(2^3)^4& Z(2^3)^4& Z(2^3)& Z(2)^0 \end{array}\right)$,

  $\left(\begin{array}{cccc} Z(2^3)& Z(2^3)^2& 0*Z(2)& Z(2)^0 \\ Z(2^3)^3& Z(2^3)^6& Z(2)^0& Z(2^3)^4 \\Z(2^3)^5& Z(2^3)^3& Z(2^3)^2& Z(2^3)^6 \\Z(2^3)^5& Z(2)^0& 0*Z(2)& Z(2^3)^3 \end{array}\right)$,
$\left(\begin{array}{cccc} Z(2^3)& Z(2^3)^5& Z(2^3)^2& Z(2^3)^2 \\
0*Z(2)& Z(2^3)^3& Z(2^3)^5& Z(2^3)\\Z(2^3)& Z(2^3)^5& Z(2^3)&
Z(2^3)^4\\Z(2^3)& Z(2^3)^2& Z(2^3)^3& Z(2^3)^3 \end{array}\right)$,

  $\left(\begin{array}{cccc} Z(2^3)^2& 0*Z(2)& Z(2^3)^2& Z(2^3)^2 \\ Z(2)^0& Z(2^3)^6& Z(2)^0& Z(2^3)^2 \\Z(2^3)^5& Z(2^3)^5& Z(2^3)^6& 0*Z(2) \\Z(2^3)^2& Z(2^3)^5& Z(2)^0& Z(2^3)^2 \end{array}\right)$,
$\left(\begin{array}{cccc} Z(2^3)^2& Z(2^3)& Z(2^3)^6& Z(2^3) \\
Z(2^3)^3& Z(2)^0& Z(2^3)^6& Z(2^3)^6 \\Z(2)^0& Z(2^3)^4& Z(2)^0&
Z(2^3)\\Z(2^3)& Z(2)^0& Z(2^3)^3& Z(2^3)^2 \end{array}\right)$,

  $\left(\begin{array}{cccc} Z(2^3)^3& Z(2^3)^2& Z(2^3)^2& Z(2)^0 \\0*Z(2)& Z(2^3)^2& Z(2^3)^5& Z(2^3)^2\\Z(2^3)^4& Z(2^3)^3& Z(2^3)^2& Z(2^3)^2 \\Z(2)^0& Z(2^3)^4& 0*Z(2)& Z(2^3)^3 \end{array}\right)$,
$\left(\begin{array}{cccc} Z(2^3)^3& Z(2^3)^4& Z(2^3)& Z(2^3)^2
\\Z(2^3)^3& Z(2^3)& Z(2^3)^5& Z(2^3)^2 \\Z(2^3)^2& Z(2^3)^5&
Z(2^3)^3& Z(2^3)^5 \\Z(2^3)& Z(2^3)& 0*Z(2)& Z(2^3)
\end{array}\right)$,

  $\left(\begin{array}{cccc} Z(2^3)^3& Z(2^3)^4& Z(2^3)^6& Z(2^3)^5 \\Z(2^3)& 0*Z(2)& Z(2^3)^3& Z(2^3)^6 \\ 0*Z(2)& Z(2^3)& 0*Z(2)& Z(2^3)^4 \\Z(2^3)^5& 0*Z(2)& Z(2^3)& Z(2^3)^3 \end{array}\right)$,
$\left(\begin{array}{cccc} Z(2^3)^3& Z(2^3)^6& Z(2^3)^4& Z(2)^0 \\
0*Z(2)& Z(2^3)^2& Z(2)^0& 0*Z(2)\\Z(2)^0& Z(2^3)^3& Z(2^3)^6&
Z(2^3)^2 \\Z(2^3)^5& Z(2^3)^5& Z(2^3)^3& Z(2^3) \end{array}\right)$,

  $\left(\begin{array}{cccc} Z(2^3)^4& Z(2)^0& 0*Z(2)& Z(2^3) \\ Z(2^3)^6& Z(2^3)^3& Z(2^3)^3& Z(2^3)^4\\Z(2)^0& Z(2^3)^4& Z(2^3)& 0*Z(2) \\Z(2^3)^2& Z(2^3)^5& Z(2^3)^4& Z(2^3)^5 \end{array}\right)$,
$\left(\begin{array}{cccc} Z(2^3)^4& Z(2^3)& Z(2^3)^2&
Z(2^3)^5\\Z(2^3)^4& Z(2^3)^6& Z(2^3)^6& Z(2^3)\\Z(2^3)^2& Z(2^3)&
Z(2^3)^2& Z(2^3)^3\\Z(2)^0& Z(2^3)& Z(2^3)^6& Z(2^3)^5
\end{array}\right)$,

  $\left(\begin{array}{cccc} Z(2^3)^4& Z(2^3)^2& 0*Z(2)& Z(2^3)^6 \\ Z(2)^0& 0*Z(2)& Z(2^3)^4& 0*Z(2)\\Z(2)^0& Z(2^3)^2& 0*Z(2)& Z(2^3)^2 \\Z(2^3)^6& Z(2)^0& Z(2)^0& Z(2^3)^4 \end{array}\right)$,
$\left(\begin{array}{cccc} Z(2^3)^4& Z(2^3)^4& Z(2)^0& Z(2^3)^3 \\
Z(2^3)^3& Z(2^3)^2& Z(2^3)& Z(2)^0 \\Z(2^3)^5& Z(2^3)^6& Z(2^3)^2&
Z(2^3)^4 \\Z(2^3)^3& Z(2^3)^5& Z(2^3)^3& Z(2^3)^4
\end{array}\right)$,

  $\left(\begin{array}{cccc} Z(2^3)^5& 0*Z(2)& Z(2^3)^4& Z(2^3) \\ Z(2^3)^4& Z(2^3)& Z(2^3)^3& 0*Z(2) \\Z(2^3)^5& Z(2^3)^4& Z(2^3)^3& Z(2)^0\\Z(2^3)^2& Z(2)^0& Z(2^3)^6& Z(2^3)^4 \end{array}\right)$,
$\left(\begin{array}{cccc} Z(2^3)^5& Z(2^3)^3& Z(2^3)&
Z(2^3)^5\\Z(2^3)^6& Z(2^3)^2& Z(2^3)^6& Z(2^3)^2 \\Z(2^3)& Z(2^3)&
Z(2^3)^6& Z(2^3) \\Z(2)^0& Z(2^3)^2& Z(2^3)^4& Z(2^3)^4
\end{array}\right)$.

The character table of $G^{s_9}$:\\
\begin{tabular}{c|cccccccccccccccc}
  & & & & & & & & & & 10& & & & & &\\\hline

$\chi_9^{(1)}$&1&1&1&1&1&1&1&1&1&1&1&1&1&1&1&1
\\$\chi_9^{(2)}$&1&-1&1&1&1&-1&-1&1&1&1&-1&-1&-1&1&-1&-1
\\$\chi_9^{(3)}$&1&1&1&-1&1&1&-1&-1&-1&1&1&-1&-1&-1&1&-1
\\$\chi_9^{(4)}$&1&-1&-1&1&1&-1&1&1&-1&-1&1&-1&1&-1&1&-1
\\$\chi_9^{(5)}$&1&1&-1&-1&1&1&1&-1&1&-1&-1&-1&1&1&-1&-1
\\$\chi_9^{(6)}$&1&-1&-1&-1&1&-1&-1&-1&1&-1&1&1&-1&1&1&1
\\$\chi_9^{(7)}$&1&1&-1&1&1&1&-1&1&-1&-1&-1&1&-1&-1&-1&1
\\$\chi_9^{(8)}$&1&-1&1&-1&1&-1&1&-1&-1&1&-1&1&1&-1&-1&1
\\$\chi_9^{(9)}$&1&1&A&-A&-1&-1&-1&A&-1&-A&-A&A&1&1&A&-A
\\$\chi_9^{(10)}$&1&1&-A&A&-1&-1&-1&-A&-1&A&A&-A&1&1&-A&A
\\$\chi_9^{(11)}$&1&-1&A&A&-1&1&-1&-A&1&-A&A&A&1&-1&-A&-A
\\$\chi_9^{(12)}$&1&-1&-A&-A&-1&1&-1&A&1&A&-A&-A&1&-1&A&A
\\$\chi_9^{(13)}$&1&1&-A&-A&-1&-1&1&A&1&A&A&A&-1&-1&-A&-A
\\$\chi_9^{(14)}$&1&1&A&A&-1&-1&1&-A&1&-A&-A&-A&-1&-1&A&A
\\$\chi_9^{(15)}$&1&-1&-A&A&-1&1&1&-A&-1&A&-A&A&-1&1&A&-A
\\$\chi_9^{(16)}$&1&-1&A&-A&-1&1&1&A&-1&-A&A&-A&-1&1&-A&A

  \end{tabular}

    \noindent where  A = -E(4)
  = -ER(-1) = -i.

The generators of $G^{s_{10}}$ are:\\
 $\left(\begin{array}{cccc} Z(2^3)^2& 0*Z(2)& Z(2^3)^2& Z(2^3)^2 \\ Z(2)^0& Z(2^3)^6& Z(2)^0& Z(2^3)^2 \\Z(2^3)^5& Z(2^3)^5& Z(2^3)^6& 0*Z(2) \\Z(2^3)^2& Z(2^3)^5& Z(2)^0& Z(2^3)^2
 \end{array}\right)$,
$\left(\begin{array}{cccc} Z(2^3)^4& Z(2^3)^4& Z(2^3)& Z(2)^0 \\
Z(2)^0& Z(2^3)^4& Z(2)^0& 0*Z(2) \\Z(2^3)& Z(2)^0& 0*Z(2)& 0*Z(2) \\
Z(2)^0& 0*Z(2)& 0*Z(2)& 0*Z(2)
  \end{array}\right)$,

  $\left(\begin{array}{cccc} Z(2^3)^3& Z(2^3)^5& Z(2^3)^6& Z(2^3) \\ Z(2^3)^2& Z(2^3)& Z(2^3)^2& Z(2)^0 \\Z(2^3)^2& Z(2)^0& Z(2^3)^6& 0*Z(2) \\Z(2)^0& 0*Z(2)& Z(2^3)^4& Z(2^3)^2
  \end{array}\right)$,
$\left(\begin{array}{cccc} Z(2^3)^4& Z(2^3)& Z(2^3)^2& Z(2^3)^5
\\Z(2^3)^4& Z(2^3)^6& Z(2^3)^6& Z(2^3) \\Z(2^3)^2& Z(2^3)& Z(2^3)^2&
Z(2^3)^3 \\Z(2)^0& Z(2^3)& Z(2^3)^6& Z(2^3)^5 \end{array}\right)$

The representatives of conjugacy classes of   $G^{s_{10}}$ are:\\
$\left(\begin{array}{cccc} Z(2)^0& 0*Z(2)& 0*Z(2)& 0*Z(2) \\ 0*Z(2)&
Z(2)^0& 0*Z(2)& 0*Z(2)\\0*Z(2)& 0*Z(2)& Z(2)^0& 0*Z(2) \\ 0*Z(2)&
0*Z(2)& 0*Z(2)& Z(2)^0 \end{array}\right)$,
$\left(\begin{array}{cccc} 0*Z(2)& 0*Z(2)& 0*Z(2)& Z(2)^0 \\ 0*Z(2)&
0*Z(2)& Z(2)^0& Z(2^3) \\0*Z(2)& Z(2)^0& Z(2^3)^4& Z(2^3)^4 \\
Z(2)^0& Z(2^3)& Z(2)^0& Z(2^3)^4 \end{array}\right)$,

  $\left(\begin{array}{cccc} Z(2)^0& 0*Z(2)& Z(2^3)^4& Z(2^3)^2 \\ Z(2^3)^4& Z(2)^0& Z(2^3)& Z(2^3)^3\\Z(2^3)^5& Z(2^3)^2& Z(2^3)^6& Z(2^3)^2\\Z(2^3)& Z(2^3)& Z(2)^0& Z(2^3)^6 \end{array}\right)$,
$\left(\begin{array}{cccc} Z(2)^0& Z(2^3)& Z(2)^0& Z(2^3)^4 \\
Z(2^3)& Z(2^3)^6& Z(2^3)^2& Z(2)^0 \\Z(2^3)^4& Z(2)^0& Z(2^3)^6&
Z(2^3) \\Z(2^3)^4& Z(2^3)^4& Z(2^3)& Z(2)^0 \end{array}\right)$,

  $\left(\begin{array}{cccc} Z(2)^0& Z(2^3)& Z(2^3)^6& Z(2^3)^5 \\Z(2^3)^4& Z(2^3)^4& Z(2^3)^6& Z(2^3)^4 \\Z(2^3)^6& Z(2^3)^6& Z(2^3)^3& Z(2^3)^5 \\Z(2^3)^3& Z(2^3)^4& Z(2^3)& Z(2^3)^2 \end{array}\right)$,
$\left(\begin{array}{cccc} Z(2)^0& Z(2^3)^2& Z(2^3)^2& Z(2^3)^3 \\
Z(2^3)& Z(2^3)& Z(2)^0& Z(2)^0\\0*Z(2)& Z(2^3)^4& Z(2^3)^4& 0*Z(2)
\\Z(2^3)^6& Z(2^3)^6& Z(2^3)^6& Z(2^3)^6 \end{array}\right)$,

  $\left(\begin{array}{cccc} Z(2)^0& Z(2^3)^4& 0*Z(2)& Z(2^3)^3 \\ Z(2^3)^2& Z(2^3)^2& Z(2^3)^2& Z(2^3) \\Z(2^3)^2& Z(2^3)^5& Z(2^3)& 0*Z(2) \\Z(2^3)^3& Z(2^3)^4& Z(2^3)^6& Z(2^3)^5 \end{array}\right)$,
$\left(\begin{array}{cccc} Z(2)^0& Z(2^3)^4& Z(2)^0& Z(2^3)^2 \\
Z(2^3)^6& Z(2^3)& Z(2)^0& Z(2^3)^4\\Z(2^3)^4& Z(2)^0& 0*Z(2)&
Z(2^3)^2 \\Z(2^3)^6& Z(2)^0& Z(2^3)^3& Z(2^3)^3 \end{array}\right)$,

  $\left(\begin{array}{cccc} Z(2^3)& Z(2^3)& 0*Z(2)& Z(2^3) \\ 0*Z(2)& Z(2^3)^3& Z(2^3)^3& Z(2^3)^3\\Z(2)^0& Z(2)^0& Z(2^3)^4& Z(2^3)^5\\Z(2^3)^4& Z(2^3)& Z(2^3)^2& Z(2^3)^5 \end{array}\right)$,
$\left(\begin{array}{cccc} Z(2^3)& Z(2^3)& Z(2)^0& Z(2^3)^6 \\
Z(2^3)^3& 0*Z(2)& Z(2^3)^5& Z(2^3)^6\\Z(2^3)^6& Z(2^3)^3& Z(2^3)^5&
Z(2^3)^6\\Z(2^3)^4& Z(2)^0& Z(2)^0& Z(2^3)^6 \end{array}\right)$,

  $\left(\begin{array}{cccc} Z(2^3)& Z(2^3)^2& 0*Z(2)& Z(2)^0 \\ Z(2^3)^3& Z(2^3)^6& Z(2)^0& Z(2^3)^4\\Z(2^3)^5& Z(2^3)^3& Z(2^3)^2& Z(2^3)^6 \\Z(2^3)^5& Z(2)^0& 0*Z(2)& Z(2^3)^3 \end{array}\right)$,
$\left(\begin{array}{cccc} Z(2^3)& Z(2^3)^5& Z(2^3)^2& Z(2^3)^2 \\
0*Z(2)& Z(2^3)^3& Z(2^3)^5& Z(2^3)\\Z(2^3)& Z(2^3)^5& Z(2^3)&
Z(2^3)^4 \\Z(2^3)& Z(2^3)^2& Z(2^3)^3& Z(2^3)^3 \end{array}\right)$,

  $\left(\begin{array}{cccc} Z(2^3)^2& 0*Z(2)& Z(2^3)^2& Z(2^3)^2 \\ Z(2)^0& Z(2^3)^6& Z(2)^0& Z(2^3)^2 \\Z(2^3)^5& Z(2^3)^5& Z(2^3)^6& 0*Z(2) \\Z(2^3)^2& Z(2^3)^5& Z(2)^0& Z(2^3)^2 \end{array}\right)$,
$\left(\begin{array}{cccc} Z(2^3)^2& Z(2^3)& Z(2^3)^2& Z(2^3)^6
\\Z(2^3)^4& Z(2^3)^4& Z(2^3)& Z(2^3)^6 \\Z(2^3)^2& Z(2^3)^3&
Z(2^3)^6& Z(2)^0 \\Z(2^3)^3& Z(2^3)^6& Z(2^3)^3& Z(2^3)^5
\end{array}\right)$,

  $\left(\begin{array}{cccc} Z(2^3)^2& Z(2^3)& Z(2^3)^6& Z(2^3) \\ Z(2^3)^3& Z(2)^0& Z(2^3)^6& Z(2^3)^6\\Z(2)^0& Z(2^3)^4& Z(2)^0& Z(2^3) \\Z(2^3)& Z(2)^0& Z(2^3)^3& Z(2^3)^2 \end{array}\right)$,
$\left(\begin{array}{cccc} Z(2^3)^2& Z(2^3)^6& Z(2^3)&
Z(2^3)^5\\Z(2)^0& Z(2^3)^3& Z(2^3)^5& Z(2^3)^3 \\Z(2^3)^4& Z(2^3)^6&
0*Z(2)& Z(2^3)^4 \\Z(2)^0& Z(2^3)^5& Z(2^3)^2& Z(2^3)^5
\end{array}\right)$,

  $\left(\begin{array}{cccc} Z(2^3)^3& Z(2^3)& Z(2^3)& Z(2^3)^4 \\ Z(2)^0& Z(2^3)^3& Z(2^3)^6& Z(2^3)^4\\0*Z(2)& Z(2^3)^4& Z(2^3)^2& Z(2^3)^6 \\Z(2^3)^6& Z(2^3)^2& Z(2^3)^3& Z(2^3)^2 \end{array}\right)$,
$\left(\begin{array}{cccc} Z(2^3)^3& Z(2^3)^2& Z(2^3)^2&
Z(2)^0\\0*Z(2)& Z(2^3)^2& Z(2^3)^5& Z(2^3)^2 \\Z(2^3)^4& Z(2^3)^3&
Z(2^3)^2& Z(2^3)^2 \\Z(2)^0& Z(2^3)^4& 0*Z(2)& Z(2^3)^3
\end{array}\right)$,

  $\left(\begin{array}{cccc} Z(2^3)^3& Z(2^3)^2& Z(2^3)^4& Z(2^3)^2 \\ Z(2^3)^3& 0*Z(2)& Z(2)^0& Z(2)^0\\Z(2)^0& Z(2)^0& Z(2^3)& Z(2^3)^4 \\Z(2^3)^6& Z(2^3)^4& Z(2^3)^6& Z(2)^0 \end{array}\right)$,
$\left(\begin{array}{cccc} Z(2^3)^3& Z(2^3)^4& Z(2^3)^6& Z(2^3)^5
\\Z(2^3)& 0*Z(2)& Z(2^3)^3& Z(2^3)^6 \\ 0*Z(2)& Z(2^3)& 0*Z(2)&
Z(2^3)^4 \\Z(2^3)^5& 0*Z(2)& Z(2^3)& Z(2^3)^3 \end{array}\right)$,

  $\left(\begin{array}{cccc} Z(2^3)^4& Z(2^3)^2& 0*Z(2)& Z(2^3)^6 \\ Z(2)^0& 0*Z(2)& Z(2^3)^4& 0*Z(2)\\Z(2)^0& Z(2^3)^2& 0*Z(2)& Z(2^3)^2 \\Z(2^3)^6& Z(2)^0& Z(2)^0& Z(2^3)^4 \end{array}\right)$,
$\left(\begin{array}{cccc} Z(2^3)^4& Z(2^3)^4& Z(2)^0& Z(2^3)^3 \\
Z(2^3)^3& Z(2^3)^2& Z(2^3)& Z(2)^0\\Z(2^3)^5& Z(2^3)^6& Z(2^3)^2&
Z(2^3)^4\\Z(2^3)^3& Z(2^3)^5& Z(2^3)^3& Z(2^3)^4
\end{array}\right)$.

The character table of $G^{s_{10}}$:\\
\begin{tabular}{c|cccccccccccccccccccccc}
  & & & & & & & & & & 10& & & & & & & & & & 20& & \\\hline

$\chi_{10}^{(1)}$&1&1&1&1&1&1&1&1&1&1&1&1&1&1&1&1&1&1&1&1&1&1
\\$\chi_{10}^{(2)}$&1&-1&1&1&1&1&-1&1&1&-1&-1&-1&1&-1&1&-1&-1&1&1&1&1&1
\\$\chi_{10}^{(3)}$&1&-1&1&1&-1&1&-1&-1&-1&-1&1&1&1&1&1&1&-1&1&-1&1&1&1
\\$\chi_{10}^{(4)}$&1&-1&-1&1&1&-1&-1&-1&1&1&-1&-1&1&1&1&1&1&1&-1&1&1&1
\\$\chi_{10}^{(5)}$&1&-1&-1&1&-1&-1&-1&1&-1&1&1&1&1&-1&1&-1&1&1&1&1&1&1
\\$\chi_{10}^{(6)}$&1&1&-1&1&-1&-1&1&1&-1&-1&-1&-1&1&1&1&1&-1&1&1&1&1&1
\\$\chi_{10}^{(7)}$&1&1&-1&1&1&-1&1&-1&1&-1&1&1&1&-1&1&-1&-1&1&-1&1&1&1
\\$\chi_{10}^{(8)}$&1&1&1&1&-1&1&1&-1&-1&1&-1&-1&1&-1&1&-1&1&1&-1&1&1&1
\\$\chi_{10}^{(9)}$&2&A&.&-2&.&.&-A&.&.&.&.&.&-2&.&2&.&.&-2&.&2&2&-2
\\$\chi_{10}^{(10)}$&2&-A&.&-2&.&.&A&.&.&.&.&.&-2&.&2&.&.&-2&.&2&2&-2
\\$\chi_{10}^{(11)}$&2&.&.&-2&A&.&.&.&-A&.&.&.&2&.&-2&.&.&2&.&-2&2&-2
\\$\chi_{10}^{(12)}$&2&.&.&-2&-A&.&.&.&A&.&.&.&2&.&-2&.&.&2&.&-2&2&-2
\\$\chi_{10}^{(13)}$&2&.&.&2&.&.&.&.&.&.&-A&A&-2&.&-2&.&.&-2&.&-2&2&2
\\$\chi_{10}^{(14)}$&2&.&.&2&.&.&.&.&.&.&A&-A&-2&.&-2&.&.&-2&.&-2&2&2
\\$\chi_{10}^{(15)}$&2&.&.&2&.&.&.&.&.&.&.&.&2&-A&2&A&.&-2&.&-2&-2&-2
\\$\chi_{10}^{(16)}$&2&.&.&2&.&.&.&.&.&.&.&.&2&A&2&-A&.&-2&.&-2&-2&-2
\\$\chi_{10}^{(17)}$&2&.&-A&2&.&A&.&.&.&.&.&.&-2&.&-2&.&.&2&.&2&-2&-2
\\$\chi_{10}^{(18)}$&2&.&A&2&.&-A&.&.&.&.&.&.&-2&.&-2&.&.&2&.&2&-2&-2
\\$\chi_{10}^{(19)}$&2&.&.&-2&.&.&.&-A&.&.&.&.&-2&.&2&.&.&2&A&-2&-2&2
\\$\chi_{10}^{(20)}$&2&.&.&-2&.&.&.&A&.&.&.&.&-2&.&2&.&.&2&-A&-2&-2&2
\\$\chi_{10}^{(21)}$&2&.&.&-2&.&.&.&.&.&A&.&.&2&.&-2&.&-A&-2&.&2&-2&2
\\$\chi_{10}^{(22)}$&2&.&.&-2&.&.&.&.&.&-A&.&.&2&.&-2&.&A&-2&.&2&-2&2
  \end{tabular}

   \noindent where   A = 2*E(4)
  = 2*ER(-1) = 2i.

The generators of $G^{s_{11}}$ are:\\
$\left(\begin{array}{cccc} Z(2^3)^4& Z(2^3)& Z(2^3)^4& Z(2^3)^6 \\
Z(2^3)^5& Z(2^3)^3& Z(2^3)^4& 0*Z(2
             ) \\ Z(2^3)^6& Z(2^3)^5& Z(2^3)^5& Z(2)^0\\Z(2^3)^6& Z(2)^0& Z(2^3)^6& Z(2^3)^3 \end{array}\right)$

The representatives of conjugacy classes of   $G^{s_{11}}$ are:\\
 $\left(\begin{array}{cccc} Z(2)^0& 0*Z(2)& 0*Z(2)& 0*Z(2) \\ 0*Z(2)& Z(2)^0& 0*Z(2)& 0*Z(2) \\0*Z(2)& 0*Z(2)& Z(2)^0& 0*Z(2) \\ 0*Z(2)& 0*Z(2)& 0*Z(2)& Z(2)^0 \end{array}\right)$,
$\left(\begin{array}{cccc} Z(2^3)^3& 0*Z(2)& Z(2^3)^4& Z(2)^0 \\
Z(2^3)^6& Z(2^3)^2& Z(2)^0& 0*Z(2) \\Z(2^3)^3& Z(2)^0& Z(2)^0&
Z(2^3)^3 \\Z(2^3)^5& Z(2^3)^5& Z(2^3)^5& Z(2^3)^5
\end{array}\right)$,

  $\left(\begin{array}{cccc} Z(2^3)^3& Z(2)^0& 0*Z(2)& Z(2^3)^6 \\Z(2^3)^6& Z(2^3)^5& Z(2^3)^4& Z(2^3)^4 \\Z(2)^0& Z(2^3)^5& Z(2^3)^3& Z(2^3)\\Z(2^3)^6& Z(2^3)^6& Z(2^3)^5& Z(2^3)^4 \end{array}\right)$,
$\left(\begin{array}{cccc} Z(2^3)^4& Z(2^3)& Z(2^3)^4& Z(2^3)^6
\\Z(2^3)^5& Z(2^3)^3& Z(2^3)^4& 0*Z(2) \\Z(2^3)^6& Z(2^3)^5&
Z(2^3)^5& Z(2)^0\\Z(2^3)^6& Z(2)^0& Z(2^3)^6& Z(2^3)^3
\end{array}\right)$,

  $\left(\begin{array}{cccc} Z(2^3)^5& Z(2^3)^3& 0*Z(2)& Z(2)^0 \\ Z(2^3)^5& Z(2)^0& Z(2)^0& Z(2^3)^4\\Z(2^3)^5& Z(2)^0& Z(2^3)^2& 0*Z(2)\\Z(2^3)^5& Z(2^3)^3& Z(2^3)^6& Z(2^3)^3 \end{array}\right)$.

The character table of $G^{s_{11}}$:\\
 \begin{tabular}{c|ccccc}
  & & & & &\\\hline

$\chi_{11}^{(1)}$&1&1&1&1&1
\\$\chi_{11}^{(2)}$&1&A&B&/B&/A
\\$\chi_{11}^{(3)}$&1&B&/A&A&/B
\\$\chi_{11}^{(4)}$&1&/B&A&/A&B
\\$\chi_{11}^{(5)}$&1&/A&/B&B&A
  \end{tabular}

  \noindent  where  A = E(5)$^2$; B = E(5)$^4$.

\vskip 0.3cm

  {\large\bf Acknowledgement}: We would like to thank Prof.
N. Andruskiewitsch and Dr. F. Fantino for suggestions and help.

\end {document}